\documentclass[twoside,11pt]{article}

\usepackage{blindtext}

\usepackage[abbrvbib, preprint]{jmlr2e}

\usepackage{titletoc}
\usepackage[utf8]{inputenc} 
\usepackage[T1]{fontenc}    
\usepackage{hyperref}       
\usepackage{url}            
\usepackage{booktabs}       
\usepackage{amsfonts}       
\usepackage{nicefrac}       

\usepackage{amsmath,amssymb,amsfonts}
\usepackage{graphicx}
\usepackage{textcomp}
\usepackage{xcolor}
\usepackage{hyperref}
\usepackage{diagbox}
\usepackage{array}%
\usepackage{makecell}
\usepackage{threeparttable}%

\usepackage{footmisc}
\usepackage{microtype}      
\usepackage{multirow}
\usepackage{algorithm, algorithmic,float}
\usepackage{graphicx} 
\usepackage{subfigure}
\usepackage{epstopdf}
\usepackage{amsmath}
\usepackage{bm}
\usepackage{graphicx}
\usepackage{caption}
\usepackage{amssymb}
\usepackage{lipsum}
\usepackage{color}
\usepackage{setspace}
\usepackage{soul}
\usepackage{colortbl}
\usepackage{wrapfig}
\usepackage{titletoc}
\setstcolor{red}
\usepackage{arydshln}
\usepackage{float}

\usepackage{xcolor}
\usepackage{colortbl}
\usepackage{makecell}
\usepackage{multirow}

\hypersetup{
colorlinks=true,
linkcolor=black,
citecolor=black,
urlcolor=black
}



\RequirePackage{mdframed}
\renewmdenv[leftmargin=0.4em, rightmargin=0.4em%
innerleftmargin=0.6em,innerrightmargin=0.6em]{quote}

\newtheorem{assumption}{Assumption}

\newcommand{\Ld}{L_{fg,x}}
\newcommand{\Lb}{L_{fg,\theta}}
\newcommand{\sigmaLG}{LG}
\newcommand{\sigmaGT}{GT}

\def \yc  #1{\textcolor{black}{#1}}
\def \nyc  #1{\textcolor{black}{#1}} 
\def \nycres  #1{\textcolor{black}{#1}}

\def \re  #1{\textcolor{black}{#1}}

\def \res #1{\textcolor{black}{#1}}

\def\BibTeX{{\rm B\kern-.05em{\sc i\kern-.025em b}\kern-.08em
    T\kern-.1667em\lower.7ex\hbox{E}\kern-.125emX}}

\definecolor{lightpurple}{RGB}{217, 217, 255}
\newcommand{\ALG}{LoPA}

\newcommand{\ALGa}{LoPA-LG}
\newcommand{\ALGb}{LoPA-GT}

\ShortHeadings{DSBO: Improved Complexity and Heterogeneity Analysis}{Niu, Xu, Sun, Huang, and Chai}

\firstpageno{1}

\begin{document}

\title{Distributed Stochastic Bilevel Optimization: Improved Complexity and Heterogeneity Analysis}

\author{\name Youcheng Niu$^\dagger$  \email ycniu@zju.edu.cn
       \AND
       \name Jinming Xu$^\dagger$* \email jimmyxu@zju.edu.cn
       \AND
       \name Ying Sun$^\ddagger$ \email ybs5190@psu.edu
        \AND
       \name Yan Huang$^\dagger$ \email huangyan5616@zju.edu.cn
        \AND
       \name Li Chai$^\dagger$ \email  chaili@zju.edu.cn\\
       \addr  $^\dagger$College of Control Science and Engineering, Zhejiang University, China\\
      $^\ddagger$School of Electrical Engineering and Computer Science, The Pennsylvania State University, USA
       }

\maketitle

\begin{abstract}
This paper consider solving a class of nonconvex-strongly-convex distributed stochastic bilevel optimization (DSBO) problems with personalized inner-level objectives. Most existing algorithms require computational loops for hypergradient estimation, leading to computational inefficiency. Moreover, the impact of data heterogeneity on convergence in bilevel problems is not explicitly characterized yet. To address these issues, we propose LoPA, a loopless personalized distributed algorithm that leverages a tracking mechanism for iterative approximation of inner-level solutions and Hessian-inverse matrices without relying on extra computation loops. Our theoretical analysis explicitly characterizes the heterogeneity across nodes (denoted by $b$), and establishes a sublinear rate of
$\mathcal{O}( {\frac{1}{{{{\left( {1 - \rho } \right)}}K}} \!+ \!\frac{{(\frac{b}{\sqrt{m}})^{\frac{2}{3}}  }}{{\left( {1 - \rho } \right)^{\frac{2}{3}} K^{\frac{2}{3}} }} \!+ \!\frac{1}{\sqrt{ K }}( {\sigma _{\operatorname{p} }}  + \frac{1}{\sqrt{m}}{\sigma _{\operatorname{c} }}  ) } )$
without the boundedness of local hypergradients, where ${\sigma _{\operatorname{p} }}$ and ${\sigma _{\operatorname{c} }}$
 represent the gradient sampling variances  associated with the inner- and  outer-level variables, respectively.  We also integrate LoPA with a gradient tracking scheme to eliminate the impact of data heterogeneity, yielding an improved rate of ${{\mathcal{O}}}(\frac{{1}}{{ (1-\rho)^2K }} \!+\! \frac{1}{{\sqrt{K}}}( \sigma_{\rm{p}}  \!+\! \frac{1}{\sqrt{m}}\sigma_{\rm{c}} ) )$. The computational complexity of  LoPA is of ${{\mathcal{O}}}({\epsilon^{-2}})$ to an $\epsilon$-stationary point, matching the communication complexity due to the loopless structure, which outperforms existing counterparts for DSBO.
 Numerical experiments validate the effectiveness of the proposed algorithm.

\end{abstract}

\begin{keywords}
distributed learning, bilevel optimization, non-convex optimization, heterogeneity analysis, gradient tracking
\end{keywords}

\section{Introduction}
Bilevel optimization is a hierarchical optimization framework that involves an outer- and  inner-level problem, where the solution of the outer-level problem depends on that of the inner-level problem. This framework has gained significant attention recently in the field of machine learning due to its wide applications in areas such as meta-learning \citep{finn2017model, rajeswaran2019meta}, neural architecture search \citep{zhang2022interpreting, xue2021rethinking}, hyperparameter selection \citep{okuno2018ell, bertrand2020implicit}, and reinforcement learning \citep{wang2016accelerating, TTSA}. With the increasing importance of large-scale machine learning, bilevel optimization has emerged as a promising approach in distributed settings, where multiple nodes with computation and communication capabilities can collaborate to improve the learning efficiency \citep{xu2015augmented, MA-DSBO, SPDB, nedic2020distributed, jiao2022asynchronous}.
Achieving this goal requires   properly coordinating the nodes. In this paper, we aim to address a class of distributed stochastic bilevel optimization (DSBO) problems consisting of $m$ nodes, each with a personalized inner-level objective as follows:
\begin{equation} \label{EQ-DPBO}
\begin{aligned}
\mathop {\min }\limits_{ x \in \mathbb{R}^n} \Phi ( x) = \frac{1}{m}\sum_{i = 1}^m \underbrace{{f_i}\left( { x,\theta _i^*( x)} \right)}_{\Phi_i(x)}, \;{\rm{s.t. \; }}\theta _i^*( x) = \arg \mathop {\min }\limits_{{\theta _i} \in \mathbb{R}^p} {g_i}( x,{\theta _i}),
\end{aligned}
\end{equation}
where $x\in \mathbb{R}^n$ and $\theta_i \in \mathbb{R}^p$ are the global and local model parameters,  respectively; $f_i: \mathbb{R}^n \times  \mathbb{R}^p \to \mathbb{R}$ denotes the outer-level objective of node $i$ which is possibly nonconvex while $g_i: \mathbb{R}^n \times  \mathbb{R}^p \to \mathbb{R}$ is the inner-level objective that is strongly convex in $\theta$ uniformly for all $x \in \mathbb{R}^n$. We consider the stochastic setting where ${f_i}(x,\theta) = {\mathbb{E}_{\varsigma _i  \sim {\mathcal{D}_{{ f_i}}}}}[ {{\hat f_i}( x,\theta, \varsigma _i)} ]$ and ${g_i}( x, \theta) = {\mathbb{E}_{\xi _i  \sim {\mathcal{D}_{{g_i}}}}}\left[ {{\hat g_i}( x,\theta,\xi _i )} \right]$, with  ${\mathcal{D}_{{ f_i}}}$ and ${\mathcal{D}_{{g_i}}}$ denoting the data distribution related to the $i$-th outer- and inner-level objective, respectively.
\\[1ex]
\textbf{Motivating Examples}. Problem \eqref{EQ-DPBO} finds a broad range of applications in practical distributed machine learning and min-max/compositional optimization problems,
ranging from few-shot learning \citep{fewshow}, adversarial learning \citep{madry2017towards}, and reinforcement learning \citep{wang2016accelerating} to fair transceiver design \citep{razaviyayn2020nonconvex}. For instance, consider the following distributed meta-learning problem:
\begin{equation} \label{EQ-APP1}
\begin{aligned}
\mathop {\min }\limits_{x \in {\mathbb{R}^n}} \frac{1}{m}\sum\limits_{i = 1}^m {\sum\limits_{t \in {{{\mathcal{T}}}_i}} {{f _i ^t}( {\theta _i^*\left( x \right)} )}},\;\operatorname{s.t.}\; \theta _i^*\left( x \right) = \arg \mathop {\min }\limits_{{\theta _i} \in {\mathbb{R}^p}} \Big\{ {\sum\limits_{t \in {{{\mathcal{T}}}_i}} \!\! {\left\langle {\theta_i ,\nabla{f _i ^t}( {x} )} \right\rangle } \!\! + \frac{\nu }{2}{{\left\| {x - {\theta _i}} \right\|}^2}} \Big\}, \nonumber
\end{aligned}
\end{equation}
where $x$ is the global model parameter to be learned,  $f_i ^t$ denotes the loss function for the $t$-th subtask  corresponding to the task set $\mathcal{T}_i$ in node $i$, and $\nu>0$  is an adjustable parameter. The objective of nodes  is to cooperatively learn a good initial global model $x$ that makes use of the knowledge obtained from past experiences among nodes to better adapt to  unseen tasks with a small number of task-specific gradient updates \citep{finn2017model, fallah2020personalized, rajeswaran2019meta}.
\renewcommand\arraystretch{0.9}

Different from the conventional single-level problem, bilevel optimization faces additional challenges due to its hierarchical structure, often leading to non-convex objectives \citep{stoBiO,BSA}.  In most cases, obtaining a closed-form expression, or computing exactly the hypergradient $\nabla \Phi_i(x)$  is difficult, due to its dependency on the  inner-level solutions $\theta_i^*(x)$  \citep{ALSET}.  For those with strongly-convex inner-level problems, the  expression of the hypergradients can be obtained by implicit function theorem and further approximated by Approximate Implicit Differentiation (AID) approaches, but it typically involves two nested loops \citep{ji2022will, BSA}: an $N$-loop to find a near-optimal solution to the inner-level function, and a $Q$-loop to approximate the  Hessian-inverse  matrices.  This is particularly challenging in large-scale machine learning applications where running these two nested loops is prohibitively computationally expensive,  and also costly in terms of the training time~\citep{BSA}.
Furthermore, in the stochastic setting, integrating the SGD method with AID  for bilevel optimization with low  computation cost and sample complexity is challenging, due to the estimation bias of the  hypergradient.
To address this issue, various approximation algorithms have been proposed to  estimate hypergradients and reduce the  bias \citep{BSA,TTSA,AmiGO, STABLE, ALSET, SUSTAIN, VRBO, FSLA, SOBA}. To deal with a wide range of large-scale machine learning tasks,  there have been some distributed algorithms recently proposed for DSBO  leveraging the distributed gradient descent approach, such as \citep{yang2022decentralized, SLAM, VRDBO, SPDB}. However, these existing approaches often require  computation loops to estimate the  inner-level solutions and Hessian-inverse matrices for solving DSBO, inducing a   high computation cost.
Therefore, the following question arises naturally: \textit{Can we design loopless decentralized learning algorithms for DSBO problems that achieve better computational complexity and training efficiency?}
Moreover, unlike the  standard distributed optimizations, the bilevel structure and the heterogeneity of inner-level solutions introduce new challenge in characterizing the node heterogeneity and obtaining sharp convergence, with its  influence on the convergence rate unclear. This leads to another important theoretical question: \textit{How to characterize the heterogeneity for DSBO problems, and how does it affect the algorithm's convergence performance?}

\textbf{Summary of Contributions.}
To address the above challenges, we propose a new loopless distributed personalized algorithm (termed {\ALG}) for solving problem \eqref{EQ-DPBO} and provide improved complexity as well as heterogeneity analysis. We summarize the key contributions as follows:
\begin{itemize}
  \item \textbf{New loopless distributed  algorithms}. We propose a new loopless distributed algorithm {\ALG} without requiring extra computation loops for personalized DSBO, which can employ either local gradient or gradient tracking scheme, termed  {\ALGa} and {\ALGb}, respectively. Different from existing distributed algorithms for DSBO problems \citep{MA-DSBO, SLAM, VRDBO, SPDB, yang2022decentralized}, {\ALG} leverages an iterative approximation approach that requires only a single SGD step per iteration to respectively track the Hessian-inverse matrix and the inner-level solution. To control the impact of the bias on the overall convergence rate caused by removing the computation loops,  {\ALG} further employs a gradient momentum step coupled with a relaxation step in the consensus update, ensuring that the bias  decays at a sufficient fast rate.

  \item \textbf{Heterogeneity  analysis}. We quantify the degree of heterogeneity among nodes in personalized DSBO by analyzing that of the inner- and outer-level functions (c.f.,  Lemma \ref{LE-hyper-heterogeneity}). Our analysis shows that  {\ALGa} is able to achieve a convergence rate  of  $\mathcal{O}( {\frac{{\kappa ^{8}}}{{{{\left( {1 - \rho } \right)}}K}} + \frac{\kappa^{\frac{16}{3}}{(\frac{b}{\sqrt{m}})^{\frac{2}{3}}  }}{{\left( {1 - \rho } \right)^{\frac{2}{3}} K^{\frac{2}{3}} }} + \frac{{\kappa^{\frac{5}{2}}}}{\sqrt{ K }}( {\sigma _{\operatorname{p} }}  + \frac{1}{\sqrt{m}}{\sigma _{\operatorname{c} }}  ) } )$ (c.f., Corollary \ref{CO-1}), where  ${\sigma _{\operatorname{p} }}={\mathcal{O}}(\kappa^{\frac{1}{2}}\sigma_{f,\theta}+\kappa^{\frac{3}{2}}\sigma_{g,\theta\theta}+\kappa^{\frac{5}{2}}\sigma_{g,\theta})$, ${\sigma _{\operatorname{c} }}={\mathcal{O}}(\sigma_{f,x}+\kappa\sigma_{g,x\theta})$ are the sampling gradient variances, $b^2={\mathcal{O}}({\kappa ^2}b_{f}^2+{\kappa ^6 b_{g}^2})$ is the heterogeneity among nodes, and $\kappa$ represents the condition number (c.f.,  Assumptions \ref{ASS-heterogeneity}, \ref{ASS-STOCHASTIC}, Equation \eqref{EQ-Lip-const}).
  \nyc{Our analysis relies on the weak assumption of the bounded inner-level heterogeneity at the optimum, while avoiding the boundedness of local hypergradients (c.f., Remark \ref{re-assumption}),}  and the resulting rate shows the clear dependence of convergence on the network connectivity $\rho$, condition number $\kappa$, sampling variances $\sigma_{\rm{p}}, \sigma_{\rm{c}}$  and node heterogeneity $b$.
The techniques used for heterogeneity analysis allow us to explicitly characterize the impact of the heterogeneity on convergence for DSBO compared to existing works \citep{MA-DSBO, SLAM, VRDBO, SPDB, yang2022decentralized} (c.f., Remarks \ref{re-assumption} and \ref{RE-het}), which is of independent interest.

  \item \textbf{Improved complexity}. The analysis of {\ALGa} reveals that the node heterogeneity affects the convergence rate by introducing a transient term of $\mathcal{O}(K^{- \frac{2}{3}})$ in DSBO, vanishing at a slower rate than the leading term, in which the
inner-level heterogeneity plays a crucial role. This motivates us further designing {\ALGb} based on gradient tracking to
 eliminate the heterogeneity. In particular, we prove that {\ALGb} improves the rate to \nyc{${{\mathcal{O}}}( {\frac{{{\kappa ^8}}}{{{{\left( {1 - \rho } \right)^2}}K}} + \frac{{ \kappa^{\frac{5}{2}} }}{ {  {\sqrt{K}} }} (\sigma_{\rm{p}}  + \frac{1}{\sqrt{m}}\sigma_{\rm{c}})   } )$}. Thanks to the loopless structure, both {\ALGa} and {\ALGb} are shown to have a computational complexity of the order of ${{\mathcal{O}}}({\epsilon^{-2}})$ to an $\epsilon$-stationary point, improving existing works on DSBO \citep{MA-DSBO, SLAM, VRDBO, SPDB, yang2022decentralized} by order of ${{\mathcal{O}}}(\log{\epsilon^{-1}})$  (See Table \ref{TA-1}). \nyc{Our analysis further reveals  that our algorithms can achieve the best-known   ${{\mathcal{O}}}({\kappa^5m^{-1}\epsilon^{-2}})$-computational complexity  w.r.t. out-level gradient evaluations  in existing  works on DSBO without mean-square smoothness (c.f.,  Remark \ref{re-improved complexity}).}
 Numerical experiments on machine learning tasks demonstrate the effectiveness of the proposed algorithm in dealing with node heterogeneity and reducing the computational cost.

\end{itemize}

\vspace{-0.5cm}
\section{Related Works}
\textbf{Bilevel optimization with SGD methods.}
There have been some efforts devoted to achieving more accurate stochastic  hypergradients and ensuring convergence in solving bilevel optimization with SGD methods, such as using  computation loops to reduce the error of approximating Hessian-inverse matrices ($Q$-loop) and increase the accuracy of inner-level solutions ($N$-loop) \citep{BSA},  increasing the batch size \citep{AmiGO}, adopting two-timescale step-sizes to eliminate steady-state stochastic variance \citep{TTSA}, incorporating additional correction terms \citep{STABLE, SUSTAIN, FSLA, VRBO, chen2024optimal}, and exploring the smoothness of objectives \citep{ALSET, SOBA}. Among these stochastic approximation algorithms, works \citep{BSA,STABLE, ALSET} achieve a computational complexity of ${ {\mathcal{O}}}\left( \epsilon ^{- 2}\log {\epsilon ^{ - 1}}\right)$ due to use of extra computation loops for estimating the hypergradients, where  $\epsilon$ represents the desired accuracy. \cite{TTSA} develop a two-timescale approximation algorithm, whereas, the nature of the two-timescale update results in a sub-optimal computational complexity of ${{\mathcal{O}}}\left( \epsilon ^{- 5/2}\log {\epsilon ^{ - 1}}\right)$ for the algorithm.  Based on warm-start strategies, a computational complexity of ${{\mathcal{O}}}\left( \epsilon ^{- 2}\right)$ is provided in \citep{AmiGO,FSLA,SOBA,chen2024optimal}. By employing  momentum accelerations  in  both outer- and inner-level optimization procedures such as STORM \citep{STORM} or SPIDER \citep{SPIDER} and imposing more strict assumptions, works \citep{SUSTAIN, VRBO} further improves the rate to ${{\mathcal{O}}}\left( \epsilon ^{- 3/2}\log {\epsilon ^{ - 1}}\right)$. While the aforementioned works provide some insights into  stochastic bilevel algorithmic design, they cannot be directly applied to distributed problem \eqref{EQ-DPBO} as considered in this work.
\begin{table}[!h]
\centering
\footnotesize{
\caption{Comparison of  existing distributed stochastic bilevel optimization algorithms with SGD methods.}
\label{TA-1}
\setlength{\tabcolsep}{0.085cm}{}
\begin{threeparttable}
\begin{tabular}{p{3cm}ccccccc} \hline
  \makecell[cc]{Algorithm}   &Setting     & \# of Loop    &Inner Step      & Complexity  &\small{HA}  & \nycres{Assumption}  &  Scheme  \\  \hline
   \makecell{MA-DSBO \\ \citep{MA-DSBO}}    &\textbf{G}   & $N$-$Q$-Loop  &SGD          &${ {\mathcal{O}}}\left( \epsilon ^{- 2}\log {\epsilon ^{ - 1}}\right)$  &No   & \textbf{$f$-BG}, \textbf{$g$-BH}  &\multirow{4}{*}{LG} \\
  \makecell{Gossip-DSBO \\ \citep{yang2022decentralized}}      &\textbf{G}   & $Q$-Loop  &SGD              &${ {\mathcal{O}}}\left( \epsilon ^{- 2}\log {\epsilon ^{ - 1}}\right)$  & No   &  \textbf{$f$-BG},  \textbf{$g$-BG} \\  \cellcolor{lightpurple}
   \makecell[cc]{ \cellcolor{lightpurple} {\ALGa}  \\   \rowcolor{lightpurple} (this work)  }  &\cellcolor{lightpurple} \textbf{P}  & \cellcolor{lightpurple} No-Loop  & \cellcolor{lightpurple} SGD               &\cellcolor{lightpurple} ${ {\mathcal{O}}}\left( \epsilon ^{- 2}\right)$
    & \cellcolor{lightpurple} Yes & \cellcolor{lightpurple} \textbf{$f$-BH},  \textbf{$g$-BHO} &  \\    \hdashline
   \makecell{SLAM \\ \citep{SLAM}}     &\textbf{G}     & $Q$-Loop  &SGD             &${ {\mathcal{O}}}\left( \epsilon ^{- 2}\log {\epsilon ^{ - 1}}\right)$ & No  & --  &\multirow{5}{*}{GT} \\
   \makecell{VRDBO \\ \citep{VRDBO}}      & \textbf{G}   & $Q$-Loop &STORM              &${ {\mathcal{O}}}( \epsilon ^{- 3/2}\log {\epsilon ^{ - 1}})$ & No  & \textbf{MSS} \\
   \makecell{SPDB \\ \citep{SPDB}}      &\textbf{P}  & $Q$-Loop  &SGD              &${ {\mathcal{O}}}( \epsilon ^{- 2}\log {\epsilon ^{ - 1}})$ & No &-- \\  \cellcolor{lightpurple}
    \makecell[cc]{  \rowcolor{lightpurple}
  {\ALGb}   \\ \rowcolor{lightpurple}  (this work)}                  & \cellcolor{lightpurple} \textbf{P}  & \cellcolor{lightpurple} No-Loop  & \cellcolor{lightpurple} SGD               & \cellcolor{lightpurple} ${ {\mathcal{O}}}\left( \epsilon ^{- 2}\right)$ & \cellcolor{lightpurple} Yes &  \cellcolor{lightpurple} --  & \\  \hline
\end{tabular}
\footnotesize{*\nycres{The algorithms above the dashed line use local gradient (LG) schemes, while those below use gradient tracking (GT) schemes}; \textbf{G} and \textbf{P} represent DSBO with global and  inner-level personalized objectives, respectively;
The complexity represents the Hessian evaluations  required to attain an $\epsilon$-stationary point; HA refers to heterogeneity analysis. \nycres{The column of `assumption'  outlines the \emph{additional} assumptions  beyond the standard Assumptions  \ref{ASS-network}–\ref{ASS-INNERLEVEL}: \textbf{$f$-BG} refers to the gradient of $f_i(x,\theta)$ being bounded for any $x$ and $\theta$, i.e., $\|\nabla f_i(x,\theta)\|^2 \leqslant C$; \textbf{$g$-BH} refers to bounded data heterogeneity of $f_i(x,\theta), i \in  \mathcal{V}$ for any $x$ and $\theta$, i.e., $\sum_{i=1}^{m}\sum_{i=j}^{m} \|\nabla g_i(x,\theta) - \nabla g_j(x,\theta)\|^2 \leqslant C$; \textbf{$g$-BHO} refers to bounded data heterogeneity of $g_i,  i \in  \mathcal{V}$ at the optimum $\theta_i^*(x)$ for any $x$;  \textbf{MSS} refers to mean-squared smoothness, i.e., $\mathbb{E}[\|\nabla\hat{f}_i(w_1,\varsigma)-\nabla\hat{f}_i(w_2,\varsigma)\|]\leqslant l \|w_1-w_2\|$ for any $w_1=(x_1,\theta_1)$ and $w_2=(x_2,\theta_2)$, which is a stronger condition than the function smoothness presented in Assumptions \ref{ASS-OUTLEVEL} and \ref{ASS-INNERLEVEL}.
}}
\end{threeparttable}
}
\end{table}
 \\[1ex]
\textbf{Distributed bilevel optimization.}
Compared to their centralized or parameter-server counterparts, distributed optimization offers several advantages on network scalability, system robustness, and privacy protection through peer-to-peer communication \citep{nedic2020distributed}. However, it also faces unique challenges, especially in dealing with data heterogeneity among nodes. In recent decades, various variants of distributed optimization algorithms have been developed, including distributed gradient descent \citep{chen2021accelerating}, gradient tracking \citep{xu2015augmented}, and alternating direction multiplier methods \citep{shi2014linear}, accompanied by theoretical advancements. Specifically, for stochastic single-level nonconvex problems, these algorithms can achieve a computational complexity of ${ {\mathcal{O}}}\left(\epsilon ^{- 2}\right)$ with SGD methods \citep{koloskova2020unified, chen2021accelerating}. However, these single-level methods  are not readily available to be adapted to tackle  the interaction between the inner and outer levels of functions in solving bilevel optimization problems due to the absence of explicit knowledge of optimal solutions to the inner-level problem.
 \\[1ex]
There have been some efforts aiming at solving distributed bilevel optimization problems, which can be generally cast into two categories: global DSBO and personalized DSBO. For global DSBO, pioneering works such as Gossip-based DSBO (termed Gossip-DSBO) \citep{yang2022decentralized}, MA-DSBO \citep{MA-DSBO}, VRDBO \citep{VRDBO}, and SLAM \citep{SLAM} have been proposed which aim to solve an inner-level problem in the finite sum structure in a distributed manner. On the other hand, few works, such as SPDB \citep{SPDB}, have been developed for solving personalized DSBO problems where each node has its own local inner-level problem. Various algorithmic frameworks have been developed in the above-mentioned works leveraging distributed optimization methods \citep{chen2021accelerating, xu2015augmented, shi2014linear} to minimize the outer- and inner-level functions and handle consensus constraints. To estimate the Hessian-inverse matrices, these frameworks utilize techniques such as Jacobian-Hessian-inverse product \citep{MA-DSBO}, Hessian-inverse-vector product \citep{rajeswaran2019meta} or Neumann series approaches \cite{BSA} to avoid explicit computation of the inverse matrices. As the estimation  of  the  Hessian-inverse matrices and inner-level solutions  obtained using SGD methods is biased \citep{ALSET}, these algorithms often incorporate  extra computation loops to reduce the bias of  the hypergradient estimation due to the approximation of the  Hessian-inverse matrices and inner-level solutions, and the computational complexity of these algorithms is typically of ${\mathcal{O}}(\epsilon^{-2}\log \epsilon^{-1})$ \citep{yang2022decentralized, MA-DSBO, SLAM, SPDB}.
 The complexity can be improved to ${ {\mathcal{O}}}\left( \epsilon ^{- 3/2}\log {\epsilon ^{ - 1}}\right)$  by utilizing variance reduced gradient estimators \citep{VRDBO} and additional mean-squared smoothness assumptions. However, it should be noted that these existing distributed algorithms for DSBO still incur high computational costs due to the extra computation loops required. Moreover,  it  remains unclear how one can properly characterize the heterogeneity among nodes as well as its detailed impact on convergence performance in DSBO.
 \\[1ex]
{\color{black}{\textbf{Parallel  works}. Recently, there have been several studies \citep{zhang2023communication, dong2023single, kong2024decentralized} concurrently exploring loopless algorithms for distributed bilevel optimization problems.  However, these works all focus on global DSBO problems where the out-level objective relies on a common inner-level solution. In particular,  \cite{zhang2023communication} employ variance reduction techniques and GT schemes to achieve an  iteration complexity of $\mathcal{O}(\epsilon^{-{3/2}})$ under the assumption of mean-square smoothness, while \cite{dong2023single} establish an iteration complexity of  $\mathcal{O}(\epsilon^{-{1}})$ with the help of GT schemes for deterministic cases.  \cite{kong2024decentralized} further provide a detailed complexity analysis for their proposed method employing LG schemes.  In contrast,  we address personalized DSBO problems where the outer-level objective relies on \emph{personalized} inner-level solutions, leading to a more challenging scenario for heterogeneity analysis. We thus provide unique analysis for quantifying the node heterogeneity and characterize its impact on the convergence for both LoPA-LG and LoPA-GT under a unified analytical framework, yielding tighter rates that clearly show the dependence of convergence on various factors, such as $\kappa$, $b$, $\sigma_{\rm{p}}, \sigma_{\rm{c}}$ and $\rho$, and improved complexity results in terms of out-level gradient evaluations under weaker  assumptions on data heterogeneity (c.f Remarks \ref{re-assumption}, \ref{re-improved complexity}, Corollaries  \ref{CO-1}, \ref{CO-2}).
}}

\section{Algorithm Design}
In this section,  we will present  the proposed {\ALG} algorithm.
Before  delving into the details of the algorithm, we first provide  some necessary  preliminaries  including relevant network models and assumptions.
\subsection{Preliminaries}

\textbf{Network models.}  We model an undirected   communication network as a weighted graph $\mathcal{G}=(\mathcal{V},\mathcal{E},W)$, where $\mathcal{V}=\{1,\cdots, m\}$ is the set of nodes, $\mathcal{E} \subset \mathcal{V} \times \mathcal{V}$ is the set of edges, and $W=[w_{ij}]_{i,j=1}^m$ is the weight  matrix. The set of neighbors of node $i$ is denoted by $\mathcal{N}_i=\{j\,|\,(i,j)\in \mathcal{E}\}$.
We make the following standard assumption on graph $\mathcal{G}$.\\[-0.5cm]
\begin{assumption}[\textbf{Network connectivity}]\label{ASS-network}
The communication network $\mathcal{G}$  is connected and  the weight matrix $W$ satisfies i)  $w_{ij}=w_{ji}>0$ if and only if $(i,j)\in \mathcal{E}$; and $w_{ij}=0$ otherwise; ii) $W$ is doubly stochastic. Consequently, we have ${\rho} \triangleq {\| { W -  {\frac{{{1_m}1_m^{\rm{T}} }}{m}}} \|^2} \in \left[ {0,1} \right)$.
\end{assumption}


In what follows, we make several assumptions on the  outer- and inner-level functions of problem~\eqref{EQ-DPBO}, which are common in the existing literature of bilevel  optimization \citep{BSA, TTSA,  MA-DSBO,yang2022decentralized,chen2024optimal,SPDB}.
\begin{assumption}[\textbf{Outer-level functions}]\label{ASS-OUTLEVEL}
Let $L_{f,x}$, ${L_{{f},\theta }}$ and  ${C_{{f},\theta }}$ be positive constants.   Each outer-level function $(x,\theta) \mapsto f_i(x, \theta)$, $i \in \mathcal{V}$  satisfies the following properties:\\[.5ex]
(i) $f_i$ is continuously differentiable;\\[.5ex]
(ii) For any $ \theta \in \mathbb{R}^p$, ${\nabla _x}{f_i}( \cdot ,\theta )$ is  $L_{f,x}$-Lipschitz-continuous in $x$; and for any $x \in \mathbb{R}^n$, ${\nabla _\theta }{f_i}(x,\cdot )$ is ${L_{{f},\theta }}$-Lipschitz-continuous in $\theta$;\\[.5ex]
(iii) For any $x \in \mathbb{R}^n$, $\|{\nabla _\theta }{f_i}(x,\theta _i^*(x) )\| \leqslant C_{f,\theta}$.
\end{assumption}

\begin{assumption}[\textbf{Inner-level functions}]\label{ASS-INNERLEVEL}
Let ${\mu _{{g}}}$, ${L_{{g},\theta }}$, ${L_{{g},x\theta }}$, ${L_{{g},\theta \theta }}$ and ${C_{{g}, x \theta}}$ be positive constants.  Each inner-level function $(x, \theta) \mapsto g_i(x,\theta)$, $i \in \mathcal{V}$ satisfies the following properties:\\[.5ex]
(i) For any $x \in \mathbb{R}^n$, $g_i(x,\cdot)$   is  ${\mu _{{g}}}$-strongly convex  in $\theta$; $g_i$  is twice continuously differentiable;\\[0.5ex]
(ii) For any $x \in \mathbb{R}^n$, ${\nabla _\theta }{g_i}\left( {x,\cdot } \right)$ is ${L_{{g},\theta }}$-Lipschitz-continuous in $\theta$;  ${\nabla _{x\theta }^2}{g_i}\left( {\cdot,\cdot } \right)$,  ${\nabla _{\theta \theta }^2}{g_i}\left( {\cdot,\cdot } \right)$  are  respectively ${L_{{g},x\theta }}$- and ${L_{{g},\theta \theta }}$-Lipschitz-continuous;\\[0.5ex]
(iii) For any $\theta \in \mathbb{R}^p$ and $x \in \mathbb{R}^n$,  $\|{\nabla _{x\theta }^2}{g_i}\left( {x,\theta} \right)\|\leqslant {C_{{g},x\theta}}$.
\end{assumption}


Next, recalling that ${f_i}(x,\theta) = {\mathbb{E}_{\varsigma _i  \sim {\mathcal{D}_{{ f_i}}}}}[ {{\hat f_i}( x,\theta, \varsigma _i)} ]$ and ${g_i}( x, \theta) = {\mathbb{E}_{\xi _i  \sim {\mathcal{D}_{{g_i}}}}}\left[ {{\hat g_i}( x,\theta,\xi _i )} \right]$, we proceed to make the following assumption regarding the data heterogeneity across nodes for problem \eqref{EQ-DPBO}, which resembles that of distributed single-level optimization \citep{chen2021accelerating, koloskova2020unified,lian2017can}.

\begin{assumption}[\textbf{Bounded    heterogeneity}]\label{ASS-heterogeneity}
Let ${\nabla _x}f(x,\theta) \triangleq \frac{1}{m}\sum\nolimits_{j = 1}^m {{\nabla _x}{f_j}(x,\theta)}$ and \\ $ {\nabla _\theta }f(x,\theta) \triangleq \frac{1}{m}\sum\nolimits_{j = 1}^m {{\nabla _\theta }{f_j}(x,\theta)}$. 
There exist positive constants $b_f^2$ and $b_g^2$ such that: \\
(i) $\sum_{i = 1}^m {{{\left\| {{\nabla _x}{f_i}(x,\theta) - {\nabla _x}f(x,\theta)} \right\|}^2}}  \leqslant b_f^2$, $\sum_{i = 1}^m {{{\left\| {{\nabla _\theta }{f_i}(x,\theta) - {\nabla _\theta }f(x,\theta)} \right\|}^2}}  \leqslant b_f^2$, for any $x$ and $\theta$;\\
(ii) $\frac{1}{m}\sum_{i = 1}^m {\sum_{j = 1}^m {{{\| {{\nabla _\theta }{g_i}(x,\theta _j^*(x))} -{  {\nabla _\theta }{g_j}(x,\theta _j^*(x))} \|}^2} \leqslant } } b_g^2$, for any $x$ and $\theta _i^*(x)$, $i \in \mathcal{V}$; \\
(iii) $\frac{1}{m}\sum_{i = 1}^m {\sum_{j = 1}^m {{{\| {\nabla _{x\theta }^2{g_i}( {x,\theta _j^*(x)} ) - \nabla _{x\theta }^2{g_j}( {x,\theta _j^*(x)} )} \|}^2}} }  \leqslant b_g^2$,\\
\phantom{(v) } $\frac{1}{m}\sum_{i = 1}^m {\sum_{j = 1}^m {{{\| {\nabla _{\theta \theta }^2{g_i}( {x,\theta _j^*(x)} ) - \nabla _{\theta \theta }^2{g_j}( {x,\theta _j^*(x)} )} \|}^2} \leqslant  b_g^2} }$, for any $x$ and $\theta _i^*(x)$, $i \in \mathcal{V}$;\\
where  $\nabla_x f$ and $\nabla_{\theta} f$ represent the partial gradient with respect to $x$ and $\theta$, respectively, while ${\nabla _{x\theta }^2}{g}$ and   ${\nabla _{\theta \theta }^2}{g}$ denote Jacobian and Hessian, respectively.
\end{assumption}
\begin{remark}[\textbf{Weaker  assumptions \nycres{on data heterogeneity}}]\label{re-assumption}
The parameters $b_{{f}}^2$ and $b_{{g}}^2$  are  introduced  to  quantify  the data heterogeneity on the outer- and inner-level functions  across nodes, respectively.
\nycres{It is worth noting that Assumption \ref{ASS-heterogeneity}(i) is weaker than that of the previous works based on LG schemes for DSBO problems, such as \citep{MA-DSBO,  yang2022decentralized} where the norm of $ \nabla _{x} f_i (x,\theta)$ is assumed to be bounded by a constant.} As for the heterogeneity in the inner-level functions (Assumption \ref{ASS-heterogeneity}(ii) and \ref{ASS-heterogeneity}(iii)), we only require that it is uniformly bounded at the optimum $\theta_i^*(x)$, $i \in \mathcal{V}$, for all $x$. This requirement is less restrictive than assuming the inner-level heterogeneity to be bounded at any $\theta$  as  in \citep{kong2024decentralized}.
\end{remark}

\subsection{The Proposed {\ALG} Algorithm}
In this section, we present our algorithm, termed {\ALG},  for problem \eqref{EQ-DPBO}. Following the standard procedures as  in distributed optimization~\citep{chen2021accelerating, nedic2009distributed}, we let each node $i$ maintain a local estimate $x_i$ for the global decision variable $x$. 
At each iteration $k$, each node $i$ alternates between a  descent step with an estimate of $\nabla {\Phi _i}\left( x_i \right)$ and average consensus ensuring the consistency of the $x_i$'s.\\[1ex]
\noindent{\bf Hypergradient construction.} By the chain rule and implicit function theorem \citep{BSA},  we can compute the Hypergradient $\nabla {\Phi _i}\left( x_i \right)$ as follows:
\[\nabla {\Phi _i}\left( x_i \right) = {\nabla _x}{f_i}(x_i,\theta _i^*(x_i)) - {\nabla _{x\theta }^2}{g_i}(x_i,\theta _i^*(x_i)){\left[ {{\nabla _{\theta \theta }^2}{g_i}(x_i,\theta _i^*(x_i))} \right]^{ - 1}}{\nabla _\theta }{f_i}(x_i,\theta _i^*(x_i)).\]

Note that computing the outer-level gradients in each iteration according to the above expression is computationally demanding.
To address this issue,  we first introduce an auxiliary  variable $\theta _i$ to approximate the inner-level solutions $\theta_i^*(x_i)$, whose update follows a simple stochastic gradient descent step (c.f.,~\eqref{EQ-ALG-a} and~\eqref{EQ-ALG-d}).
As such, the Hessian-inverse-vector products  ${\left[ {{\nabla _{\theta \theta }^2}{g_i}(x_i,\theta _i^*(x_i))} \right]^{ - 1}}{\nabla _\theta }{f_i}(x_i,\theta _i^*(x))$ (abbreviated  Hv) and \nyc{local hypergradient thus can be approximately computed as follows \citep{dagreou2022framework, AmiGO}}:
\begin{align}
  &v_i = {\left[ {{\nabla _{\theta \theta }^2}{g_i}(x_i,\theta_i )} \right]^{ - 1}}{\nabla _\theta }{f_i}(x_i,\theta_i ) \label{EQ-surrogate-v},\\
  &s_i = {\nabla _x}{f_i}(x_i,\theta_i ) - {\nabla _{x\theta }^2}{g_i}(x_i,\theta_i )v_i.\label{EQ-surrogate-s}
\end{align}

In the presence of stochasticity, steps~\eqref{EQ-surrogate-v} and~\eqref{EQ-surrogate-s} still cannot be computed directly as the gradient and Hessian are unknown.
To overcome this issue, for step~\eqref{EQ-surrogate-s}, we  replace ${\nabla _x}{f_i}$ and $\nabla _{x\theta }^2{g_i}$ with their stochastic estimates, respectively (c.f.,,~\eqref{EQ-ALG-f}).
As for step~\eqref{EQ-surrogate-v}, notice that $v_i$ can be regarded as the solution of the following strongly convex problem: $$v_i = \arg\min_{v} \big\{ \frac{1}{2} v^{\rm T}  \nabla _{\theta \theta }^2{g_i}(x_i,\theta_i )v - {\nabla _\theta }{f_i}(x_i,\theta_i ) v \big\}.$$

Instead of directly computing the solution using stochastic approximation methods, we propose to approximate it by performing only one stochastic gradient iteration which warm-starts with the value of $v_i$ initialized to its value from the previous iteration. By doing so, we aim to approximate the solution more efficiently while taking advantage of the progress made in the previous iteration. Further details are provided in Remark \ref{rek-NS-CG} while the specific updates are given by~\eqref{EQ-ALG-b} and~\eqref{EQ-ALG-e}.

Putting all the ingredients together,  {\ALG}  designs the following  updates for the local inner-level variable $\theta _i^{k+1}$, the local  Hv estimate variable $v_i^{k+1}$ and the local hypergradient $s_i^{k+1}$ at iteration $k+1$:
\begin{align}
  &   \theta _i^{k + 1} = \theta _i^k - \beta d_i^k, \label{EQ-ALG-a} \\
  & v_i^{k + 1} = v_i^k - \lambda h_i^k.  \label{EQ-ALG-b}
\end{align}
where $\beta>0$ and $\lambda>0$ are  step-sizes. The directions $ d_i^k$ and $ h_i^k$ are further updated as:
 \begin{align}
  &  d_i^{k + 1} = {\nabla _\theta }{\hat g_i}(x_i^{k + 1},\theta _i^{k + 1};\xi _{i,1}^{k + 1}),  \label{EQ-ALG-d} \\
  & h_i^{k + 1} = \nabla _{\theta \theta }^2{\hat g_i}(x_i^{k + 1},\theta _i^{k + 1};\xi _{i,2}^{k + 1})v_i^{k + 1} - {\nabla _\theta }{\hat f_i}( {x_i^{k + 1},\theta _i^{k + 1};\varsigma _{i,1}^{k + 1}} ). \label{EQ-ALG-e}
\end{align}
With $\theta_i^{k+1}$ and $v_i^{k+1}$ at hand, the local hypergradient is approximated by
\begin{equation}\label{EQ-ALG-f}
\begin{aligned}
&\nyc{s_i^{k + 1}}= {\nabla _x}{\hat f_i}( {x_i^{k + 1},\theta _i^{k + 1};\varsigma _{i,2}^{k + 1}} ) - \nabla _{x\theta }^2{\hat g_i}( {x_i^{k + 1},\theta _i^{k + 1};\xi _{i,3}^{k + 1}} )v_i^{k + 1}.
\end{aligned}
\end{equation}
Here,  $\nabla_\theta \hat{g}_i$ (resp. $\nabla _{\theta \theta }^2{\hat g_i}$, ${\nabla _\theta }{\hat f_i}$, ${\nabla _x}{\hat f_i}$, $\nabla _{x\theta }^2{\hat g_i}$) denotes  a stochastic gradient estimate of $\nabla_\theta {g}_i$ (resp. $\nabla _{\theta \theta }^2{ g_i}$, ${\nabla _\theta }{ f_i}$, ${\nabla _x}{ f_i}$, $\nabla _{x\theta }^2{ g_i}$) depending on the random variable $\xi _{i,1}^{k+1}$ (resp. $\xi _{i,2}^{k+1}$, $\xi _{i,3}^{k+1}$, $\varsigma _{i,1}^{k+1}$, $\varsigma _{i,2}^{k+1}$).

However, using the above stochastic estimators will lead to steady-state  errors under Assumption \ref{ASS-network}-\ref{ASS-heterogeneity} and single-timescale step-sizes \citep{ALSET, TTSA}, unless an increasing number of batch sizes is used \citep{AmiGO} or extra smoothness conditions are imposed on Hv variables \citep{SOBA} (refer to Section \ref{SE-diff} for more details).  To address this issue, we further introduce a gradient momentum step to reduce the impact  as follows:

\begin{equation}\label{EQ-ALG-g}
\begin{aligned}
   \nyc{ z_i^{k+1} = s_i^{k} + (1 - \gamma )(z_i^k - s_i^k),}
\end{aligned}
\end{equation}
where $0<\gamma<1$. The gradient momentum step helps control the order of the sampling variance by maintaining a moving average of the past  approximated hypergradients \citep{MA-DSBO, ghadimi2020single}.\\[1ex]
\noindent{\bf Distributed gradient descent/tracking.} The update of local copy $x_i$ follows standard distributed gradient method as:
\begin{equation}\label{EQ-ALG-c}
\begin{aligned}
x_i^{k + 1} = (1-\tau)x_i^k+ \tau (\sum\nolimits_{j \in {{\mathcal{N}}_i}} {{w_{ij}}x_j^k}  - {\alpha}y_i^k).
\end{aligned}
\end{equation}
where $\alpha > 0$ denotes the step size, and a relaxation step characterized by a parameter $\tau \in (0, 1)$ is employed to smooth both the consensus and the local gradient processes.  Then,  two alternative choices of the direction $y_i^k$ are considered:
 \begin{align}
\text{(Local gradient scheme)} \;\; &y_i^{k+1}=z_i^{k + 1}.  \label{EQ-ALG-h1} \\
\text{(Gradient tracking scheme)} \;\; &y_i^{k + 1} = \sum\nolimits_{j \in {{\mathcal{N}}_i}} {{w_{ij}}y_j^k}  + z_i^{k + 1} - z_i^k.\label{EQ-ALG-h2}
\end{align}
Eq.~\eqref{EQ-ALG-h1} yields a distributed gradient descent type algorithm where $y_i^k$ estimates the local hypergradient $\nabla f_i$; whereas~\eqref{EQ-ALG-h2} employs the gradient tracking technique so that $y_i^k$ estimates the average hypergradient $\nabla \Phi$.


The overall {\ALG} algorithm is summarized in Algorithm~\ref{alg:1}, where we refer to the versions with local gradient and gradient tracking schemes as {\ALGa} and {\ALGb}, respectively.
\begin{algorithm}[ht]
	\caption{{\ALG}}
	\label{alg:1}
	\begin{algorithmic}[1]
		\STATE \textbf{Require}: Initialize  $\theta_i^0$,  $v_i^0$, $x_i^0$, $s_i^0$,  $z_i^0$, $y_i^0$, $i \in \mathcal{V}$ and set step-sizes $\{\alpha, \beta, \lambda, \gamma, \tau \}$.
            \FOR {{$ k=0,1,2..., K$}, each node $i \in \mathcal{V}$ in parallel}
             \STATE Sample batch $\zeta_i^{k+1}= \{\xi_{i,1}^{k+1}, \xi_{i,2}^{k+1}, \xi_{i,3}^{k+1}, \varsigma_{i,1}^{k+1}, \varsigma_{i,2}^{k+1}\}$.\\
             \STATE Communicate with  neighboring node  $j \in \mathcal{N}_i$.\\
             \STATE  {Update state variables  ${\theta}_i^{k+1}$, ${v}_i^{k+1}$, $x_i^{k+1}$ according to  \eqref{EQ-ALG-a}, \eqref{EQ-ALG-b}}, \eqref{EQ-ALG-c};\\
             \STATE  {Update local gradient estimates $d_i^{k+1}$, ${h}_i^{k+1}$, $\nyc{s_i^{k + 1}}$, $z_i^{k+1}$  according to \eqref{EQ-ALG-d}, \eqref{EQ-ALG-e}, \eqref{EQ-ALG-f}, \eqref{EQ-ALG-g}};\\
             \STATE{Update the descent direction of outer-level variables  $y_i^{k+1}$ as follows:}
             \IF  {gradient tracking scheme is not used}
              \STATE{{\ALGa}: Update $y_i^{k+1}$  according to  \eqref{EQ-ALG-h1};}
             \ELSE{}
              \STATE{{\ALGb}: Update $y_i^{k+1}$  according to \eqref{EQ-ALG-h2}. }
             \ENDIF
       \ENDFOR
	\end{algorithmic}
\end{algorithm}
\begin{remark}[\textbf{Iterative approximation approach  for Hv}] \label{rek-NS-CG}
In estimating hypergradients, {\ALG} takes a different approach from existing distributed algorithms   that use Neumann Series (NS) methods and \nycres{Summed Hessian Inverse Approximation (SHIA) methods}
to directly approximate the Hessian-inverse matrices to a high-precision  at each iteration $k$. To be more specific, the key idea of the NS \citep{BSA, TTSA} and \nycres{SHIA methods \citep{dagreou2022framework,rajeswaran2019meta}} is to approximate the Hessian-inverse matrices and Hv in multiple iterations, respectively. The approximation process of these two methods can be summarized as follows:
\begin{align}
& \textbf{ ({NS}):} \;Q\lambda \prod_{t = 0}^Q {{{\left( {I - \lambda {\nabla _{\theta \theta }^2}{ g_i}\left( {{x_i},{\theta _i}} \right)} \right)}}}   \approx  {\left[ {{\nabla _{\theta \theta }^2}{g_i}\left( {{x_i},{\theta _i}} \right)} \right]^{ - 1}},  \label{EQ-NSCG}\\
& \nycres{\textbf{({SHIA})}}: \;\nyc{ \lambda \sum_{t = 0}^Q\Big({\prod_{j = 0}^t {{{{\left( {I - \lambda {\nabla _{\theta \theta }^2}{ g_i}( {{x_i},{\theta _i}} )} \right)}}}} }\Big) {\nabla _\theta }{ f_i}( {{x_i},{\theta _i}} ) \!\approx \! {\left[ {{\nabla _{\theta \theta }^2}{g_i}\left( {{x_i},{\theta _i}} \right)} \right]^{ - 1}}{\nabla _\theta }{f_i} ( {{x_i},{\theta _i}} ). }\notag
\end{align}
We can know from the above expressions that  the high-precision approximation generally requires a large $Q$, which leads to extra computation loops  at each iteration $k$. For examples, the state-of-the-art  works \citep{SPDB,SLAM,VRDBO} require $Q$ obeying $\mathcal{O}(\log{\epsilon ^{-1}})$.   Unlike these methods,  {\ALG}  adopts an iterative approximation approach with one stochastic gradient iteration for tracking the states of Hessian-inverse matrices and inner-level solutions. Thus, {\ALG}  enjoys a loopless structure in the algorithmic design and  achieves  a  computational complexity of $\mathcal{O}(1)$ with respect to the number of Hessian evaluations at each iteration $k$ while maintaining  the same complexity for outer- and inner-level gradient and Jacobian evaluations at each iteration $k$.
\end{remark}
\vspace{-0.2cm}

\section{Convergence Results}
In this section, we respectively analyze the performance of  {\ALGa} and   {\ALGb} for  nonconvex-strongly-convex cases.

\subsection{Preliminaries}
We make the following assumption on the  stochastic gradients used for estimating the gradients of the outer- and inner-level functions.

Let
\begin{equation}
\begin{aligned}
{\mathcal{F}^k} = \sigma \left\{ {\bigcup\nolimits_{i = 1}^m {(x_i^0,\theta _i^0,v_i^0,z_i^0,y_i^0, \cdots ,x_i^k,\theta _i^k,v_i^k,z_i^k,y_i^k)} } \right\}
\end{aligned}
\end{equation}
be the $\sigma$-algebra generated by the random variables up to the $k$-th iteration.

\begin{assumption}[\textbf{Stochastic gradient estimates}]\label{ASS-STOCHASTIC}
We assume the random variables  $\xi _{i,1}^{k}$, $\xi _{i,2}^{k}$, $\xi _{i,3}^{k}$, $\varsigma _{i,1}^{k}$, $\varsigma _{i,2}^{k}$ are mutually independent  for any iteration $k$; and also independent across all the iterations. Furthermore, for any $x\in \mathbb{R}^n$, $\theta \in  \mathbb{R}^p$ and $k \geqslant 0$,  the followings hold:\\[0.5ex]
(i) Unbiased estimators:
  \begin{align*}
      & \mathbb{E}[{\nabla _\theta }{\hat g_i}(x,\theta ;\!\xi _{i,1}^{k} )]\!\!=\!\!{\nabla _\theta }{ g_i}(x,\theta), \,
      \mathbb{E}[{\nabla _{\theta \theta }^2}{\hat g_i}(x,\theta ;\!\xi _{i,2}^{k})]\!\!=\!\!{\nabla _{\theta \theta }^2}{g_i}(x,\theta ),
      \mathbb{E}[{\nabla _{x\theta }^2}{\hat g_i}(x,\theta ;\!\xi _{i,3}^{k})]\!\!=\!\!{\nabla _{x \theta }^2}{g_i}(x,\theta ),\\
      & \mathbb{E}[{\nabla _\theta}{\hat f_i}(x,\theta ; \varsigma _{i,1}^{k})]={\nabla _\theta}{f_i}(x,\theta),\, \mathbb{E}[{\nabla _x}{\hat f_i}(x,\theta ;\varsigma _{i,2}^{k})]={\nabla _x}{ f_i}(x,\theta).
  \end{align*}
(ii) Bounded stochastic variances:
  \begin{align*}
 & \mathbb{E}[ {{{\| {{\nabla _\theta }{g_i}(x,\theta ;\xi _{i,1}^{k}) - {\nabla _\theta }{g_i}(x,\theta )} \|}^2}} ]\! \leqslant\! {\sigma _{g,\theta}^2}, \mathbb{E}[ {{{\| {\nabla _{\theta \theta }^2{{\hat g}_i}(x,\theta ;\xi _{i,2}^{k}) - \nabla _{\theta \theta }^2{g_i}(x,\theta )} \|}^2}} ] \!\leqslant \!{\sigma _{g,\theta\theta}^2}, \\
 &\mathbb{E}[ {{{\| {\nabla _{x\theta }^2{{\hat g}_i}(x,\theta ;\xi _{i,3}^{k}) - \nabla _{x\theta }^2{g_i}(x,\theta )} \|}^2}} ] \!\!\leqslant \!\! {\sigma _{g,x\theta}^2}, \mathbb{E}[ {{{\| {{\nabla _\theta }{{\hat f}_i}(x,\theta ;\varsigma _{i,1}^{k}) - {\nabla _\theta }{f_i}(x,\theta )} \|}^2}} ] \!\!\leqslant\!\!{\sigma _{f,\theta}^2}, \\
 & \mathbb{E}[ {{{\| {{\nabla _x}{{\hat f}_i}(x,\theta ;\varsigma _{i,2}^{k }) - {\nabla _x}{f_i}(x,\theta )} \|}^2}}] \leqslant {\sigma _{f,x}^2}.
\end{align*}
\end{assumption}

The following two propositions provide the smoothness property of $\nabla \Phi(x)$, $\theta_i^*(x)$, and $v_i^*(x)$, as well as the boundedness of $v_i^*(x)$, with the first proposition being adapted from \citep{BSA} and the second being derived from the properties of certain gradients. For completeness, we provide the complete proof in Section \ref{sec-Pro}.

\begin{proposition}[\textbf{Smoothness property}]\label{PR-smooth}
Suppose Assumptions \ref{ASS-OUTLEVEL} and \ref{ASS-INNERLEVEL}  hold. Let $\bar \nabla {f_i}( {x,\theta})$ $\triangleq\nabla_x f(x,\theta)-\nabla_{x\theta}^2g_i(x,\theta)v_i(x,\theta)$ be a surrogate of the local  hypergradient $ \nabla {f_i}(x,\theta_i^*(x))$ and denote ${v_i}\left( {x,\theta } \right)\triangleq {\left[ {\nabla _{\theta \theta }^2{g_i}\left( {x,\theta } \right)} \right]^{ - 1}}{\nabla _\theta }{f_i}\left( {x,\theta } \right),  v_i^*\left( x \right) \triangleq {v_i}\left( {x,\theta _i^*(x)} \right)$.  Then given any $x,x' \in \mathbb{R}^n$,
it holds that: $\forall i  \in \mathcal{V}$,
\[\begin{gathered}
  \!\| {\theta _i^*(x) \!-\!\theta _i^*(x')} \| \!\leqslant\! {L_{\theta ^*}}\!\| {x \!- \!x'} \|, \| {{v_i}( {x,\theta _i^*(x) } ) \!- \!{v_i}\!( {x',\theta _i^*(x')} )} \| \! \leqslant \! {L_{{v}}}( {\| {x \!-\! x'} \|\! +\! \| {\theta _i^*(x)  \!-\! \theta _i^*(x')} \|} ), \hfill \\
  \!\| {v_i^*\!\left( x \right) \!-\! v_i^*(x')} \|\! \leqslant\! {L_{v^*}}\!\| {x \!-\! x'} \|, \| {\bar \nabla \!{f_i}( {x,\theta _i^*(x)} ) \!-\! \bar \nabla \!{f_i}( {x',\theta _i^*(x')} \!)} \| \!\!\leqslant \! {L_f}( {\| {x - x'} \| \!+ \!\| {\theta _i^*(x)  \!- \!\theta _i^*(x')} \|} ), \hfill \\
  \!\| {\nabla \Phi (x) - \nabla \Phi (x')} \| \leqslant L\| {x - x'} \|, \hfill
\end{gathered} \]
\end{proposition}
where the Lipschitz  constants are provided  as follows:
\begin{align}
  &{L_{{\theta ^*}}} \!\triangleq \!\frac{{{C_{g, x\theta}}}}{{{\mu _g}}}\!=\! {{\mathcal{O}}}\!\left( {\kappa } \right),  {L_v} \triangleq{\frac{{{L_{f,\theta }}}}{{{\mu _g}}}\!+\! \frac{{{C_{f,\theta }}{L_{g,\theta \theta }}}}{{{\mu _g ^2}}}}\!=\!{{\mathcal{O}}}\!( {\kappa ^2} )\!, {L_{{v^*}}} \!\triangleq \!( {\frac{{{L_{f,\theta }}}}{{{\mu _g}}} \!+\! \frac{{{C_{f,\theta }}{L_{g,\theta \theta }}}}{{{\mu _g ^2}}}} )\!( {1 \!+ \!{L_{{\theta ^*}}}}\! )\!= \!{{\mathcal{O}}}\! \left( {\kappa ^3} \right),\hfill \nonumber \\
  &{L_f} \!\triangleq\! {L_{f,x}} \!+ \!{C_{g,x\theta }}{L_v} \!+\! \frac{{{C_{f,\theta }}{L_{g,x\theta }}}}{{{\mu _g}}}=\nyc{{{\mathcal{O}}}\left( {\kappa ^2} \right)}, L \!\triangleq\! ( {{L_{f,x}} \!+ \!{C_{g,x\theta }}{L_v} \!+ \!\frac{{{C_{f,\theta }}{L_{g,x\theta }}}}{{{\mu _g}}}} )( {1\! +\! {L_{{\theta ^*}}}} )=\nyc{{{\mathcal{O}}}\left( {\kappa ^3} \right)}, \hfill \label{EQ-Lip-const}
\end{align}
with \nyc{$\kappa=\frac{\max\{L_{f,x}, L_{f,\theta}, L_{g,\theta}, L_{g,x\theta}, L_{g,\theta\theta}\}}{\mu_g}$} denoting the condition number.
\nyc{
\begin{proposition}[\textbf{Boundness property}]\label{PR-boundness}
Suppose Assumptions \ref{ASS-OUTLEVEL} and \ref{ASS-INNERLEVEL}  hold. Then, there exists a constant $M=\frac{{{C_{f,\theta }}}}{{{\mu _g}}}$ such that the following holds for any $x$: $\forall i  \in \mathcal{V}$,
\begin{equation}\label{EQ-V-Boundness}
\begin{aligned}
\| {v_i^*{(x)}} \| \leqslant  M.
\end{aligned}
\end{equation}
where $v_i^*\left( x \right)= {v_i}\left( {x,\theta _i^*(x)} \right)$ as defined in Proposition \ref{PR-smooth}.
\end{proposition}
}


\subsection{Convergence of {\ALGa} and {\ALGb}}
\textbf{Convergence of {\ALGa}.} We first analyze {\ALGa} that uses local gradient scheme \eqref{EQ-ALG-h1}. To derive the convergence results of {\ALGa}, one key step is to explore the heterogeneity on  overall hypergradients. The following lemma shows the boundness and composition  of the  heterogeneity on overall hypergradients.

\begin{lemma}[\textbf{Bounded  heterogeneity on  overall hypergradients}]\label{LE-hyper-heterogeneity}
Suppose Assumptions \ref{ASS-OUTLEVEL}, \ref{ASS-INNERLEVEL} and \ref{ASS-heterogeneity} hold. Let $\nabla {\Phi _i}\left( x \right)$ be the local hypergradient of node $i$ evaluated at $x$. Then, we have
\begin{equation}
\begin{aligned}
\sum\limits_{i = 1}^m {{{\left\| {\nabla {\Phi _i}\left( x \right) - \nabla \Phi \left( x \right)} \right\|}^2}} \leqslant {b^2},
\end{aligned}
\end{equation}
where $b^2\triangleq C_1(\mu_g,C_{g,x\theta})b_{{f}}^2 +  C_2(\mu_g,L_{f,x},L_{f,\theta},L_{g,x\theta},L_{g,\theta\theta},C_{f,\theta},C_{g,x\theta})b_{{g }}^2$ with $C_1(\mu_g,C_{g,x\theta})$  and $C_2(\mu_g,L_{f,x},L_{f,\theta},L_{g,x\theta},L_{g,\theta\theta},C_{f,\theta},C_{g,x\theta})$ being the constants defined in Appendix \ref{sec-poof-heterogeneity}.
\end{lemma}
The proof of Lemma \ref{LE-hyper-heterogeneity} is deferred to Appendix \ref{sec-poof-heterogeneity}. {\hfill $\blacksquare$}
\\

 Note that the heterogeneity of the overall hypergradients consists of two main components:  the inner-level heterogeneity $b_{{f}}^2$ and  the outer-level heterogeneity  $b_{{g}}^2$. To further quantify the effect of  heterogeneity in each level,  we   characterize its dependency on  the condition number $\kappa$. It is not difficult to show  that $C_1(\mu_g,C_{g,x\theta})=\mathcal{O}(\kappa^2)$  and $C_2(\mu_g,L_{f,x},L_{f,\theta},L_{g,x\theta},L_{g,\theta\theta},C_{f,\theta},C_{g,x\theta})=\mathcal{O}(\kappa^6)$.  According to the definition of $b^2$, we can observe that  $b^2$ is of ${\mathcal{O}}({\kappa^2}b_{f}^2+{\kappa^6 b_{g}^2})$. This observation  suggests that the heterogeneity of inner-level objective functions plays a crucial role in determining the  heterogeneity of overall hypergradients.  Lemma \ref{LE-hyper-heterogeneity} provides a novel characterization on the heterogeneity in personalized DSBO, \nyc{distinguishing itself from previous works \citep{MA-DSBO,chen2021accelerating} where
  the out-level objective depends on a common inner-level solution and thus the out-level heterogeneity is  independent of the inner-level heterogeneity.} Now, we are thus ready to present the convergence of {\ALGa}.

\begin{theorem}\label{TH-1}
Suppose Assumptions \ref{ASS-network}, \ref{ASS-OUTLEVEL}, \ref{ASS-INNERLEVEL}, \ref{ASS-heterogeneity} and \ref{ASS-STOCHASTIC} hold.  Consider the sequence  $\{x_i^k, \theta_i^k, $ $v_i^k, z_i^k, y_i^k\}$ generated by Algorithm 1 employing  local gradient scheme as depicted in \eqref{EQ-ALG-h1}. Let $\bar x^k=(1/m)\sum_{i=1}^{m} x_i^k$,  \nycres{${L_{fg,x}} = 2L_{f,x}^2 + 4{M^2}L_{g,x\theta }^2$ and  ${L_{fg,\theta}}=2L_{f,\theta }^2+ 4{M^2}L_{g,\theta \theta }^2$} with $M=\frac{{{C_{f,\theta }}}}{{{\mu _g}}}$.  There exists a proper choice of step-sizes $\alpha, \beta, \lambda, \gamma, \tau$ such that, \yc{for any total number of iterations $K$}, we have
\begin{equation}\label{EQ-TH-1-INQ}
\begin{aligned}
\frac{1}{{K + 1}}\sum\limits_{k = 0}^K {\mathbb{E}[ {{{\| {\nabla \Phi ({{\bar x}^k})} \|}^2}} ]}  \leqslant \frac{{4({V^0} - {V^K})}}{\tau\alpha(K+1)} + 4\alpha \nycres{\sigma_{\sigmaLG}^2} + \frac{4\vartheta}{m}\alpha^2 {b^2},
\end{aligned}
\end{equation}
where  $ \nycres{\sigma_{\sigmaLG} ^2}=(\frac{1}{m}\frac{d_3}{\tau } + \frac{d_4}{\tau } )(\sigma _{f,x}^2 + 2{M^2}\sigma _{g,x\theta }^2)\frac{{{\gamma ^2}}}{{{\alpha ^2}}} + 2\frac{d_1}{\tau } ( {\sigma _{f,\theta }^2 + 2{M^2}\sigma _{g,\theta\theta }^2} )\frac{{{\lambda ^2}}}{{{\alpha ^2}}} + 2\frac{d_2}{\tau } \sigma _{g,\theta }^2\frac{{{\beta ^2}}}{{{\alpha ^2}}}$ and $\vartheta   =\frac{{24 \alpha \varphi }}{{{{(1 - \rho )}^2}\gamma }}$ with $\varphi  = ({L_{fg,x}} + \frac{{32C_{g,x\theta }^2{L_{fg,\theta }}}}{{{\mu _g ^2}}})(1 + \frac{{4L_{g,\theta }^2}}{{\omega _\theta ^2}})$ and $\omega_{\theta}  = \frac{{{\mu _g}{L_{g,\theta }}}}{{2\left( {{\mu _g} + {L_{g,\theta }}} \right)}}$. The coefficients $d_0$, $d_1$, $d_2$, $d_3$, $d_4$, $d_5$, $d_6$ within the unified Lyapunov function $V^k$ as depicted in  \eqref{EQ-V} are defined as follows: ${d_0} = 1$, ${d_1} = \frac{{8C_{g,x\theta }^2\tau \alpha }}{{{\mu _g}\lambda }}$, ${d_2} = ( {L_{fg,x}}+\frac{{32C_{g,x\theta }^2{L_{fg,\theta }}}}{{\mu _g^2}})\frac{{\tau \alpha }}{{{\omega _\theta }\beta }}$, ${d_3} = \frac{{\tau \alpha }}{{2\gamma }}$, ${d_4} = \frac{{24\varphi \tau {\alpha ^4}}}{{{{(1 - \rho )}^2}\gamma }}$, ${d_5} = \frac{{2\varphi \alpha }}{{1 - \rho }}$, $d_6=0$, and the network spectral gap $\rho$ is defined in Assumption~\ref{ASS-network}.
\end{theorem}
The proof sketch and supporting lemmas of Theorem \ref{TH-1} are provided in Section \ref{sec-proof} and  the complete proof is deferred to Appendix \ref{sec-proof-TH-1}. {\hfill $\blacksquare$}

\begin{corollary}\label{CO-1}
Consider the same setting as Theorem \ref{TH-1}.  If the step-sizes are properly chosen such that $\alpha  = \min \left\{ {u,{{\left( {\frac{{{a_0}}}{{{a_1}\left( {K + 1} \right)}}} \right)}^{\frac{1}{2}}}, {{\left( {\frac{{{a_0}}}{{{a_2}\left( {K + 1} \right)}}} \right)}^{\frac{1}{3}}}} \right\} \leqslant \frac{1}{m}$ for a large value of $K$, where $u$, $a_0, a_1$ and $ a_2$ are given in Appendix \ref{sec-proof-CO-1}, and {$\gamma=\mathcal{O}(\alpha)$, $\lambda=\mathcal{O}(\alpha)$, $\beta=\mathcal{O}(\alpha)$},  then  we have\footnote{The symbol $\mathcal{O}$ hides both the constants and parameters associated with the properties of functions.}
\begin{equation}\label{EQ-rate-1}
\begin{aligned}
\frac{1}{{K + 1}}\sum\limits_{k = 0}^K {{\mathbb{E}}[ {{{\| \nabla {\Phi ({{\bar x}^k})} \|}^2}} ]}  = {{\mathcal{O}}}\Big( {\frac{{\kappa ^{8}}}{{{{\left( {1 - \rho } \right)}}K}} + \frac{\kappa^{\frac{16}{3}}{(\frac{b}{\sqrt{m}})^{\frac{2}{3}}  }}{{\left( {1 - \rho } \right)^{\frac{2}{3}} K^{\frac{2}{3}} }} + \frac{{1}}{\sqrt{ K }}\kappa^{\frac{5}{2}}( {\sigma _{\operatorname{p} }}  + \frac{1}{\sqrt{m}}{\sigma _{\operatorname{c} }}  ) } \Big).
\end{aligned}
\end{equation}

\noindent where \nyc{${\sigma _{\operatorname{p} }}={\mathcal{O}}(\kappa^{\frac{1}{2}}\sigma_{f,\theta}+\kappa^{\frac{3}{2}}\sigma_{g,\theta\theta}+\kappa^{\frac{5}{2}}\sigma_{g,\theta})$ and ${\sigma _{\operatorname{c} }}={\mathcal{O}}(\sigma_{f,x}+\kappa\sigma_{g,x\theta})$} are the  gradient sampling variances associated with the inner-level and the out-level variables, respectively.


\end{corollary}
The proof of Corollary \ref{CO-1} is deferred to Appendix  \ref{sec-proof-CO-1}. {\hfill $\blacksquare$}

\begin{remark}[\textbf{Heterogeneity analysis}] \nyc{\label{RE-het}  Corollary \ref{CO-1} characterizes the detailed impact of the  heterogeneity at each level on the convergence rate in DSBO by unique heterogeneity analysis under weaker Assumption 4 on data heterogeneity, which does not require the boundedness of local hypergradients or the boundedness of inner-level heterogeneity at any point (c.f., Assumption \ref{ASS-heterogeneity}).}
The above Corollary \ref{CO-1} also shows that   {\ALGa}  has a convergence rate of  ${\mathcal{O}}(\frac{1}{{\sqrt K }})$, where $K$ denotes the total number of iterations. In particular,  the analysis of {\ALGa}  reveals that the data heterogeneity affects the convergence rate by introducing a transient term of $\mathcal{O}(K^{- \frac{2}{3}})$ in DSBO, vanishing at a slower rate than the leading term, \nyc{in which the inner-level heterogeneity $b_g ^2$ plays a crucial role in personalized DSBO}.  To the best of our knowledge, these results are  first presented in our work.
\end{remark}
\textbf{Convergence of {\ALGb}.} Now, we move on to present the main results for {\ALGb} that employs gradient tracking scheme \eqref{EQ-ALG-h2} in the following Theorem \ref{TH-2}.

\begin{theorem}\label{TH-2}
Suppose Assumptions \ref{ASS-network}, \ref{ASS-OUTLEVEL}, \ref{ASS-INNERLEVEL}, and \ref{ASS-STOCHASTIC} hold.  Consider the sequence   $\{x_i^k, \theta_i^k, v_i^k,$ $ z_i^k, y_i^k\}$  generated by Algorithm \ref{alg:1} employing  gradient tracking scheme \eqref{EQ-ALG-h2}. Let $\bar x^k=(1/m)\sum_{i=1}^{m} x_i^k$,  \nycres{${L_{fg,x}} = 2L_{f,x}^2 + 4{M^2}L_{g,x\theta }^2$ and ${L_{fg,\theta}}=2L_{f,\theta }^2+ 4{M^2}L_{g,\theta \theta }^2$}  with $M=\frac{{{C_{f,\theta }}}}{{{\mu _g}}}$.  There exists a proper choice of the step-sizes $\alpha, \beta, \lambda, \gamma, \tau$ such that, for any total number of iterations $K$, we have
\begin{equation}\label{EQ-TH-2-INQ}
\begin{aligned}
\frac{1}{K+1}\sum\limits_{k = 0}^K  {{\mathbb{E}}[ {{{\| \nabla{\Phi ({{\bar x}^k})} \|}^2}} ]} \leqslant \frac{{2\left( {{V^0} - {V^K}} \right)}}{{ \tau\alpha (K+1)}} + {2{\alpha \nycres{{\sigma_{\sigmaGT} ^2}}}},
 \end{aligned}
\end{equation}
where $\nycres{{\sigma_{\sigmaGT} ^2}}= (\frac{1}{m}\frac{{{d_3}}}{\tau } + \frac{{{d_4}}}{\tau } + \frac{{2{d_6}}}{\tau(1-\rho) })(\sigma _{f,x}^2 +2 {M^2}\sigma _{g,x\theta }^2)\frac{{{\gamma ^2}}}{{{\alpha ^2}}} + 2\frac{{{d_1}}}{\tau }( {\sigma _{f,\theta }^2 + 2{M^2}\sigma _{g,\theta \theta }^2} )\frac{{{\lambda ^2}}}{{{\alpha ^2}}} + 2\frac{{{d_2}}}{\tau }\sigma _{g,\theta }^2\frac{{{\beta ^2}}}{{{\alpha ^2}}}$ with $\varphi  = ({L_{fg,x}} + \frac{{32C_{g,x\theta }^2{L_{fg,\theta }}}}{{{\mu _g ^2}}})(1 + \frac{{4L_{g,\theta }^2}}{{\omega _\theta ^2}})$ and $\omega_{\theta}  = \frac{{{\mu _g}{L_{g,\theta }}}}{{2\left( {{\mu _g} + {L_{g,\theta }}} \right)}}$. The coefficients $d_0$, $d_1$, $d_2$, $d_3$, $d_4$, $d_5$, $d_6$ within the unified Lyapunov function $V^k$ as depicted in  \eqref{EQ-V} are defined as follows: ${d_0} = 1$, ${d_1} = \frac{{8C_{g,x\theta }^2\tau \alpha }}{{{\mu _g}\lambda }}$, ${d_2} = ( {L_{fg,x}}+\frac{{32C_{g,x\theta }^2{L_{fg,\theta }}}}{{\mu _g^2}})\frac{{\tau \alpha }}{{{\omega _\theta }\beta }}$, ${d_3} = \frac{{\tau \alpha }}{{2\gamma }}$, ${d_4} = \frac{{64\varphi \tau \gamma {\alpha ^2}}}{{{{(1 - \rho )}^4}}}$, ${d_5} = \frac{{4\varphi \alpha }}{{(1 - \rho )}}$, $d_6=\frac{{16\varphi \tau {\alpha ^2}}}{{{{(1 - \rho )}^3}}}$, and  the network spectral gap $\rho$ is  defined in Assumption \ref{ASS-network}.
\end{theorem}
The proof sketch and supporting lemmas of Theorem \ref{TH-2} are provided in Section \ref{sec-proof} and the complete proof of Theorem \ref{TH-2} is deferred to Appendix \ref{sec-proof-TH-2}. {\hfill $\blacksquare$}

\begin{corollary}\label{CO-2}
Consider the same setting as Theorem \ref{TH-2}.    If the step-sizes are properly chosen such that $\alpha  = \min \left\{ {u',{{\left( {\frac{{{a'_0}}}{{{a'_1}(K + 1)}}} \right)}^{\frac{1}{2}}}} \right\} \leqslant \frac{1}{m}$ for  a large value of $K$, where  $u'$, ${a'_0}$ and $ {a'_1}$  are given in Appendix \ref{sec-proof-CO-2},   and {$\gamma=\mathcal{O}(\alpha)$, $\lambda=\mathcal{O}(\alpha)$, $\beta=\mathcal{O}(\alpha)$}, then we have
\begin{equation}\label{EQ-rate-2}
\begin{aligned}
\frac{1}{{K + 1}}\sum\limits_{k = 0}^K {{\mathbb{E}}[ {{{\| \nabla {\Phi ({{\bar x}^k})} \|}^2}} ]}  ={{\mathcal{O}}}\Big( {\frac{{{{{\kappa ^8}}}}}{{{{\left( {1 - \rho } \right)^2}}K}} + \frac{{ {\kappa ^{\frac{5}{2}}} }}{ {  {\sqrt{K}} }} (\sigma_{\operatorname{p}}  + \frac{1}{\sqrt{m}}\sigma_{\operatorname{c}})   } \Big).
\end{aligned}
\end{equation}
\end{corollary}
\noindent where ${\sigma _{\operatorname{p} }}={\mathcal{O}}(\kappa^{\frac{1}{2}}\sigma_{f,\theta}+\kappa^{\frac{3}{2}}\sigma_{g,\theta\theta}+\kappa^{\frac{5}{2}}\sigma_{g,\theta})$ and ${\sigma _{\operatorname{c} }}={\mathcal{O}}(\sigma_{f,x}+\kappa\sigma_{g,x\theta})$.
The proof of Corollary \ref{CO-2} is deferred to Appendix \ref{sec-proof-CO-2}. {\hfill $\blacksquare$}
\begin{remark}[\textbf{Linear speed-up}]
Corollaries \ref{CO-1} and \ref{CO-2} demonstrate that, in the context of DSBO problems with personalized inner-level objectives, the  term $\sigma_{\rm{p}}$ associated with the inner-level variable $\theta$  does not decay with $m$.  In contrast,  the term $\sigma_{\rm{c}}$ induced by  the stochastic gradients   with respect to the out-level variable $x$  decays at order $O(\frac{1}{\sqrt{mK}})$. This lack of linear speed-up in $\sigma_{\rm{p}}$ is due to the fact that the inner-level solutions $\theta_i^*(x)$ are different for each node $i$ such that neither  the local inner-level variables or the local Hv variables  do not need to reach a consensus in personalized DSBO, and thus the related variances can only be reduced via temporal averaging, yielding an order of $O(\frac{1}{\sqrt{K}})$. In contrast, the outer-level variable $x$ is shared across the nodes and can achieve a speed up with respect to the network size $m$.
Our result also implies that to achieve the right balance of the two terms, the batch size should be increased by $m$ times  when estimating the local partial gradients $\nabla_{\theta}\hat{f}_i$,  $\nabla_{\theta}\hat{g}_i$, $\nabla_{\theta\theta}\hat{g}_i$. In this case, both {\ALGa} and  {\ALGb}  can obtain a convergence rate of $\mathcal{O}(\frac{1}{\sqrt{mK}})$, leading to an iteration complexity of $\mathcal{O}(m^{-1}\epsilon^{-2})$.
\end{remark}

\begin{remark}[\textbf{Improved complexity}] \label{re-improved complexity}
Corollary \ref{CO-2}  shows that {\ALGb} has a convergence rate of ${\mathcal{O}}(\frac{1}{{\sqrt K }})$. It is also noted from the above result that, the effect of the heterogeneity is eliminated by gradient tracking scheme \eqref{EQ-ALG-h2} in personalized DSBO even when the inner-level variables lack communication. Thanks to the loopless structure, both {\ALGa} and {\ALGb} can achieve a computational  complexity of ${\mathcal{O}}(\epsilon^{-2})$ (w.r.t. the number of Hessian evaluations) to attain an $\epsilon$-stationary point.  This computational complexity improves  existing state-of-the-art works on DSBO problems by the order of  ${\mathcal{O}}(\log(\epsilon ^{-1}))$. Note that the computational complexity for inner- and outer-level gradient and Jacobian evaluations is also  of the order ${\mathcal{O}}(\epsilon^{-2})$ for {\ALGa} and {\ALGb}. \nyc{Our analysis for both LoPA-LG and LoPA-GT provides tighter rates that clearly show the dependence of convergence on
the condition number, heterogeneity at each level, sampling variances at each level, and network connectivity, in which  the detailed effects of sampling variances (c.f., \eqref{EQ-rate-1} and \eqref{EQ-rate-2}) enable us to improve the rate by properly optimizing batch sizes. In particular, when the batch sizes of the gradient evaluations $\nabla_{\theta}\hat f$, $\nabla_{\theta}\hat g$,  $\nabla_{\theta \theta}^2\hat g$, and $\nabla_{x \theta}^2\hat g$ are respectively chosen as $\mathcal{O}(m\kappa)$, $\mathcal{O}({m}\kappa^{5})$,   $\mathcal{O}({m}\kappa^3)$, and $\mathcal{O}(\kappa^2)$ times relative to $\nabla_{x}\hat f$, an improved iteration rate of $\mathcal{O}(\frac{\kappa^{\frac{5}{2}}}{\sqrt{mK}})$  can be obtained, which matches the result in centralized counterparts \citep{ALSET}. In this case,  a computational complexity of $\mathcal{O}(\kappa^5m^{-1}\epsilon^{-2})$ w.r.t.  out-level gradient evaluations (also $\mathcal{O}(\kappa^8\epsilon^{-2})$  w.r.t. Hessian evaluations\footnote{Here, the absence of linear speed-up of Hessian complexity is due to the fact that, different from global DSBO problems, the Hessian-inverse-vector product variables do not need to reach a consensus in personalized DSBO.}) can be achieved for both LoPA-LG and LoPA-GT due to the loopless structure, which are the best-known result among existing works on DSBO without mean-square smoothness assumption.}
\end{remark}
\nyc{\textbf{Dependency on network connectivity}}. \nyc{Note that it follows from the results \eqref{EQ-TH-1-INQ}, \eqref{EQ-TH-2-INQ} and the corresponding definitions of $d_0, \cdots ,d_6$ in  \eqref{EQ-the1-dddd} and \eqref{EQ-the2-dddd2} that Theorems \ref{TH-1} and \ref{TH-2} have respective dependencies of $\mathcal{O}(\frac{1}{(1-\rho)^2})$ and $\mathcal{O}(\frac{1}{(1-\rho)^4})$ on  network connectivity  without proper choice of step-sizes.  These results of dependency are aligned with that of existing single-level distributed methods  \citep{lian2017can,alghunaim2022unified}  using LG schemes and GT schemes. By optimizing the step-size selection, we can obtain a tighter dependence of   $\mathcal{O}(\frac{1}{1-\rho})$ and $\mathcal{O}(\frac{1}{(1-\rho)^2})$ for LoPA-LG and LoPA-GT (c.f.,  Corollary \ref{CO-1} and Corollary \ref{CO-2}), respectively, where the dependence of  $\mathcal{O}(\frac{1}{1-\rho})$ for LG schemes is superior to  $\mathcal{O}(\frac{1}{(1-\rho)2})$ provided in the work \citep{yang2022decentralized}.  Notably, our result also reveals that LoPA-GT has a higher dependence on the network spectral gap compared to LoPA-LG, due to the introduction of GT schemes to eliminate heterogeneity. }
\vspace{0.5cm}
\\
\nycres{\textbf{Convergence of LoPA-LG and LoPA-GT in  deterministic case.} In what follows, we present the convergence results in Corollaries \ref{CO-3} and \ref{CO-4} for special cases where the inner- and outer-level functions in problem \eqref{EQ-DPBO} are deterministic. These results can be directly derived by following a similar proof to that of Corollaries \ref{CO-1} and \ref{CO-2}, with the sampling variances set to zero.}

\begin{corollary}\label{CO-3}
\nycres{Consider the deterministic case for Problem \eqref{EQ-DPBO} under the same setting as Theorem \ref{TH-1} for LoPA-LG. Suppose Assumptions \ref{ASS-network}-\ref{ASS-heterogeneity} hold.
If the step-sizes are properly chosen such that $\alpha  = \min \left\{ {u, {{\left( {\frac{{{a_0}}}{{{a_2}\left( {K + 1} \right)}}} \right)}^{\frac{1}{3}}}} \right\}$ for a large value of $K$, where $u$, $a_0$ and $ a_2$ are the same as that of Corollary \ref{CO-1},
and {$\gamma=\mathcal{O}(\alpha)$, $\lambda=\mathcal{O}(\alpha)$, $\beta=\mathcal{O}(\alpha)$},  then  we have}
\begin{equation}\label{EQ-COR-3}
\begin{aligned}
\frac{1}{{K + 1}}\sum\limits_{k = 0}^K {{\mathbb{E}}[ {{{\| \nabla {\Phi ({{\bar x}^k})} \|}^2}} ]}  = {{\mathcal{O}}}\Big( {\frac{{\kappa ^{8}}}{{{{\left( {1 - \rho } \right)}}K}} + \frac{\kappa^{\frac{16}{3}}{(\frac{b}{\sqrt{m}})^{\frac{2}{3}}  }}{{\left( {1 - \rho } \right)^{\frac{2}{3}} K^{\frac{2}{3}} }}  } \Big).\nonumber
\end{aligned}
\end{equation}
\end{corollary}

\begin{corollary}\label{CO-4}
\nycres{Consider the deterministic case for Problem \eqref{EQ-DPBO} under the same setting as Theorem \ref{TH-2} for LoPA-GT. Suppose Assumptions \ref{ASS-network}-\ref{ASS-INNERLEVEL}  hold.
If the step-sizes are chosen such that $\alpha \leqslant u'$, where $u'$ is the same as that of Corollary \ref{CO-2},
and {$\gamma=\mathcal{O}(\alpha)$, $\lambda=\mathcal{O}(\alpha)$, $\beta=\mathcal{O}(\alpha)$}, then we have}
\begin{equation}\label{EQ-COR-4}
\begin{aligned}
\frac{1}{{K + 1}}\sum\limits_{k = 0}^K {{\mathbb{E}}[ {{{\| \nabla {\Phi ({{\bar x}^k})} \|}^2}} ]}  ={{\mathcal{O}}}\Big( {\frac{{{{{\kappa ^8}}}}}{{{{\left( {1 - \rho } \right)^2}}K}}} \Big). \nonumber
\end{aligned}
\end{equation}
\end{corollary}

\vspace{0.5cm}
\nycres{From Corollaries \ref{CO-3} and \ref{CO-4}, we observe that the convergence rates of LoPA-LG and LoPA-GT for the deterministic case can be improved to $\mathcal{O}(\frac{1}{K^{2/3}})$ and $\mathcal{O}(\frac{1}{K})$, respectively.}

\subsection{Proof Sketch and Supporting Lemmas for Theorems \ref{TH-1} and \ref{TH-2} } \label{sec-proof}
Before presenting the proof sketch, we first introduce some essential notations as follows:
\begin{equation}
\begin{aligned}
  {v_i}\left( {x,\theta } \right)&\triangleq {\left[ {\nabla _{\theta \theta }^2{g_i}\left( {x,\theta } \right)} \right]^{ - 1}}{\nabla _\theta }{f_i}\left( {x,\theta } \right),\; v_i^*\left( x \right) \triangleq {v_i}\left( {x,\theta _i^*(x)} \right),  \\
\bar x^k&\triangleq(1/m)\sum\nolimits_{i=1}^{m} x_i^k, \; {x^{k}}\triangleq{\rm{col}}{\{x_i^k\}}_{i=1}^m, \;{\theta ^*}({{\bar x}^k}) \triangleq \operatorname{col} \{ \theta _i^*({{\bar x}^k})\} _{i = 1}^m, \\
{v^*}({{\bar x}^k}) & \triangleq \operatorname{col} \{ v_i^*({{\bar x}^k})\} _{i = 1}^m, \nabla \tilde \Phi ({{\bar x}^{k }}) =  \operatorname{col}\{ \nabla {\Phi _i}({{\bar x}^{k }})\} _{i = 1}^m.
\end{aligned}
\end{equation}
The notations $z^k$, $y^k$, $\theta^k$, $s^k$, $\bar y^k$, $\bar z^k$, $\bar s^k$ share the similar definitions. The key idea of the proof for Theorems \ref{TH-1} and \ref{TH-2} is to characterize the dynamics of the following unified Lyapunov function with properly selected coefficients $d_0$, $d_1$, $d_2$, $d_3$, $d_4$, $d_5$, $d_6$:
\begin{equation}\label{EQ-V}
\begin{aligned}
  {V^k} = &{d_0}\Phi ({{\bar x}^k}) + {d_1}\frac{1}{m}\underbrace {{{\| {{v^k} - {v^*}({{\bar x}^k})} \|}^2}}_{\operatorname{Hv} \:\operatorname{errors} } + {d_2}\frac{1}{m}\underbrace {{{\| {{\theta ^k} - {\theta ^*}({{\bar x}^k})} \|}^2}}_{\operatorname{inner - level\:errors}} + \nyc{{d_3}\underbrace {{{\| { {\nabla \Phi ({{\bar x}^k}) - {{\bar z}^k}}} \|}^2}}_{\operatorname{ave-variance} \: \operatorname{errors}}} \\
   &  +  {d_4}\nyc{\frac{1}{m}\underbrace {{{\| { {\nabla \tilde \Phi ({{\bar x}^k}) - {z^k}}} \|}^2}}_{\operatorname{variance} \: \operatorname{errors}}}+  {d_5}\frac{1}{m}\underbrace {{{\| {{x^k} - {1_m} \otimes {{\bar x}^k}} \|}^2}}_{\operatorname {consensus\:errors}} + {d_6}\frac{1}{m}\underbrace {{{\| {{y^k} - {1_m} \otimes {{\bar y}^k}} \|}^2}}_{\operatorname{gradient\:errors}}, \\
 \end{aligned}
\end{equation}
where the detailed definitions of these coefficients can be founded in \eqref{EQ-the1-dddd} and \eqref{EQ-the2-dddd2}, corresponding to Theorems 6 and 9, respectively. To this end, we proceed to derive iterative evolution for each term of $V^k$  in  expectation according to the following four key steps:
\\
\\
\textbf{Step 1 (Quantifying the descent of the overall objective function)}:  We begin by quantifying the descent of the overall objective function $\Phi ({{\bar x}^k})$ evaluated at the average point  by using its smoothness and  the tracking property of $\bar y^k$ for $\bar z^k$. This  descent    is controlled by the hypergradient approximation errors $\mathbb{E}[ {{{\| {\nabla \Phi ({{\bar x}^k}) - {{\bar z}^k}} \|}^2}} ]$ in  Lemma \ref{LE-descent}.
 \begin{lemma}[\textbf{Descent lemma}] \label{LE-descent}
 Consider the sequence $\{x_i^k, \theta_i^k, v_i^k, z_i^k, y_i^k\}$ generated by Algorithm \ref{alg:1}. Suppose Assumptions \ref{ASS-network}, \ref{ASS-OUTLEVEL}, \ref{ASS-INNERLEVEL} and \ref{ASS-STOCHASTIC}  hold. Then, we have:
 \begin{equation}\label{EQ-descentlemma}
\begin{aligned}
  \mathbb{E}[ {\Phi ({{\bar x}^{k + 1}})} ] \leqslant & \mathbb{E}[ {\Phi ({{\bar x}^k})} ] - \frac{\tau\alpha }{2}\mathbb{E}[ {{{\|\nabla {\Phi ({{\bar x}^k})} \|}^2}} ] - \frac{\tau\alpha }{2}\left( {1 - \tau\alpha L} \right)\mathbb{E}[ {{{\| {{{\bar y}^k}} \|}^2}} ] \\
   &+ \nyc{\frac{\tau \alpha}{2} \mathbb{E}[ {{{\| \nabla{\Phi ({{\bar x}^k}) - {{\bar z}^k}} \|}^2}} ]}.
 \end{aligned}
\end{equation}
\end{lemma}
The proof of Lemma \ref{LE-descent} is provided in  Section \ref{sec-proof-LE-descent}.  {\hfill $\blacksquare$}
\\
\\
\textbf{Step 2 (Characterizing the average variance errors and  hypergradient errors)}: We then deal with the average variance errors $\mathbb{E}[ {{{\|  { \nabla \Phi (\bar{x}^k)} -{{ \bar z}^k} \|}^2}} ]$ according to the bounded variances of different stochastic gradients  and the updates of ${{ z}^k}$. However,  bounding the hypergradient errors  $\mathbb{E}[ {{{\| {\nabla \Phi ({{\bar x}^k}) - {{\bar s}^k}} \|}^2}} ]$  is more challenging  as it requires an investigation into  how the iterative approximation strategies with one stochastic gradient iteration influence  the evolution of Hv errors and inner-level errors. Lemma \ref{LE-HypergradientE} shows that the Hv product errors $\mathbb{E}[ {{{\| {{v^k} - {v^*}({{\bar x}^k})} \|}^2}} ]$, inner-level errors $\mathbb{E}[ {{{\| {{\theta ^k} - {\theta ^*}({{\bar x}^k})} \|}^2}} ]$ and the consensus errors $\mathbb{E}[ {{{\| {{x^k} - {1_m} \otimes {{\bar x}^k}} \|}^2}} ]$ jointly control the hypergradient approximation errors $\mathbb{E}[ {{{\| {\nabla \Phi ({{\bar x}^k}) - \mathbb{E}[{{\bar{s}^k}}]} \|}^2}} ]$, while the term $\mathbb{E}[ {{{\| {\nabla \Phi ({{\bar x}^k}) - \mathbb{E}[{{\bar{s}^k}}}] \|}^2}} ]$ controls the average variance errors $\mathbb{E}[ {{{\|  { \nabla \Phi (\bar{x}^k)} -{{ \bar z}^k} \|}^2}} ]$.
  In what follows, we focus on quantifying these four error terms and establishing their recursions in Lemmas  \ref{LE-Vstar}, \ref{LE-ThetaStar}, \ref{LE-CE}.

\begin{lemma} [\textbf{Hypergradient approximation errors and average variance errors}]  \label{LE-HypergradientE}
$\;\;\;$
Consider the sequence $\{x_i^k, \theta_i^k, v_i^k, z_i^k, y_i^k\}$  generated by Algorithm \ref{alg:1}. Suppose Assumptions \ref{ASS-network}, \ref{ASS-OUTLEVEL}, \ref{ASS-INNERLEVEL} and \ref{ASS-STOCHASTIC} hold. If the step-size $\gamma$  satisfies $0<\gamma<1$,  then we have:
\nyc{
\begin{equation}\label{LE-stochastic-error}
\begin{aligned}
\mathbb{E}[ {{{\| {\nabla \Phi ({{\bar x}^{k + 1}}) - {{\bar z}^{k + 1}}} \|}^2}} ] \leqslant& (1 - \gamma )\mathbb{E}[ {{{\| {\nabla \Phi ({{\bar x}^k}) - {{\bar z}^k}} \|}^2}} ] + {r_z}\alpha \mathbb{E}[ {{{\| {\nabla \Phi ({{\bar x}^k}) - \mathbb{E}[{{\bar{s}^k}}]} \|}^2}} ] \\
&+ {r_y}{\tau ^2}{\alpha ^2}\mathbb{E}[ {{{\| {{{\bar y}^k}} \|}^2}} ] + \frac{1}{m}\sigma _{\bar z}^2{\alpha ^2},
\end{aligned}
\end{equation}
}
and
\begin{equation}\label{EQ-hypergradient-error}
\begin{aligned}
\mathbb{E}[ {{{\| {\nabla \Phi ({{\bar x}^k}) - \mathbb{E}[{{\bar{s}^k}}]} \|}^2}} ] \leqslant& \frac{\yc{\Ld}}{m}\mathbb{E}[ {{{\| {{x^k} - {1_m} \otimes {{\bar x}^k}} \|}^2}} ] + \frac{\yc{\Ld}}{m}\mathbb{E}[ {{{\| {{\theta ^k} - {\theta ^*}({{\bar x}^k})} \|}^2}} ] \\
&+ \frac{{4C_{g,x\theta }^2}}{m}\mathbb{E}[ {{{\| {{v^k} - {v^*}({{\bar x}^k})} \|}^2}} ].
 \end{aligned}
\end{equation}
where ${r_z} \triangleq \frac{{2\gamma }}{\alpha }$, ${r_y} \triangleq \frac{{2{L^2}}}{\gamma }$, and $\sigma _{\bar z}^2\triangleq \res{({\sigma _{f,x}^2}+\nyc{\hat M^2}{\sigma _{g,x\theta}^2})}{\frac{\gamma^2}{\alpha^2}}$ with \nyc{$\hat M^2 \triangleq 2M^2+ 2\frac{1}{m}\|v^k-v^*(\bar x^k)\|^2$} and  $\Ld= 2L_{f,x}^2 + 4{M^2}L_{g,x\theta }^2$.
\end{lemma}
The proof of Lemma \ref{LE-HypergradientE} is provided in  Section \ref{sec-proof-LE-HypergradientE}.  {\hfill $\blacksquare$}
\begin{lemma}[\textbf{Hessian-inverse-vector product errors}] \label{LE-Vstar}
Consider the sequence $\{x_i^k, \theta_i^k,$ $ v_i^k, z_i^k, y_i^k\}$  generated by Algorithm \ref{alg:1}. Suppose Assumptions \ref{ASS-network}, \ref{ASS-OUTLEVEL}, \ref{ASS-INNERLEVEL} and \ref{ASS-STOCHASTIC} hold.  If the step-size $\lambda$  satisfies  \begin{equation}\label{EQ-lambda}
\begin{aligned}
\lambda<\frac{1}{{\mu _g}},
\end{aligned}
\end{equation}
then we have:
\begin{equation}\label{EQ-v*-v}
\begin{aligned}
  \mathbb{E}[ {{{\| {{v^{k + 1}} - {v^*}({{\bar x}^{k + 1}})} \|}^2}} ] \leqslant& ( 1 -{ {\mu _g}{\lambda}} )\mathbb{E}[ {{{\| {{v^k} - {v^*}({{\bar x}^k})} \|}^2}} ] + {q_x}\alpha \mathbb{E}[ {{{\| {{x^k} - {1_m} \otimes {{\bar x}^k}} \|}^2}} ] \\
     &+ {q_x}\alpha \mathbb{E}[ {{{\| {{\theta ^k} - {\theta ^*}({{\bar x}^k})} \|}^2}} ] + m{q_s}{\tau^2 \alpha ^2}\mathbb{E}[ {{{\| {{{\bar y}^k}} \|}^2}} ] + m{\sigma _v^2 }{\alpha ^2}, \\
\end{aligned}
\end{equation}
where ${q_{x  }}\triangleq \frac{{4\yc{\Lb} {\lambda }}}{{{\mu _g}}\alpha}$, ${\sigma _v^2 }\triangleq 2(\sigma _{f,\theta}^2+\nyc{\hat M^2}\sigma _{g,\theta\theta}^2){\frac{\lambda^2}{\alpha^2}}$, ${q_{s}} \triangleq \frac{{2L_{{v^*}}^2 }}{{\varpi {\lambda }}}$ with
\yc{$\Lb=2L_{f,\theta }^2+ 4{M^2}L_{g,\theta \theta }^2 $} and $\varpi  =\frac{{{\mu _g}}}{3}$.
\end{lemma}
The proof of Lemma \ref{LE-Vstar} is provided in  Section \ref{sec-proof-LE-Vstar}.  {\hfill $\blacksquare$}
\begin{lemma}[\textbf{Inner-level errors}] \label{LE-ThetaStar}
Consider the sequence $\{x_i^k, \theta_i^k, v_i^k, z_i^k, y_i^k\}$ generated by Algorithm \ref{alg:1}. Suppose Assumptions \ref{ASS-network}, \ref{ASS-OUTLEVEL}, \ref{ASS-INNERLEVEL} and \ref{ASS-STOCHASTIC} hold. If the step-size $\beta$ satisfies
\begin{equation}\label{EQ-beta}
\begin{aligned}
{\beta} <  \min \left\{  \frac{2}{{{\mu _g} + {L_{g,\theta }}}}, \frac{{{\mu _g} + {L_{g,\theta }}}}{{2{\mu _g}{L_{g,\theta }}}}  \right \},
\end{aligned}
\end{equation}
then we have:
\begin{equation}\label{LE-theta-k-k-1}
\begin{aligned}
  \mathbb{E}[ {{{\| {{\theta ^{k + 1}} - {\theta ^*}({{\bar x}^{k + 1}})} \|}^2}} ] \leqslant& ( {1 -  \frac{{{\mu _g}{L_{g,\theta }}}}{{ {{\mu _g} + {L_{g,\theta }}} }}}{\beta} )\mathbb{E}[ {{{\| {{\theta ^k} - {\theta ^*}({{\bar x}^k})} \|}^2}} ] \\
  & + {p_x}\alpha \mathbb{E}[ {{{\| {{x^k} - {1_m} \otimes {{\bar x}^k}} \|}^2}} ] + m{p_s}{\tau ^2\alpha ^2}\mathbb{E}[ {{{\| {{{\bar y}^k}} \|}^2}} ] + m{\sigma _{\theta}^2}{\alpha ^2}, \\
\end{aligned}
\end{equation}
where \yc{$p_{x}  \triangleq  \frac{4L_{g,\theta }^2{\beta }}{\omega _\theta \alpha}$}, ${\sigma _{\theta}^2} \triangleq 2 \res{\sigma _{g,\theta}^2} \frac{\beta ^2}{\alpha^2}$, ${p_{s}} \triangleq \frac{{2L_{{\theta ^*}}^2 }}{{{\omega _\theta }{\beta}}}$ with $\omega_{\theta}  = \frac{{{\mu _g}{L_{g,\theta }}}}{{2\left( {{\mu _g} + {L_{g,\theta }}} \right)}}$.
\end{lemma}
The proof of Lemma \ref{LE-ThetaStar} is provided in  Section \ref{sec-proof-LE-ThetaStar}.  {\hfill $\blacksquare$}

\begin{lemma}[\textbf{Consensus errors}] \label{LE-CE}
Consider the sequence $\{x_i^k, \theta_i^k, v_i^k, z_i^k, y_i^k\}$ generated by Algorithm \ref{alg:1}. Suppose Assumptions \ref{ASS-network}, \ref{ASS-OUTLEVEL}, \ref{ASS-INNERLEVEL} and \ref{ASS-STOCHASTIC} hold.
Then, we have:
\begin{equation}\label{LE-consensus}
\begin{aligned}
&\mathbb{E}[ {{{\| {{x^{k + 1}} - {1_m} \otimes {{\bar x}^{k + 1}}} \|}^2}} ]
\leqslant   (1-\re{\tau\frac{{1 - \rho }}{2}})\mathbb{E}[ {{{\| {{x^k} - {1_m} \otimes {{\bar x}^k}} \|}^2}} ] + \re{\frac{{2\tau{\alpha ^2}}}{{1 - \rho }}}\mathbb{E}[ {{{\| {{y^k} - {1_m} \otimes \bar y^k} \|}^2}} ],
\end{aligned}
\end{equation}
where $0<\tau<1$, $\rho  = {\| {\mathcal{W} - \mathcal{J}}  \|^2} \in \left[ {0,1} \right)$ with $\mathcal{J}\triangleq \frac{1_m 1_m^{\rm T}}{m} \otimes  I_n$.
\end{lemma}
The proof of Lemma \ref{LE-CE} is provided in  Section \ref{sec-proof-CE}.  {\hfill $\blacksquare$}
\\
\\
 \textbf{Step 3 ({Characterizing} the gradient errors)}: The next  step is to upperbound the  gradient error parts ${\mathbb{E}}[ {{{\| {{y^k} - {1_m} \otimes {{\bar y}^k}} \|}^2}} ]$ in consensus errors caused by the data heterogeneity across nodes. Hence, we respectively demonstrate how the gradient errors change in {\ALGa} with $y^{k+1}=z^{k+1}$ and  {\ALGb} with $y^{k+1}={\mathcal{W}}y^{k}+z^{k+1}-z^k$ in Lemma \ref{LE-heterogeity} and Lemma  \ref{LE-TE}. In particular, it is shown that the gradient errors are impacted by the  heterogeneity  for {\ALGa} with  local gradient scheme \eqref{EQ-ALG-h1}, whereas these error terms will decay  as the iteration progresses  for {\ALGb} with  gradient-tracking scheme \eqref{EQ-ALG-h2}. Lemma \ref{LE-stochastic-error-2} further provides the evolution for the variance errors $\mathbb{E}[ {{{\| {\nabla \tilde{\Phi}(\bar{x}^k) - {z^k}} \|}^2}}]$ induced by the gradient errors ${\mathbb{E}}[ {{{\| {{y^k} - {1_m} \otimes {{\bar y}^k}} \|}^2}} ]$.
 \begin{lemma}[\textbf{Gradient errors for {\ALGa}}]  \label{LE-heterogeity}
Consider the sequence $\{x_i^k, \theta_i^k, v_i^k, z_i^k, y_i^k\}$ generated by Algorithm \ref{alg:1} employing the local gradient scheme \eqref{EQ-ALG-h1}. Suppose Assumptions \ref{ASS-network}, \ref{ASS-OUTLEVEL}, \ref{ASS-INNERLEVEL}, \ref{ASS-heterogeneity} and \ref{ASS-STOCHASTIC} hold.  Then,  we have
\begin{equation}\label{EQ-y-heterogeity}
\begin{aligned}
&\mathbb{E}[{\| {{y^k} - {1_m} \otimes {{\bar y}^k}} \|^2}]   \\
   \leqslant&  3{b^2} + 3m\mathbb{E}[{\| {\nabla \Phi ( {{{\bar x}^k}} )} \|^2}] +3\mathbb{E}[ {{{\| {\nabla \tilde \Phi ({{\bar x}^k}) - {z^k}} \|}^2}} ],  \\
\end{aligned}
\end{equation}
where $\nabla \tilde{\Phi} (\bar{x}^k)={\rm col}\{\nabla \Phi _i(\bar{x}^k)\}_{i=1}^{m}$, and ${b^2}$ is the heterogeneity on overall hypergradients  denoted in Lemma \ref{LE-hyper-heterogeneity}.
\end{lemma}
The proof of Lemma \ref{LE-heterogeity} is provided in  Section \ref{sec-proof-LE-heterogeity}.  {\hfill $\blacksquare$}
\begin{lemma}[\textbf{Gradient errors for {\ALGb}}] \label{LE-TE}
Consider the sequence $\{x_i^k, \theta_i^k, v_i^k, z_i^k, y_i^k\}$  generated by Algorithm \ref{alg:1} employing
the gradient tracking scheme \eqref{EQ-ALG-h2}. Suppose Assumptions \ref{ASS-network}, \ref{ASS-OUTLEVEL}, \ref{ASS-INNERLEVEL} and \ref{ASS-STOCHASTIC} hold.  Then,  we have:
\begin{equation}
\begin{aligned}
 & {\mathbb{E}} [{\| {{y^{k + 1}} - {1_m} \otimes {{\bar y}^{k + 1}}} \|}^2] \\
   \leqslant & \frac{{1 + \rho }}{2}{\mathbb{E}}[{\| {{y^k} - {1_m} \otimes {{\bar y}^k}} \|^2}] + \nyc{\frac{4}{{1 - \rho }}\gamma ^2{\mathbb{E}}[{\| {\nabla \tilde \Phi ({{\bar x}^k}) - {\mathbb{E}}[{s^k}]} \|^2}]} \\
   &+ \frac{4}{{1 - \rho }}{{\gamma^2}}\nyc{{\mathbb{E}}[{{\| \nabla \tilde \Phi ({{\bar x}^k}) - {z^k}\|}^2}]}  + \frac{2}{{1 - \rho }}m\sigma_y^2{\alpha ^2},
\end{aligned}
\end{equation}
where $\sigma _y^2 \triangleq ({\sigma _{f,x}^2}+\nyc{\hat M^2}{\sigma _{g,x\theta}^2})\frac{\gamma ^2}{\alpha^2}$.
\end{lemma}
The proof of Lemma \ref{LE-TE} is provided in  Section \ref{sec-proof-LE-TE}.  {\hfill $\blacksquare$}

\begin{lemma}[\textbf{Variance errors}] \label{LE-stochastic-error-2}
Consider the sequence $\{x_i^k, \theta_i^k, v_i^k, z_i^k, y_i^k\}$ generated by Algorithm \ref{alg:1}. Suppose Assumptions  \ref{ASS-network}, \ref{ASS-OUTLEVEL}, \ref{ASS-INNERLEVEL} and \ref{ASS-STOCHASTIC} hold. Recall that $\nabla \tilde{\Phi} (\bar{x}^k)={\rm col}\{\nabla \Phi _i(\bar{x}^k)\}_{i=1}^{m}$. Then, we have
\nyc{
\begin{equation}\label{EQ-stochastic-error-2}
\begin{aligned}
\mathbb{E}[ {{{\| {\nabla \tilde{\Phi} ({{\bar x}^{k + 1}}) - {{ z}^{k + 1}}} \|}^2}} ] \leqslant& (1 - \gamma )\mathbb{E}[ {{{\| {\nabla \tilde{\Phi} ({{ \bar x}^k}) - {{ z}^k}} \|}^2}} ] + {r_z}\alpha \mathbb{E}[ {{{\| {\nabla \tilde{\Phi} ({{\bar x}^k}) - \mathbb{E}[{{{s}^k}}]} \|}^2}} ] \\
&+ m{r_y}{\tau ^2}{\alpha ^2}\mathbb{E}[ {{{\| {{{\bar y}^k}} \|}^2}} ] + m\sigma _{ z}^2{\alpha ^2},
\end{aligned}
\end{equation}
}
and
\nyc{
\begin{equation}\label{EQ-hypergradient-error-2}
\begin{aligned}
\mathbb{E}[ {{{\| {\nabla \tilde{\Phi} ({{\bar x}^k}) - \mathbb{E}[{{{s}^k}}]} \|}^2}} ] \leqslant& {L_{fg,x}}\mathbb{E}[ {{{\| {{x^k} - {1_m} \otimes {{\bar x}^k}} \|}^2}} ] + {L_{fg,x}}\mathbb{E}[ {{{\| {{\theta ^k} - {\theta ^*}({{\bar x}^k})} \|}^2}} ] \\
&+ {{4C_{g,x\theta }^2}}\mathbb{E}[ {{{\| {{v^k} - {v^*}({{\bar x}^k})} \|}^2}} ].
 \end{aligned}
\end{equation}
}
 where $\sigma _{z}^2 \triangleq ({\sigma _{f,x}^2}+\nyc{\hat M^2}{\sigma _{g,x\theta}^2}){\frac{\gamma^2}{\alpha^2}}$.
\end{lemma}
The proof of Lemma \ref{LE-stochastic-error-2} is similar to  Lemma \ref{LE-HypergradientE} and we omit it here.  {\hfill $\blacksquare$}
\\
\\
 \textbf{Step 4 (Integrating Steps 1, 2, 3 to obtain the overall dynamics with the unified  Lyapunov function)}: Finally, by integrating the obtained results  and employing small-gain-like techniques, we establish the dynamics of $V^k$ in Sections \ref{sec-proof-TH-1} and Section \ref{sec-proof-TH-2} for {\ALGa} and {\ALGb}, respectively, using  a set of carefully chosen  coefficients $d_0$, $d_1$, $d_2$, $d_3$, $d_4$, $d_5$, $d_6$.

\subsection{Further Discussions on Convergence Analysis}\label{SE-diff}

\textbf{Unique heterogeneity analysis and weaker  assumption on heterogeneity.} We adopt unique heterogeneity analysis approaches to investigate the impact of heterogeneity. Unlike single-level optimization and global DSBO, personalized DSBO involves distinct inner-level solutions across nodes, making the overall heterogeneity depend intricately on both outer-level and inner-level heterogeneity.  This results in a fundamentally different and more challenging scenario for heterogeneity analysis. Moreover, we adopt a weaker Assumption \ref{ASS-heterogeneity} on heterogeneity than that of works on DSBO, where our assumption  requires only the boundedness of inner-level heterogeneity at the optimum $\theta_i^*(x)$ (instead of being at any point), while avoiding the boundedness of local hypergradients (c.f., {Remark \ref{re-assumption}}). This weaker assumption on heterogeneity, together with the different inner-level solutions, {presents a unique challenge in explicitly characterizing the heterogeneity} in personalized DSBO, as it requires new techniques to analyze the heterogeneity of the inner-level solutions. To tackle these challenges,  based on the structure of the hypergradient (c.f., Assumption \ref{ASS-heterogeneity} (ii) and (iii)), we introduce the heterogeneity of the Hessian and Jacobian matrices to characterize the inner-level heterogeneity in DSBO. Furthermore, by establishing the intricate relationships between inner-level and outer-level heterogeneity, we  explicitly characterize the overall heterogeneity (c.f., {Lemma \ref{LE-hyper-heterogeneity}}) that clearly shows the impact of inner- and outer-level heterogeneity and provides new insights into heterogeneity in personalized DSBO (cf. Remark~\ref{RE-het}).
 \\[1ex]
\textbf{Different  hypergradient estimation   and unified analytical framework.}  Our analysis for hypergradient approximation errors  differs from the existing literature \citep{MA-DSBO, SLAM, VRDBO, SPDB, yang2022decentralized} due to the additional errors from loopless approximation and network dynamic (c.f., Lemma \ref{LE-HypergradientE}). It poses an  additional challenge for technique analysis in controlling  the variances and analyzing heterogeneity.  In contrast, previous works \citep{BSA, VRBO, SPDB} utilize NS and HSIA methods \eqref{EQ-NSCG} along with extra  loops with increasing number of iterations $Q$ to obtain highly accurate hypergradient approximations $\tilde \nabla \Phi ({{\bar x}^k})$ (corresponding to the term $z^k$ in \eqref{EQ-descentlemma}), i.e., $\mathbb{E}[ {{{\| {\nabla \Phi ({{\bar x}^k}) - \tilde \nabla \Phi ({{\bar x}^k}) } \|^2}}} ] \leqslant \frac{{{C_{g,x\theta }}{C_{f,\theta }}}}{{{\mu _g}}}{(1 - \frac{{{\mu _g}}}{{{L_{g,\theta }}}})^Q}$.
\nyc{Besides, we adopt a different technique  to asymptotically control the bound of  $\|v_i^k\|$ without requiring  the boundedness of the stochastic gradients or the boundedness of $\|v_i^k\|$ \citep{VRDBO,yang2022decentralized} (c.f., Lemmas \ref{LE-HypergradientE}, \ref{LE-Vstar}). On the other hand, we also establish a unified analytical framework for both LoPA-LG and LopA-GT with tighter convergence rates. In particular, we design a unified Lyapunov function and employ small-gain-like techniques to establish a overall error evolution dynamic for both LoPA-LG and LoPA-GT. Besides, leveraging unique heterogeneity analysis technique, we are able to bound the gradient tracking errors via  variance-related errors (c.f., Lemma \ref{LE-stochastic-error-2}) for LoPA-LG and LoPA-GT in a unified manner (c.f. Lemma \ref{LE-hyper-heterogeneity}, \ref{LE-heterogeity}-\ref{LE-stochastic-error-2}), while
characterizing the effect of the heterogeneity. This above analysis, indeed, provides an unified analytical framework that enables us to  establish intricate relationships among the errors as well as their detailed effects on the convergence performance  (c.f., {Theorems \ref{TH-1}, \ref{TH-2}}),
leading to  tighter convergence rates (c.f., Corollaries  \ref{CO-1}, \ref{CO-2}) and  improved complexity (c.f., Remark  \ref{re-improved complexity}).
The above analysis and obtained results distinguish our work significantly from existing works on DSBO.}
\\[1ex]
\textbf{The necessity of   gradient momentum steps.} The gradient momentum step \eqref{EQ-ALG-g} is introduced to
control  stochastic variances in hypergradient estimation  by leveraging the moving average of historical hypergradient  information and  selecting proper step-size $\gamma$.  This  step allows us to  bound only the term $\mathbb{E}[ \| \nabla \Phi (\bar{x}^k) - \mathbb{E}[\bar{s}^k] \|^2 ]$ instead of $\mathbb{E}[ \| \nabla \Phi (\bar{x}^k) - \bar{s}^k \|^2 ]$ in Lemma \ref{LE-HypergradientE}, thereby ensuring the convergence of {\ALG}.
In particular, let us consider the case where we exclude the recursion \eqref{EQ-ALG-g}, i.e., $\gamma=1$ and  $z^k = s^k$.  In this scenario,  the average variance error term $\tau \alpha\mathbb{E}[ {{{\| \nabla{\Phi ({{\bar x}^k}) - {{\bar z}^k}} \|}^2}} ]$
reduces to  the term $\tau \alpha\mathbb{E}[ {{{\| \nabla{\Phi ({{\bar x}^k}) - {{\bar s}^k}} \|}^2}} ]$,  which can only be bounded by a sampling variance of order ${\mathcal{O}}(\alpha)$.   As a result,  this case   either leads to non-decaying variance errors or requires the proper choice of two-timescale step-sizes to control the decay of the term  $\mathbb{E}[{| {{{\bar y}^k}} |^2}]$ in \eqref{EQ-v*-v} and  \eqref{LE-theta-k-k-1},
where the later will lead to a sub-optimal rate of ${\mathcal{O}}(K^{-2/5})$ \citep{TTSA}. Instead,  with the  recursion \eqref{EQ-ALG-g},  when the  step-size $\gamma$ is properly taken as $\gamma=\mathcal{O}({\alpha})$, the  term  $\mathbb{E}[ \| \nabla \Phi (\bar{x}^k) - \bar{z}^k \|^2 ]$ can be asymptotically controlled by  the term $\mathbb{E}[ \| \nabla \Phi (\bar{x}^k) - \mathbb{E}[\bar{s}^k] \|^2 ]$ and the sampling variance of an order of ${\mathcal{O}}(\alpha^2)$  (c.f., the last term in Lemma \ref{LE-HypergradientE}).   \nyc{  Similar  techniques have been employed in  \citep{MA-DSBO, ghadimi2020single, chen2024optimal}. There are other kinds of gradient momentum methods based on double gradient evaluations \citep{STORM,SPIDER}  or  unbiased averaging stochastic gradients \citep{SAGA} that can help achieve a faster rate under some extra conditions \citep{VRDBO,SOBA}}.
\section{Numerical Experiments}
In this section, we  present two numerical experiments to test the performance of the proposed {\ALG} algorithm  and verify the theoretical findings. \nycres{The first experiment considers classification problems for a 10-class classification task, while the second focuses on hyperparameter optimization in $l_2$-regularized binary logistic regression problems for a two-class classification task.}

\subsection{Distributed  Classification} In this subsection, we evaluate the effectiveness of our algorithms on a class of distributed classification  problems involving heterogeneous datasets. We employ MNIST datasets to train $m$ personalized  classifiers \nycres{for a 10-class classification task.} Specifically, we construct a classifier in node $i$ that consists of a hidden layer  with parameters $x$ shared across all nodes followed by sigmoid activation functions and a linear layer with parameters $\theta_i$ adapted to node-specific samples. The cross-entropy loss is used  as the outer-level objective  $f_i$.  Regarding the inner-level objective function $g_i$, we include a quadratic regularization term to the parameters $\theta_i$ based on the cross-entropy loss to avoid overfitting to local samples. Specifically, $f_i$ and $g_i$
 take the following form:
 \begin{equation}
\begin{aligned}
&\underset{x\in \mathbb{R} ^n}{\min}\frac{1}{m}\sum_{i=1}^m{f_i(x,\theta_i^*(x))=\frac{1}{m}\sum_{i=1}^m{\sum \nolimits_{(s_{ij},b_{ij})\in \mathcal{D} _i}{b_{ij}\ln\mathrm{(}y_j(x,\theta _{i}^{*}\left( x \right) ;s_{ij}))}}},\;
\\
&\;\mathrm{s}.\mathrm{t}. \; \theta _{i}^{*}\left( x \right) =arg\underset{\theta _i\in \mathbb{R} ^p}{\min}\,\,g_i(x,\theta _i)=\sum \nolimits_{(s_{ij},b_{ij})\in \mathcal{D} _i}{b_{ij}\ln\mathrm{(}y_j(x,\theta _i;s_{ij}))}+\frac{\mu}{2}\left\| \theta _i \right\| ^2,
\end{aligned}
\end{equation}
where $({s_{ij}},{b_{ij}})\in {\mathcal{D}_{i}} $ denotes the $j$-th sample assigned to the node $i$ with $s_{ij}$ being the $j$-th feature and $b_{ij}$ being $j$-th label in one-shot form; $\mu$ denotes the penalty coefficient; \nyc{${y_j}(x,\theta ;{s_{ij}})$ represents the output of the neural network parameterized  by the layer weight $x$ and $\theta$.}

\begin{figure}[ht]
\centering
\subfigure[Loss.]{
\begin{minipage}[ht]{0.30\linewidth}
\centering
\includegraphics[width=0.95\textwidth,height=0.18\textheight]{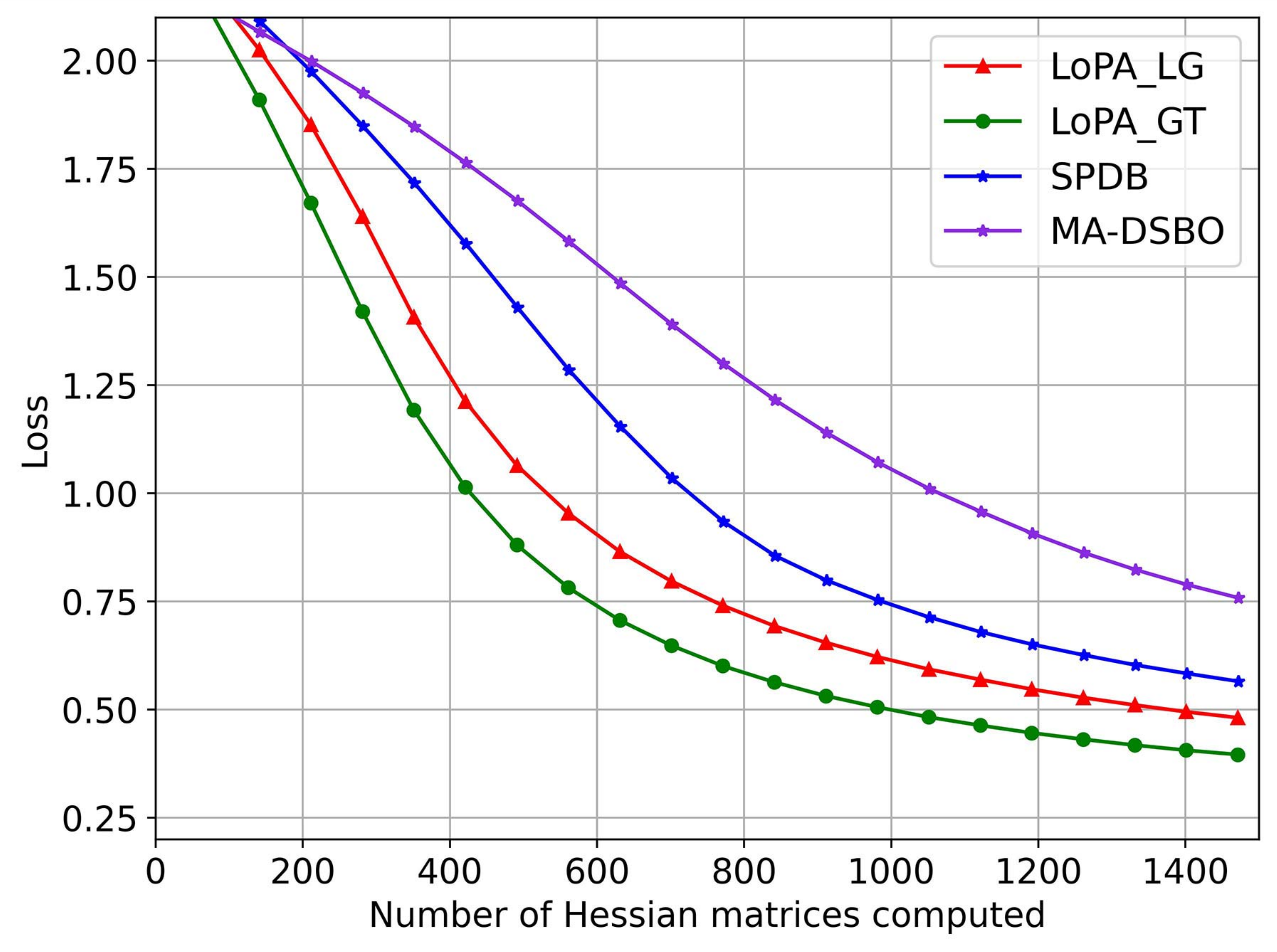} 
\end{minipage}%
}%
\subfigure[Training accuracy.]{
\begin{minipage}[ht]{0.30\linewidth}
\centering
\includegraphics[width=0.95\textwidth,height=0.18\textheight]{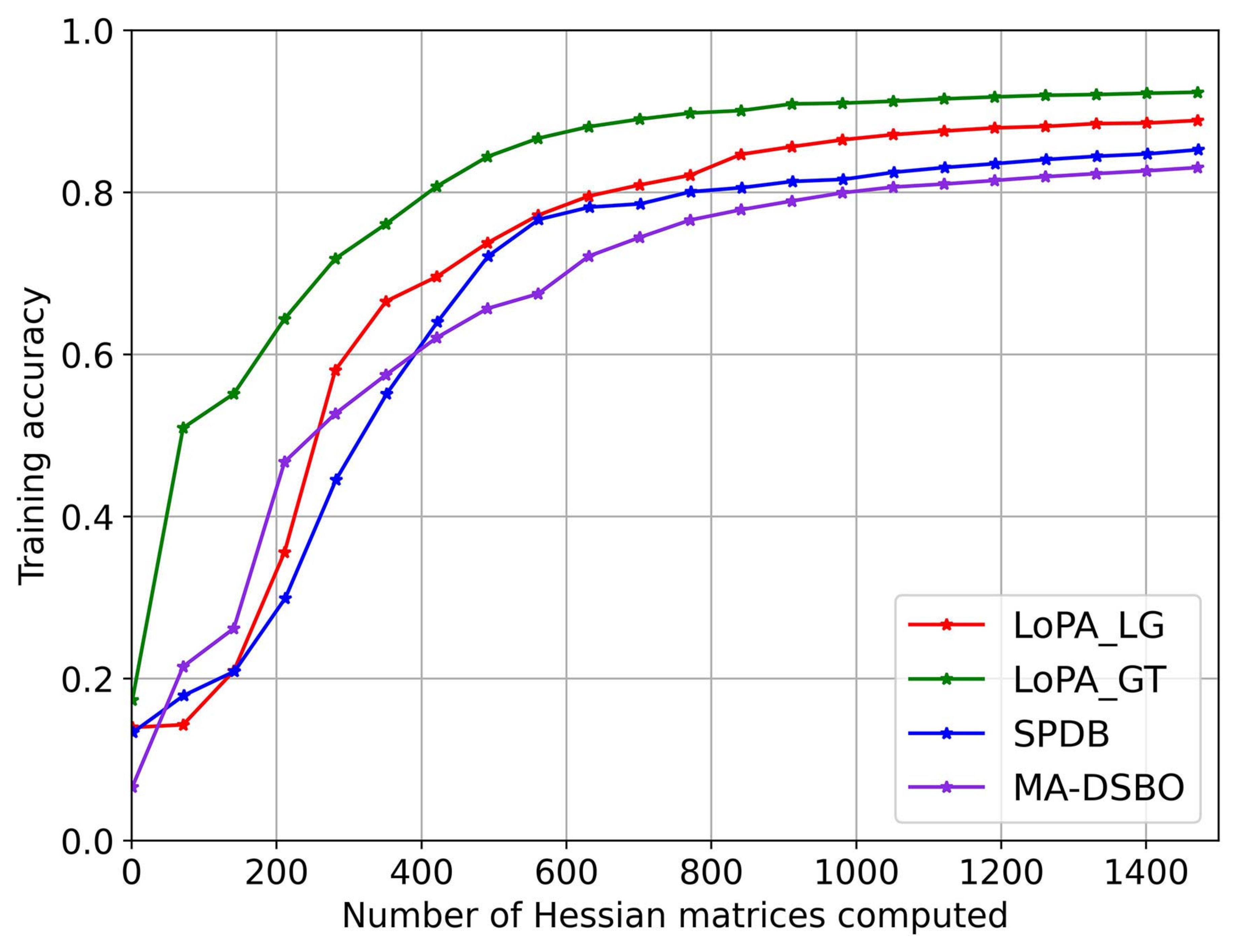}
\end{minipage}%
}%
\subfigure[Testing accuracy.]{
\begin{minipage}[ht]{0.30\linewidth}
\centering
\includegraphics[width=0.95\textwidth,height=0.18\textheight]{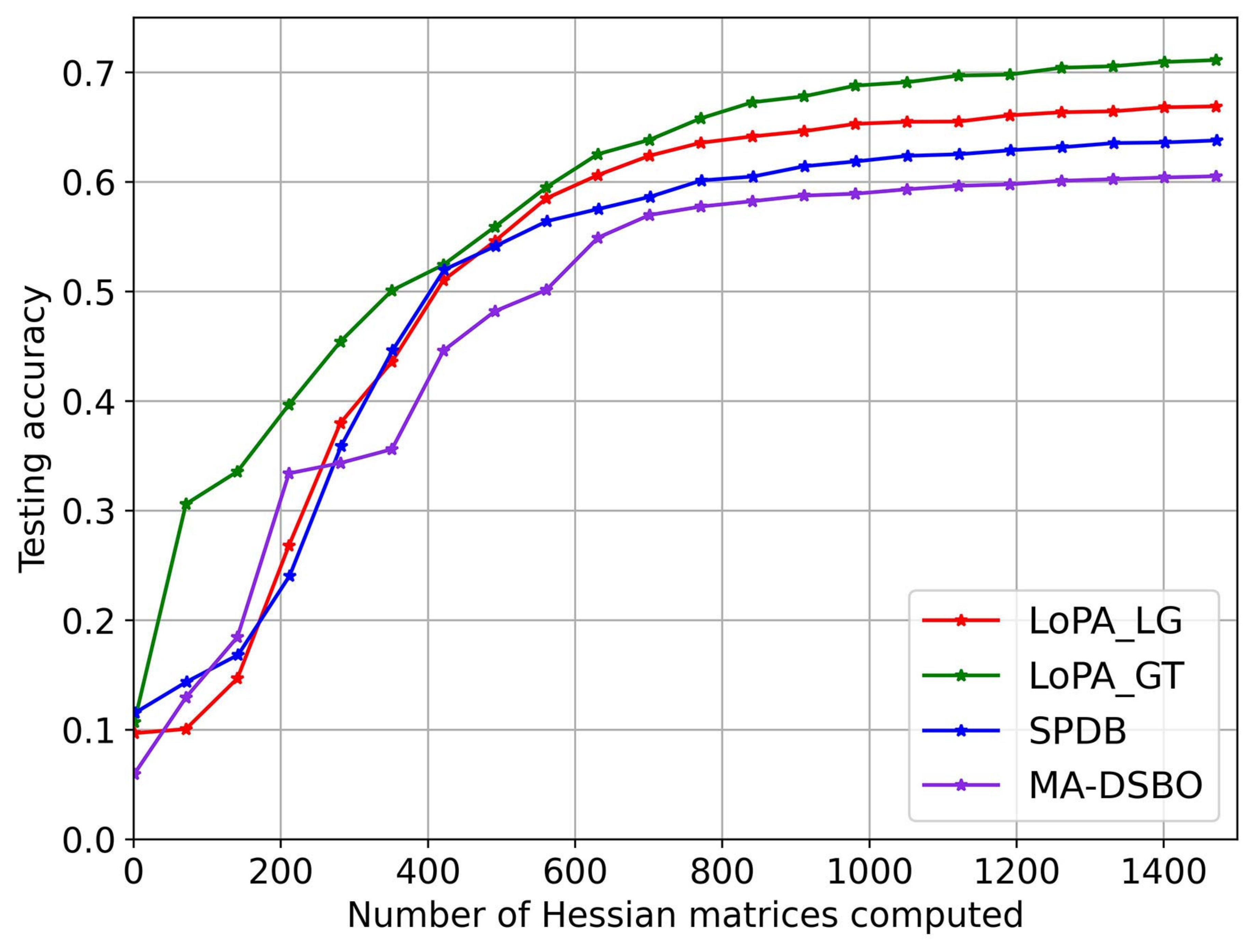}
\end{minipage}%
}%
\centering
\caption{\nycres{Performance comparison of SPDB, MA-DSBO and our LoPA-LG and LoPA-GT algorithms over $4$ nodes for a 10-class
classification task using MNIST dataset.}}
\label{fig:node4}
\end{figure}

\begin{figure}[ht]
\centering
\subfigure[Loss.]{
\begin{minipage}[ht]{0.30\linewidth}
\centering
\includegraphics[width=0.95\textwidth,height=0.18\textheight]{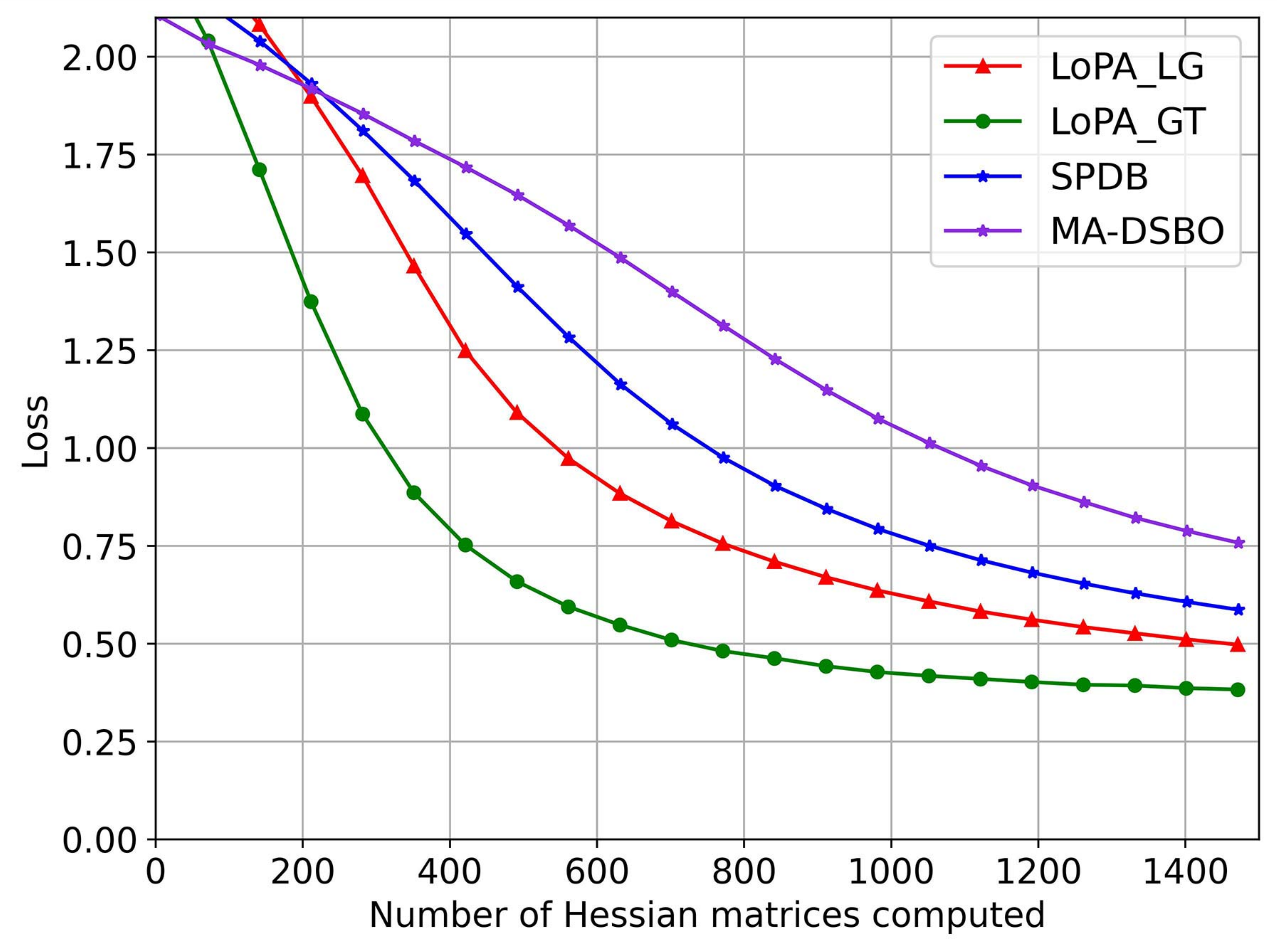} 
\end{minipage}%
}%
\subfigure[Training accuracy.]{
\begin{minipage}[ht]{0.30\linewidth}
\centering
\includegraphics[width=0.95\textwidth,height=0.18\textheight]{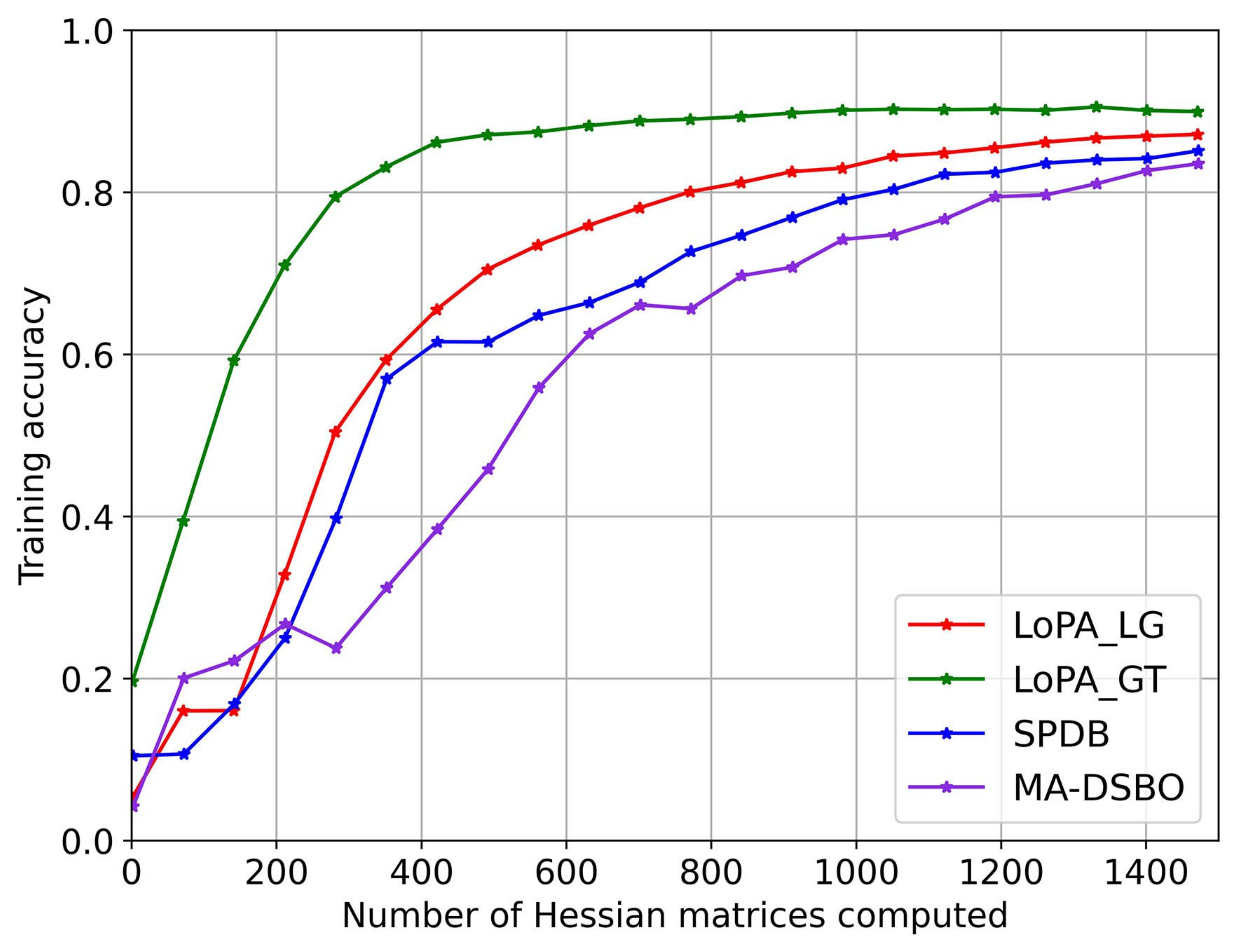}
\end{minipage}%
}%
\subfigure[Testing accuracy.]{
\begin{minipage}[ht]{0.30\linewidth}
\centering
\includegraphics[width=0.95\textwidth,height=0.18\textheight]{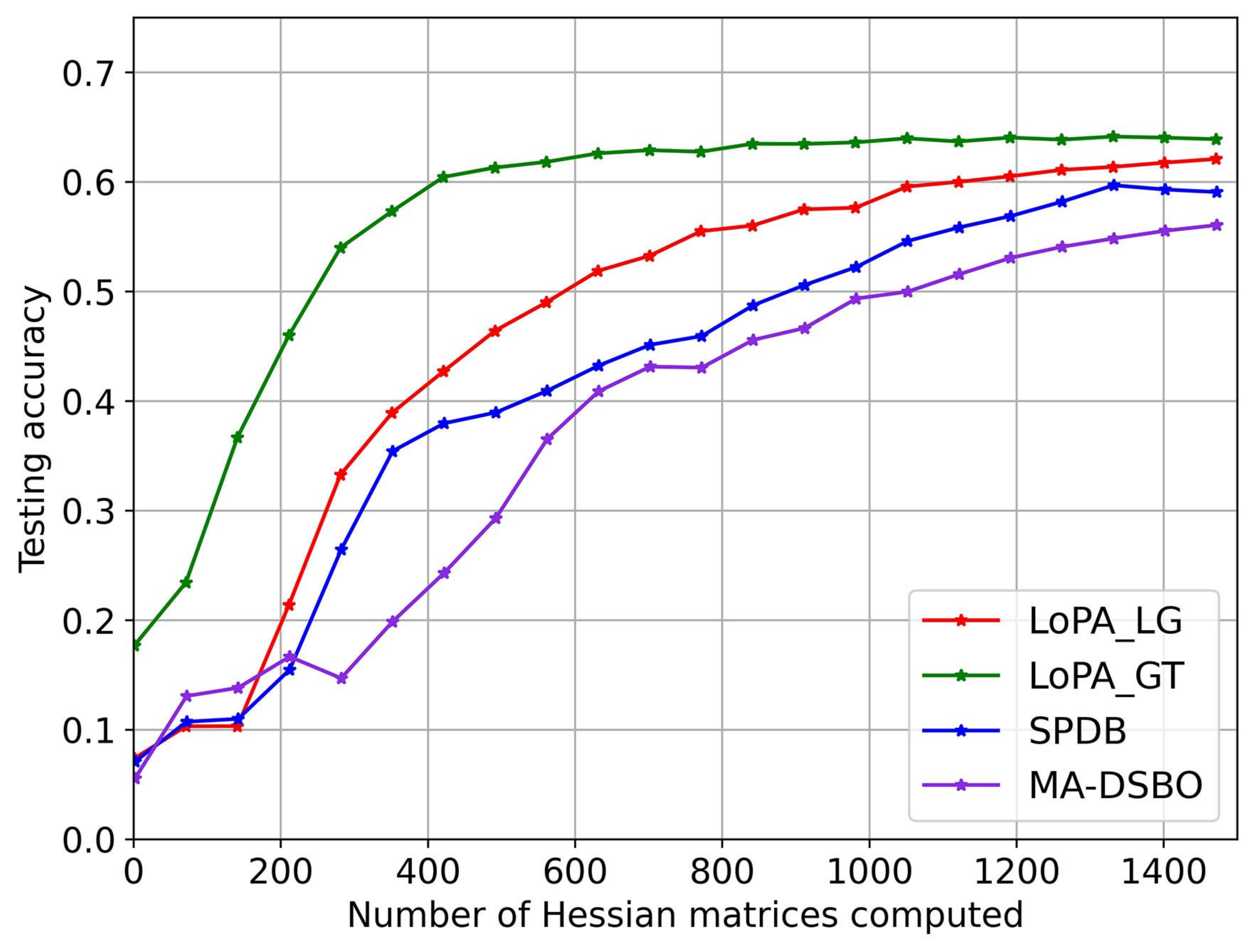}
\end{minipage}%
}%
\centering
\caption{\nycres{Performance comparison of SPDB, MA-DSBO and our LoPA-LG and LoPA-GT algorithms over $8$ nodes for a 10-class
classification task using MNIST dataset.}}
\label{fig:node8}
\end{figure}

\nycres{The experiments on the 10-class classification task are conducted under two distinct scenarios:} i) $m=4$ with each node having 14000 samples; and ii) $m=8$ with each node 6500 samples, where communication networks are generated by  random {Erd\H{o}s–R\'enyi} graphs, and $\mu=0.52$. Each node $i$ is assigned with a random subset  of overall 10 classes such that  each node has different label distributions and high data heterogeneity. \nycres{For the classifier in each node, the input is a 784-dimensional vector. The first layer of the classifier contains 28 neurons, while the second layer contains 10 neurons.} The step-sizes  are set as $\alpha=0.01$, $\beta=0.01$, $\lambda=0.008$, $\gamma=0.4$, $\tau=0.4$ both for   {\ALGa}  and   {\ALGb}.  \nycres{We compare our algorithms with the state-of-the-art SPDB  \citep{SPDB} and MA-DSBO \citep{MA-DSBO} algorithms}, where the mini-batch sizes are set to 50 for all algorithms. These two algorithms respectively utilize the NS and SHIA methods for estimating Hessian inverse matrices. To adapt the MA-DSBO algorithm for personalized settings, we modify the updates of inner-level variables and Hessian-inverse-vector variables in MA-DSBO so that they exclusively perform local gradient updates.

The experiment results for loss, training accuracy, and testing accuracy  are presented in Figures \ref{fig:node4} and  \ref{fig:node8}. It can be observed from Figure \ref{fig:node4}  that the proposed {\ALGa} and {\ALGb} algorithms have lower computational complexity than SPDB and MA-DSBO in terms of the number of Hessian matrices to achieve the same desired accuracy. When the number of nodes increases to 8, we can observe similar results as shown in Figure \ref{fig:node4}. This further demonstrates the scalability  of the proposed algorithms.
Moreover, the difference between {\ALGa} and {\ALGb} highlights the advantages of gradient tracking.

 \begin{wrapfigure}[15]{r}{0.33\textwidth}
\begin{minipage}{0.33\textwidth}
\centering
\includegraphics[width=1\textwidth,height=0.18\textheight]{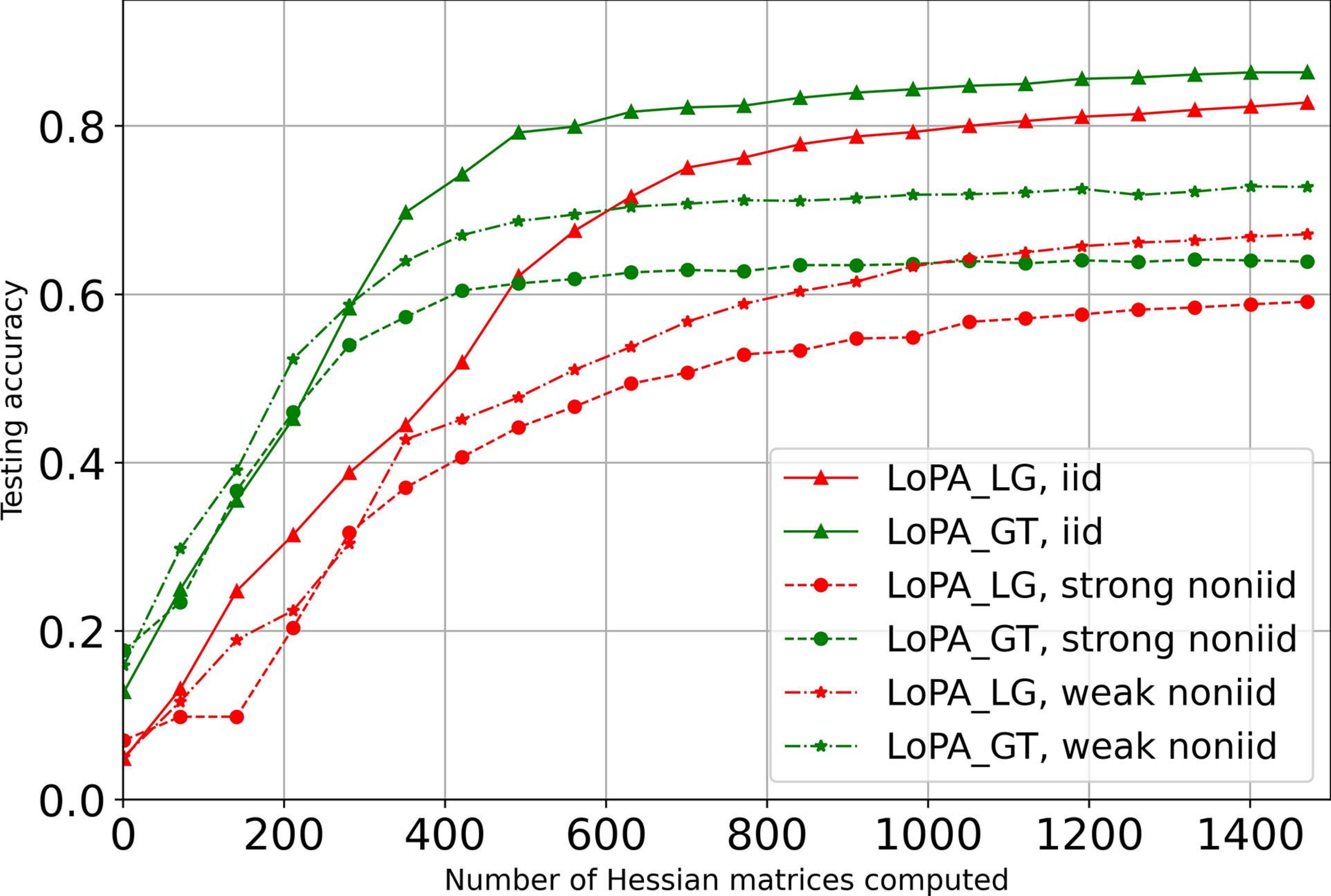} 
\vspace{-0.4cm}
\caption{\nycres{Testing accuracy of {\ALGa} and {\ALGb} under different data heterogeneity for a 10-class
classification task using MNIST dataset.}}
\label{fig:heter}
\end{minipage}
\end{wrapfigure}
\begin{figure}[!htpb]
\centering
\subfigure[IID dataset.]{
\begin{minipage}[ht]{0.30\linewidth}
\centering
\includegraphics[width=1\textwidth,height=0.15\textheight]{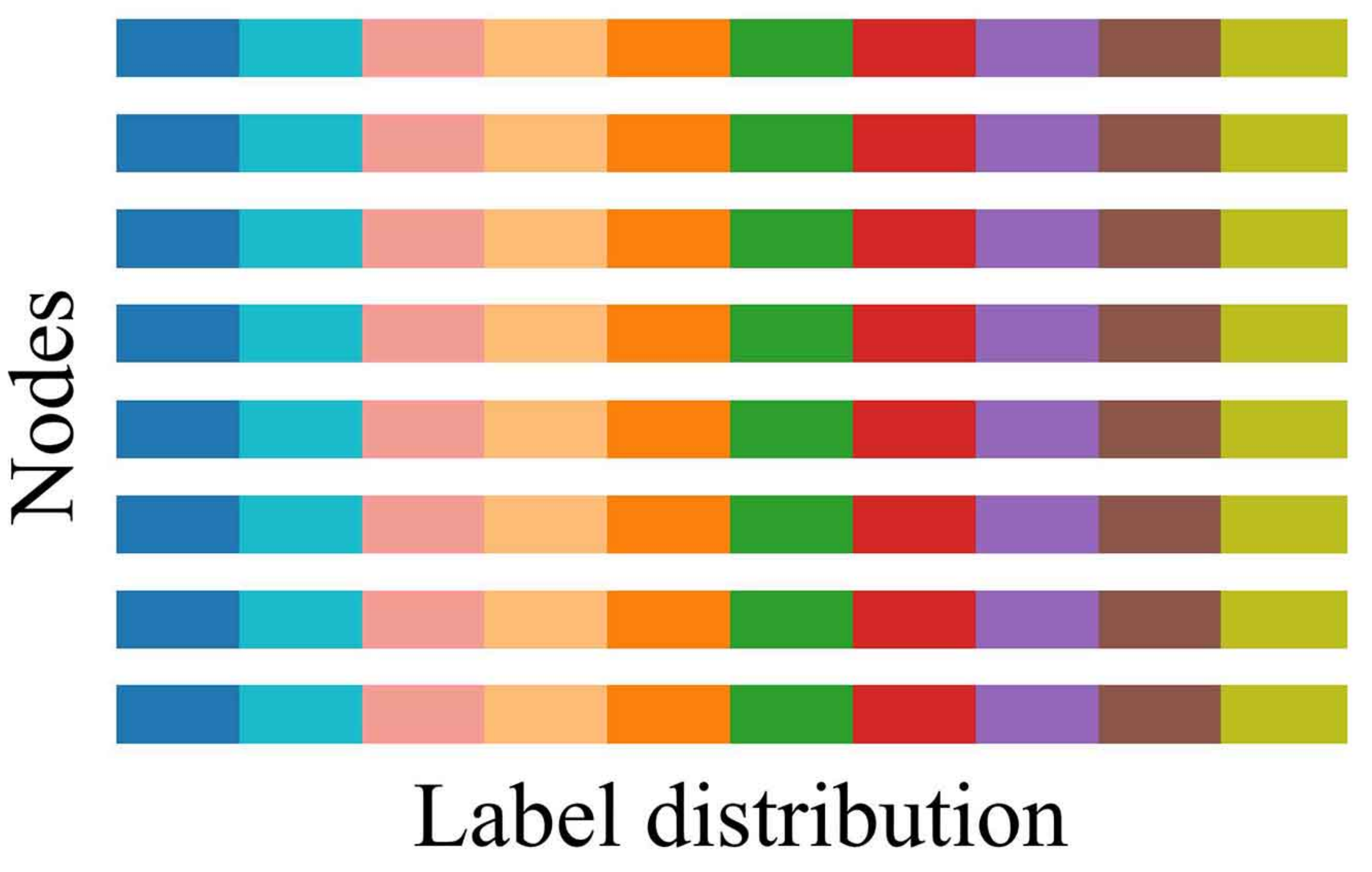} 
\end{minipage}%
}%
\subfigure[Weak non-IID dataset.]{
\begin{minipage}[ht]{0.30\linewidth}
\centering
\includegraphics[width=1\textwidth,height=0.15\textheight]{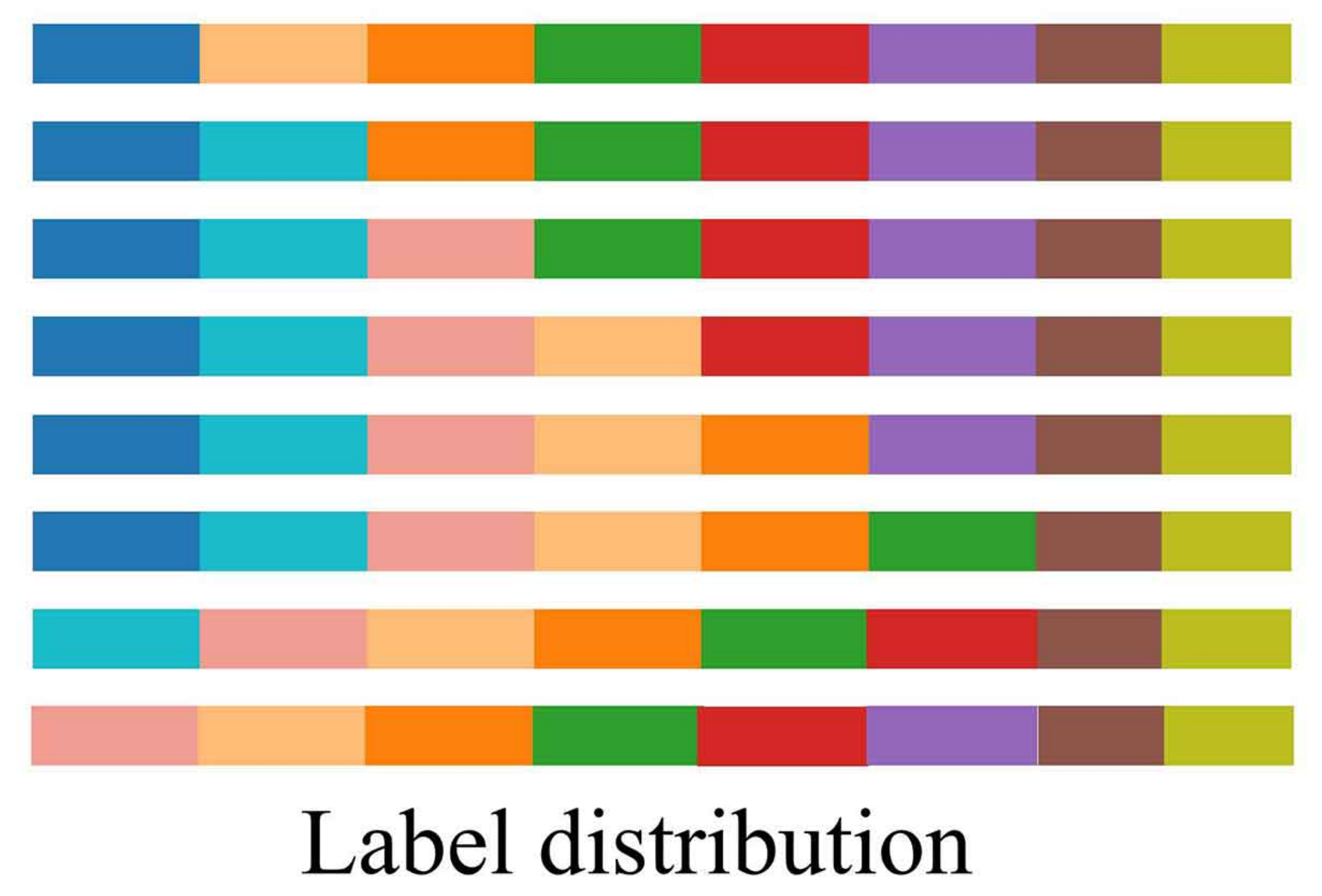}
\end{minipage}
}%
\subfigure[Strong non-IID dataset.]{
\begin{minipage}[ht]{0.30\linewidth}
\centering
\includegraphics[width=1\textwidth,height=0.15\textheight]{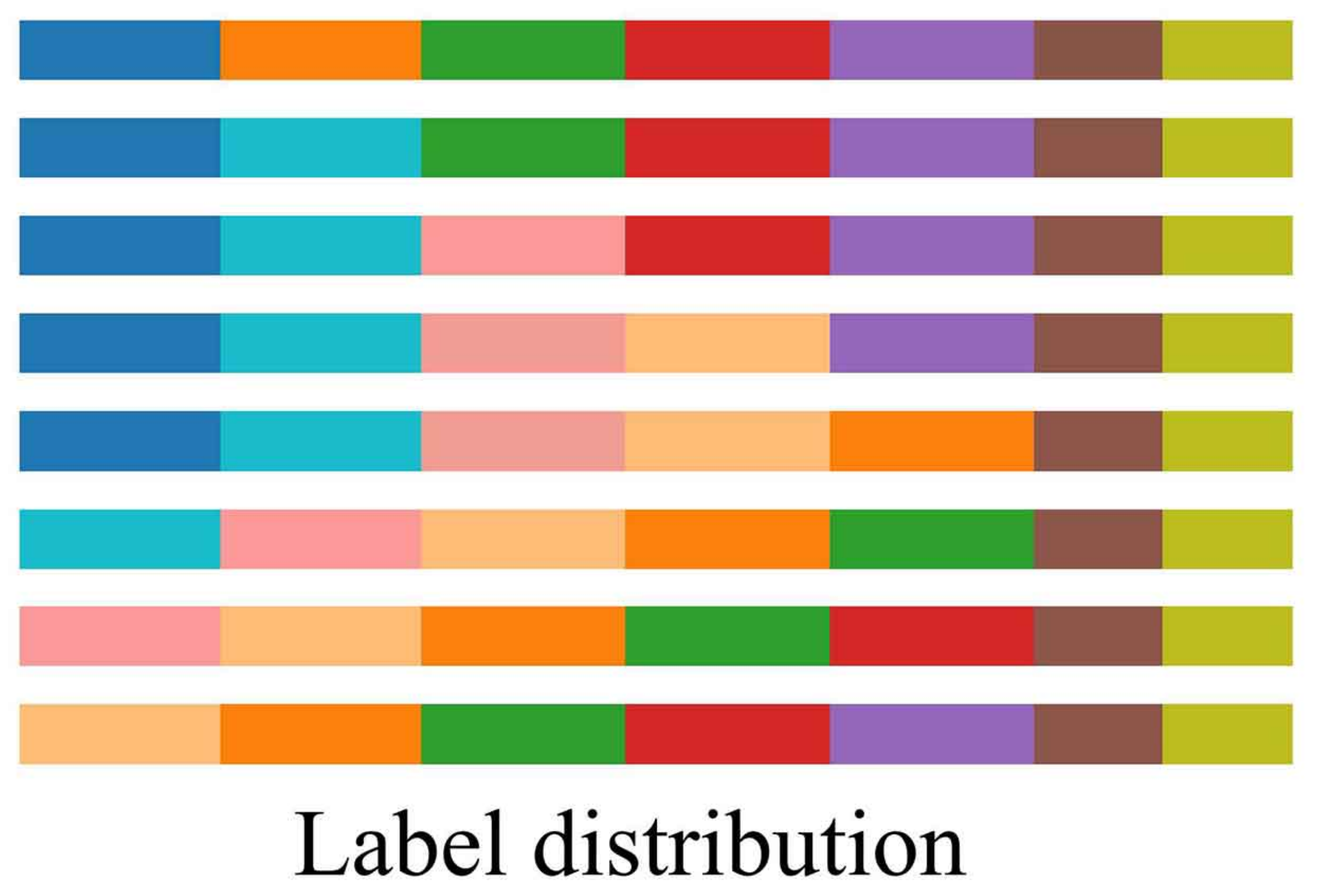}
\end{minipage}%
}%
\centering
\caption{Synthetic label distributions with different levels of data heterogeneity across nodes. The label classes are represented with different colors.}

\label{fig:distribution}
\end{figure}

\noindent\textbf{Impact of heterogeneity. } We conduct an addtional experiment to evaluate the performance of our proposed algorithms  under different settings of heterogeneous label distributions. Specifically, we consider a network of 8 nodes with each holding  7500 samples, and  generate the following three different label distributions among nodes:  i) independent and identically distributed (IID) datasets; ii) non-IID datasets with strong heterogeneity; iii) non-IID dataset with weak heterogeneity. The considered  three label distributions  are depicted  in  Figure \ref{fig:distribution}.


We plot the testing accuracy of {\ALGa} and {\ALGb}  with the same parameter settings as the previous experiments in Figure \ref{fig:heter}. The experiment results, as shown in Figure \ref{fig:heter}, demonstrate that {\ALGb} can maintain a relatively higher accuracy compared to {\ALGa} as the level of data heterogeneity increases. This suggests that {\ALGb} is more robust against data heterogeneity, verifying the theoretical results and the effectiveness of {\ALGb} in scenarios with heterogeneous datasets.

\subsection{Hyperparameter Optimization}
\begin{figure}[ht]
\centering
\subfigure{
\begin{minipage}[ht]{0.30\linewidth}
\centering
\includegraphics[width=0.98\textwidth,height=0.18\textheight]
{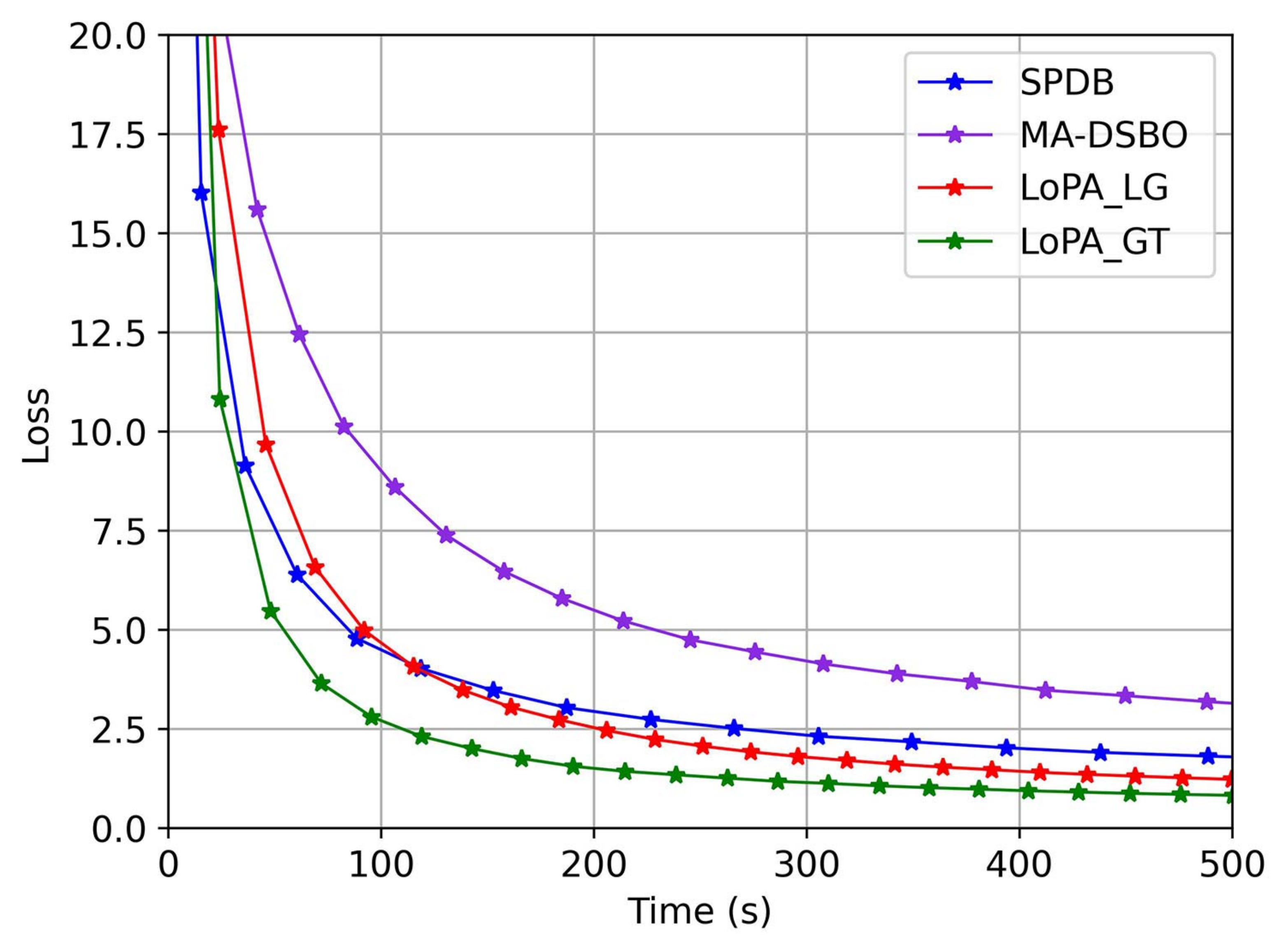} 
\end{minipage}%
}%
\subfigure{
\begin{minipage}[ht]{0.30\linewidth}
\centering
\includegraphics[width=0.98\textwidth,height=0.18\textheight]{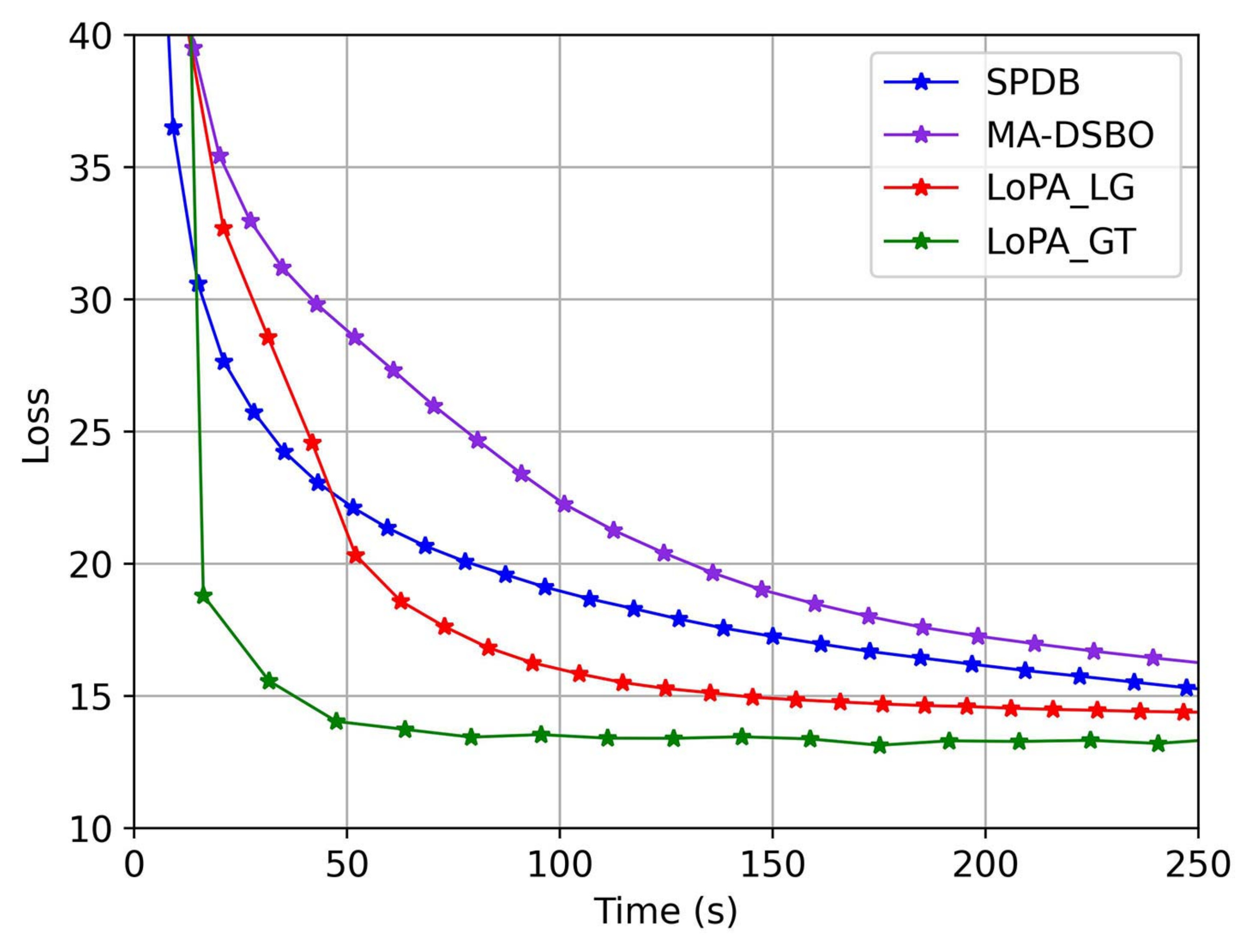} 
\end{minipage}%
}%
\subfigure{
\begin{minipage}[ht]{0.30\linewidth}
\centering
\includegraphics[width=0.98\textwidth,height=0.18\textheight]
{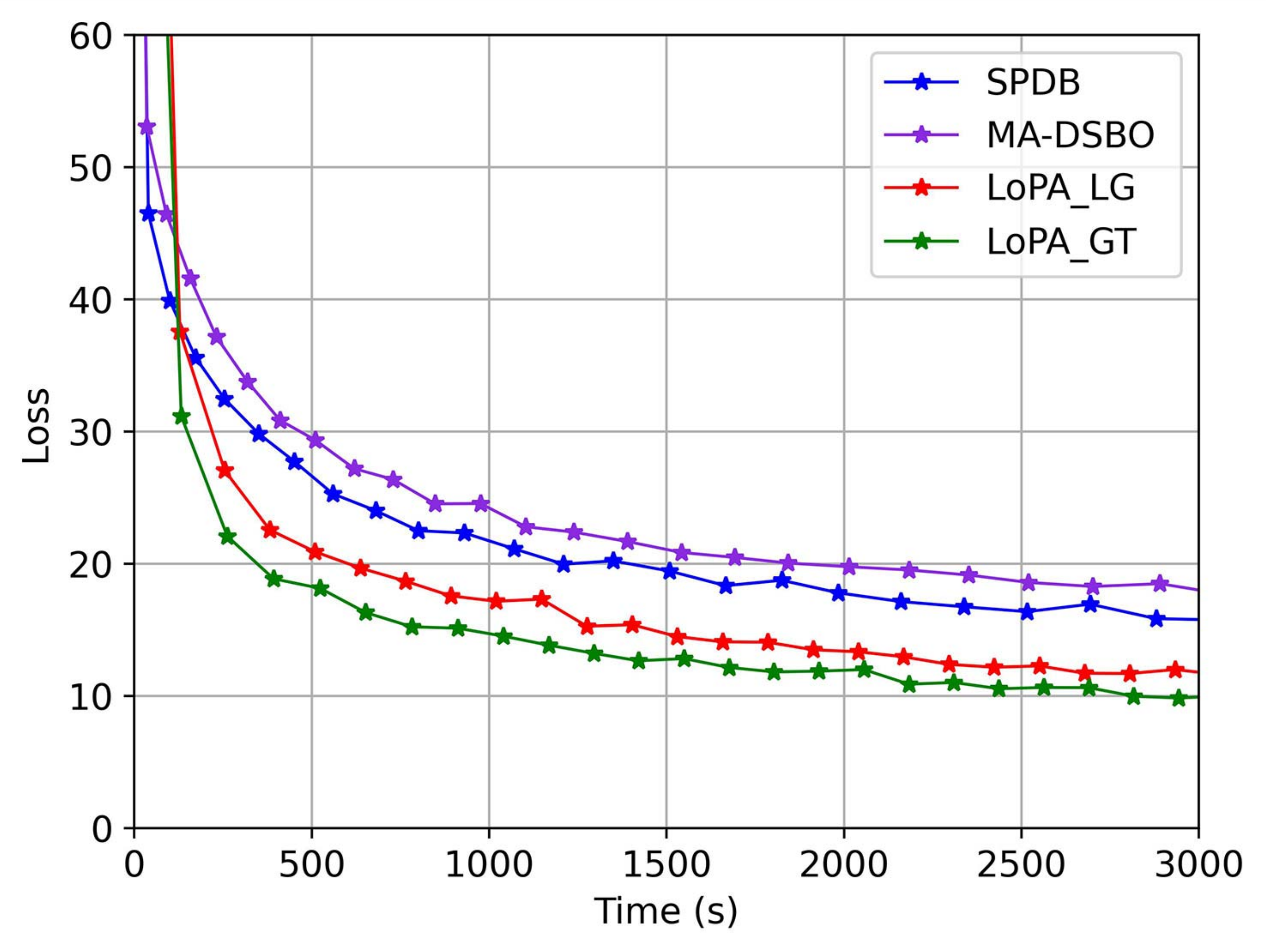} 
\end{minipage}%
}%
\vspace{-0.4cm}
\\
\subfigure{
\begin{minipage}[ht]{0.30\linewidth}
\centering
\includegraphics[width=0.98\textwidth,height=0.18\textheight]
{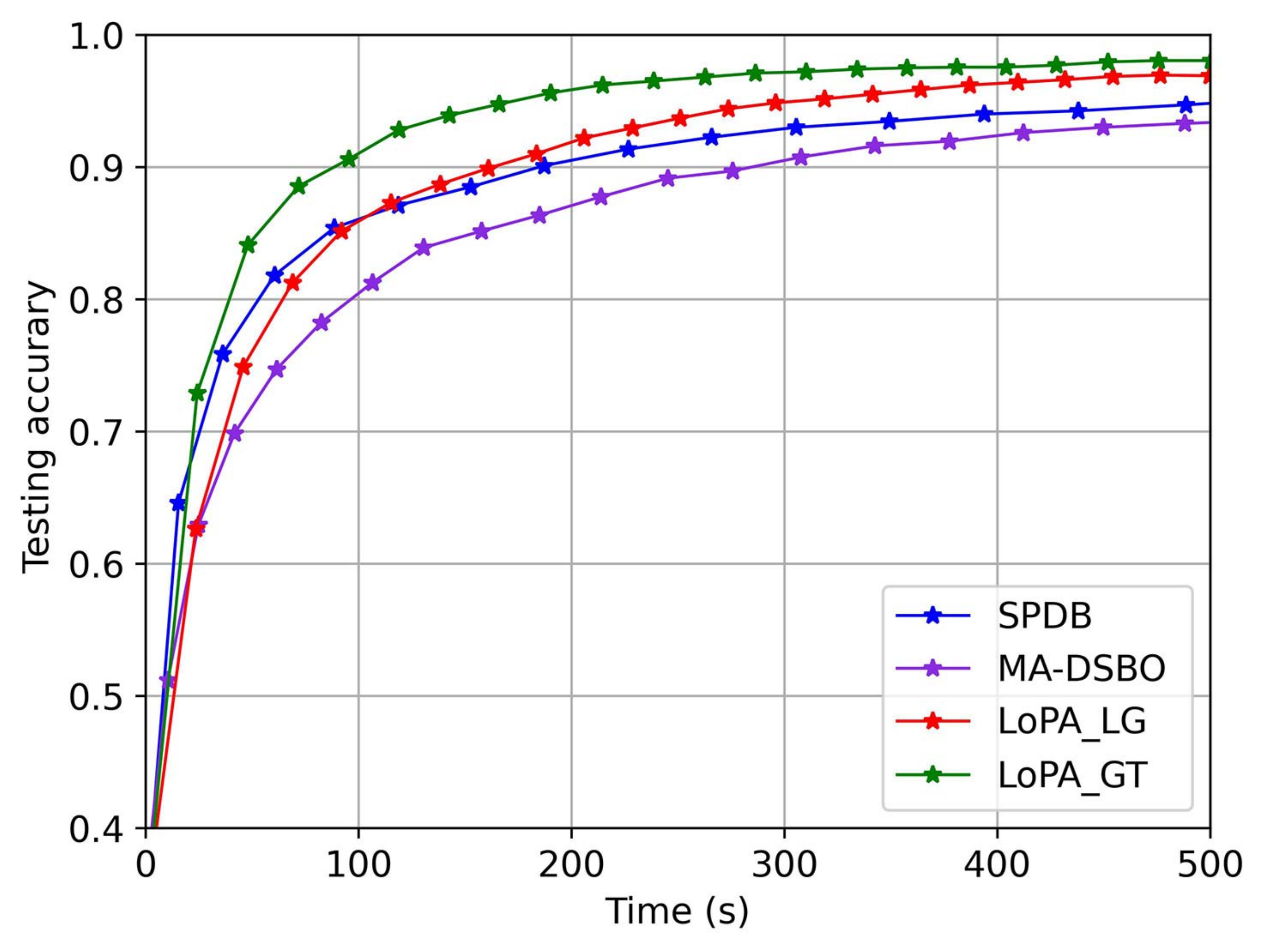} 
\end{minipage}%
}%
\subfigure{
\begin{minipage}[ht]{0.30\linewidth}
\centering
\includegraphics[width=0.98\textwidth,height=0.18\textheight]{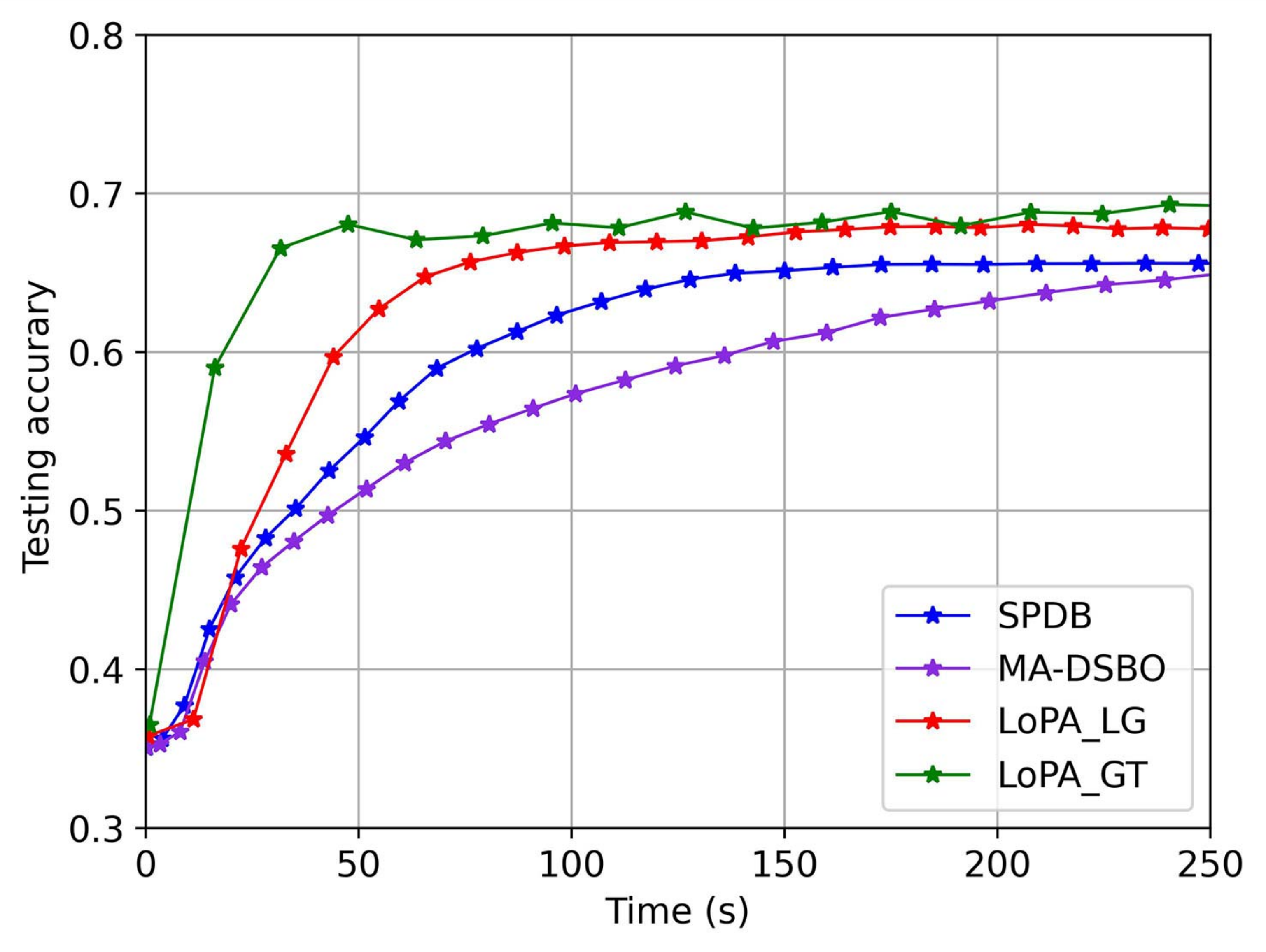} 
\end{minipage}%
}%
\subfigure{
\begin{minipage}[ht]{0.30\linewidth}
\centering
\includegraphics[width=0.98\textwidth,height=0.18\textheight]
{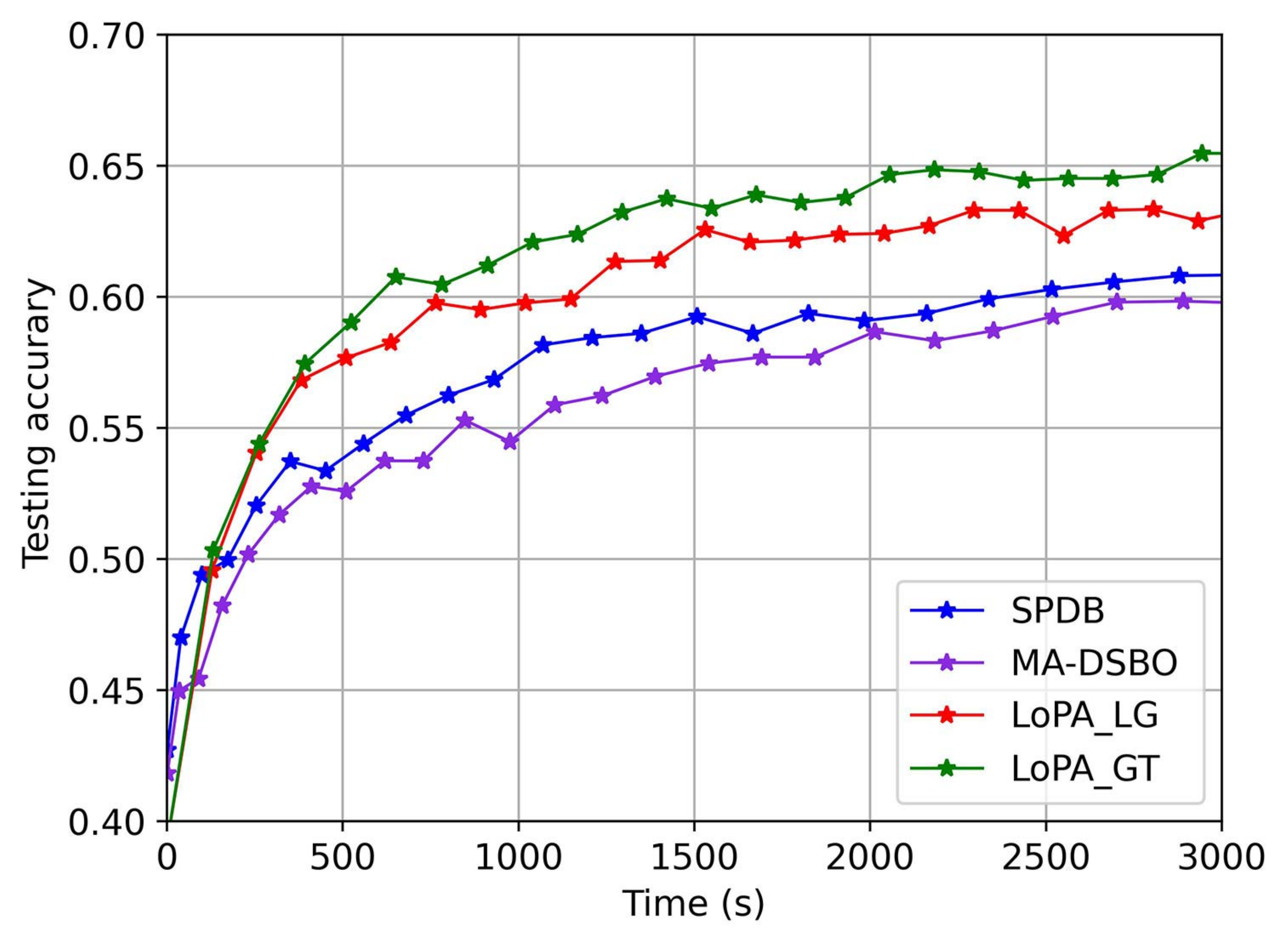} 
\end{minipage}%
}%
\vspace{-0.2cm}
\centering
\caption{Performance comparison of SPDB, MA-DSBO and our LoPA-LG and LoPA-GT algorithms  w.r.t. the computational time for  hyperparameter optimization on \nycres{binary logistic regression problems under different datasets: i) MNIST (first column); ii) covtype (second column); iii) cifar10 (third column).}}
\label{fig:P_hessian}
\end{figure}

\begin{figure}[ht]
\centering
\subfigure{
\begin{minipage}[ht]{0.30\linewidth}
\centering
\includegraphics[width=0.98\textwidth,height=0.18\textheight]
{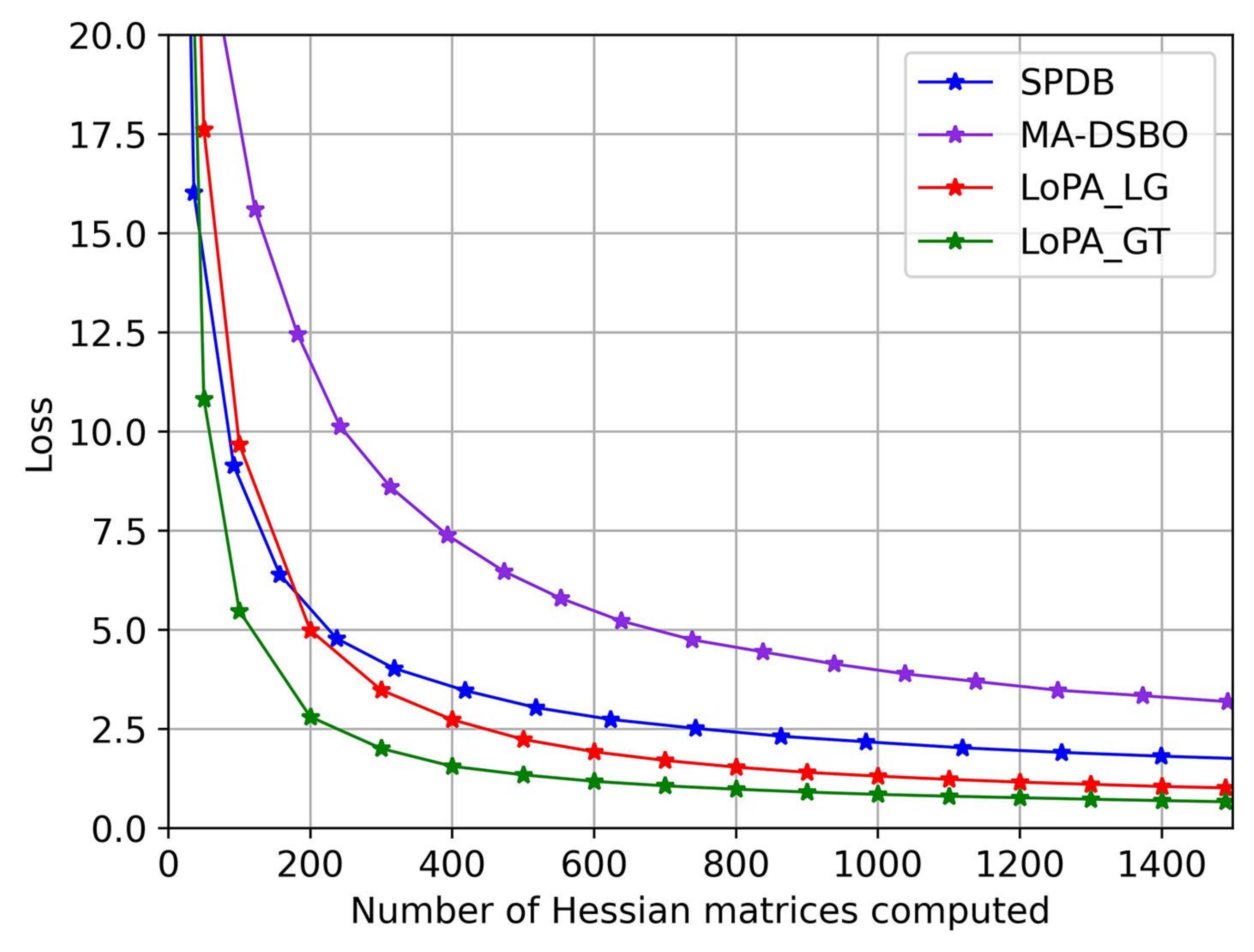} 
\end{minipage}%
}%
\subfigure{
\begin{minipage}[ht]{0.30\linewidth}
\centering
\includegraphics[width=0.98\textwidth,height=0.18\textheight]{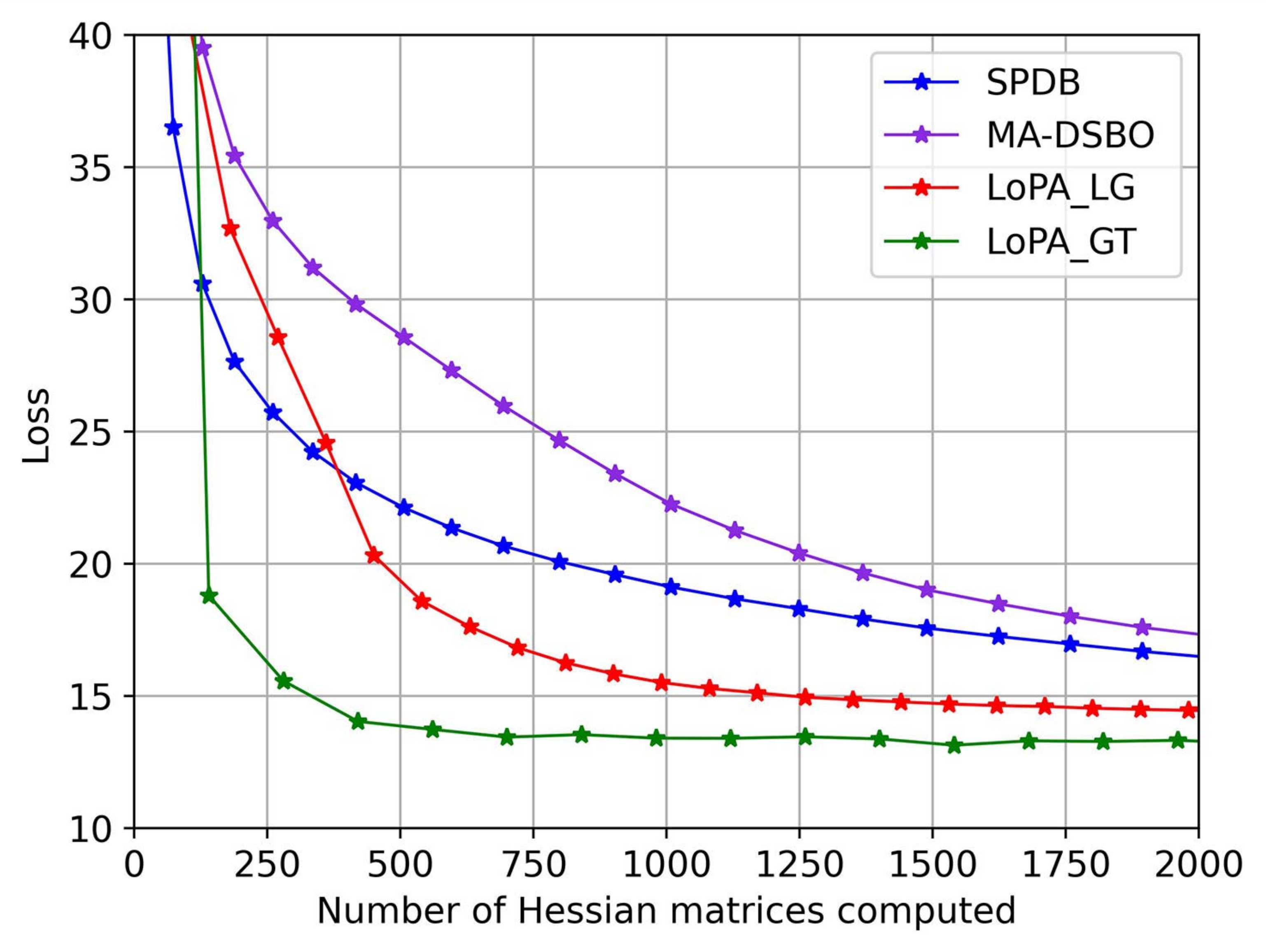} 
\end{minipage}%
}%
\subfigure{
\begin{minipage}[ht]{0.30\linewidth}
\centering
\includegraphics[width=0.98\textwidth,height=0.18\textheight]
{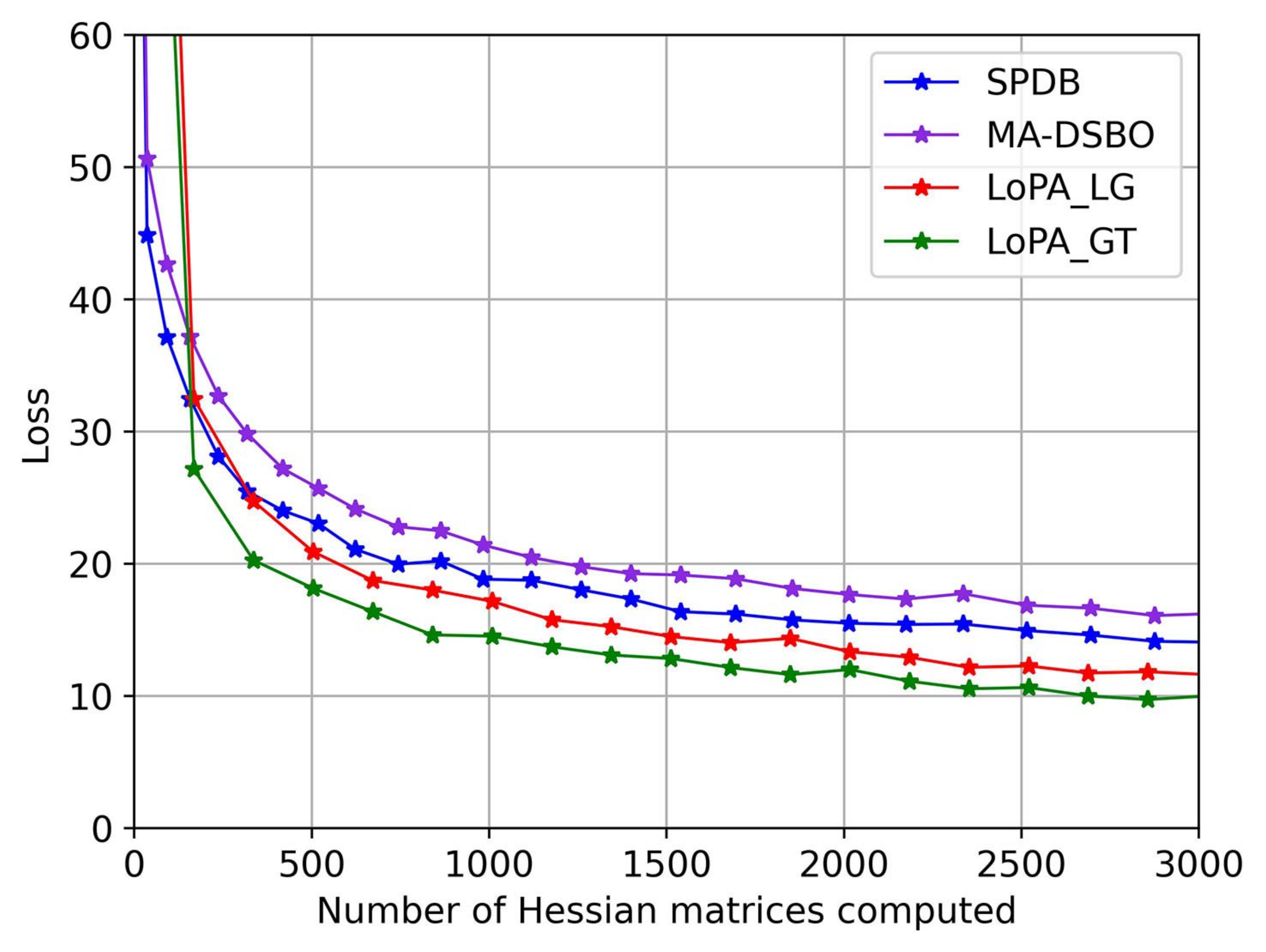} 
\end{minipage}%
}%
\vspace{-0.2cm}
\\
\subfigure{
\begin{minipage}[ht]{0.30\linewidth}
\centering
\includegraphics[width=0.98\textwidth,height=0.18\textheight]
{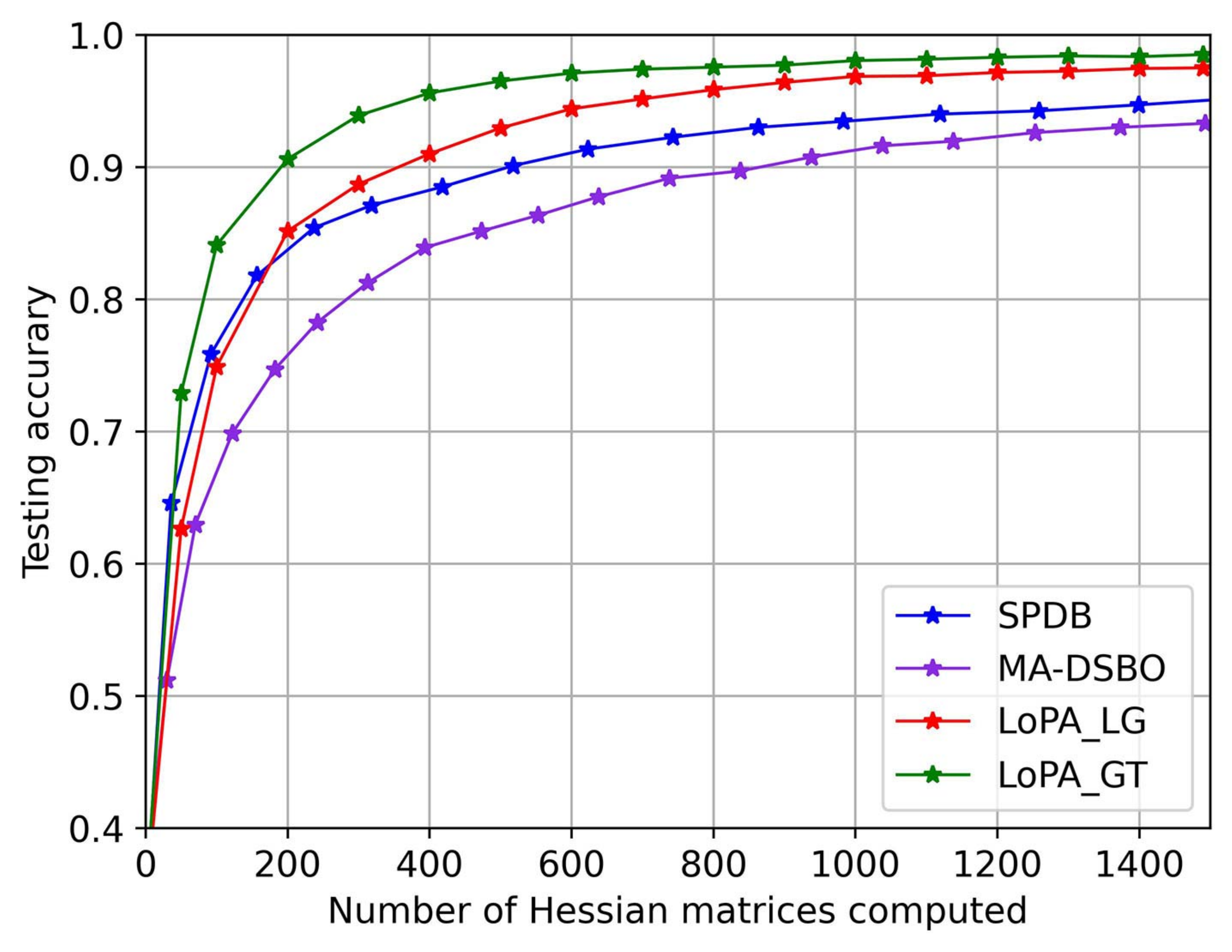} 
\end{minipage}%
}%
\subfigure{
\begin{minipage}[ht]{0.30\linewidth}
\centering
\includegraphics[width=0.98\textwidth,height=0.18\textheight]{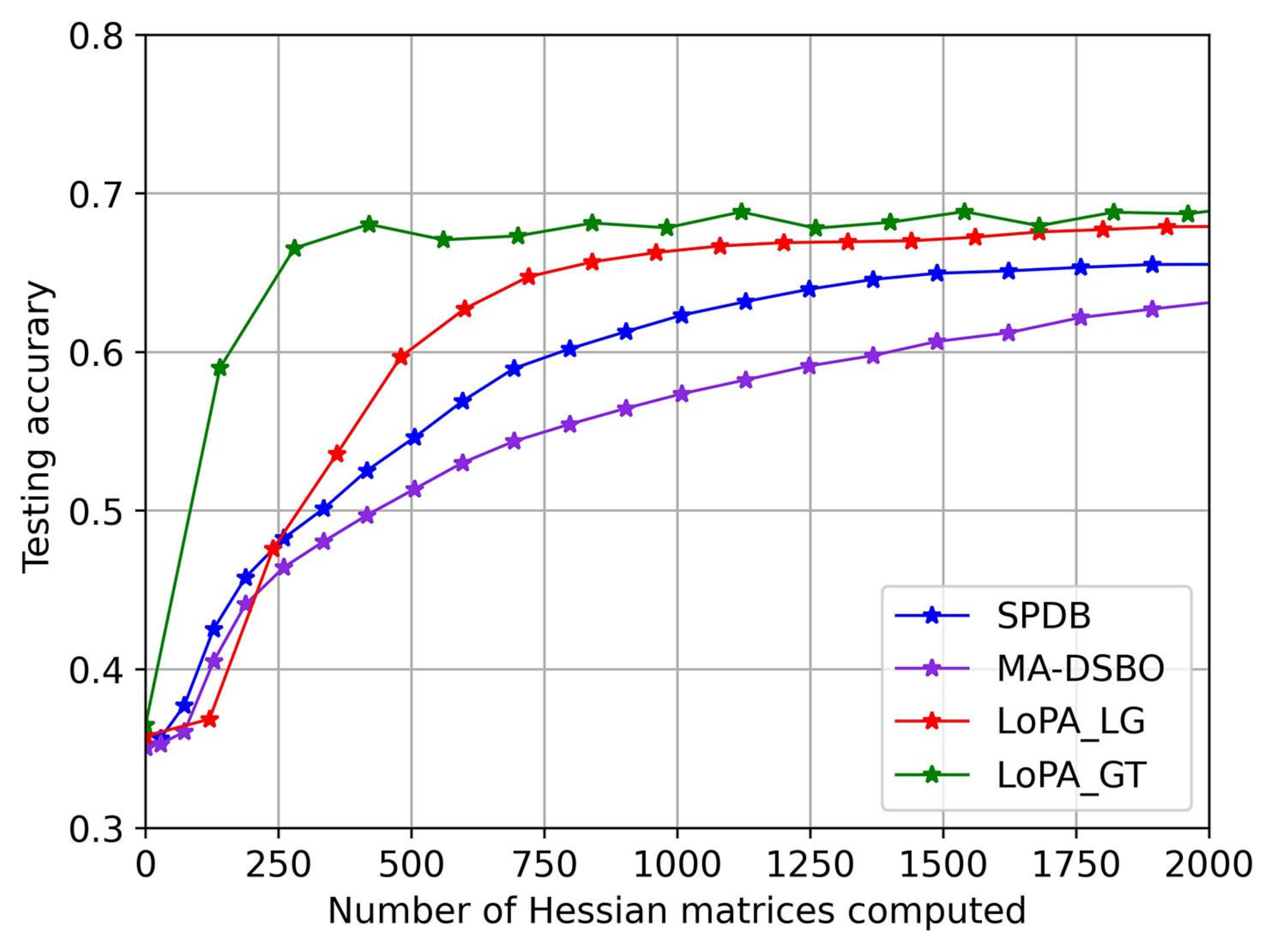} 
\end{minipage}%
}%
\subfigure{
\begin{minipage}[ht]{0.30\linewidth}
\centering
\includegraphics[width=0.98\textwidth,height=0.18\textheight]
{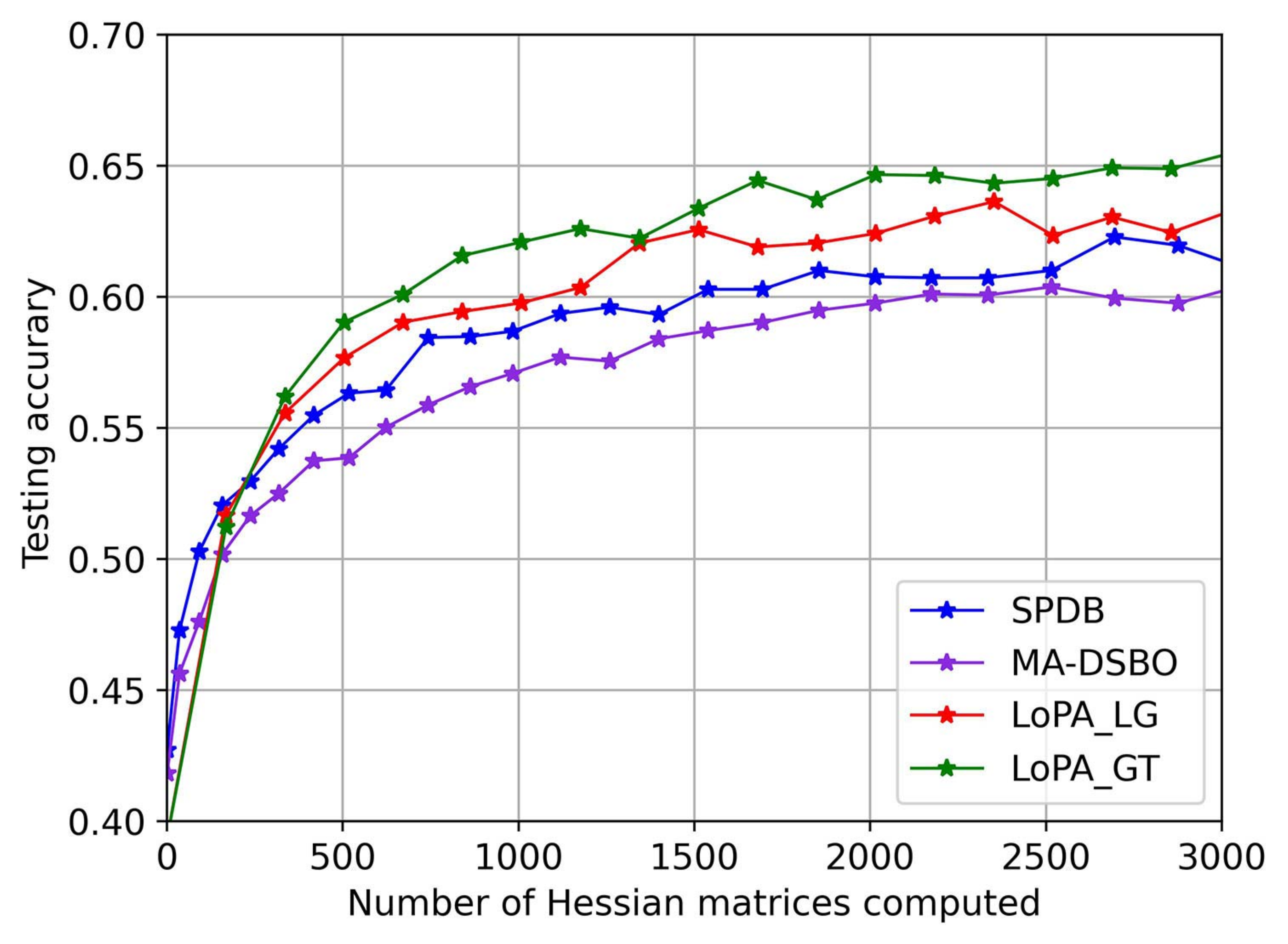} 
\end{minipage}%
}%
\vspace{-0.3cm}
\centering
\caption{Performance comparison  of SPDB, MA-DSBO and our LoPA-LG and LoPA-GT algorithms w.r.t. the number of Hessian matrices  computed for  hyperparameter optimization on \nycres{binary logistic regression problems under different datasets: i) MNIST (first column); ii) covtype (second column); iii) cifar10 (third column).} }
\label{fig:P_time}
\end{figure}

In this subsection, we delve into \nycres{ hyperparameter optimization  in  $l_2$-regularized binary logistic regression problems for a two-class classification scenario}.  The hyperparameter optimization problem can be framed as a class of bilevel optimization.  At the inner level, the goal is to minimize the logistic regression loss on the training data, taking into account the $l_2$-regularization term for a given hyperparameter. Simultaneously, at the outer level, the objective is to maximize the performance of the $l_2$-regularized logistic regression model on a validation set by optimizing the hyperparameter. In this experiment,  we consider the case that each node $i$ has local  validation and training sets, i.e., ${{\mathcal {D}}_i^{\operatorname{val}}}$ and ${{\mathcal {D}}_i^{\operatorname{train}}}$. The goal  of nodes is to cooperatively determine an optimal regularization strength $\lambda$ that can enhance overall performance, while optimizing their personalized  model parameter, i.e.,
\begin{equation}\label{EQ-exp2}
\begin{aligned}
  &\mathop {\min }\limits_{\lambda \in \mathbb{R}^n}  \frac{1}{m}\sum\limits_{i = 1}^m {{f_i}(\lambda ,\theta _i^*(\lambda ))}  =  \frac{1}{m}\sum\limits_{i = 1}^m {\sum\nolimits_{({s_{ij}},{b_{ij}}) \in {\mathcal {D}}_i^{\operatorname{val} }} {\log (1 + {e^{ - ({b_{ij}}s_{ij}^{\operatorname{T}} \theta _i^*(\lambda ))}})} }  \hfill \\
  &\operatorname{s.t.} {\text{   }}\theta _i^*(\lambda ) = \arg \mathop {\min }\limits_{{\theta _i} \in \mathbb{R}^p} {g_i}(\lambda ,{\theta _i}) = \sum\nolimits_{({s_{ij}},{b_{ij}}) \in {\mathcal {D}}_i^{\operatorname{train} }} {\log (1 + {e^{ - ({b_{ij}}s_{ij}^{\operatorname{T}} {\theta _i})}})}  + \theta _i^{\operatorname{T}} \operatorname{diag} \{{{\text{e}}^\lambda }{\} }{\theta _i}.
\end{aligned}
\end{equation}
where ${{\text{e}}^\lambda } = \operatorname{col} \{ {e^{{\lambda _t}}}\} _{t = 1}^p$ with $p=n$,  $\theta_i^*(\lambda)$ is the optimal model in node $i$ given the hyperparameter $\lambda$, and \nycres{$(s_{ij},b_{ij})$ represents the $j$-th sample in node $i$ with $s_{ij} \in \mathbb{R}^{p}$ being the feature and  $b_{ij} \in \mathbb{R}$ being the corresponding label}.

Firstly, we conduct the  experiment
across  various datasets including MNIST (784 features, 12000 samples for digits `0' and `1'), covtype (54 features, 90000 samples for the `Lodgepole' and `Ponderosa' pine classes), and cifar10 (3072 features, 6000 samples for  the `dog' and `horse' classes). The experiment  is implemented in a connected network with  $m=10$, where  the validation and training sets for each node are randomly assigned with a uniform number of samples. To further validate the efficacy of our proposed algorithms, we comprehensively compare them with state-of-the-art SPDB \citep{SPDB} and MA-DSBO \citep{MA-DSBO} algorithms in terms of computational time and the number of Hessian estimates.
The mini-batch sizes are set to 40 for MNIST, 200 for covtype, and 15 for cifar10 across all algorithms. The step-sizes in both the proposed algorithms and the baseline algorithms are manually adjusted for optimal performance. The results are shown in  Figures \ref{fig:P_hessian} and \ref{fig:P_time}.   It can be seen from  Figure \ref{fig:P_hessian} that the proposed algorithms achieve a certain training and loss accuracy with fewer number of Hessian matrices compared to the SPDB \citep{SPDB} and MA-DSBO \citep{MA-DSBO} algorithms. Besides, the proposed algorithms  demonstrate  a reduced time requirement, as depicted in Figure \ref{fig:P_time}. These results further verify the superiority of the proposed algorithms in terms of computational complexity and running efficiency.
\begin{figure}[ht]
\centering
\subfigure[$\rho=0.92$]{
\begin{minipage}[ht]{0.3\linewidth}
\centering
\includegraphics[width=0.9\textwidth,height=0.15\textheight]
{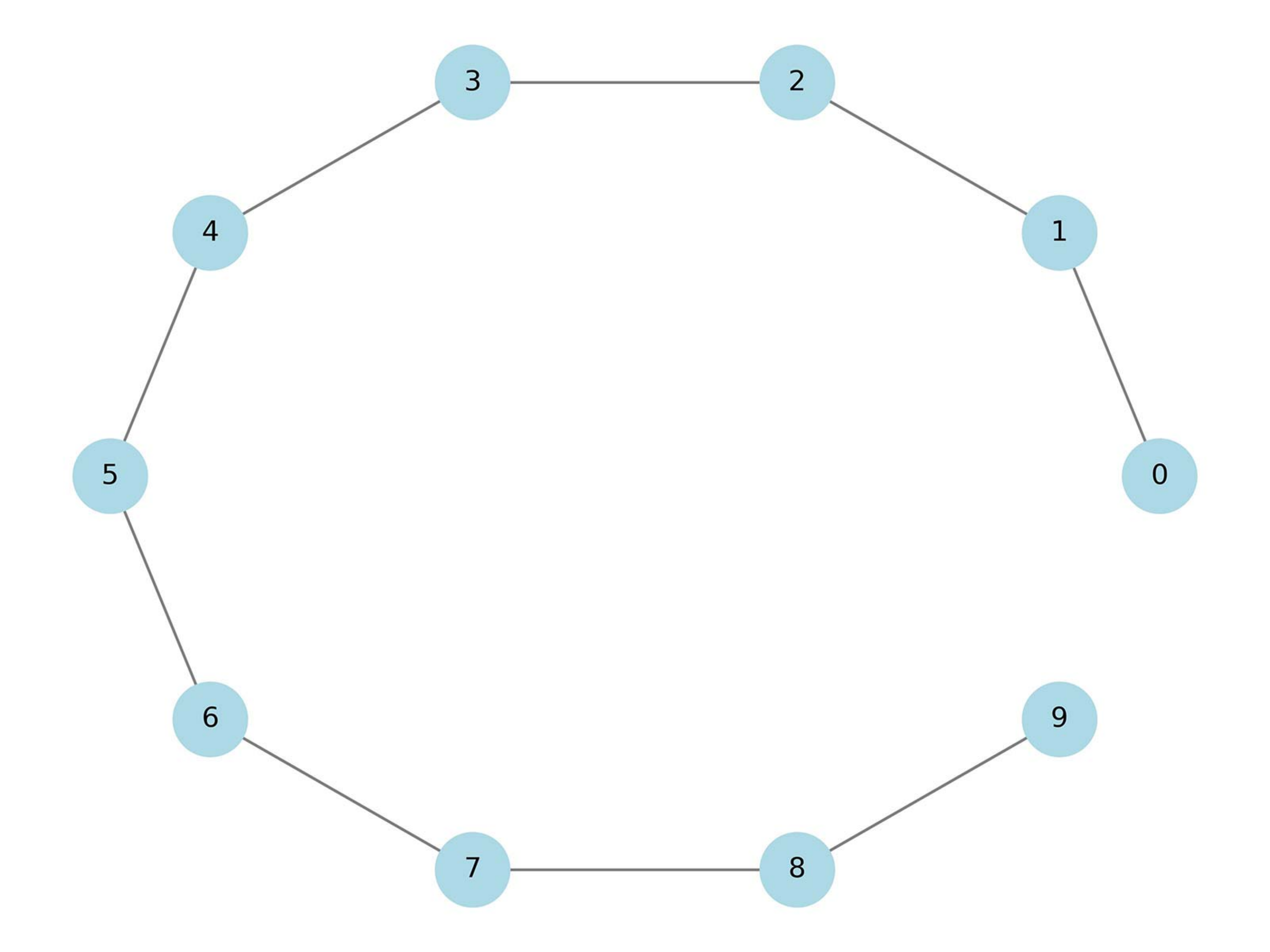} 
\end{minipage}%
}%
\subfigure[$\rho=0.79$]{
\begin{minipage}[ht]{0.25\linewidth}
\centering
\includegraphics[width=0.9\textwidth,height=0.15\textheight]
{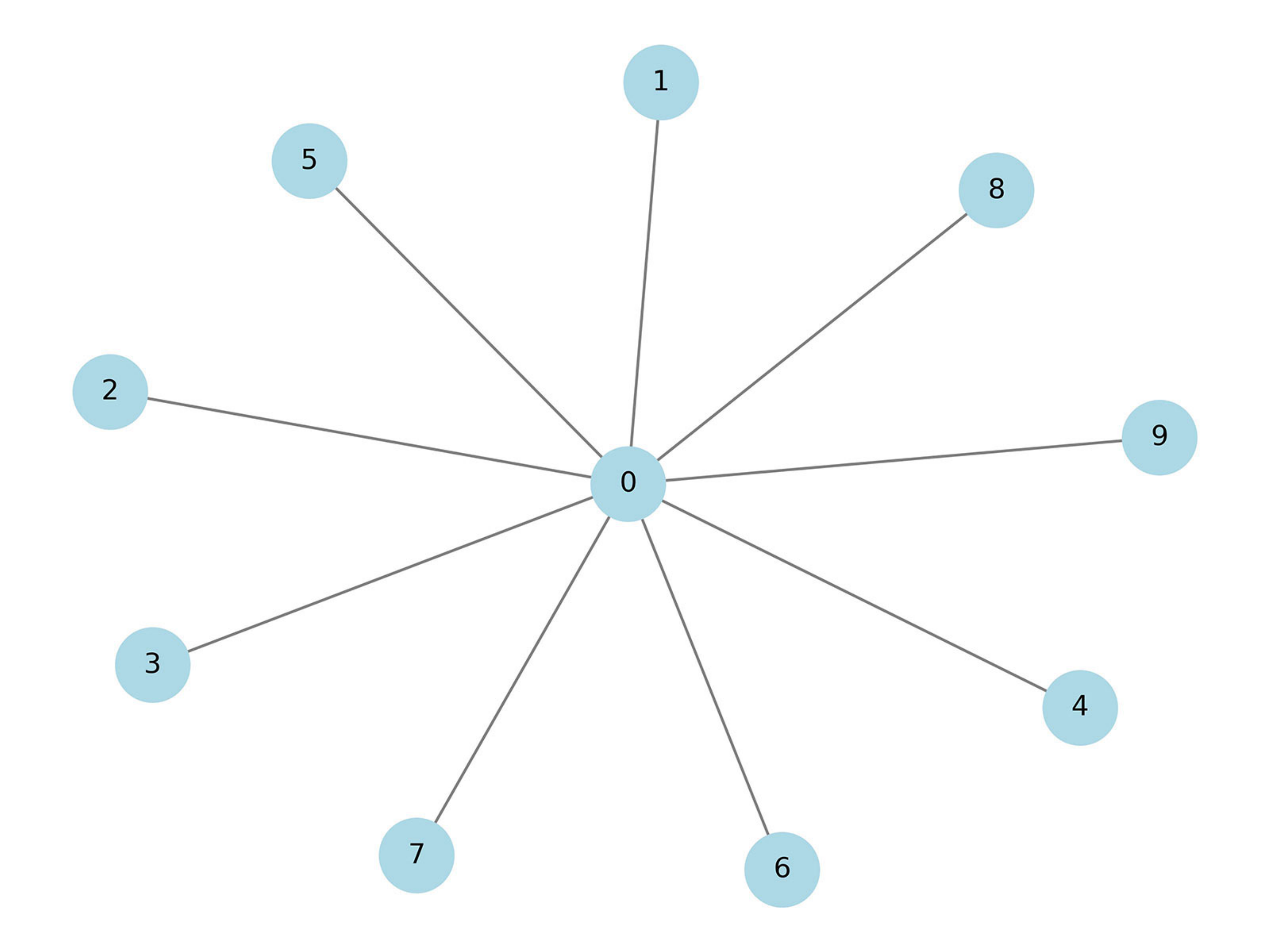} 
\end{minipage}%
}%
\subfigure[$\rho=0.42$]{
\begin{minipage}[ht]{0.25\linewidth}
\centering
\includegraphics[width=0.9\textwidth,height=0.15\textheight]
{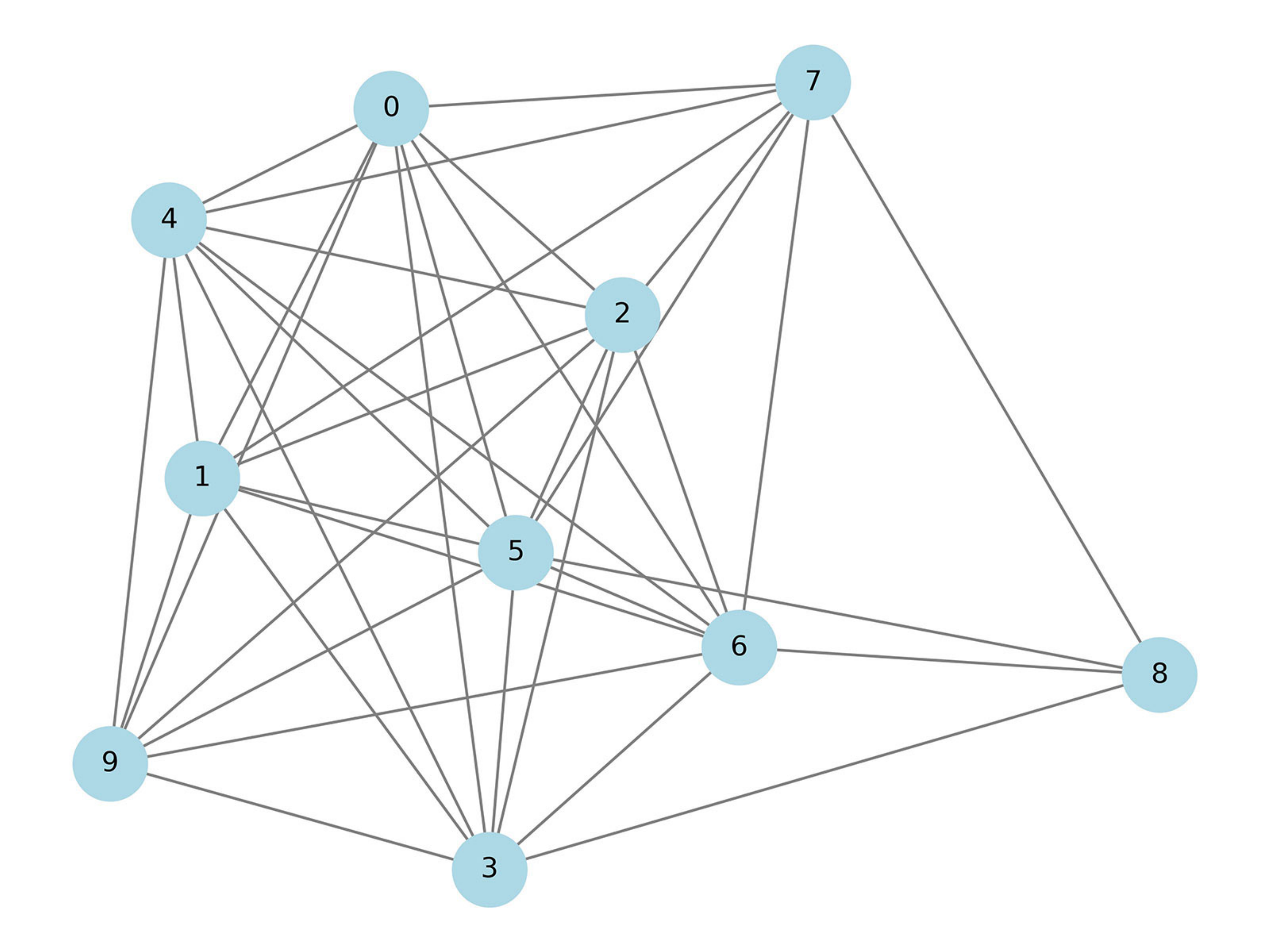} 
\end{minipage}%
}%
\\
\subfigure{
\begin{minipage}[ht]{0.3\linewidth}
\centering
\includegraphics[width=0.98\textwidth,height=0.18\textheight]
{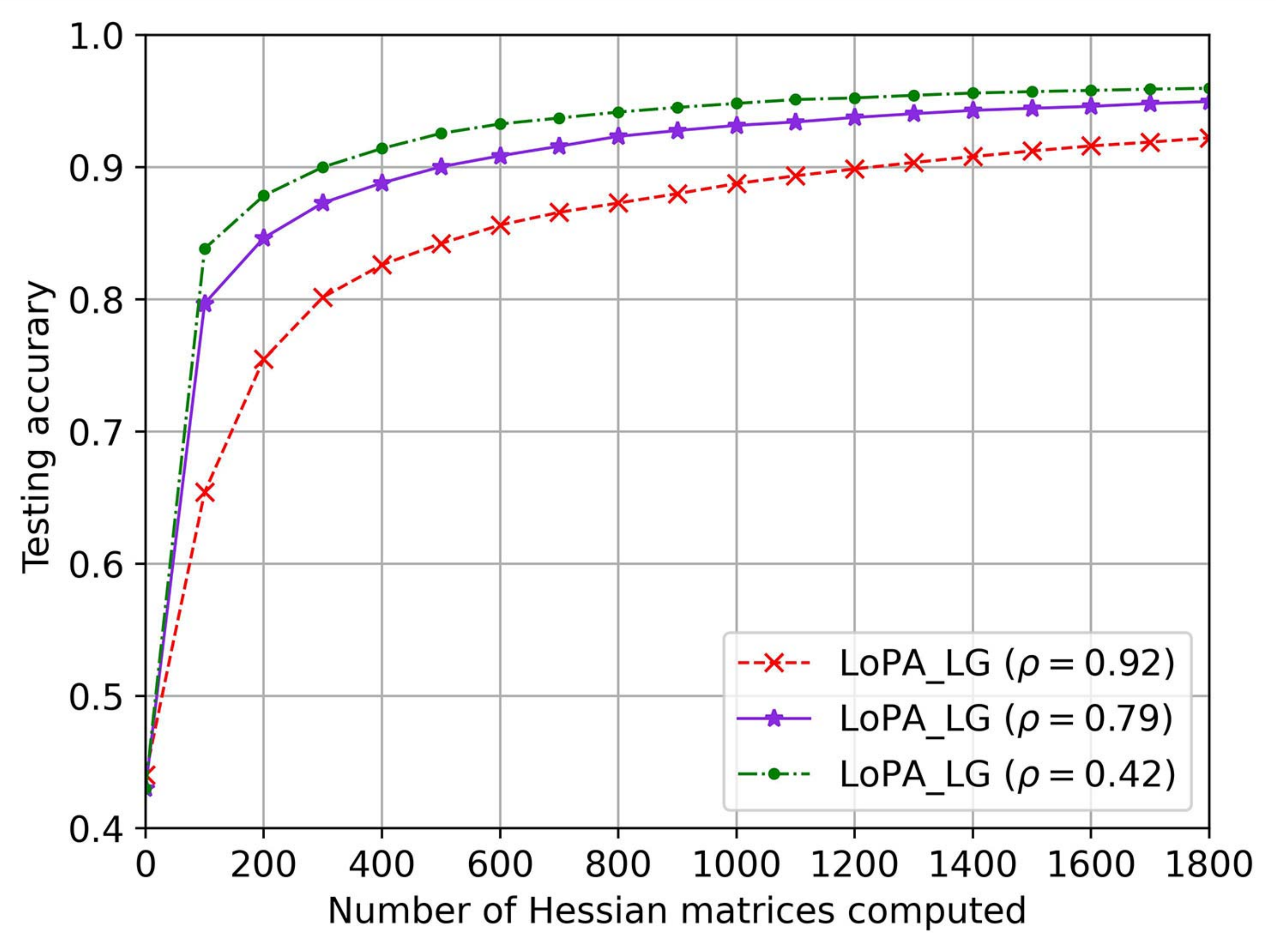} 
\end{minipage}%
}%
\subfigure{
\begin{minipage}[ht]{0.3\linewidth}
\centering
\includegraphics[width=0.98\textwidth,height=0.18\textheight]
{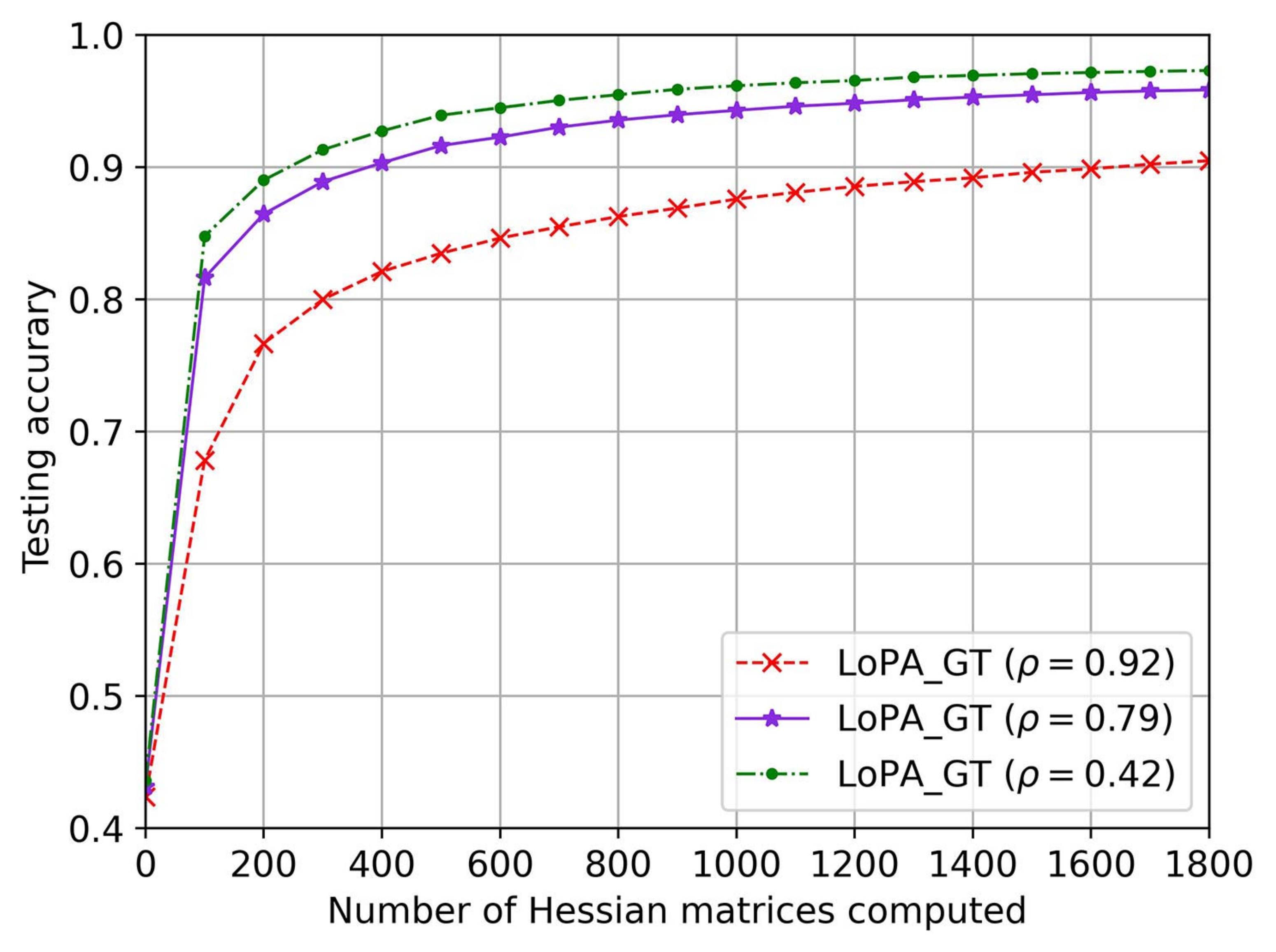} 
\end{minipage}%
}%
\vspace{-0.3cm}
\centering
\caption{\nycres{Performance of  LoPA-LG and LoPA-GT across different topologies  w.r.t.   the number of Hessian matrices  computed for  hyperparameter optimization on  binary logistic regression problems using MNIST dataset. Left: testing accuracy of LG schemes w.r.t. the network spectral gap $\rho$; Right: testing accuracy of GT schemes  w.r.t. the network spectral gap $\rho$. The figure shows GT is more sensitive to the variation of network connectivity.}}
\label{fig:net_diff}
\end{figure}

 \nycres{Additionally, we evaluate the performance of our LoPA-LG and LoPA-GT algorithms across three  types of networks with the network spectral gap $\rho = 0.92$, $\rho = 0.79$, and $\rho = 0.42$ for $m=10$, as shown in Figure \ref{fig:net_diff}(a)-\ref{fig:net_diff}(c). The doubly stochastic matrix is generated using the Metropolis rule. This experiment also considers hyperparameter optimization on the $l_2$-regularized binary logistic regression problem in \eqref{EQ-exp2}  utilizing 6000 samples from the digits `0' and `1' of the MNIST dataset, with validation and training sets randomly assigned to each node in uniform quantities. The mini-batch sizes are set to 25 for both LoPA-LA and LoPA-GT. The results in Figure \ref{fig:net_diff} indicate that as the network spectral gap  $\rho$ increases (i.e., a decrease in network connectivity), both the LoPA-LG and LoPA-GT algorithms exhibit a degradation in performance. Notably, when $\rho$ shifts from $0.79$ to $0.92$, the LoPA-GT algorithm experiences a more significant performance drop compared to LoPA-LG, demonstrating  the GT scheme is more sensitive to the variation of network connectivity. This observation corroborates our theoretical findings in Corollaries \ref{CO-1} and \ref{CO-2}, which
suggest that the GT scheme exhibits greater dependence on network connectivity. }

\begin{figure}[ht]
\subfigure{
\begin{minipage}[ht]{0.30\linewidth}
\centering
\includegraphics[width=0.98\textwidth,height=0.18\textheight]
{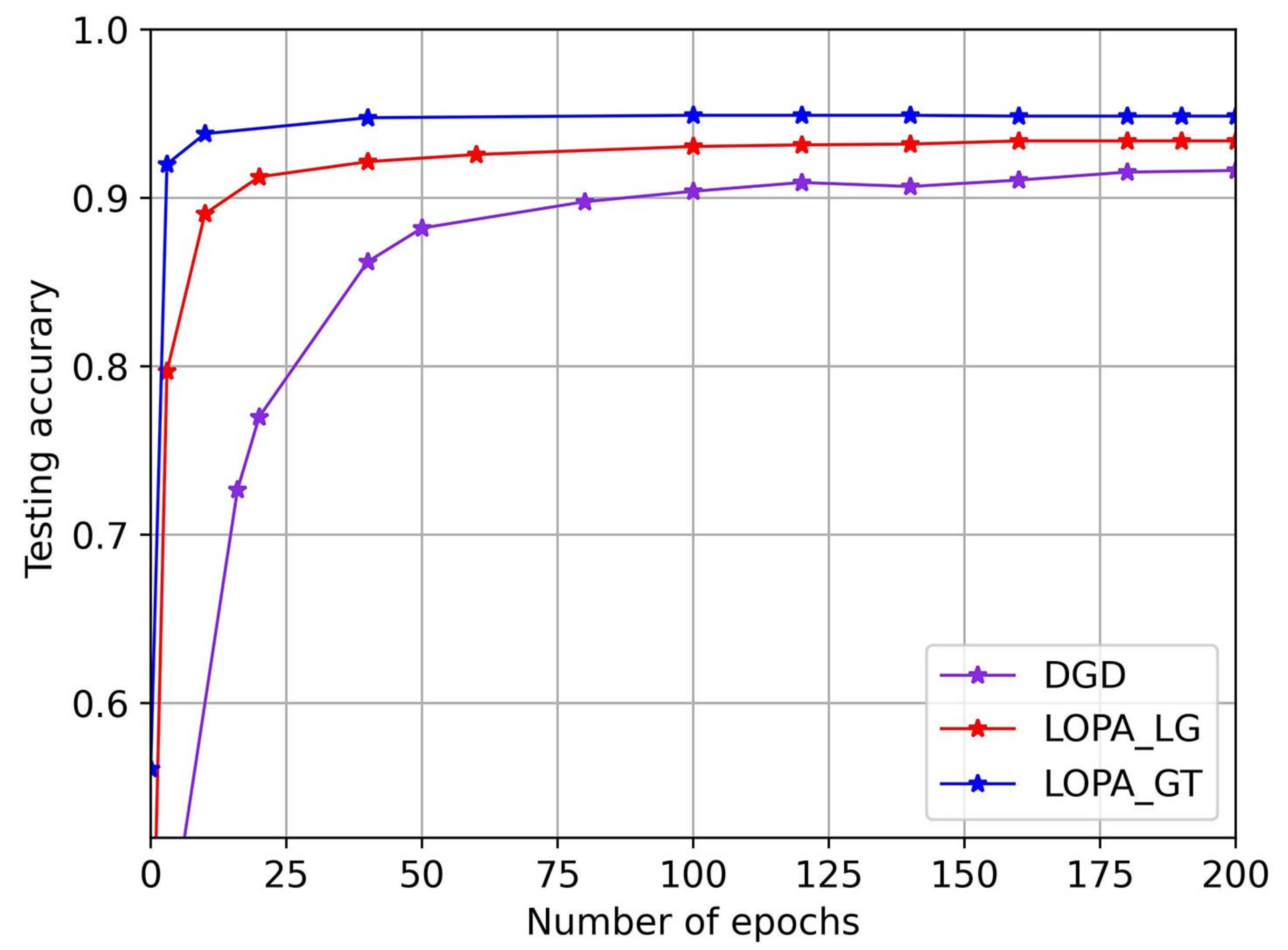} 
\end{minipage}%
}%
\subfigure{
\begin{minipage}[ht]{0.30\linewidth}
\centering
\includegraphics[width=0.98\textwidth,height=0.18\textheight]{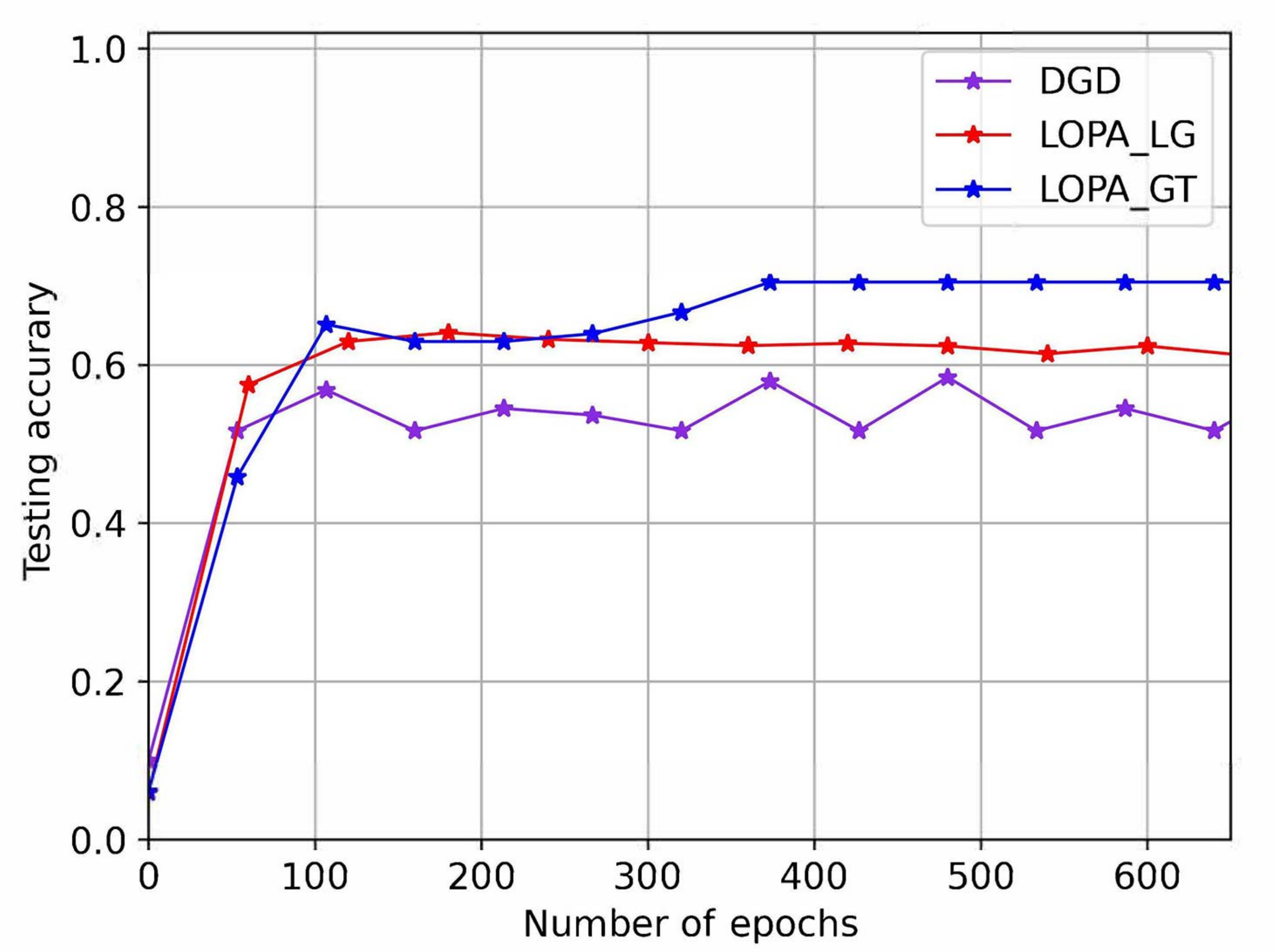} 
\end{minipage}%
}%
\subfigure{
\begin{minipage}[ht]{0.30\linewidth}
\centering
\includegraphics[width=0.98\textwidth,height=0.18\textheight]
{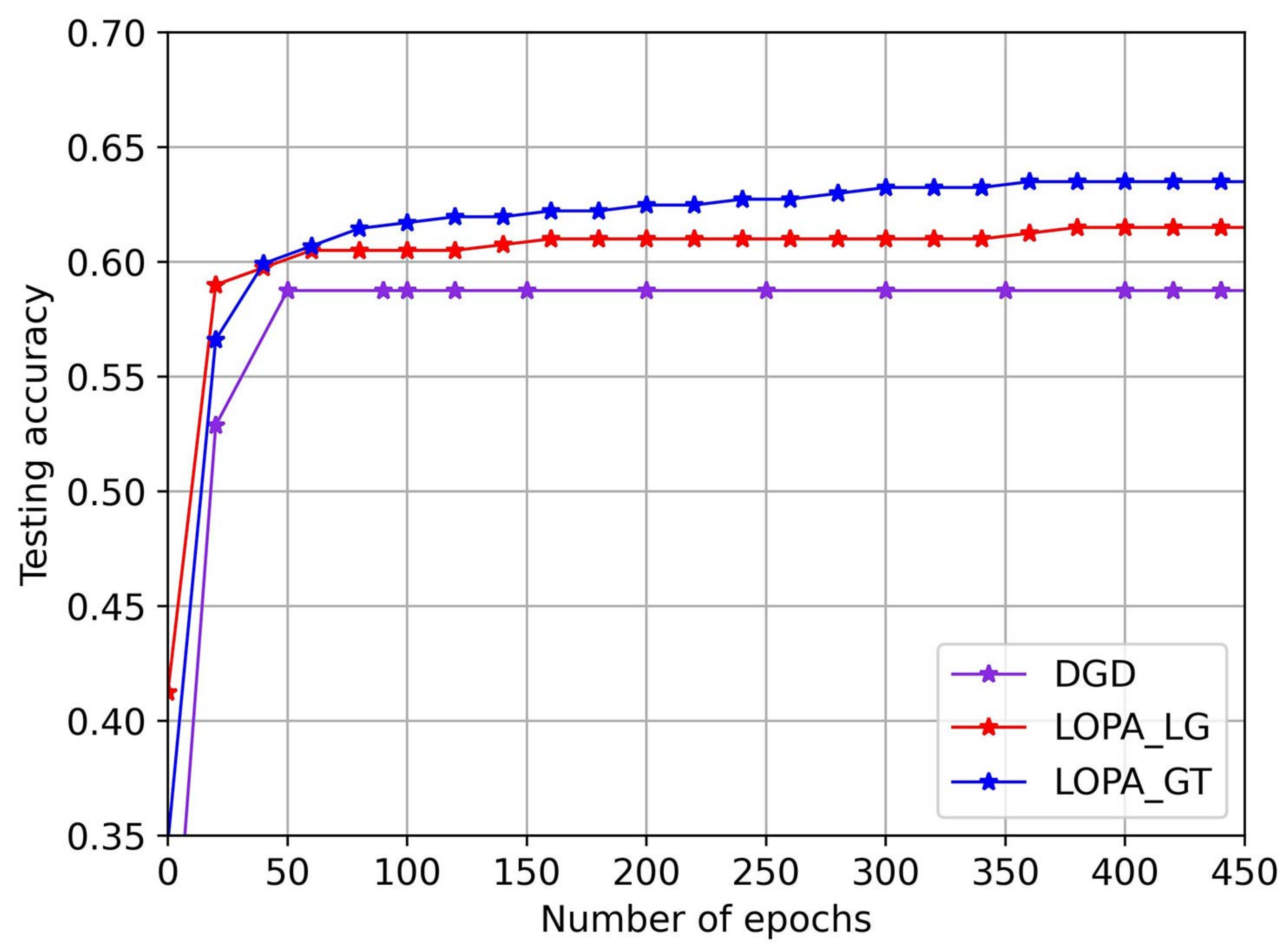} 
\end{minipage}%
}%
\centering
\caption{\nyc{Performance comparison of DGD and our LoPA-LG and  LoPA-GT algorithms for distributed binary logistic regression problems. i) MNIST (first column); ii) covtype (second column); iii) cifar10 (third column).} }
\label{fig:comparsion_DGD}
\end{figure}

 Furthermore, we compare the proposed algorithms with the single-level distributed DGD algorithm \citep{lian2017can} \nycres{under a heterogeneous label distribution scenario} with $m=10$.  Specifically, we conduct the experiment using three types of datasets: MNIST (with 4000 samples for the digits `1' and `3'), covtype (with 30000 samples for the `cottonwood' and `aspen'  classes), and cifar10 (with 2000 samples for the `plane' and `car' classes).    In this case, the proportion of positive and negative labels in half of the nodes is set to be approximately $0.8 \pm a$ with a random $a \in [-0.15, 0.15]$, and in the remaining nodes, it is around $0.35 \pm a$ with a random $a \in [-0.15, 0.15]$.   The DGD algorithm is implemented to train a common regression model by solving the $l_2$-regularized binary logistic regression problem with a fixed $l_2 $-regularization term $\frac{\mu}{2}\|\theta\|^2$, where the regularization coefficient $\mu $ is  manually adjusted to  $0.47, 0.25, 0.15$ for the MINIST, covetype and cifar10 datasets, respectively. Our algorithms aim to train personalized regression models through the hyperparameter optimization on $l_2 $-regularized binary logistic regression problems as outlined in \eqref{EQ-exp2}.
 The mini-batch sizes are set to 20 for MNIST, 150 for covtype, and 10 for cifar10 across all algorithms. The results regarding the average testing accuracy among nodes  are provided in Figure
\ref{fig:comparsion_DGD}. It can be observed that the proposed LoPA-LG and LoPA-GT algorithms achieve superior performance under the heterogeneous scenario without manually adjusting the regularization coefficient.

\section{Conclusion}
In this paper, we have proposed a new loopless  algorithm {\ALG} for solving nonconvex-strongly-convex DSBO problems with personalized inner-level objectives. The proposed {\ALG} algorithm is shown to converge sublinearly, while significantly reducing the computational complexity of gradient evaluation. We have also explicitly characterized the detailed impact of  each-level data heterogeneity on the convergence under  local gradient schemes, and have shown how the heterogeneity is eliminated by employing  gradient tracking schemes. Moreover, we have introduced a unified analytical framework that enables us to obtain a tighter convergence rate and the best known computational complexity for out-level gradient evaluations in DSBO problems, while explicitly revealing the crucial role of the inner-level heterogeneity. Numerical experiments were conducted  to verify the effectiveness of {\ALG} and demonstrate the impact of the heterogeneity.


\bibliography{ref}

\begin{thebibliography}{45}
\providecommand{\natexlab}[1]{#1}
\providecommand{\url}[1]{\texttt{#1}}
\expandafter\ifx\csname urlstyle\endcsname\relax
  \providecommand{\doi}[1]{doi: #1}\else
  \providecommand{\doi}{doi: \begingroup \urlstyle{rm}\Url}\fi

\bibitem[Alghunaim and Yuan(2022)]{alghunaim2022unified}
S.~A. Alghunaim and K.~Yuan.
\newblock A unified and refined convergence analysis for non-convex
  decentralized learning.
\newblock \emph{IEEE Transactions on Signal Processing}, 70:\penalty0
  3264--3279, 2022.

\bibitem[Arbel and Mairal(2022)]{AmiGO}
M.~Arbel and J.~Mairal.
\newblock Amortized implicit differentiation for stochastic bilevel
  optimization.
\newblock \emph{The Tenth International Conference on Learning
  Representations}, 2022.

\bibitem[Bertrand et~al.(2020)Bertrand, Klopfenstein, Blondel, Vaiter,
  Gramfort, and Salmon]{bertrand2020implicit}
Q.~Bertrand, Q.~Klopfenstein, M.~Blondel, S.~Vaiter, A.~Gramfort, and
  J.~Salmon.
\newblock Implicit differentiation of lasso-type models for hyperparameter
  optimization.
\newblock In \emph{International Conference on Machine Learning}, pages
  810--821. PMLR, 2020.

\bibitem[Chen et~al.(2021{\natexlab{a}})Chen, Sun, and Yin]{ALSET}
T.~Chen, Y.~Sun, and W.~Yin.
\newblock Closing the gap: Tighter analysis of alternating stochastic gradient
  methods for bilevel problems.
\newblock \emph{Advances in Neural Information Processing Systems},
  34:\penalty0 25294--25307, 2021{\natexlab{a}}.

\bibitem[Chen et~al.(2022{\natexlab{a}})Chen, Sun, Xiao, and Yin]{STABLE}
T.~Chen, Y.~Sun, Q.~Xiao, and W.~Yin.
\newblock A single-timescale method for stochastic bilevel optimization.
\newblock In \emph{International Conference on Artificial Intelligence and
  Statistics}, pages 2466--2488. PMLR, 2022{\natexlab{a}}.

\bibitem[Chen et~al.(2022{\natexlab{b}})Chen, Huang, Ma, and
  Balasubramanian]{MA-DSBO}
X.~Chen, M.~Huang, S.~Ma, and K.~Balasubramanian.
\newblock Decentralized stochastic bilevel optimization with improved
  per-iteration complexity.
\newblock \emph{arXiv preprint arXiv:2210.12839}, 2022{\natexlab{b}}.

\bibitem[Chen et~al.(2024)Chen, Xiao, and Balasubramanian]{chen2024optimal}
X.~Chen, T.~Xiao, and K.~Balasubramanian.
\newblock Optimal algorithms for stochastic bilevel optimization under relaxed
  smoothness conditions.
\newblock \emph{Journal of Machine Learning Research}, 25\penalty0
  (151):\penalty0 1--51, 2024.

\bibitem[Chen et~al.(2021{\natexlab{b}})Chen, Yuan, Zhang, Pan, Xu, and
  Yin]{chen2021accelerating}
Y.~Chen, K.~Yuan, Y.~Zhang, P.~Pan, Y.~Xu, and W.~Yin.
\newblock Accelerating gossip {SGD} with periodic global averaging.
\newblock In \emph{International Conference on Machine Learning}, pages
  1791--1802. PMLR, 2021{\natexlab{b}}.

\bibitem[Cutkosky and Orabona(2019)]{STORM}
A.~Cutkosky and F.~Orabona.
\newblock Momentum-based variance reduction in non-convex {SGD}.
\newblock \emph{Advances in neural information processing systems}, 32, 2019.

\bibitem[Dagr{\'e}ou et~al.(2022{\natexlab{a}})Dagr{\'e}ou, Ablin, Vaiter, and
  Moreau]{SOBA}
M.~Dagr{\'e}ou, P.~Ablin, S.~Vaiter, and T.~Moreau.
\newblock A framework for bilevel optimization that enables stochastic and
  global variance reduction algorithms.
\newblock \emph{arXiv preprint arXiv:2201.13409}, 2022{\natexlab{a}}.

\bibitem[Dagr{\'e}ou et~al.(2022{\natexlab{b}})Dagr{\'e}ou, Ablin, Vaiter, and
  Moreau]{dagreou2022framework}
M.~Dagr{\'e}ou, P.~Ablin, S.~Vaiter, and T.~Moreau.
\newblock A framework for bilevel optimization that enables stochastic and
  global variance reduction algorithms.
\newblock \emph{arXiv preprint arXiv:2201.13409}, 2022{\natexlab{b}}.

\bibitem[Defazio et~al.(2014)Defazio, Bach, and Lacoste-Julien]{SAGA}
A.~Defazio, F.~Bach, and S.~Lacoste-Julien.
\newblock Saga: A fast incremental gradient method with support for
  non-strongly convex composite objectives.
\newblock \emph{Advances in neural information processing systems}, 27, 2014.

\bibitem[Dong et~al.(2023)Dong, Ma, Yang, and Yin]{dong2023single}
Y.~Dong, S.~Ma, J.~Yang, and C.~Yin.
\newblock A single-loop algorithm for decentralized bilevel optimization.
\newblock \emph{arXiv preprint arXiv:2311.08945}, 2023.

\bibitem[Fallah et~al.(2020)Fallah, Mokhtari, and
  Ozdaglar]{fallah2020personalized}
A.~Fallah, A.~Mokhtari, and A.~Ozdaglar.
\newblock Personalized federated learning with theoretical guarantees: A
  model-agnostic meta-learning approach.
\newblock \emph{Advances in Neural Information Processing Systems},
  33:\penalty0 3557--3568, 2020.

\bibitem[Fang et~al.(2018)Fang, Li, Lin, and Zhang]{SPIDER}
C.~Fang, C.~J. Li, Z.~Lin, and T.~Zhang.
\newblock Spider: Near-optimal non-convex optimization via stochastic
  path-integrated differential estimator.
\newblock \emph{Advances in Neural Information Processing Systems}, 31, 2018.

\bibitem[Finn et~al.(2017)Finn, Abbeel, and Levine]{finn2017model}
C.~Finn, P.~Abbeel, and S.~Levine.
\newblock Model-agnostic meta-learning for fast adaptation of deep networks.
\newblock In \emph{International conference on machine learning}, pages
  1126--1135. PMLR, 2017.

\bibitem[Gao et~al.(2022)Gao, Gu, and Thai]{VRDBO}
H.~Gao, B.~Gu, and M.~T. Thai.
\newblock Stochastic bilevel distributed optimization over a network.
\newblock \emph{arXiv preprint arXiv:2206.15025}, 2022.

\bibitem[Ghadimi and Wang(2018)]{BSA}
S.~Ghadimi and M.~Wang.
\newblock Approximation methods for bilevel programming.
\newblock \emph{arXiv preprint arXiv:1802.02246}, 2018.

\bibitem[Ghadimi et~al.(2020)Ghadimi, Ruszczynski, and Wang]{ghadimi2020single}
S.~Ghadimi, A.~Ruszczynski, and M.~Wang.
\newblock A single timescale stochastic approximation method for nested
  stochastic optimization.
\newblock \emph{SIAM Journal on Optimization}, 30\penalty0 (1):\penalty0
  960--979, 2020.

\bibitem[Goldblum et~al.(2020)Goldblum, Fowl, and Goldstein]{fewshow}
M.~Goldblum, L.~Fowl, and T.~Goldstein.
\newblock Adversarially robust few-shot learning: A meta-learning approach.
\newblock \emph{Advances in Neural Information Processing Systems},
  33:\penalty0 17886--17895, 2020.

\bibitem[Hong et~al.(2023)Hong, Wai, Wang, and Yang]{TTSA}
M.~Hong, H.-T. Wai, Z.~Wang, and Z.~Yang.
\newblock A two-timescale stochastic algorithm framework for bilevel
  optimization: Complexity analysis and application to actor-critic.
\newblock \emph{SIAM Journal on Optimization}, 33\penalty0 (1):\penalty0
  147--180, 2023.

\bibitem[Ji et~al.(2021)Ji, Yang, and Liang]{stoBiO}
K.~Ji, J.~Yang, and Y.~Liang.
\newblock Bilevel optimization: Convergence analysis and enhanced design.
\newblock In \emph{International conference on machine learning}, pages
  4882--4892. PMLR, 2021.

\bibitem[Ji et~al.(2022)Ji, Liu, Liang, and Ying]{ji2022will}
K.~Ji, M.~Liu, Y.~Liang, and L.~Ying.
\newblock Will bilevel optimizers benefit from loops.
\newblock \emph{arXiv preprint arXiv:2205.14224}, 2022.

\bibitem[Jiao et~al.(2022)Jiao, Yang, Wu, Song, and Jian]{jiao2022asynchronous}
Y.~Jiao, K.~Yang, T.~Wu, D.~Song, and C.~Jian.
\newblock Asynchronous distributed bilevel optimization.
\newblock \emph{arXiv preprint arXiv:2212.10048}, 2022.

\bibitem[Khanduri et~al.(2021)Khanduri, Zeng, Hong, Wai, Wang, and
  Yang]{SUSTAIN}
P.~Khanduri, S.~Zeng, M.~Hong, H.-T. Wai, Z.~Wang, and Z.~Yang.
\newblock A near-optimal algorithm for stochastic bilevel optimization via
  double-momentum.
\newblock \emph{Advances in neural information processing systems},
  34:\penalty0 30271--30283, 2021.

\bibitem[Koloskova et~al.(2020)Koloskova, Loizou, Boreiri, Jaggi, and
  Stich]{koloskova2020unified}
A.~Koloskova, N.~Loizou, S.~Boreiri, M.~Jaggi, and S.~Stich.
\newblock A unified theory of decentralized {SGD} with changing topology and
  local updates.
\newblock In \emph{International Conference on Machine Learning}, pages
  5381--5393. PMLR, 2020.

\bibitem[Kong et~al.(2024)Kong, Zhu, Lu, Huang, and
  Yuan]{kong2024decentralized}
B.~Kong, S.~Zhu, S.~Lu, X.~Huang, and K.~Yuan.
\newblock Decentralized bilevel optimization over graphs: Loopless algorithmic
  update and transient iteration complexity.
\newblock \emph{arXiv preprint arXiv:2402.03167}, 2024.

\bibitem[Li et~al.(2022)Li, Gu, and Huang]{FSLA}
J.~Li, B.~Gu, and H.~Huang.
\newblock A fully single loop algorithm for bilevel optimization without
  {Hessian} inverse.
\newblock In \emph{Proceedings of the AAAI Conference on Artificial
  Intelligence}, volume~36, pages 7426--7434, 2022.

\bibitem[Lian et~al.(2017)Lian, Zhang, Zhang, Hsieh, Zhang, and
  Liu]{lian2017can}
X.~Lian, C.~Zhang, H.~Zhang, C.-J. Hsieh, W.~Zhang, and J.~Liu.
\newblock Can decentralized algorithms outperform centralized algorithms? a
  case study for decentralized parallel stochastic gradient descent.
\newblock \emph{Advances in neural information processing systems}, 30, 2017.

\bibitem[Lu et~al.(2022{\natexlab{a}})Lu, Cui, Squillante, Kingsbury, and
  Horesh]{SPDB}
S.~Lu, X.~Cui, M.~S. Squillante, B.~Kingsbury, and L.~Horesh.
\newblock Decentralized bilevel optimization for personalized client learning.
\newblock In \emph{ICASSP 2022-2022 IEEE International Conference on Acoustics,
  Speech and Signal Processing (ICASSP)}, pages 5543--5547. IEEE,
  2022{\natexlab{a}}.

\bibitem[Lu et~al.(2022{\natexlab{b}})Lu, Zeng, Cui, Squillante, Horesh,
  Kingsbury, Liu, and Hong]{SLAM}
S.~Lu, S.~Zeng, X.~Cui, M.~S. Squillante, L.~Horesh, B.~Kingsbury, J.~Liu, and
  M.~Hong.
\newblock A stochastic linearized augmented lagrangian method for decentralized
  bilevel optimization.
\newblock In \emph{Advances in Neural Information Processing Systems},
  2022{\natexlab{b}}.

\bibitem[Madry et~al.(2017)Madry, Makelov, Schmidt, Tsipras, and
  Vladu]{madry2017towards}
A.~Madry, A.~Makelov, L.~Schmidt, D.~Tsipras, and A.~Vladu.
\newblock Towards deep learning models resistant to adversarial attacks.
\newblock \emph{arXiv preprint arXiv:1706.06083}, 2017.

\bibitem[Nedic(2020)]{nedic2020distributed}
A.~Nedic.
\newblock Distributed gradient methods for convex machine learning problems in
  networks: Distributed optimization.
\newblock \emph{IEEE Signal Processing Magazine}, 37\penalty0 (3):\penalty0
  92--101, 2020.

\bibitem[Nedic and Ozdaglar(2009)]{nedic2009distributed}
A.~Nedic and A.~Ozdaglar.
\newblock Distributed subgradient methods for multi-agent optimization.
\newblock \emph{IEEE Transactions on Automatic Control}, 54\penalty0
  (1):\penalty0 48--61, 2009.

\bibitem[Okuno et~al.(2018)Okuno, Takeda, Kawana, and Watanabe]{okuno2018ell}
T.~Okuno, A.~Takeda, A.~Kawana, and M.~Watanabe.
\newblock On $l_p$-hyperparameter learning via bilevel nonsmooth optimization.
\newblock \emph{arXiv preprint arXiv:1806.01520}, 2018.

\bibitem[Rajeswaran et~al.(2019)Rajeswaran, Finn, Kakade, and
  Levine]{rajeswaran2019meta}
A.~Rajeswaran, C.~Finn, S.~M. Kakade, and S.~Levine.
\newblock Meta-learning with implicit gradients.
\newblock \emph{Advances in neural information processing systems}, 32, 2019.

\bibitem[Razaviyayn et~al.(2020)Razaviyayn, Huang, Lu, Nouiehed, Sanjabi, and
  Hong]{razaviyayn2020nonconvex}
M.~Razaviyayn, T.~Huang, S.~Lu, M.~Nouiehed, M.~Sanjabi, and M.~Hong.
\newblock Nonconvex min-max optimization: Applications, challenges, and recent
  theoretical advances.
\newblock \emph{IEEE Signal Processing Magazine}, 37\penalty0 (5):\penalty0
  55--66, 2020.

\bibitem[Shi et~al.(2014)Shi, Ling, Yuan, Wu, and Yin]{shi2014linear}
W.~Shi, Q.~Ling, K.~Yuan, G.~Wu, and W.~Yin.
\newblock On the linear convergence of the {ADMM} in decentralized consensus
  optimization.
\newblock \emph{IEEE Transactions on Signal Processing}, 62\penalty0
  (7):\penalty0 1750--1761, 2014.

\bibitem[Wang et~al.(2016)Wang, Liu, and Fang]{wang2016accelerating}
M.~Wang, J.~Liu, and E.~Fang.
\newblock Accelerating stochastic composition optimization.
\newblock \emph{Advances in Neural Information Processing Systems}, 29, 2016.

\bibitem[Xu et~al.(2015)Xu, Zhu, Soh, and Xie]{xu2015augmented}
J.~Xu, S.~Zhu, Y.~C. Soh, and L.~Xie.
\newblock Augmented distributed gradient methods for multi-agent optimization
  under uncoordinated constant stepsizes.
\newblock In \emph{2015 54th IEEE Conference on Decision and Control (CDC)},
  pages 2055--2060. IEEE, 2015.

\bibitem[Xue et~al.(2021)Xue, Wang, Yan, Hu, Yang, and Sun]{xue2021rethinking}
C.~Xue, X.~Wang, J.~Yan, Y.~Hu, X.~Yang, and K.~Sun.
\newblock Rethinking bi-level optimization in neural architecture search: a
  {Gibbs} sampling perspective.
\newblock In \emph{Proceedings of the AAAI Conference on Artificial
  Intelligence}, volume~35, pages 10551--10559, 2021.

\bibitem[Yang et~al.(2021)Yang, Ji, and Liang]{VRBO}
J.~Yang, K.~Ji, and Y.~Liang.
\newblock Provably faster algorithms for bilevel optimization.
\newblock \emph{Advances in Neural Information Processing Systems},
  34:\penalty0 13670--13682, 2021.

\bibitem[Yang et~al.(2022)Yang, Zhang, and Wang]{yang2022decentralized}
S.~Yang, X.~Zhang, and M.~Wang.
\newblock Decentralized gossip-based stochastic bilevel optimization over
  communication networks.
\newblock \emph{arXiv preprint arXiv:2206.10870}, 2022.

\bibitem[Zhang et~al.(2022)Zhang, Huang, and Yang]{zhang2022interpreting}
M.~Zhang, W.~Huang, and B.~Yang.
\newblock Interpreting operation selection in differentiable architecture
  search: A perspective from influence-directed explanations.
\newblock \emph{Advances in Neural Information Processing Systems},
  35:\penalty0 31902--31914, 2022.

\bibitem[Zhang et~al.(2023)Zhang, Thai, Wu, and Gao]{zhang2023communication}
Y.~Zhang, M.~T. Thai, J.~Wu, and H.~Gao.
\newblock On the communication complexity of decentralized bilevel
  optimization.
\newblock \emph{arXiv preprint arXiv:2311.11342}, 2023.

\end{thebibliography}

\clearpage
\appendix

 \noindent{\LARGE{\textbf{Appendix}}}
 \begin{spacing}{0.85}
 {\small  \tableofcontents}
 \end{spacing}


\newpage
\section{Technical Preliminaries}
\textbf{Notation.}
%
For notional convenience, we  define  some compact notations as follows:\\
${\nabla _\theta }\hat G({x^k},{\theta ^k};\xi _1^{k}) \triangleq \operatorname{col} \{ {{\nabla _\theta }{\hat g_i}(x_i^k,\theta _i^k;\xi_{i,1}^{k })} \}_{i = 1}^m$, $ \nabla _{\theta \theta }^2\hat G({x^k},{\theta ^k};\xi _2^{k}) \triangleq \operatorname{diag} \{ {\nabla _{\theta \theta }^2 \hat g_i(x_i^k,\theta _i^k;\xi _{i,2}^{k})} \}_{i = 1}^m$, \\
$\nabla _{x\theta }^2\hat G( {{x^k},{\theta ^k};\xi _3^{k}} ) \triangleq \operatorname{diag} \{ {\nabla _{x\theta }^2\hat g_i(x_i^k,\theta _i^k;\xi _{i,3}^{k})} \}_{i = 1}^m$, $ {\nabla _\theta }\hat F( {{x^k},{\theta ^k};\varsigma _1^{k }} ) \triangleq \operatorname{col} \{ {{\nabla _\theta }{\hat f_i}(x_i^k,\theta _i^k;\varsigma _{i,1}^{k})} \}_{i = 1}^m$, \\
 ${\nabla _x}\hat F( {{x^k},{\theta ^k};\varsigma _2^{k }} ) \triangleq \operatorname{col} \{ {{\nabla _x}{\hat f_i}(x_i^k,\theta _i^k;\varsigma _{i,2}^{k})} \}_{i = 1}^m$.
 \renewcommand\arraystretch{1.2}
 \\
 \\
\textbf{The proposed {\ALG}  algorithm  in a compact form}. For the sake of subsequent analysis, let $\mathcal{W} = W \otimes I_n$, and we can then rewrite the {\ALG} algorithm in a more compact form as follows:
\begin{subequations}\label{ALG-LPDBO}
\begin{align}
 &{\theta ^{k + 1}} = {\theta ^k} - \beta {d^k}, \hfill \label{ALG-LPDBO-a}\\
  &{v^{k + 1}} = {v^k} - \gamma {h^k}, \hfill \label{ALG-LPDBO-b} \\
  &{x^{k + 1}} = (1-\tau)x^k+\tau(\mathcal{W}{x^k} - \alpha {y^k}), \hfill \label{ALG-LPDBO-c}\\
  &{d^{k + 1}} = {\nabla _\theta }\hat G({x^{k + 1}},{\theta ^{k + 1}};\xi _1^{k + 1}), \label{ALG-LPDBO-d}\\
  &{h^{k + 1}} = \nabla _{\theta \theta }^2\hat G({x^{k + 1}},{\theta ^{k + 1}};\xi _2^{k + 1}){v^{k+1}} - {\nabla _\theta }\hat F( {{x^{k + 1}},{\theta ^{k + 1}};\varsigma _1^{k + 1}} ), \label{ALG-LPDBO-e}\\
   & {s^{k + 1}}= {\nabla _x}\hat F( {{x^{k + 1}},{\theta ^{k + 1}};\varsigma _2^{k + 1}} ) - \nabla _{x\theta }^2\hat G( {{x^{k + 1}},{\theta ^{k + 1}};\xi _3^{k + 1}} ){v^{k + 1}}, \hfill  \label{ALG-LPDBO-f}\\
  &\nyc{{z^{k + 1}} = {s^{k}} + (1 - \gamma )( {{z^k} - {s^k}} ),} \hfill \label{ALG-LPDBO-g}
\end{align}
\end{subequations}
with {\ALGa} updating $y^{k+1}$ as:
\begin{equation}\label{ALG-LPDBO-h1}
\begin{aligned}
 & {y^{k + 1}} = z^{k + 1},
\end{aligned}
\end{equation}
while {\ALGb} updating $y^{k+1}$ as:
\begin{equation}\label{ALG-LPDBO-h2}
\begin{aligned}
 & {y^{k + 1}} = \mathcal{W}{y^k} + {z^{k + 1}} - {z^k}.
\end{aligned}
\end{equation}
\\
\textbf{Basic inequalities.}  In the subsequent analysis, we will utilize a set of fundamental inequalities and equalities to simplify the analysis as follows:

$\bullet$  Young's inequality with  parameter $\eta>0$: ${\| {a + b} \|^2} \le  ( {1 + \frac{1}{\eta }} ){\| a \|^2} + ( {1 + \eta } ){\| b \|^2},\forall a,b.$

$\bullet$  Jensen’s inequality with $l_2$-norm for any vectors $x_1, \cdots, x_m$:  ${\| {\frac{1}{m}\sum\nolimits_{i = 1}^m {{x_i}} } \|^2} \le  \frac{1}{m}\sum\nolimits_{i = 1}^m {{{\| {{x_i}} \|}^2}}. $

$\bullet$   Variance decomposition for a stochastic vector $x$: $\mathbb{E}[ {\| {x - \mathbb{E}[ x ]} \|^2} ] = \mathbb{E}[ {{{\| x \|}^2}} ] - {\| {\mathbb{E}[ x ]} \|^2}.$
\section{Proof of Theorems and Corollaries}

\subsection{Proof of Theorem \ref{TH-1}} \label{sec-proof-TH-1}
To analyze the convergence of {\ALGa}, we need to properly select the  coefficients $d_0$, $d_1$, $d_2$, $d_3$, $d_4$, $d_5$, $d_6$ of the  Lyapunov function \eqref{EQ-V} and  establish the  dynamic of the  function based on the results in Section \ref{sec-proof}.
To this end,  we first set the coefficients $d_0$, $d_1$ and $d_2$   as follows:
\begin{equation}\label{EQ-the1-d0-3}
\begin{aligned}
{d_0} = 1,{d_1} = \frac{{8C_{g,x\theta }^2\tau \alpha }}{{{\mu _g}\lambda }},{d_2} = (\frac{{8C_{g,x\theta }^2\tau \alpha }}{{{\mu _g}\lambda }}{q_x} + {L_{fg,x}}\tau )\frac{\alpha }{{{\omega _\theta }\beta }},{d_3} = \frac{{\tau \alpha }}{{2\gamma }},
\end{aligned}
\end{equation}
where the parameters  $\Ld=2L_{f,x}^2 + 4{M^2}L_{g,x\theta }^2$,  \yc{${q_{x  }}= \frac{{4\yc{\Lb}}\lambda}{{{\mu _g}}\alpha}$} and  $\omega_{\theta}  =\frac{{{\mu _g}{L_{g,\theta }}}}{{2\left( {{\mu _g} + {L_{g,\theta }}} \right)}}$ are defined in Lemmas  \ref{LE-HypergradientE},    \ref{LE-Vstar} and \ref{LE-ThetaStar}, respectively. Then, considering above coefficients and combining Lemmas \ref{LE-descent}-\ref{LE-ThetaStar}, we can reach the following inequality:
\begin{equation}\label{EQ-the1-basisV}
\begin{aligned}
  &{d_0}\mathbb{E}[ {\Phi ({{\bar x}^{k + 1}})} ] \!+ \!{d_1}\frac{1}{m}\mathbb{E}[ {{{\| {{v^{k + 1}} \!-\! {v^*}({{\bar x}^{k + 1}})} \|}^2}} ] \!\\
  &+\! {d_2}\frac{1}{m}\mathbb{E}[ {{{\| {{\theta ^{k + 1}} - {\theta ^*}({{\bar x}^{k + 1}})} \|}^2}} ]\!+\!{d_3}\frac{1}{m}\mathbb{E}[ {{{\| {\nabla \Phi ({{\bar x}^{k + 1}}) \!-\! {{\bar z}^{k + 1}}}\|}^2}} ]  \hfill \\
   \leqslant& {d_0}\mathbb{E}[ {\Phi ({{\bar x}^k})} ] + {d_1}\frac{1}{m}\mathbb{E}[ {{{\| {{v^k} - {v^*}({{\bar x}^k})} \|}^2}} ] + {d_2}\frac{1}{m}\mathbb{E}[ {{{\| {{\theta ^k} - {\theta ^*}({{\bar x}^k})} \|}^2}} ]+ {d_3}\frac{1}{m}\mathbb{E}[ {{{\| {\nabla \Phi ({{\bar x}^{k }}) - {{\bar z}^{k }}}\|}^2}} ] \hfill \\
   &- \frac{{{d_0}}}{2}\tau\alpha \mathbb{E}[ {{{\| \nabla {\Phi ({{\bar x}^k})} \|}^2}} ]- ( {\frac{{{d_0}}}{2}\tau\alpha ( {1 - \tau\alpha L} )  - {d_1}{q_s}{\tau ^2}{\alpha ^2} - {d_2}{p_s}{\tau ^2}{\alpha ^2} - {d_3}{r_y}{\tau ^2}{\alpha ^2}} )\mathbb{E}[ {{{\| {{{\bar y}^k}} \|}^2}} ] \hfill \\
   &- {d_1}\frac{{{\mu _g}\lambda }}{2}\frac{1}{m}\mathbb{E}[ {{{\| {{v^k} - {v^*}({{\bar x}^k})} \|}^2}} ] - {d_2}{\omega _\theta }\beta  \frac{1}{m}\mathbb{E}[ {{{\| {{\theta ^k} - {\theta ^*}({{\bar x}^k})} \|}^2}} ] \\
   &+( {{L_{fg,x}}{d_3}{r_z} + {d_1}{q_x} + {d_2}{p_x}} )\frac{1}{m}\alpha \mathbb{E}[ {{{\| {{x^k} - {1_m} \otimes {{\bar x}^k}} \|}^2}} ] \hfill \\
   &+ \frac{1}{m}{d_3}\sigma _{\bar z}^2{\alpha ^2} + {d_1}\sigma _v^2{\alpha ^2} + {d_2}\sigma _\theta ^2{\alpha ^2}. \hfill
 \end{aligned}
\end{equation}
We next deal with the gradient error term  $\mathbb{E}[ {{{\| {{y^k} - {1_m} \otimes {{\bar y}^k}} \|}^2}} ]$ induced by the consensus errors $\mathbb{E}[ {{{\| {{x^k} - {1_m} \otimes {{\bar x}^k}} \|}^2}} ]$ under local gradient scheme \eqref{EQ-ALG-h1}. Specifically, we let
\[\begin{gathered}
  {d_4} = \left( {{L_{fg,x}}{d_3}{r_z} + {d_1}{q_x} + {d_2}{p_x}} \right)\frac{{24}}{{{{(1 - \rho )}^2}}}\frac{\alpha }{\gamma }{\alpha ^3}, \hfill \\
  {d_5} = 2\left( {{L_{fg,x}}{d_3}{r_z} + {d_1}{q_x} + {d_2}{p_x}} \right)\frac{{2\alpha }}{{\tau (1 - \rho )}}. \hfill \\
\end{gathered} \]
Then, employing  Lemmas \ref{LE-heterogeity} and \ref{LE-stochastic-error-2} gives us:
\begin{equation}\label{EQ-V-d0-d3-xxx}
\begin{aligned}
  &{d_4}\frac{1}{m}\mathbb{E}[ {{{\| {\nabla \tilde \Phi ({{\bar x}^{k + 1}}) - {z^{k + 1}}} \|}^2}} ] + {d_5}\frac{1}{m}\mathbb{E}[ {{{\| {{x^{k + 1}} - {1_m} \otimes {{\bar x}^{k + 1}}} \|}^2}} ] \hfill \\
   \leqslant& {d_4}\frac{1}{m}\mathbb{E}\left[ {{{\| {\nabla \tilde \Phi ({{\bar x}^k}) - {z^k}} \|}^2}} \right] + {d_5}\frac{1}{m}\mathbb{E}[ {{{\| {{x^k} - {1_m} \otimes {{\bar x}^k}} \|}^2}} ] - {d_5}\tau \frac{{1 - \rho }}{2}\frac{1}{m}\mathbb{E}[ {{{\| {{x^k} - {1_m} \otimes {{\bar x}^k}} \|}^2}} ] \hfill \\
   &+ \! {d_5}\frac{{6\tau {\alpha ^2}}}{{1 - \rho }}\frac{{{b^2}}}{m} \!+\! {d_5}\frac{{6\tau {\alpha ^2}}}{{1 - \rho }}[ {{{\| {\nabla \Phi ({{\bar x}^k})} \|}^2}} ] \!+\! {d_4}{r_z}\alpha \frac{1}{m}\mathbb{E}[ {{{\| {\nabla \tilde \Phi ({{\bar x}^k}) - {s^k}} \|}^2}} ]\! +\! {d_4}{r_y}{\tau ^2}{\alpha ^2}\mathbb{E}[ {{{\| {{{\bar y}^k}} \|}^2}} ] \!+\! {d_4}\sigma _z^2{\alpha ^2} \hfill \\
   \leqslant &{d_4}\frac{1}{m}\mathbb{E}[ {{{\| {\nabla \tilde \Phi ({{\bar x}^k}) - {z^k}} \|}^2}} ] + {d_5}\frac{1}{m}\mathbb{E}[ {{{\| {{x^k} - {1_m} \otimes {{\bar x}^k}} \|}^2}} ] - {d_5}\tau \frac{{1 - \rho }}{2}\frac{1}{m}\mathbb{E}[ {{{\| {{x^k} - {1_m} \otimes {{\bar x}^k}} \|}^2}} ] \hfill \\
   &+ {d_5}\frac{{6\tau {\alpha ^2}}}{{1 - \rho }}[ {{{\| {\nabla \Phi ({{\bar x}^k})} \|}^2}} ] \\
   &+ {d_4}{r_y}{\tau ^2}{\alpha ^2}\mathbb{E}[ {{{\| {{{\bar y}^k}} \|}^2}} ] + {d_4}\sigma _z^2{\alpha ^2} + ( {{L_{fg,x}}{d_3}{r_z} + {d_1}{q_x} + {d_2}{p_x}} )\frac{{24\alpha }}{{{{(1 - \rho )}^2}\gamma }}\frac{{{b^2}{\alpha ^3}}}{m} \hfill \\
   &+ 4C_{g,x\theta }^2{d_4}{r_z}\alpha \frac{1}{m}\mathbb{E}[ {{{\| {{v^k} - {v^*}({{\bar x}^k})} \|}^2}} ] + {L_{fg,x}}{d_4}{r_z}\alpha \frac{1}{m}\mathbb{E}[ {{{\| {{\theta ^k} - {\theta ^*}({{\bar x}^k})} \|}^2}} ] \\
   &+ {L_{fg,x}}{d_4}{r_z}\alpha \frac{1}{m}\mathbb{E}[ {{{\| {{x^k} - {1_m} \otimes {{\bar x}^k}} \|}^2}} ], \hfill \\
 \end{aligned}
\end{equation}
where the last step uses the boundedness of the term $\mathbb{E}[ {{{\| {\nabla \tilde \Phi ({{\bar x}^k}) - {s^k}} \|}^2}} ]$ in Lemma \ref{LE-stochastic-error-2}. We proceed in eliminating the term $\mathbb{E}[ {{{\| {{x^k} - {1_m} \otimes {{\bar x}^k}} \|}^2}} ]$ in \eqref{EQ-the1-basisV}. It is noted that for the term $\mathbb{E}[ {{{\| {{x^k} - {1_m} \otimes {{\bar x}^k}} \|}^2}} ]$ in \eqref{EQ-V-d0-d3-xxx} we have $-{d_5}\tau \frac{{1 - \rho }}{2} = -2\left( {{L_{fg,x}}{d_3}{r_z} + {d_1}{q_x} + {d_2}{p_x}} \right)\alpha $. Then,  incorporating  the inequality \eqref{EQ-V-d0-d3-xxx} into the inequality \eqref{EQ-the1-basisV} and letting $d_6 =0$, we obtain:
\begin{equation}\label{EQ-the1-basisV2}
\begin{aligned}
  &\mathbb{E}[ {{V^{k + 1}}} ] \\
  \leqslant& \mathbb{E}[ {{V^k}} ] \hfill - \frac{{{d_0}}}{4}\tau\alpha \mathbb{E}[ {{{\| \nabla{\Phi ({{\bar x}^k})} \|}^2}} ]  - ( {\frac{{{d_0}}}{4}\tau \alpha  - {d_5}\frac{{6\tau {\alpha ^2}}}{{1 - \rho }}})\mathbb{E}[ {{{\| \nabla{\Phi ({{\bar x}^k})} \|}^2}} ]\\
   &- ({( {\frac{{{d_0}}}{2}\tau \alpha \left( {1 - \tau \alpha L} \right) - {d_1}{q_s}{\tau ^2}{\alpha ^2} - {d_2}{p_s}{\tau ^2}{\alpha ^2} - {d_3}{r_y}{\tau ^2}{\alpha ^2} - {d_4}{r_y}{\tau ^2}{\alpha ^2}} )} )\mathbb{E}[ {{{\| {{{\bar y}^k}} \|}^2}} ] \hfill \\
  & - \!({{d_1}\frac{{{\mu _g}\lambda }}{2}\! -\! 4C_{g,x\theta }^2{d_4}{r_z}\alpha })\frac{1}{m}\mathbb{E}[ {{{\| {{v^k} \!-\! {v^*}({{\bar x}^k})} \|}^2}} ]\!-\!({{d_2}{\omega _\theta }\beta  - {L_{fg,x}}{d_4}{r_z}\alpha })\frac{1}{m}\alpha \mathbb{E}[ {{{\| {{\theta ^k} \!-\! {\theta ^*}({{\bar x}^k})} \|}^2}} ] \hfill \\
  & - ({{d_5}\frac{{\tau( 1 - \rho) }}{4} - {L_{fg,x}}{d_4}{r_z}\alpha })\frac{1}{m}\mathbb{E}[ {{{\| {{x^k} - {1_m} \otimes {{\bar x}^k}} \|}^2}} ] \hfill \\
  &+ (\frac{1}{m}{d_3}\sigma _{\bar z}^2 + {d_1}\sigma _v^2 + {d_2}\sigma _\theta ^2 + {d_4}\sigma _z^2){\alpha ^2} + \left( {{L_{fg,x}}{d_3}{r_z} + {d_1}{q_x} + {d_2}{p_x}} \right)\frac{{24\alpha }}{{{{(1 - \rho )}^2}\gamma }}\frac{{{b^2}{\alpha ^3}}}{m} \hfill \\
    =& \mathbb{E}[ {{V^k}} ] \hfill - \frac{{{d_0}}}{4}\tau\alpha \mathbb{E}[ {{{\| \nabla{\Phi ({{\bar x}^k})} \|}^2}} ]  - ( {\frac{{{d_0}}}{4}\tau \alpha  - {d_5}\frac{{6\tau {\alpha ^2}}}{{1 - \rho }}})\mathbb{E}[ {{{\| \nabla{\Phi ({{\bar x}^k})} \|}^2}} ]\\
   &- ({( {\frac{{{d_0}}}{2}\tau \alpha \left( {1 - \tau \alpha L} \right) - {d_1}{q_s}{\tau ^2}{\alpha ^2} - {d_2}{p_s}{\tau ^2}{\alpha ^2} - {d_3}{r_y}{\tau ^2}{\alpha ^2} - {d_4}{r_y}{\tau ^2}{\alpha ^2}} )} )\mathbb{E}[ {{{\| {{{\bar y}^k}} \|}^2}} ] \hfill \\
  & - ({{d_1}\frac{{{\mu _g}\lambda }}{2} - 4C_{g,x\theta }^2{d_4}{r_z}\alpha }-2{(\frac{1}{m}{d_3}\frac{{{\gamma ^2}}}{{{\alpha ^2}}} + {d_1}\frac{{{\lambda ^2}}}{{{\alpha ^2}}} + {d_4}\frac{{{\gamma ^2}}}{{{\alpha ^2}}})\sigma _{g,x\theta }^2{\alpha ^2})})\frac{1}{m}\mathbb{E}[ {{{\| {{v^k} - {v^*}({{\bar x}^k})} \|}^2}} ]\\
  &-({{d_2}{\omega _\theta }\beta  - {L_{fg,x}}{d_4}{r_z}\alpha })\frac{1}{m}\alpha \mathbb{E}[ {{{\| {{\theta ^k} - {\theta ^*}({{\bar x}^k})} \|}^2}} ] \hfill \\
  & - ({{d_5}\frac{{\tau( 1 - \rho) }}{4} - {L_{fg,x}}{d_4}{r_z}\alpha })\frac{1}{m}\mathbb{E}[ {{{\| {{x^k} - {1_m} \otimes {{\bar x}^k}} \|}^2}} ] \hfill \\
  &+ ((\frac{1}{m}{d_3} + {d_4})(\sigma _{f,x}^2 + 2{M^2}\sigma _{g,x\theta }^2)\frac{{{\gamma ^2}}}{{{\alpha ^2}}} + 2{d_1}\left( {\sigma _{f,\theta }^2 + 2{M^2}\sigma _{g,\theta \theta }^2} \right)\frac{{{\lambda ^2}}}{{{\alpha ^2}}} + 2{d_2}\sigma _{g,\theta }^2\frac{{{\beta ^2}}}{{{\alpha ^2}}}){\alpha ^2} \\
  &+ \left( {{L_{fg,x}}{d_3}{r_z} + {d_1}{q_x} + {d_2}{p_x}} \right)\frac{{24\alpha }}{{{{(1 - \rho )}^2}\gamma }}\frac{{{b^2}{\alpha ^3}}}{m}. \hfill \\
 \end{aligned}
\end{equation}
where  the last equality uses the fact that $\sigma_{\bar z}^2=\sigma_{z}^2$ in Lemmas \ref{LE-HypergradientE} and \ref{LE-stochastic-error-2} and the definitions of $\sigma_{\bar z}^2, \sigma_v^2, \sigma_{\theta}^2$ in Lemmas \ref{LE-HypergradientE}, \ref{LE-Vstar}, \ref{LE-ThetaStar}.
If the step-size $\alpha$ satisfies the following conditions:
\begin{align}
\frac{{{d_0}}}{2}\tau \alpha \left( {1 - \tau \alpha L} \right) - {d_1}{q_s}{\tau ^2}{\alpha ^2} - {d_2}{p_s}{\tau ^2}{\alpha ^2} - {d_3}{r_y}{\tau ^2}{\alpha ^2} - {d_4}{r_y}{\tau ^2}{\alpha ^2} \geqslant 0, \hfill\label{EQ-the1-C1} \\
\frac{{{d_0}}}{4}\alpha  - {d_5}\frac{{6\tau {\alpha ^2}}}{{1 - \rho }} \geqslant 0, \label{EQ-the1-C2} \\
{d_1}\frac{{{\mu _g}\lambda }}{2} - 4C_{g,x\theta }^2{d_4}{r_z}\alpha  - 2(\frac{1}{m}{d_3}\frac{{{\gamma ^2}}}{{{\alpha ^2}}} + {d_1}\frac{{{\lambda ^2}}}{{{\alpha ^2}}} + {d_4}\frac{{{\gamma ^2}}}{{{\alpha ^2}}})\sigma _{g,x\theta }^2{\alpha ^2} \geqslant 0, \hfill \label{EQ-the1-C3} \\
{d_2}{\omega _\theta }\beta  - {L_{fg,x}}{d_4}{r_z}\alpha  \geqslant 0, \hfill \label{EQ-the1-C4} \\
{d_5}\frac{{\tau (1 - \rho) }}{4} - {L_{fg,x}}{d_4}{r_z}\alpha  \geqslant 0, \hfill \label{EQ-the1-C5}
\end{align}
then we further have:
\begin{align}
\mathbb{E}[ {{V^{k + 1}}} ] \leqslant& \mathbb{E}[ {{V^k}} ] - \frac{{{d_0}}}{4}\tau\alpha \mathbb{E}[ {{{\| \nabla{\Phi ({{\bar x}^k})} \|}^2}} ] +\left( {{L_{fg,x}}{d_3}{r_z} + {d_1}{q_x} + {d_2}{p_x}} \right)\frac{{24\alpha }}{{{{(1 - \rho )}^2}\gamma }}\frac{{{b^2}{\alpha ^3}}}{m} \label{EQ-1-inq} \\
&+((\frac{1}{m}{d_3} + {d_4})(\sigma _{f,x}^2 + 2{M^2}\sigma _{g,x\theta }^2)\frac{{{\gamma ^2}}}{{{\alpha ^2}}} + 2{d_1}\left( {\sigma _{f,\theta }^2 + 2{M^2}\sigma _{g,\theta \theta }^2} \right)\frac{{{\lambda ^2}}}{{{\alpha ^2}}} + 2{d_2}\sigma _{g,\theta }^2\frac{{{\beta ^2}}}{{{\alpha ^2}}}){\alpha ^2}, \nonumber
\end{align}
where
the  coefficients $d_0$, $d_1$, $d_2$, $d_3$, $d_4$, $d_5$, $d_6$ of the  Lyapunov function \eqref{EQ-V} are  as follows:
\begin{align}
 &{d_0} = 1,
 {d_1} = \frac{{8C_{g,x\theta }^2\tau \alpha }}{{{\mu _g}\lambda }},
 {d_2} = (\frac{{8C_{g,x\theta }^2\tau \alpha }}{{{\mu _g}\lambda }}{q_x} + {L_{fg,x}}\tau )\frac{\alpha }{{{\omega _\theta }\beta }},
 {d_3} = \frac{{\tau \alpha }}{{2\gamma }},\; \label{EQ-the1-dddd} \\
 &{d_4} = \!\left( {{L_{fg,x}}{d_3}{r_z} \!+\! {d_1}{q_x}\! +\!{d_2}{p_x}} \right)\frac{{24}}{{{{(1 - \rho )}^2}}}\frac{\alpha }{\gamma }{\alpha ^3},{d_5} = 2\left( {{L_{fg,x}}{d_3}{r_z} \!+\! {d_1}{q_x} \!+ \!{d_2}{p_x}} \right)\frac{{2\alpha }}{{\tau (1\! - \!\rho )}}, d_6=0. \hfill \nonumber
\end{align}
Next, we proceed to find the sufficient conditions for the step-sizes to satisfy the conditions \eqref{EQ-the1-C1} to \eqref{EQ-the1-C5}. To address the conditions \eqref{EQ-the1-C1}-\eqref{EQ-the1-C5}, we start by simplifying the term ${L_{fg,x}}{d_3}{r_z} + {d_1}{q_x} + {d_2}{p_x}$ in $d_4$ and $d_5$  as:
\begin{equation}\label{EQ-varphi}
\begin{aligned}
 {L_{fg,x}}{d_3}{r_z}  + {d_1}{q_x} + {d_2}{p_x} \hfill
   =  &({L_{fg,x}} + \frac{{32C_{g,x\theta }^2{L_{fg,\theta }}}}{{{\mu _g {\mu _g}}}} + \left( {\frac{{32C_{g,x\theta }^2{L_{fg,\theta }}}}{{{\mu _g^2}}} + {L_{fg,x}}} \right)\frac{{4L_{g,\theta }^2}}{{\omega _\theta ^2}})\tau  \hfill \\
   =& \underbrace{({L_{fg,x}} + \frac{{32C_{g,x\theta }^2{L_{fg,\theta }}}}{{{\mu _g^2}}})(1 + \frac{{4L_{g,\theta }^2}}{{\omega _\theta ^2}}))}_{\triangleq \varphi  }\tau.
\end{aligned}
\end{equation}
\nycres{Then,  $d_2$, $d_4$ and $d_5$ can be  simplified as follows}:
\begin{align}
{d_2} = ({L_{fg,x}}+\frac{{32C_{g,x\theta }^2{L_{fg,\theta }}}}{{\mu _g^2}} )\frac{{\tau \alpha }}{{{\omega _\theta }\beta }}, {d_4} = \frac{{24\varphi \tau {\alpha ^4}}}{{{{(1 - \rho )}^2}\gamma }}, {d_5} = \frac{{2\varphi \alpha }}{{(1 - \rho )}}.  \nonumber
\end{align}

To ensure that condition \eqref{EQ-the1-C1} holds, a sufficient selection condition for the step-sizes $\alpha$, $\beta$, $\lambda$, $\gamma$, $\tau$ is given as:
\begin{align}
\alpha  \leqslant& {u_1} \triangleq  {\frac{1}{{2\tau L}} }, \label{EQ-1-u1} \\
\lambda  \geqslant& \frac{{48{L_{{v^*}}}{C_{g,x\theta }}}}{{{\mu _g}}}\tau \alpha, \label{EQ-1-lambda}\\
\beta  \geqslant& \frac{{6{L_{{\theta ^*}}}}}{{{\omega _\theta }}}{(\frac{{32C_{g,x\theta }^2{L_{fg,\theta }}}}{{{\mu _g ^2}}} + {L_{fg,x}})^{1/2}}\tau \alpha, \label{EQ-1-beta}\\
\gamma  \geqslant&  \max \{ 4L\tau \alpha ,\frac{{32L{\varphi ^{1/2}}}}{{(1 - \rho )}}\tau {\alpha ^2}\}, \label{EQ-1-gamma}
\end{align}
with $0<\tau<1$.
Furthermore, we can derive the following sufficient selection condition for the step-sizes $\alpha$, $\gamma$ and $\lambda$ to satisfy the conditions \eqref{EQ-the1-C2}-\eqref{EQ-the1-C5}:
\begin{align}
&\alpha \leqslant {u_2} \triangleq \min \Bigg\{ \frac{{1 - \rho }}{{12{\varphi ^{1/2}}}},\frac{{{{(1 - \rho )}^{2/3}}}}{{4{\varphi ^{1/3}}}},\frac{{{{(\frac{{32C_{g,x\theta }^2{L_{fg,\theta }}}}{{{\mu _g ^2}}} + {L_{fg,x}})}^{1/3}}{{(1 - \rho )}^{2/3}}}}{{4{\varphi ^{1/3}}}}, \\
&\quad \quad \quad \quad  \quad \quad \quad \quad  \quad    \frac{{{{(1 - \rho )}^{1/2}}}}{{3{{({L_{fg,x}})}^{1/4}}}}, {\frac{{C_{g,x\theta }^{2/3}{{(1 - \rho )}^{2/3}}}}{{4{\varphi ^{1/3}}\sigma _{g,x\theta }^{2/3}}}} \Bigg\}, \nonumber \\
&\gamma  \leqslant \frac{{C_{g,x\theta }^2}}{{\sigma _{g,x\theta }^2}}, \;\lambda \leqslant \frac{{{\mu _g}}}{{16\sigma _{g,x\theta }^2}}, \label{EQ-1-u2}
\end{align}
where $L_{fg,x}$, $L_{fg,\theta}$ and  $\varphi$ are  given by     \eqref{LE-HypergradientE}, \eqref{LE-Vstar} and \eqref{EQ-varphi}, respectively.  It is noted that   $\frac{{{\varphi ^{1/2}}\alpha }}{{1 - \rho }} \leqslant \frac{1}{{12}}$ by the condition \eqref{EQ-1-u2}, which implies that the term $\frac{{32{L^{1/2}}{\varphi ^{1/2}}}}{{(1 - \rho )}}\tau {\alpha ^2}$ in \eqref{EQ-1-gamma}  can be bounded by
\[\frac{{32{L^{1/2}}{\varphi ^{1/2}}}}{{(1 - \rho )}}\tau {\alpha ^2} \leqslant \frac{{32}}{{12}}\tau L\alpha  < 4L\tau \alpha . \]
\nycres{By combining Lemmas \ref{LE-descent}-\ref{LE-ThetaStar} and the above inequality,  we have the following     condition for the  step-sizes $\lambda$, $\beta$ and  $\gamma$ on the basis of \eqref{EQ-1-lambda}-\eqref{EQ-1-u2}:}
\begin{align}
\min\{\frac{1}{{{\mu _g}}}, \frac{{{\mu _g}}}{{8\sigma _{g,x\theta }^2}}\}> \lambda  \geqslant& \frac{{48{L_{{v^*}}}{C_{g,x\theta }}}}{{{\mu _g}}}\tau \alpha, \label{EQ-1-lambda1}
\\
\min \left\{  \frac{2}{{{\mu _g} + {L_{g,\theta }}}}, \frac{{{\mu _g} + {L_{g,\theta }}}}{{2{\mu _g}{L_{g,\theta }}}}  \right \} > \beta  \geqslant& \frac{{6{L_{{\theta ^*}}}}}{{{\omega _\theta }}}{(\frac{{32C_{g,x\theta }^2{L_{fg,\theta }}}}{{{\mu _g ^2}}} + {L_{fg,x}})^{1/2}}\tau \alpha, \label{EQ-1-beta1}
\\
\min\{1, \frac{{2C_{g,x\theta }^2}}{{\sigma _{g,x\theta }^2}}\}> \gamma  \geqslant&   4L\tau \alpha, \label{EQ-1-gamma1}
\end{align}
\nycres{where the selection for $\lambda$, $\beta$ and $\gamma$ can be guaranteed to be non-empty when the step-size $\alpha$ satisfies the following condition:}
\begin{align}
\alpha \leqslant u_3\triangleq \min \Big\{ \frac{\min \{1,\frac{\mu _{g}^{2}}{8\sigma _{g,x\theta}^{2}}\}}{96\tau L_{v^*}C_{g,x\theta}},\frac{\min \left\{ \frac{2}{\mu _g+L_{g,\theta}},\frac{\mu _g+L_{g,\theta}}{2\mu _gL_{g,\theta}} \right\} w_{\theta}}{12\tau (\frac{32C_{g,x\theta}^{2}L_{fg,\theta}}{\mu _{g}^{2}}+L_{fg,x})^{1/2}L_{\theta ^*}},\frac{\min \{1,\frac{2C_{g,x\theta}^{2}}{\sigma _{g,x\theta}^{2}}\}}{8L\tau} \Big\}.\label{EQ-1-u3}
\end{align}
\nycres{As a result, to ensure that the conditions \eqref{EQ-the1-C1}-\eqref{EQ-the1-C5} hold, with $0<\tau<1$, we can properly select the step-size $\alpha$ as follows:
\begin{equation}\label{EQ-1-u}
\begin{aligned}
\nycres{\alpha \leqslant u \triangleq \min \left\{ {{u_1},{u_2}},u_3 \right\},}
\end{aligned}
\end{equation}
where  $u_1$, $u_2$, and $u_3$ are defined by \eqref{EQ-1-u1}, \eqref{EQ-1-u2}, and \eqref{EQ-1-u3}, respectively, and are clearly independent of $\lambda$, $\beta$, and $\gamma$.
Subsequently, we can determine the step sizes $\lambda$, $\beta$ and $\gamma$ according to  \eqref{EQ-1-lambda1}-\eqref{EQ-1-gamma1}.}
Then, by combining the definitions of $d_0$ in \eqref{EQ-the1-dddd} and $\varphi$ in \eqref{EQ-varphi} and considering condition \eqref{EQ-1-lambda1}-\eqref{EQ-1-u},  the inequality \eqref{EQ-1-inq}  can be derived as
\begin{equation}
\begin{aligned}
\frac{\mathbb{E}[ {{V^{k + 1}}} ]}{\tau} \leqslant\frac{\mathbb{E}[ {{V^{k }}} ]}{\tau} - \frac{d_0}{4}\alpha \mathbb{E}[ {{{\| \nabla{\Phi ({{\bar x}^k})} \|}^2}} ] +{\alpha ^2}\sigma _{\sigmaLG}^2 + \frac{\vartheta }{m}{\alpha ^3}{b^2},
\end{aligned}
\end{equation}
where  $\vartheta$  and $\sigma_{\sigmaLG}^2$ are given by:
\begin{align}
\vartheta  \triangleq &\frac{{24\alpha \varphi }}{{{{(1 - \rho )}^2}\gamma }},\label{EQ-sigma-th1}
\\
  \sigma _{\sigmaLG}^2 \triangleq  &(\frac{1}{m}\frac{d_3}{\tau } + \frac{d_4}{\tau } )(\sigma _{f,x}^2 + 2{M^2}\sigma _{g,x\theta }^2)\frac{{{\gamma ^2}}}{{{\alpha ^2}}} + 2\frac{d_1}{\tau } \left( {\sigma _{f,\theta }^2 + 2{M^2}\sigma _{g,\theta\theta }^2} \right)\frac{{{\lambda ^2}}}{{{\alpha ^2}}} + 2\frac{d_2}{\tau } \sigma _{g,\theta }^2\frac{{{\beta ^2}}}{{{\alpha ^2}}} \hfill , \label{EQ-sigma-the1}
\end{align}
where $d_0=1$, ${d_2} = ({L_{fg,x}}+\frac{{32C_{g,x\theta }^2{L_{fg,\theta }}}}{{\mu _g^2}} )\frac{{\tau \alpha }}{{{\omega _\theta }\beta }}$, $d_3 = \frac{{\tau \alpha }}{{2\gamma }}$, ${d_4} = \frac{{24\varphi \tau {\alpha ^4}}}{{{{(1 - \rho )}^2}\gamma }}$,  ${d_5} = \frac{{2\varphi \alpha }}{{(1 - \rho )}}$, $d_6=0$ with $\varphi={({L_{fg,x}} + \frac{{32C_{g,x\theta }^2{L_{fg,\theta }}}}{{{\mu _g^2}}})(1 + \frac{{4L_{g,\theta }^2}}{{\omega _\theta ^2}})}$.
 Now, summing up and telescoping the above inequality  from $k=0$ to $K$ yields:
\begin{equation}\label{EQ-summing}
\begin{aligned}
\frac{1}{{K + 1}}\sum\limits_{k = 0}^K {\mathbb{E}[ {{{\| {\nabla \Phi ({{\bar x}^k})} \|}^2}} ]}  \leqslant \frac{{4({V^0} - {V^K})}}{\tau\alpha(K+1)} + 4\alpha \sigma _{\sigmaLG}^2 + \frac{{4\vartheta }}{m}\alpha^2 {b^2},
\end{aligned}
\end{equation}
where we use the fact that $d_0=1$.
This completes the proof.   {\hfill $\blacksquare$}

\subsection{Proof of Corollary \ref{CO-1}} \label{sec-proof-CO-1}

In the subsequent analysis,  based  on Theorem \ref{TH-1} we will select the step-sizes $\alpha$, $\gamma$,  $\beta$, $\lambda$ in term of the number of iterations $K$.  When the conditions \eqref{EQ-1-lambda1}-\eqref{EQ-1-gamma1} and \eqref{EQ-1-u} hold and the step-sizes $\gamma, \lambda, \beta$ are taken as $\gamma  = {c_\gamma }\alpha ,\lambda  = {c_\lambda }\alpha ,\beta  = {c_\beta }\alpha $ with the positive parameters ${c_\gamma }, {c_\lambda }, {c_\beta}$ being  independent of the terms $K$, $1-\rho$ and satisfying
\begin{equation}\label{EQ-cccc1}
\begin{aligned}
 c_{\lambda} \geqslant \frac{{48{L_{{v^*}}}{C_{g,x\theta }}}}{{{\mu _g}}}\tau,
 {c_\beta } \geqslant \frac{{6{L_{{\theta ^*}}}}}{{{\omega _\theta }}}{(\frac{{32C_{g,x\theta }^2{L_{fg,\theta }}}}{{{\mu _g^2}}} + {L_{fg,x}})^{1/2}}\tau,
 {c_\gamma } \geqslant 4L\tau,
\end{aligned}
\end{equation}
then from  \eqref{EQ-sigma-th1} we have that ${{ d}_0} = 1,{{ d}_1} = \frac{{8C_{g,x\theta }^2\tau }}{{{\mu _g}{c_\lambda }}},{d_2} = (\frac{{32C_{g,x\theta }^2{L_{fg,\theta }}}}{{{\mu _g ^2}}} + {L_{fg,x}})\frac{\tau }{{{\omega _\theta }{c_\beta }}}$, ${ d_3} = \frac{\tau }{{2{c_\gamma }}}$ are independent of the step-size $\alpha$. When the parameters ${c_\gamma }, {c_\lambda }, {c_\beta}$ satisfy \eqref{EQ-cccc1} and their values are properly selected,   the conditions \eqref{EQ-1-lambda1}-\eqref{EQ-1-gamma1} always hold. For examples,  a sufficient  condition for the upper bounds of the parameters  ${c_\gamma }, {c_\lambda }, {c_\beta}$ could be
$ \frac{96L_{v^*}C_{g,x\theta}}{\mu _g}\tau
 > c_{\lambda}$, $\frac{12L_{\theta ^*}}{\omega _{\theta}}(\frac{32C_{g,x\theta}^{2}L_{fg,\theta}}{\mu _{g}^{2}}+L_{fg,x})^{1/2}\tau>c_{\beta}$, and $8L\tau  >c_{\gamma}$, if we consider a sufficient case $\alpha \leqslant u_3$ under the condition \eqref{EQ-1-u} and substitute  it  into the conditions \eqref{EQ-1-lambda1}-\eqref{EQ-1-gamma1}.
Then it follows from the definition of $\varphi$ in \eqref{EQ-varphi}  that $\vartheta$ in \eqref{EQ-sigma-th1} can be rewritten as follows:
\begin{align}
  \vartheta  = \frac{{24  \varphi }}{{{{(1 - \rho )}^2}c_ \gamma }} = {\mathcal{O}}(\frac{1}{{{{(1 - \rho )}^2}}}). \label{EQ-summing-vartheta}
\end{align}
Furthermore, with the above-mentioned selection condition for the step-sizes  $\lambda$, $\beta$, $\gamma$,  we also have that the terms
$\sigma_v$, $\sigma_{\theta}$, $\sigma _{\bar z}$ and ${\sigma _z}$ are  independent  of the step-size $\alpha$. Therefore, by combining  the results that ${d_1}$,  ${d_2}$,  ${d_3}$ are  independent  of the step-size $\alpha$ as well as $d_4={\mathcal{O}(\alpha^3 d_3)}$, the variance related term $\sigma_{\sigmaLG}$ in \eqref{EQ-sigma-the1} can be further derived as:
\begin{align}
  {\sigma _{\sigmaLG}} = &\sqrt {(\frac{1}{m}\frac{d_3}{\tau } + \frac{d_4}{\tau } )(\sigma _{f,x}^2 + 2{M^2}\sigma _{g,x\theta }^2)\frac{{{\gamma ^2}}}{{{\alpha ^2}}} + 2\frac{d_1}{\tau } \left( {\sigma _{f,\theta }^2 + 2{M^2}\sigma _{g,\theta\theta }^2} \right)\frac{{{\lambda ^2}}}{{{\alpha ^2}}} + 2\frac{d_2}{\tau } \sigma _{g,\theta }^2\frac{{{\beta ^2}}}{{{\alpha ^2}}}}  \nonumber \hfill \\
   \leqslant & \sqrt {\frac{2}{m}\frac{d_3}{\tau }(\sigma _{f,x}^2 + 2{M^2}\sigma _{g,x\theta }^2)\frac{{{\lambda ^2}}}{{{\alpha ^2}}} + 2\frac{d_1}{\tau } \left( {\sigma _{f,\theta }^2 + 2{M^2}\sigma _{g,\theta\theta }^2} \right)\frac{{{\gamma ^2}}}{{{\alpha ^2}}} + 2\frac{d_2}{\tau } \sigma _{g,\theta }^2\frac{{{\beta ^2}}}{{{\alpha ^2}}}}  \label{EQ-sigma-r} \hfill \\
   \leqslant &\sqrt{2} {\underbrace {{(\sqrt {\frac{d_1}{\tau}} (\sigma _{f,\theta} + 2{M}\sigma _{g,\theta\theta }){c_\lambda} + \sqrt {\frac{d_2}{\tau}} \sigma _{g,\theta }{c_\beta})}}_{ \triangleq \sigma _{\operatorname{p}} }} + \frac{\sqrt 2}{{\sqrt m }}\underbrace{\sqrt {\frac{d_3}{\tau}}  (\sigma _{f,x} +2 {M}\sigma _{g,x\theta }){c_\gamma}}_{ \triangleq {\sigma _c}} \hfill \triangleq  {{\hat \sigma }_{\sigmaLG}}, \nonumber
\end{align}

\noindent
where the first inequality  holds as
  $d_4$ is high-order term w.r.t $\alpha$ such  that  ${d_4} \ll \frac{1}{m}{d_3}$ under a large number of iterations $K$, and the second inequality is derived by  the   triangle inequality. Meanwhile, from \eqref{EQ-sigma-r},  it follows that
\begin{equation}\label{EQ-sigma-detail}
\begin{aligned}
 & {\sigma _{\operatorname{p} }} =\mathcal{O}({\sigma _{f,\theta }} + {\sigma _{g,\theta \theta }} + {\sigma _{g,\theta }}), \\
  &{\sigma _{\operatorname{c} }}  = \mathcal{O}({\sigma _{f,x}} + {\sigma _{g,x\theta }}). \\
\end{aligned}
\end{equation}
By the condition \eqref{EQ-1-u}, we have that  $\alpha\leqslant u=\mathcal{O}(1-\rho)$. Then, combining the definition of $\varphi$ in \eqref{EQ-varphi} and the condition that ${c_\gamma }, {c_\lambda }, {c_\beta}$ are  independent of the term  $1-\rho$, it follows from \eqref{EQ-sigma-th1} that  the coefficients $d_4$ and $d_5$ can be bounded as follows:
\begin{align}
  {d_4} =&  \frac{{24\varphi \tau \alpha ^4 }}{{{{(1 - \rho )}^2}\gamma }}\leqslant  \frac{{24\varphi\tau  }}{{{{(1 - \rho )}^2}c_\gamma }}{u^3} \triangleq {{\hat d}_4} = \mathcal{O}(1 - \rho ), \label{EQ-the1-hatd4} \hfill \\
  {d_5} =&  \frac{{2\varphi }}{{(1 - \rho )}}\alpha  \leqslant \frac{{2\varphi }}{{(1 - \rho )}}u \triangleq {{\hat d}_5} = \mathcal{O}(1), \hfill \label{EQ-the1-hatd5}
\end{align}
which implies that
\begin{equation} \label{EQ-the1-V0}
\begin{aligned}
   \frac{{V^0} - {V^K}}{\tau}  \leqslant& \frac{V^0}{\tau} \\
   \leqslant& \frac{ d_0}{\tau}\Phi ({{\bar x}^0}) + \frac{d_1}{\tau}\frac{1}{m}{\left\| {{v^0} - {v^*}({{\bar x}^0})} \right\|^2} + \frac{ d_2}{\tau}\frac{1}{m}{\left\| {{\theta ^0} - {\theta ^*}({{\bar x}^0})} \right\|^2} + \frac{ d_3}{\tau}{\left\| {\nabla \Phi ({{\bar x}^0}) - {{\bar z}^0}} \right\|^2} \\
  &+ \frac{\hat d_4}{\tau}\frac{1}{m}{\left\| {\nabla \tilde \Phi ({{\bar x}^0}) - {z^0}} \right\|^2} + \frac{\hat d_5}{\tau}\frac{1}{m}{\left\| {{x^0} - {1_m} \otimes {{\bar x}^0}} \right\|^2} \\
   \triangleq& {{\hat V}^0} =\mathcal{O}(1 - \rho ) ,
\end{aligned}
\end{equation}
where the last equality holds due to the fact that the coefficients $ d_0$, $ d_1$, $ d_2$, and $d_3$ are also independent of the term  $1-\rho$.

For the sake of simplicity, let us consider the following notations:
\nyc{
\begin{equation}\label{EQ-a123}
\begin{aligned}
{a_0} \triangleq {4 \hat V^0},{a_1} \triangleq {{4}}{\hat \sigma_{\sigmaLG}^2},{a_2} \triangleq \frac{{4  \vartheta }}{m}{b^2}.
\end{aligned}
\end{equation}
}
Then, combining  \eqref{EQ-summing-vartheta},
\eqref{EQ-sigma-r} and \eqref{EQ-the1-V0},  the inequality  \eqref{EQ-summing} becomes:
\begin{equation}
\begin{aligned}
\frac{1}{{K + 1}}\sum\limits_{k = 0}^K {{\mathbb{E}}[ {{{\| \nabla{\Phi ({{\bar x}^k})} \|}^2}} ]} \leqslant {a_0}\frac{1}{{\alpha \left( {K + 1} \right)}} + {a_1}\alpha  + {a_2}\alpha^2. \end{aligned}
\end{equation}
When the step-size ${\alpha}$ is taken as  $\alpha  = \min \big\{ {u,{{( {\frac{{{a_0}}}{{{a_1}\left( {K + 1} \right)}}} )}^{\frac{1}{2}}}, {( {\frac{{{a_0}}}{{{a_2}\left( {K + 1} \right)}}} )}^{\frac{1}{3}}} \big\}$ and \res{the step-sizes $\gamma, \lambda, \beta$ are taken as $\gamma  = {c_\gamma }\alpha ,\lambda  = {c_\lambda }\alpha ,\beta  = {c_\beta }\alpha$ with ${c_\gamma }, {c_\lambda }, {c_\beta}$ being  independent of  the terms $K$ and $1-\rho$}, we can proceed with the following discussion: \\
i) When ${( {\frac{{{a_0}}}{{{a_1}\left( {K + 1} \right)}}} )^{\frac{1}{2}}}$ is smallest, we set $\alpha  = {( {\frac{{{a_0}}}{{{a_1}\left( {K + 1} \right)}}} )^{\frac{1}{2}}}$. According to the fact that ${( {\frac{{{a_0}}}{{{a_1}( {K + 1} )}}} )^{\frac{1}{2}}} \leqslant {( {\frac{{{a_0}}}{{{a_2}( {K + 1} )}}} )^{\frac{1}{3}}}$, we have:
  \begin{equation}\label{EQ-discussion1}
\begin{aligned}
\frac{1}{{K + 1}}\sum\limits_{k = 0}^K {{\mathbb{E}}[ {{{\| \nabla{\Phi ({{\bar x}^k})} \|}^2}} ]} \leqslant & {( {\frac{{{a_1}{a_0}}}{{ {K + 1} }}} )^{\frac{1}{2}}} + {( {\frac{{{a_1}{a_0}}}{{ {K + 1} }}} )^{\frac{1}{2}}} + {a_2}{( {\frac{{{a_0}}}{{{a_1}\left( {K + 1} \right)}}} )} \\
\leqslant & 2{( {\frac{{{a_1}{a_0}}}{{ {K + 1} }}} )^{\frac{1}{2}}} + {a_2^{\frac{1}{3}}} {( {\frac{{{a_0}}}{{ {K + 1} }}} )^{\frac{2}{3}}}.
\end{aligned}
\end{equation}
ii)  When ${( {\frac{{{a_0}}}{{{a_2}( {K + 1} )}}} )^{\frac{1}{3}}}$ is smallest, we set $\alpha  = {( {\frac{{{a_0}}}{{{a_2}( {K + 1} )}}} )^{\frac{1}{3}}}$. According to the fact that ${( {\frac{{{a_0}}}{{{a_2}( {K + 1} )}}} )^{\frac{1}{3}}} \leqslant {( {\frac{{{a_0}}}{{{a_1}( {K + 1} )}}} )^{\frac{1}{2}}}$, we have:
  \begin{equation}\label{EQ-discussion2}
\begin{aligned}
 \frac{1}{{K + 1}}\sum\limits_{k = 0}^K {{\mathbb{E}}[ {{{\| \nabla{\Phi ({{\bar x}^k})} \|}^2}} ]}
  \leqslant &2{a_2^{\frac{1}{3}}}{( {\frac{{{a_0}}}{{ {K + 1} }}} )^{\frac{2}{3}}} + {( {\frac{{{a_1}{a_0}}}{{{K + 1} }}} )^{\frac{1}{2}}}.
  \end{aligned}
  \end{equation}
iii)  When $u$ is smallest, we set $\alpha=u$. According to  $u \leqslant {( {\frac{{{a_0}}}{{{a_1}( {K + 1} )}}} )^{\frac{1}{2}}},u \leqslant {( {\frac{{{a_0}}}{{{a_2}( {K + 1} )}}} )^{\frac{1}{3}}}$, we have:
  \begin{equation}\label{EQ-discussion3}
\begin{aligned}
\frac{1}{{K + 1}}\sum\limits_{k = 0}^K {{\mathbb{E}}[ {{{\|\nabla {\Phi ({{\bar x}^k})} \|}^2}} ]}
\leqslant & \frac{{{a_0}}}{{u( {K + 1} )}} + {( {\frac{{{a_1}{a_0}}}{{ {K + 1} }}} )^{\frac{1}{2}}} +{a_2^{\frac{1}{3}}} {( {\frac{{{a_0}}}{{ {K + 1} }}} )^{\frac{2}{3}}}.
\end{aligned}
\end{equation}
According to the above discussion regarding \eqref{EQ-discussion1}, \eqref{EQ-discussion2}, and \eqref{EQ-discussion3}, we can conclude that:
\begin{equation}\label{EQ-a0-a2}
\begin{aligned}
\frac{1}{{K + 1}}\sum\limits_{k = 0}^K {{\mathbb{E}}[ {{{\| \nabla {\Phi ({{\bar x}^k})} \|}^2}} ]}
\leqslant & \frac{{{a_0}}}{{u( {K + 1} )}} + 2{( {\frac{{{a_1}{a_0}}}{{ {K + 1} }}} )^{\frac{1}{2}}} + 2{a_2^{\frac{1}{3}}}{( {\frac{{{a_0}}}{{{K + 1}}}} )^{\frac{2}{3}}}\\
\leqslant &
\frac{{4( {{\hat V^0}} )}}{{ u\left( {K + 1} \right)}} + \frac{{8\hat \sigma_{\sigmaLG}  {( {{\hat V^0} } )}^{\frac{1}{2}} }}{{  {(K + 1)} }^{\frac{1}{2}}} + \frac{ 8(\frac{\vartheta b^2}{m})^{\frac{1}{3}} ( {{\hat V^0}} )^{\frac{2}{3}}} {({ {K + 1} })^{\frac{2}{3}}}.
\end{aligned}
\end{equation}
\nyc{In what follows, we will examine how the parameters in \eqref{EQ-a0-a2} depend on the terms $\kappa$ and $1-\rho$. To this end, we first note that with  the Lipschitz constants defined in \eqref{EQ-Lip-const} and the auxiliary coefficients $L_{fg,x}$, $L_{fg,\theta}$, $\omega_{\theta}$ defined in Lemma \ref{LE-HypergradientE},
 \ref{LE-Vstar}, \ref{LE-ThetaStar} respectively. Then it follows from \eqref{EQ-varphi} that  $\varphi=\mathcal{O}(\kappa^{6})$ and ${\vartheta } ={\mathcal{O}}(\frac{\kappa^{6}}{{{{(1 - \rho )}^2}}})$ in \eqref{EQ-summing-vartheta}. To obtain the best $\kappa$-dependence, we select  $c_{\gamma}=\mathcal{O}(1), c_{\beta} = \mathcal{O}(1), c _{\lambda} =\mathcal{O}(1)$, $\tau = \mathcal{O}(\kappa ^{-4})$ and $u = {{\mathcal{O}}}\left( \kappa^{-3}{(1 - \rho) } \right)$ in terms of $\kappa$ and $1-\rho$. When we initialize the outer-level variables as  $x_i^0=x_j^0, \forall i,j \in {{\mathcal{V}}}$,  it follows that  $\|{{x^0} - {1_m} \otimes {{\bar x}^0}}\|^2=0$ holds. Furthermore, from \eqref{EQ-the1-V0} we know that ${\hat V}^0=\mathcal{O}(\kappa ^{5})$ and ${\hat V}^0$ is  independent of the term $\frac{1}{1-\rho}$. Therefore, recalling the definition of  $\hat \sigma _{\sigmaLG} = {{\mathcal{O}}}( {\sigma _{\operatorname{p} }}  + \frac{1}{\sqrt{m}}{\sigma _{\operatorname{c} }} )$ in \eqref{EQ-sigma-r}, we have
\begin{equation}\label{EQ-a0-a2-noV}
\begin{aligned}
\frac{1}{{K + 1}}\sum\limits_{k = 0}^K {{\mathbb{E}}[ {{{\| \nabla {\Phi ({{\bar x}^k})} \|}^2}} ]}  = {{\mathcal{O}}}\Big( {\frac{{\kappa ^{8}}}{{{{\left( {1 - \rho } \right)}}K}} + \frac{\kappa^{\frac{16}{3}}{(\frac{b}{\sqrt{m}})^{\frac{2}{3}}  }}{{\left( {1 - \rho } \right)^{\frac{2}{3}} K^{\frac{2}{3}} }} + \frac{{1}}{\sqrt{ K }}\kappa^{\frac{5}{2}}( {\sigma _{\operatorname{p} }}  + \frac{1}{\sqrt{m}}{\sigma _{\operatorname{c} }}  ) } \Big),
\end{aligned}
\end{equation}
where  we use the fact that ${\sigma _{\operatorname{p} }}={\mathcal{O}}(\kappa^{\frac{1}{2}}\sigma_{f,\theta}+\kappa^{\frac{3}{2}}\sigma_{g,\theta\theta}+\kappa^{\frac{5}{2}}\sigma_{g,\theta})$ and ${\sigma _{\operatorname{c} }}={\mathcal{O}}(\sigma_{f,x}+\kappa\sigma_{g,x\theta})$
 by \eqref{EQ-sigma-r}. }
This completes the proof.   {\hfill $\blacksquare$}

\subsection{Proof of Theorem \ref{TH-2}} \label{sec-proof-TH-2}
The proof of Theorem \ref{TH-2} follows the similar steps to that of Theorem \ref{TH-1}. Particularly, in addition to $d_0$, $d_1$, $d_2$, $d_3$, $d_4$, $d_5$, we  need to properly select the coefficient  $d_6$ in order to establish the dynamic of  the Lyapunov function \eqref{EQ-V} for {\ALGb} based on the results of Section \ref{sec-proof}. To be specific, by letting $${d_0} = 1,{d_1} = \frac{{8C_{g,x\theta }^2\tau \alpha }}{{{\mu _g}\lambda }},{d_2} = (\frac{{8C_{g,x\theta }^2\tau \alpha }}{{{\mu _g}\lambda }}{q_x} + {L_{fg,x}}\tau )\frac{\alpha }{{{\omega _\theta }\beta }},{d_3} = \frac{{\tau \alpha }}{{2\gamma }}$$
 and combining Lemmas \ref{LE-descent}-\ref{LE-ThetaStar},  the  inequality \eqref{EQ-the1-basisV} can also be derived for {\ALGb}. In the  inequality \eqref{EQ-the1-basisV}, we recall that
the  related parameters are defined in  Lemmas \ref{LE-descent}-\ref{LE-ThetaStar}. In what follows, we focus on dealing  with the consensus errors $\mathbb{E}[ {{{\| {{x^k} - {1_m} \otimes {{\bar x}^k}} \|}^2}} ]$. Note that the term $\mathbb{E}[ {{{\| {{x^k} - {1_m} \otimes {{\bar x}^k}} \|}^2}} ]$ is controlled by the evolution of the gradient errors $\mathbb{E}[{\| {{y^k} - {1_m} \otimes {{\bar y}^k}} \|^2}] $ under  gradient tracking scheme \eqref{EQ-ALG-h2}, while the  errors $\mathbb{E}[{\| {{y^k} - {1_m} \otimes {{\bar y}^k}} \|^2}]$  are further influenced  by the variance errors  $\mathbb{E}[{\| {\nabla \tilde \Phi ({{\bar x}^k}) - {s^k}} \|^2}]$. Motivated by this fact, we first let
\[\begin{gathered}
  {d_4} = \left( {{L_{fg,x}}{d_3}{r_z} + {d_1}{q_x} + {d_2}{p_x}} \right)\frac{{64\gamma {\alpha ^2}}}{{{{(1 - \rho )}^4}}}, \hfill \\
  {d_5} = 2\left( {{L_{fg,x}}{d_3}{r_z} + {d_1}{q_x} + {d_2}{p_x}} \right)\frac{{2\alpha }}{{\tau (1 - \rho )}}, \hfill \\
  {d_6} = \left( {{L_{fg,x}}{d_3}{r_z} + {d_1}{q_x} + {d_2}{p_x}} \right)\frac{{16{\alpha ^2}}}{{{{(1 - \rho )}^3}}}. \hfill \\
\end{gathered} \]
Then, by integrating Lemmas \ref{LE-HypergradientE}, \ref{LE-TE},  \ref{LE-stochastic-error-2},   the  following inequality \eqref{EQ-the1-basisV} can be derived:
\begin{align}
  &{d_4}\mathbb{E}[ {{{\| {\nabla \tilde \Phi ({{\bar x}^{k + 1}}) - {z^{k + 1}}} \|}^2}} ] + {d_5}\mathbb{E}[ {{{\| {{x^{k + 1}} - {1_m} \otimes {{\bar x}^{k + 1}}} \|}^2}} ] + {d_6}{\mathbb{E}}[ {{{\| {{y^{k + 1}} - {1_m} \otimes {{\bar y}^{k + 1}}} \|}^2}} ] \hfill \nonumber \\
   \leqslant& {d_4}\mathbb{E}\left[ {{{\| {\nabla \tilde \Phi ({{\bar x}^k}) - {z^k}} \|}^2}} \right] + {d_5}(1 - \tau \frac{{1 - \rho }}{2})\mathbb{E}[ {{{\| {{x^k} - {1_m} \otimes {{\bar x}^k}} \|}^2}} ] + {d_6}{\mathbb{E}}[ {{{\| {{y^k} - {1_m} \otimes {{\bar y}^k}} \|}^2}} ] \hfill \nonumber \\
   &\!+ \!({d_4}{r_z}\alpha  + {d_6}\frac{4}{{1 - \rho }}{\gamma ^2})\mathbb{E}[ {{{\| {\nabla \tilde \Phi ({{\bar x}^k}) \!- \!{s^k}} \|}^2}} ] \! \nonumber \\
   &+ \!m{d_4}{r_y}{\tau ^2}{\alpha ^2}\mathbb{E}[ {{{\| {{{\bar y}^k}} \|}^2}} ] \!+\! m{d_4}\sigma _z^2{\alpha ^2} \! + \!m\frac{2d_6}{1-\rho}\sigma _y^2{\alpha ^2}. \hfill  \label{EQ-the2-d4-d6}
 \end{align}
Now, noting that \[{d_5}\tau \frac{{1 - \rho }}{2} = 2\left( {{L_{fg,x}}{d_3}{r_z} + {d_1}{q_x} + {d_2}{p_x}} \right)\alpha. \]
and combining  the  inequalities \eqref{EQ-the1-basisV} and \eqref{EQ-the2-d4-d6}, we can establish the following dynamic:
\begin{equation}
\begin{aligned}
  \mathbb{E}[ {{V^{k + 1}}} ] \leqslant& \mathbb{E}[ {{V^k}} ] - \frac{1}{2}{d_0}\tau\alpha \mathbb{E}[ {{{\| {\nabla \Phi ({{\bar x}^k})} \|}^2}} ] \hfill \\
   &-(\frac{{{d_0}}}{2}\tau \alpha \left( {1 - \tau \alpha L} \right) - {d_1}{q_s}{\tau ^2}{\alpha ^2} - {d_2}{p_s}{\tau ^2}{\alpha ^2} - {d_3}{r_y}{\tau ^2}{\alpha ^2} - {d_4}{r_y}{\tau ^2}{\alpha ^2})\mathbb{E}[ {{{\| {{{\bar y}^k}} \|}^2}} ] \hfill \\
   &- ({d_1}\frac{{{\mu _g}\lambda }}{2} - 4C_{g,x\theta }^2({d_4}{r_z}\alpha  + {d_6}\frac{4}{{1 - \rho }}{\gamma ^2}))\frac{1}{m}\mathbb{E}[ {{{\| {{v^k} - {v^*}({{\bar x}^k})} \|}^2}} ] \hfill
\end{aligned}
\end{equation}
\vspace{-2cm}
\begin{equation}
\begin{aligned}
   &- ({d_2}{\omega _\theta }\beta  - {L_{fg,x}}({d_4}{r_z}\alpha  + {d_6}\frac{4}{{1 - \rho }}{\gamma ^2}))\frac{1}{m}\mathbb{E}[ {{{\| {{\theta ^k} - {\theta ^*}({{\bar x}^k})} \|}^2}} ] \hfill \\
   &- ({d_5}\tau \frac{{1 - \rho }}{4} - {L_{fg,x}}({d_4}{r_z}\alpha  + {d_6}\frac{4}{{1 - \rho }}{\gamma ^2}))\frac{1}{m}\mathbb{E}[ {{{\| {{x^k} - {1_m} \otimes {{\bar x}^k}} \|}^2}} ] \hfill \\
   \quad \quad \quad \; \; &- ({d_4}{\gamma}-d_6\frac{4}{{1 - \rho }}{\res{\gamma^2} })\frac{1}{m}\mathbb{E}[ {{{\| {{s^k} - {z^k}} \|}^2}} ] \nonumber \\
   &+ \frac{1}{m}{d_3}\sigma _{\bar z}^2{\alpha ^2} + {d_1}\sigma _v^2{\alpha ^2} + {d_2}\sigma _\theta ^2{\alpha ^2} + {d_4}\sigma _z^2{\alpha ^2} + \frac{2d_6}{1-\rho}\sigma _y^2{\alpha ^2} \hfill  \\
   =& \mathbb{E}[{V^k}] - \frac{{1}}{2}{d_0}\tau \alpha \mathbb{E}[ {{{\| {\Phi ({{\bar x}^k})} \|}^2}} ] \hfill \\
    & - (\frac{{{d_0}}}{2}\tau \alpha \left( {1 - \tau \alpha L} \right) - {d_1}{q_s}{\tau ^2}{\alpha ^2} - {d_2}{p_s}{\tau ^2}{\alpha ^2} - {d_3}{r_y}{\tau ^2}{\alpha ^2} - {d_4}{r_y}{\tau ^2}{\alpha ^2})\mathbb{E}[ {{{\| {{{\bar y}^k}} \|}^2}} ] \\
     &- ({d_1}\frac{{{\mu _g}\lambda }}{2} - 4C_{g,x\theta }^2({d_4}{r_z}\alpha  + {d_6}\frac{4}{{1 - \rho }}{\gamma ^2}))\frac{1}{m}\mathbb{E}[ {{{\| {{v^k} - {v^*}({{\bar x}^k})} \|}^2}} ] \hfill \\
   &+2(\frac{1}{m}{d_3}\frac{{{\gamma ^2}}}{{{\alpha ^2}}} + {d_1}\frac{{{\lambda ^2}}}{{{\alpha ^2}}} + {d_4}\frac{{{\gamma ^2}}}{{{\alpha ^2}}} + \frac{2d_6}{1-\rho}\frac{{{\gamma ^2}}}{{{\alpha ^2}}})\sigma _{g,\theta \theta }^2{\alpha ^2}\frac{1}{m}\mathbb{E}[ {{{\| {{v^k} - {v^*}({{\bar x}^k})} \|}^2}} ]\\ \nonumber
   &- ({d_2}{\omega _\theta }\beta  - {L_{fg,x}}({d_4}{r_z}\alpha  + {d_6}\frac{4}{{1 - \rho }}{\gamma ^2}))\frac{1}{m}\alpha \mathbb{E}[ {{{\| {{\theta ^k} - {\theta ^*}({{\bar x}^k})} \|}^2}} ] \hfill \\
   &- ({d_5}\tau \frac{{1 - \rho }}{4} - {L_{fg,x}}({d_4}{r_z}\alpha  + {d_6}\frac{4}{{1 - \rho }}{\gamma ^2}))\frac{1}{m}\mathbb{E}[ {{{\| {{x^k} - {1_m} \otimes {{\bar x}^k}} \|}^2}} ] \hfill \\
   &+ (\frac{1}{m}{d_3} + {d_4} + \frac{2d_6}{1-\rho})(\sigma _{f,x}^2 + 2{M^2}\sigma _{g,x\theta }^2)\frac{{{\gamma ^2}}}{{{\alpha ^2}}}{\alpha ^2} \\
   &+ 2{d_1}\left( {\sigma _{f,\theta }^2 + 2{M^2}\sigma _{g,\theta \theta }^2} \right)\frac{{{\lambda ^2}}}{{{\alpha ^2}}}{\alpha ^2} + 2{d_2}\sigma _{g,\theta }^2\frac{{{\beta ^2}}}{{{\alpha ^2}}}{\alpha ^2},
 \end{aligned}
\end{equation}
where  the last equality uses the fact that $\sigma_{\bar z}^2=\sigma_{z}^2=\sigma_{y}^2$ in Lemmas \ref{LE-HypergradientE}, \ref{LE-stochastic-error-2}, \ref{LE-TE} and the definitions of $\sigma_{\bar z}^2, \sigma_v^2, \sigma_{\theta}^2$ in Lemmas \ref{LE-HypergradientE}, \ref{LE-Vstar}, \ref{LE-ThetaStar}.
Furthermore, when the following conditions hold:
\begin{align}
 \frac{{{d_0}}}{2}\tau \alpha \left( {1 - \tau \alpha L} \right) - {d_1}{q_s}{\tau ^2}{\alpha ^2} - {d_2}{p_s}{\tau ^2}{\alpha ^2} - {d_3}{r_y}{\tau ^2}{\alpha ^2} - {d_4}{r_y}{\tau ^2}{\alpha ^2} \geqslant& 0, \label{EQ-the2-C1} \hfill \\
 {d_1}\frac{{{\mu _g}\lambda }}{2} \!- 4C_{g,x\theta }^2({d_4}{r_z}\alpha \! +\! {d_6}\frac{4}{{1 - \rho }}{\gamma ^2}) \!- \!2(\frac{1}{m}{d_3}\frac{{{\gamma ^2}}}{{{\alpha ^2}}} \!+ \!{d_1}\frac{{{\lambda ^2}}}{{{\alpha ^2}}} \!+ \!{d_4}\frac{{{\gamma ^2}}}{{{\alpha ^2}}} \!+ \!\frac{2d_6}{1-\rho}\frac{{{\gamma ^2}}}{{{\alpha ^2}}})\sigma _{g,\theta \theta }^2{\alpha ^2} \geqslant &0, \label{EQ-the2-C2}\hfill \\
  {d_2}{\omega _\theta }\beta  - {L_{fg,x}}({d_4}{r_z}\alpha  + {d_6}\frac{4}{{1 - \rho }}{\gamma ^2}) \geqslant& 0, \label{EQ-the2-C3}\hfill \\
  {d_5}\tau \frac{{1 - \rho }}{4} - {L_{fg,x}}({d_4}{r_z}\alpha  + {d_6}\frac{4}{{1 - \rho }}{\gamma ^2}) \geqslant& 0, \label{EQ-the2-C4}
\end{align}
we  derive that:
\begin{equation}\label{EQ-inq-1}
\begin{aligned}
\mathbb{E}[ {{V^{k + 1}}} ] \leqslant&  \mathbb{E}[ {{V^k}} ] - \frac{{{d_0}}}{2}\tau\alpha \mathbb{E}[ {{{\| {\nabla \Phi ({{\bar x}^k})} \|}^2}} ] \\
& + (\frac{1}{m}{d_3} \!+\! {d_4}\!+\! \frac{2d_6}{1-\rho})(\sigma _{f,x}^2 + 2{M^2}\sigma _{g,x\theta }^2)\frac{{{\gamma ^2}}}{{{\alpha ^2}}}{\alpha ^2} \! \\
&+\ 2{d_1}\left( {\sigma _{f,\theta }^2 \!+ \!2{M^2}\sigma _{g,\theta \theta }^2} \right)\frac{{{\lambda ^2}}}{{{\alpha ^2}}}{\alpha ^2} \!+ \!2{d_2}\sigma _{g,\theta }^2\frac{{{\beta ^2}}}{{{\alpha ^2}}}{\alpha ^2},
\end{aligned}
\end{equation}
where we recall that the  coefficients $d_0$, $d_1$, $d_2$, $d_3$, $d_4$, $d_5$, $d_6$  of the  Lyapunov function \eqref{EQ-V} are given by:
\begin{equation}\label{EQ-the2-dddd2}
\begin{aligned}
  & {d_0} = 1,{d_1} = \frac{{8C_{g,x\theta }^2\tau \alpha }}{{{\mu _g}\lambda }},{d_2} = (\frac{{8C_{g,x\theta }^2\tau \alpha }}{{{\mu _g}\lambda }}{q_x} + {L_{fg,x}}\tau )\frac{\alpha }{{{\omega _\theta }\beta }},{d_3} = \frac{{\tau \alpha }}{{2\gamma }}, \hfill \\
& {d_4} = \left( {{L_{fg,x}}{d_3}{r_z} + {d_1}{q_x} + {d_2}{p_x}} \right)\frac{{64\gamma {\alpha ^2}}}{{{{(1 - \rho )}^4}}},
  {d_5} = 2\left( {{L_{fg,x}}{d_3}{r_z} + {d_1}{q_x} + {d_2}{p_x}} \right)\frac{{2\alpha }}{{\tau (1 - \rho )}}, \\
  &{d_6} = \left( {{L_{fg,x}}{d_3}{r_z} + {d_1}{q_x} + {d_2}{p_x}} \right)\frac{{16{\alpha ^2}}}{{{{(1 - \rho )}^3}}}. \hfill
\end{aligned}
\end{equation}
Next, we proceed to find the sufficient conditions for the step-sizes to make the conditions \eqref{EQ-the2-C1}-\eqref{EQ-the2-C4} hold.   To this end, we first
recall that like \eqref{EQ-varphi} the term ${{L_{fg,x}}{d_3}{r_z} + {d_1}{q_x} + {d_2}{p_x}}$ in $d_4$ can be simplified as follows:
\begin{equation}\label{EQ-varphi-the2}
\begin{aligned}
{{L_{fg,x}}{d_3}{r_z} + {d_1}{q_x} + {d_2}{p_x}}
   =\underbrace{({L_{fg,x}} + \frac{{32C_{g,x\theta }^2{L_{fg,\theta }}}}{{{\mu _g ^2}}})(1 + \frac{{4L_{g,\theta }^2}}{{\omega _\theta ^2}}))}_{\triangleq \varphi  }\tau.
\end{aligned}
\end{equation}
\nycres{Then,  $d_2$, $d_4$, $d_5$ and $d_6$ can be simplified as follows:}
\begin{align}
{d_2} = ({L_{fg,x}} + \frac{{32C_{g,x\theta }^2{L_{fg,\theta }}}}{{\mu _g^2}})\frac{{\tau \alpha }}{{{\omega _\theta }\beta }},  {d_4} = \frac{{64\varphi \tau \gamma {\alpha ^2}}}{{{{(1 - \rho )}^4}}}, {d_5} = \frac{{4\varphi \alpha }}{{(1 - \rho )}}, {d_6} = \frac{{16\varphi \tau {\alpha ^2}}}{{{{(1 - \rho )}^3}}}.   \nonumber
\end{align}
Now, we can derive  that the condition \eqref{EQ-the2-C1} holds if
\begin{align}
    \alpha  \leqslant & {u'_1} \triangleq \min \Big\{ {\frac{1}{{2\tau L}},\frac{{{{(1 - \rho )}^{2}}}}{{32{\varphi ^{1/2}}(\tau L)^{1/2} }}} \Big\}, \label{EQ-u3-1} \\
  \lambda  \geqslant& \frac{{48{L_{{v^*}}}{C_{g,x\theta }}}}{{{\mu _g}}}\tau \alpha,  \label{EQ-u3-2} \\
  \beta  \geqslant &\frac{{6{L_{{\theta ^*}}}}}{\omega_{\theta} }{(\frac{{32C_{g,x\theta }^2{L_{fg,\theta }}}}{{\mu _g^2}} + {L_{fg,x}})^{1/2}}\tau \alpha , \label{EQ-u3-3} \\
  \gamma  \geqslant& 4L\tau \alpha, \label{EQ-u3-4}
\end{align}
where $\varphi$ is  given by  \eqref{EQ-varphi-the2} and $0<\tau<1$.
Besides, based on the condition  \eqref{EQ-u3-4}, a sufficient condition to make the inequalities \eqref{EQ-the2-C2}-\eqref{EQ-the2-C4} hold is:
\begin{equation}\label{EQ-u3-gamma}
\begin{aligned}
  \gamma  \leqslant & {u' _2} = \min \Big\{ \frac{4}{5}\frac{{C_{g,x\theta }^2}}{{\sigma _{g,\theta \theta }^2}},\frac{{{{(1 - \rho )}^{4/3}}{{(\tau L)}^{1/3}}}}{{6{\varphi ^{1/3}}}},\frac{{{{(1 - \rho )}^{4/3}}{{(\tau L)}^{1/3}}C_{g,x\theta }^{2/3}}}{{3{\varphi ^{1/3}}\sigma _{g,\theta \theta }^{2/3}}},  \\
  &\quad\quad\quad \frac{{{{(1 - \rho )}^{4/3}}{{(L\tau )}^{1/3}}}}{{4{\varphi ^{1/3}}L_{fg,x}^{1/3}}}{(\frac{{32C_{g,x\theta }^2{L_{fg,\theta }}}}{{\mu _g^2}} + {L_{fg,x}})^{1/3}},\frac{{{{(1 - \rho )}^{4/3}}{{(L\tau )}^{1/3}}}}{{4L_{fg,x}^{1/3}}}\Big\},  \hfill \\
  \lambda  \leqslant &\frac{{{\mu _g}}}{{20\sigma _{g,\theta \theta }^2}}. \hfill \\
\end{aligned}
\end{equation}
\nycres{By combining Lemmas \ref{LE-descent}-\ref{LE-ThetaStar} and the conditions \eqref{EQ-u3-2}-\eqref{EQ-u3-gamma}, we  have  the following condition for the step-sizes $\lambda, \beta,  \gamma$:}
\begin{align}
 \min\ \left\{\frac{1}{{{\mu _g}}},\frac{{{\mu _g}}}{{20\sigma _{g,\theta \theta }^2}}\right\} >  \lambda  \geqslant& \frac{{48{L_{{v^*}}}{C_{g,x\theta }}}}{{{\mu _g}}}\tau \alpha,  \label{EQ-u3-2-upper} \\
  \min \left\{ {\frac{2}{{{\mu _g} + {L_{g,\theta }}}},\frac{{{\mu _g} + {L_{g,\theta }}}}{{2{\mu _g}{L_{g,\theta }}}}} \right\} > \beta  \geqslant &\frac{{6{L_{{\theta ^*}}}}}{\omega_{\theta} }{(\frac{{32C_{g,x\theta }^2{L_{fg,\theta }}}}{{\mu _g^2}} + {L_{fg,x}})^{1/2}}\tau \alpha , \label{EQ-u3-3-upper} \\
 \min\{1, u'_2\}> \gamma  \geqslant& 4L\tau \alpha, \label{EQ-u3-4-upper}
\end{align}
where the selection for $\lambda$, $\beta$ and $\gamma$ can be guaranteed to be nonempty, if  the step-size $\alpha$ satisfy the following condition:
\begin{align} \label{EQ-u3-3-3}
\alpha \leqslant \frac{\min \{1,u'_2\}}{8L\tau}, \alpha \leqslant u'_3\triangleq \min \Big\{ \frac{\min \{1,\frac{\mu _{g}^{2}}{8\sigma _{g,x\theta}^{2}}\}}{96\tau L_{v^*}C_{g,x\theta}},\frac{\min \left\{ \frac{2}{\mu _g+L_{g,\theta}},\frac{\mu _g+L_{g,\theta}}{2\mu _gL_{g,\theta}} \right\} w_{\theta}}{12\tau (\frac{32C_{g,x\theta}^{2}L_{fg,\theta}}{\mu _{g}^{2}}+L_{fg,x})^{1/2}L_{\theta ^*}} \Big\}.
\end{align}
\nycres{As a result, to ensure that the conditions \eqref{EQ-the2-C1}-\eqref{EQ-the2-C4} hold, with $0<\tau<1$,  we can properly select the step-size $\alpha$ as follows:}
\begin{equation}\label{EQ-condition-step-size-2}
\begin{aligned}
\alpha \leqslant u' \triangleq   \min \left\{ {{u'_1},\frac{\min \{1,u'_2\}}{8\tau L}}, u'_3 \right\}.
\end{aligned}
\end{equation}
\nycres{where $u'_1$, $u'_2$, and $u'_3$ are defined by \eqref{EQ-u3-1}, \eqref{EQ-u3-gamma}, and \eqref{EQ-u3-3}, respectively, and are clearly independent of $\lambda$, $\beta$, and $\gamma$. Subsequently, we can determine the step-sizes $\lambda$, $\beta$ and $\gamma$ according to \eqref{EQ-u3-2-upper}-\eqref{EQ-u3-4-upper}.} Therefore, under the condition \eqref{EQ-u3-2-upper}-\eqref{EQ-condition-step-size-2},  it holds that
\begin{equation}\label{EQ-inq}
\begin{aligned}
\frac{\mathbb{E}[ {{V^{k + 1}}} ]}{\tau} \leqslant & \frac{\mathbb{E}[ {{V^{k }}} ]}{\tau} - \frac{{{d_0}}}{2}\alpha \mathbb{E}[ {{{\| {\nabla \Phi ({{\bar x}^k})} \|}^2}} ]+ {\alpha ^2}\sigma_{\sigmaGT}^2,
\end{aligned}
\end{equation}
where ${\sigma_{\sigmaGT} ^2}$ is denoted as:
\begin{equation}\label{EQ-inq-variance-2}
\begin{aligned}
\!\! \!\! \sigma _{\sigmaGT}^2 \triangleq &
 (\frac{1}{m}\frac{{d_3}}{\tau} +\! \frac{{d_4}}{\tau} \!+ \!\frac{{2d_6}}{\tau(1-\rho)})(\sigma _{f,x}^2 + 2{M^2}\sigma _{g,x\theta }^2)\frac{{{\gamma ^2}}}{{{\alpha ^2}}} \!+ \!2\frac{{d_1}}{\tau}\left( {\sigma _{f,\theta }^2 + 2{M^2}\sigma _{g,\theta \theta }^2} \right)\frac{{{\lambda ^2}}}{{{\alpha ^2}}} \!+\! 2\frac{{d_2}}{\tau}\sigma _{g,\theta }^2\frac{{{\beta ^2}}}{{{\alpha ^2}}}.
\end{aligned}
\end{equation}
where ${d_0} = 1$, ${d_1} = \frac{{8C_{g,x\theta }^2\tau \alpha }}{{{\mu _g}\lambda }}$, ${d_2} = ( {L_{fg,x}}+\frac{{32C_{g,x\theta }^2{L_{fg,\theta }}}}{{\mu _g^2}})\frac{{\tau \alpha }}{{{\omega _\theta }\beta }}$, ${d_3} = \frac{{\tau \alpha }}{{2\gamma }}$, ${d_4} = \frac{{64\varphi \tau \gamma {\alpha ^2}}}{{{{(1 - \rho )}^4}}}$, ${d_5} = \frac{{4\varphi \alpha }}{{(1 - \rho )}}$, $d_6=\frac{{16\varphi \tau {\alpha ^2}}}{{{{(1 - \rho )}^3}}}$ with  $\varphi={({L_{fg,x}} + \frac{{32C_{g,x\theta }^2{L_{fg,\theta }}}}{{{\mu _g^2}}})(1 + \frac{{4L_{g,\theta }^2}}{{\omega _\theta ^2}})}$.  In what follows, by summing up and telescoping the inequality \eqref{EQ-inq} from $k=0$ to $K$ and combining the fact   $d_0=1$, it follows that:
\begin{equation}\label{EQ-phik}
\begin{aligned}
\frac{1}{K+1}\sum\limits_{k = 0}^K  {{\mathbb{E}}[ {{{\| \nabla{\Phi ({{\bar x}^k})} \|}^2}} ]} \leqslant \frac{{2\left( {{V^0} - {V^K}} \right)}}{{ \tau\alpha (K+1)}} + {2{\alpha {\sigma_{\sigmaGT} ^2}}}.
 \end{aligned}
\end{equation}
This completes the proof.   {\hfill $\blacksquare$}

\subsection{Proof of Corollary \ref{CO-2}} \label{sec-proof-CO-2}
Similar to the proof of Corollary \ref{CO-1}, when the conditions \eqref{EQ-u3-2-upper}-\eqref{EQ-u3-4-upper} and \eqref{EQ-condition-step-size-2} hold and the step-sizes $\gamma, \lambda, \beta$ are taken as $\gamma  = {c_\gamma }\alpha ,\lambda  = {c_\lambda }\alpha ,\beta  = {c_\beta }\alpha$ with positive parameters ${c_\gamma }, {c_\lambda }, {c_\beta}$ being  independent of the terms $K$ and $1-\rho$
and satisfying
\begin{equation}\label{EQ-cccc2}
\begin{aligned}
 c_{\lambda} \geqslant \frac{{48{L_{{v^*}}}{C_{g,x\theta }}}}{{{\mu _g}}}\tau,
 {c_\beta } \geqslant \frac{{6{L_{{\theta ^*}}}}}{{{\omega _\theta }}}{(\frac{{32C_{g,x\theta }^2{L_{fg,\theta }}}}{{{\mu _g^2}}} + {L_{fg,x}})^{1/2}}\tau,
 {c_\gamma } \geqslant 4L\tau,
\end{aligned}
\end{equation}
it follows from \eqref{EQ-inq-variance-2} that ${{ d}_0} = 1,{{ d}_1} = \frac{{8C_{g,x\theta }^2\tau }}{{{\mu _g}{c_\lambda }}},{d_2} = (\frac{{32C_{g,x\theta }^2{L_{fg,\theta }}}}{{{\mu _g ^2}}} + {L_{fg,x}})\frac{\tau }{{{\omega _\theta }{c_\beta }}}$, ${ d_3} = \frac{\tau }{{2{c_\gamma }}}$ are independent of the step-size $\alpha$ and $1-\rho$. When the parameters ${c_\gamma }, {c_\lambda }, {c_\beta}$ satisfy \eqref{EQ-cccc2} and their values are properly selected,   the conditions \eqref{EQ-u3-2-upper}-\eqref{EQ-u3-4-upper} always hold. For examples,  a sufficient  condition for the upper bounds of the parameters  ${c_\gamma }, {c_\lambda }, {c_\beta}$ could be
$ \frac{96L_{v^*}C_{g,x\theta}}{\mu _g}\tau
 > c_{\lambda}$, $\frac{12L_{\theta ^*}}{\omega _{\theta}}(\frac{32C_{g,x\theta}^{2}L_{fg,\theta}}{\mu _{g}^{2}}+L_{fg,x})^{1/2}\tau>c_{\beta}$, and $8L\tau \alpha >c_{\gamma}$.
In addition, from \eqref{EQ-the2-dddd2} we note that ${d_4} = O({\gamma ^2}\alpha {d_3}),{d_6} = O(\gamma \alpha {d_3})$, which ensures that ${d_4}  \ll  \frac{1}{m}{d_3},{d_6}  \ll  \frac{1}{m}{d_3}$ under a large number of iterations $K$. Thus,  the variance related term $\sigma_{\sigmaGT}$ in \eqref{EQ-inq-variance-2} can be further derived as:
\begin{equation}\label{EQ-sigma-r-dot}
\begin{aligned}
  {\sigma _{\sigmaGT}}
   \leqslant&  {\sqrt {\frac{4}{m}\frac{{{d_3}}}{\tau }(\sigma _{f,x}^2 + 2{M^2}\sigma _{g,x\theta }^2)\frac{{{\gamma ^2}}}{{{\alpha ^2}}} + 2\frac{{{d_1}}}{\tau }\left( {\sigma _{f,\theta }^2 + 2{M^2}\sigma _{g,\theta \theta }^2} \right)\frac{{{\lambda ^2}}}{{{\alpha ^2}}} + 2\frac{{{d_2}}}{\tau }\sigma _{g,\theta }^2\frac{{{\beta ^2}}}{{{\alpha ^2}}}} }\\
   \leqslant& \sqrt{2} {\underbrace {{(\sqrt {\frac{d_1}{\tau}} (\sigma _{f,\theta} + 2{M}\sigma _{g,\theta\theta }){c_\lambda} + \sqrt {\frac{d_2}{\tau}} \sigma _{g,\theta }{c_\beta})}}_{ \triangleq \sigma _{\operatorname{p}} }} + \frac{2}{{\sqrt m }}\underbrace{\sqrt {\frac{d_3}{\tau}}  (\sigma _{f,x} +2 {M}\sigma _{g,x\theta }){c_\gamma}}_{ \triangleq {\sigma _c}}
   \triangleq {{\hat \sigma }_{\sigmaGT}},
 \end{aligned}
\end{equation}
where the last inequality  the  triangle inequality. By \eqref{EQ-sigma-detail}, it is known that ${\sigma _{\operatorname{p} }}  = \mathcal{O}({\sigma _{f,\theta }} + {\sigma _{g,\theta \theta }} + {\sigma _{g,\theta }})$ and ${\sigma _{\operatorname{c} }}  = \mathcal{O}({\sigma _{f,x}} + {\sigma _{g,x\theta }})$. By  \eqref{EQ-u3-1},
\eqref{EQ-u3-gamma}, \eqref{EQ-u3-3}, and \eqref{EQ-condition-step-size-2}, we have that $u'=\mathcal{O}({(1 - \rho )^2})$. Then, by the this fact and the condition that $\gamma<1$ in \eqref{EQ-u3-4-upper},  the coefficients $d_4$, $d_5$, $d_6$ in \eqref{EQ-inq-variance-2} can be bounded as follows:
\begin{equation}
\begin{aligned}
  {d_4} =& \frac{{64\varphi \tau \gamma {\alpha ^2}}}{{{{(1 - \rho )}^4}}} \leqslant \frac{{64\varphi \tau {{(u')}^2}}}{{{{(1 - \rho )}^4}}} \triangleq {{\hat d}_4} = O(1), \\
  {d_5} =&  \frac{{4\varphi  \alpha }}{{(1 - \rho )}} \leqslant \frac{{4\varphi  u'}}{{(1 - \rho )}} \triangleq {{\hat d}_5} = O(1 - \rho ), \\
  {d_6} =& \frac{{16\varphi \tau {\alpha ^2}}}{{{{(1 - \rho )}^3}}} \leqslant \frac{{16\varphi \tau {(u')^2}}}{{{{(1 - \rho )}^3}}} \triangleq {{\hat d}_6} = O(1 - \rho ). \\
\end{aligned}
\end{equation}
Then, for the term $\frac{{{V^0} - {V^K}}}{\tau }$ in \eqref{EQ-phik}, it follows that
\begin{equation}\label{EQ-the2-vvvv}
\begin{aligned}
\frac{{{V^0} - {V^K}}}{\tau } \leqslant & \frac{{{V^0}}}{\tau } \\
   \leqslant& \frac{1}{\tau }{{ d}_0}\Phi ({{\bar x}^0}) + \frac{1}{\tau }{{ d}_1}\frac{1}{m}{\left\| {{v^0} - {v^*}({{\bar x}^0})} \right\|^2}\! +\! \frac{1}{\tau }{{ d}_2}\frac{1}{m}{\left\| {{\theta ^0} \!-\! {\theta ^*}({{\bar x}^0})} \right\|^2} \!+\! \frac{1}{\tau }{{ d}_3}{\left\| {\nabla \Phi ({{\bar x}^0}) - {{\bar z}^0}} \right\|^2} \\
  &+ \frac{1}{\tau }{{\hat d}_4}\frac{1}{m}{\left\| {\nabla \tilde \Phi ({{\bar x}^0}) - {z^0}} \right\|^2} + \frac{1}{\tau }{{\hat d}_5}\frac{1}{m}{\left\| {{x^0} - {1_m} \otimes {{\bar x}^0}} \right\|^2} + \frac{1}{\tau }{{\hat d}_6}\frac{1}{m}{\left\| {{y^0} - {1_m} \otimes {{\bar y}^0}} \right\|^2} \\
   \triangleq &{{\hat V}^0} = O(1 - \rho ). \\
\end{aligned}
\end{equation}
In what follows, letting
\begin{equation}
\begin{aligned}
{a'_0} \triangleq 2{{\hat V}^0},{a'_1} \triangleq 2{\hat \sigma_{\sigmaGT}^2 },
\end{aligned}
\end{equation}
and combining the inequality \eqref{EQ-sigma-r-dot}, the inequality \eqref{EQ-phik} can be further derived as:
\nyc{
\begin{equation}
\begin{aligned}\frac{1}{{K + 1}}\sum\limits_{k = 0}^K {{\mathbb{E}}[ {{{\| \nabla{\Phi ({{\bar x}^k})} \|}^2}} ]}  \leqslant {a' _0}\frac{1}{{\alpha (K + 1)}} + {a' _1}\alpha.
\end{aligned}
\end{equation}
}

When we take the step-size $\alpha$ as $\alpha  = \min \Big\{ {u',{{( {\frac{{{a'_0}}}{{{a'_1}(K + 1)}}} )}^{\frac{1}{2}}}} \Big\}$ and  the step-sizes $\gamma, \lambda, \beta$  as $\gamma  = {c_\gamma }\alpha ,\lambda  = {c_\lambda }\alpha ,\beta  = {c_\beta }\alpha$ with ${c_\gamma }, {c_\lambda }, {c_\beta}$ being independent of the terms $K$ and $1-\rho$, we have:
\begin{equation}\label{EQ-a0-a1-dot}
\begin{aligned}
\frac{1}{{K + 1}}\sum\limits_{k = 0}^K {{\mathbb{E}}[ {{{\| {\Phi ({{\bar x}^k})} \|}^2}} ]}  \leqslant &\frac{{{a'_0}}}{{u'(K + 1)}} + 2{( {\frac{{{a'_1}{a'_0}}}{{K + 1}}} )^{\frac{1}{2}}}\\
\leqslant  &\frac{{2{\hat V}^0}}{{ u'(K + 1)}} + \frac{{4{\hat \sigma _{\sigmaGT}} \sqrt { {{\hat V}^0}} }}{{ \sqrt {K + 1} }}.
 \end{aligned}
\end{equation}
Then, we proceed to analyze how the parameters in \eqref{EQ-a0-a2} depend on the terms $\kappa$ and $1-\rho$. Firstly, combining the fact that $\varphi  = O({\kappa ^6})$ by \eqref{EQ-varphi-the2}, we set  $\tau  = \mathcal{O}({\kappa ^{ - 4}}),u' = \mathcal{O}({\kappa ^{ - 3}}(1-\rho)^2)$ and $c_\lambda=\mathcal{O}(1), c_\beta =\mathcal{O}(1), c_\gamma =\mathcal{O}(1)$ in terms of $\kappa$ and $1-\rho$. Furthermore, when we initialize the outer-level variables as  $x_i^0=x_j^0, \forall i,j \in {{\mathcal{V}}}$, i.e., $\|{{x^0} - {1_m} \otimes {{\bar x}^0}}\|^2=0$,  we can derive that ${{\hat V}^0} = O({\kappa ^5})$, and ${{\hat V}^0}$ is independent of the term $\frac{1}{1-\rho}$ by \eqref{EQ-the2-vvvv}. Thus, with the result  $\hat \sigma _{\sigmaGT} = {{\mathcal{O}}}( {{{\sigma_{\operatorname{p}}  + \frac{1}{\sqrt{m}}\sigma_{\operatorname{c}}}}} )$ in \eqref{EQ-sigma-r-dot}, it follows from \eqref{EQ-a0-a1-dot}  that:
\begin{equation}\label{EQ-CO2-phi-kapper}
\begin{aligned}
\frac{1}{{K + 1}}\sum\limits_{k = 0}^K {{\mathbb{E}}[ {{{\| \nabla {\Phi ({{\bar x}^k})} \|}^2}} ]}  ={{\mathcal{O}}}\Big( {\frac{{{{{\kappa ^8}}}}}{{{{\left( {1 - \rho } \right)^2}}K}} + \frac{{ {\kappa ^{\frac{5}{2}}} }}{ {  {\sqrt{K}} }} (\sigma_{\operatorname{p}}  + \frac{1}{\sqrt{m}}\sigma_{\operatorname{c}})   } \Big).
 \end{aligned}
\end{equation}
where similar to  Corollary \ref{CO-1},  the terms $\sigma_{\operatorname{p}}$ and $\sigma_{\operatorname{c}}$ follow that ${\sigma _{\operatorname{p} }}={\mathcal{O}}(\kappa^{\frac{1}{2}}\sigma_{f,\theta}+\kappa^{\frac{3}{2}}\sigma_{g,\theta\theta}+\kappa^{\frac{5}{2}}\sigma_{g,\theta})$ and ${\sigma _{\operatorname{c} }}={\mathcal{O}}(\sigma_{f,x}+\kappa\sigma_{g,x\theta})$ by \eqref{EQ-sigma-r-dot}.
This completes the proof.   {\hfill $\blacksquare$}

\section{Proof of Supporting Propositions} \label{sec-Pro}
\subsection{Proof of Proposition \ref{PR-smooth} }

\textbf{Lipschitz continuity of $\theta _i^*(x)$}. Recalling the definition of $\theta _i^*(x)$,  the expression of $\theta _i^*(x)$ is given by \citep{BSA}:
\[\nabla \theta _i^*(x) = -{\nabla _{ x\theta}^2}{g_i}\left( {x,\theta _i^*(x)} \right){\left[ {{\nabla _{\theta \theta }^2}{g_i}\left( {x,\theta _i^*(x)} \right)} \right]^{ - 1}}.\]
Following the strong convexity of $g_i$  and bounded Jacibian of $\nabla_{x\theta}^2g_i$ in  Assumption \ref{ASS-STOCHASTIC}, we have
\begin{equation}\label{EQ-vixo}
\begin{aligned}
\left\| {\nabla \theta _i^*(x)} \right\| = \left\| {{\nabla _{x\theta}^2}{g_i}\left( {x,\theta _i^*(x)} \right){{\left[ {{\nabla _{\theta \theta }^2}{g_i}\left( {x,\theta _i^*(x)} \right)} \right]}^{ - 1}}} \right\| \leqslant \frac{{{C_{{g}, x\theta}}}}{{{\mu _{{g}}}}},
 \end{aligned}
\end{equation}
which implies that for any $x$, $x'$:
\begin{equation}
\begin{aligned}
\left\| {\theta _i^*(x) - \theta _i^*(x')} \right\| \leqslant \underbrace{\frac{{{C_{g, x\theta}}}}{{{\mu _g}}}}_{\triangleq L_{\theta ^*}}\left\| {x - x'} \right\|.
 \end{aligned}
\end{equation}
\textbf{Lipschitz continuity of $v_i^*\left( x \right)$}.
Recalling the definition of ${v_i}\left( {x,\theta } \right)$, we known that ${v_i}\left( {x,\theta } \right)$ admits the following expression:
\[{v_i}(x,\theta ) = {\left[ {\nabla _{\theta \theta }^2{g_i}\left( {x,\theta } \right)} \right]^{ - 1}}{\nabla _x}{f_i}(x,\theta ).\]
Letting $x$ and $x'$ be any two points in $\mathbb{R}^n$ and considering the case $\theta=\theta_i^*(x)$ and $\theta'=\theta_i^*(x')$, it follows that
\begin{equation}\label{EQ-v-xo}
\begin{aligned}
  &\| {{v_i}(x,\theta ) - {v_i}(x',\theta ')} \| \hfill \\
   \mathop  \leqslant \limits^{({\rm a})} & \| {{{[ {\nabla _{\theta \theta }^2{g_i}( {x,\theta } )} ]}^{ - 1}}{\nabla _\theta }{f_i}(x,\theta ) - {{[ {\nabla _{\theta \theta }^2{g_i}( {x,\theta } )} ]}^{ - 1}}{\nabla _\theta }{f_i}(x',\theta ')} \| \hfill \\
   &+ \| {{{[ {\nabla _{\theta \theta }^2{g_i}\left( {x,\theta } \right)} ]}^{ - 1}}{\nabla _\theta }{f_i}(x',\theta ') - {{[ {\nabla _{\theta \theta }^2{g_i}( {x',\theta '} )} ]}^{ - 1}}{\nabla _\theta }{f_i}(x',\theta ')} \| \hfill \\
   \mathop  \leqslant \limits^{({\rm b})}  &\frac{1}{{{\mu _g}}}\| {{\nabla _\theta }{f_i}(x,\theta ) - {\nabla _\theta }{f_i}(x',\theta ')} \| \hfill \\
   &+ {C_{f,\theta }}\| {{{[ {\nabla _{\theta \theta }^2{g_i}( {x',\theta '} )} ]}^{ - 1}}( {\nabla _{\theta \theta }^2{g_i}( {x,\theta } ) - \nabla _{\theta \theta }^2{g_i}( {x',\theta '} )} ){{[ {\nabla _{\theta \theta }^2{g_i}( {x,\theta } )} ]}^{ - 1}}} \| \hfill \\
   \mathop  \leqslant \limits^{({\rm c})} & \frac{1}{{{\mu _g}}}{L_{f,\theta }}( {\| {x - x'} \| + \| {\theta  - \theta '} \|} ) + \frac{{{C_{f,\theta }}{L_{g,\theta \theta }}}}{{{\mu _g ^2}}}( {\| {x - x'} \| + \| {\theta  - \theta '} \|} ) \hfill \\
   =& \underbrace {( {\frac{{{L_{f,\theta }}}}{{{\mu _g}}} + \frac{{{C_{f,\theta }}{L_{g,\theta \theta }}}}{{{\mu _g ^2}}}} )}_{ \triangleq {L_v}}( {\| {x - x'} \| + \| {\theta  - \theta '} \|} ), \hfill \\
 \end{aligned}
\end{equation}
where step (a) follows from the triangle inequality; step (b) uses the upper bound of the Hessian-inverse matrix with the parameter $\frac{1}{\mu _g}$ related to the strong convexity constant as well as the boundness of $\nabla_{\theta}f_i(x,\theta_i^*(x))$ under
the case $\theta=\theta_i^*(x)$ and $\theta'=\theta_i^*(x')$ with  Assumption \ref{ASS-OUTLEVEL}; step (c) follows from   $\nabla_{\theta\theta}^2g_i$ in Assumption \ref{ASS-STOCHASTIC}.
Then, recalling the definition of $v_i^*\left( x \right)$, we know that $v_i^*\left( x \right)=v_i\left( x, \theta_i^*(x) \right)$. Thus, given any two $x$, $x'$, by taking $\theta=\theta_i^*(x)$ and $\theta'=\theta_i^*(x')$ in \eqref{EQ-v-xo}, it follows  that
\begin{equation}
\begin{aligned}
\left\| {v_i^*(x) - v_i^*(x')} \right\| \hfill  = &\left\| {{v_i}(x,\theta _i^*(x)) - {v_i}(x',\theta _i^*(x'))} \right\| \hfill \\
   \leqslant& {L_v}\left( {\left\| {x - x'} \right\| + \left\| {\theta _i^*(x) - \theta _i^*(x')} \right\|} \right) \hfill \\
   \leqslant& {L_v}\left( {1 + {L_{{\theta ^*}}}} \right)\left\| {x - x'} \right\|, \hfill \\
 \end{aligned}
\end{equation}
where the last inequality is obtained employing the Lipschitz continuity of $\theta _i^*(x)$ as mentioned  earlier.
\\
\\
\textbf{Lipschitz continuity of $\nabla \Phi (x)$}.
Letting $x$ and $x'$ be any two points in $\mathbb{R}^n$ as well as $\theta=\theta_i^*(x)$ and $\theta'=\theta_i^*(x')$ and following the definition of $\bar \nabla {f_i}\left( {x,\theta } \right)$ and  the triangle inequality yield that
\begin{equation}\label{EQ-bar-nablaf}
\begin{aligned}
  \| {\bar \nabla {f_i}( {x,\theta } ) - \bar \nabla {f_i}( {x',\theta '} )} \|
   \leqslant & \| {{\nabla _x}{f_i}( {x,\theta } ) - {\nabla _x}{f_i}( {x',\theta '} )} \| \hfill \\
   &+ \| {\nabla _{x\theta }^2{g_i}( {x,\theta } )( {{v_i}\left( {x,\theta } \right) - {v_i}( {x',\theta '} )} )} \| \hfill \\
   &+ \| {( {\nabla _{x\theta }^2{g_i}( {x,\theta } ) - \nabla _{x\theta }^2{g_i}( {x',\theta '} )} ){v_i}( {x',\theta '} )} \|. \hfill \\
 \end{aligned}
\end{equation}
Utilizing Lipschitz continuity of $\nabla_xf_i$, the boundness of $\nabla_{x\theta}^2g_i$ and Lipschitz continuity of $v_i(x,\theta)$, the first two terms on right hand of \eqref{EQ-bar-nablaf} can be bounded by:
\begin{equation}
\begin{aligned}
  &\left\| {{\nabla _x}{f_i}\left( {x,\theta } \right) - {\nabla _x}{f_i}\left( {x',\theta '} \right)} \right\| + \left\| {\nabla _{x\theta }^2{g_i}\left( {x,\theta } \right)\left( {{v_i}\left( {x,\theta } \right) - {v_i}\left( {x',\theta '} \right)} \right)} \right\| \hfill \\
   \leqslant & ({L_{f,x}} + {C_{g,x\theta }}{L_v})\left( {\left\| {x - x'} \right\| + \left\| {\theta  - \theta '} \right\|} \right). \hfill \\
\end{aligned}
\end{equation}
Note that for  the case $\theta=\theta_i^*(x)$ and $\theta'=\theta_i^*(x')$, we have $\left\| {{v_i}\left( {x',\theta '} \right)} \right\| \leqslant \frac{{{C_{f,\theta }}}}{{{\mu _g}}}$  by employing the  boundness of the
Hessian-inverse matrix $\nabla _{\theta\theta}^2 g_i$ and the gradient $\nabla_{\theta}f_i(x,\theta_i^*(x))$. Then, combining the Lipschitz continuity of $\nabla_{x\theta}^2g_i$, one can get
\[\left\| {\left( {\nabla _{x\theta }^2{g_i}\left( {x,\theta } \right) - \nabla _{x\theta }^2{g_i}\left( {x',\theta '} \right)} \right){v_i}\left( {x',\theta '} \right)} \right\| \leqslant \frac{{{C_{f,\theta }}{L_{g,x\theta }}}}{{{\mu _g}}}\left( {\left\| {x - x'} \right\| + \left\| {\theta  - \theta '} \right\|} \right).\]
Then the inequality $\left\| {\bar \nabla {f_i}\left( {x,\theta } \right) - \bar \nabla {f_i}\left( {x',\theta '} \right)} \right\| \leqslant {L_f}\left( {\left\| {x - x'} \right\| + \left\| {\theta  - \theta '} \right\|} \right)$
can be derived by integrating above inequalities, with $L_f$ defined in \eqref{EQ-Lip-const}.
It is noted from  the expression of $\nabla \Phi _i (x)$ that $\nabla \Phi_i(x)=\bar \nabla {f_i}\left( {x,\theta_i^*(x) } \right)$. Then, given any two $x$, $x'$, by taking $\theta=\theta_i^*(x)$ and $\theta'=\theta_i^*(x')$ in \eqref{EQ-bar-nablaf}, we can employ the result \eqref{EQ-bar-nablaf}  to derive the following inequality:
\begin{equation}
\begin{aligned}
  \left\| {\nabla {\Phi _i}(x) - \nabla {\Phi _i}(x')} \right\|
   \leqslant & ( {{L_{f,x}} + {C_{g,x\theta }}{L_v} + \frac{{{C_{f,\theta }}{L_{g,x\theta }}}}{{{\mu _g}}}} )\left( {\left\| {x - x'} \right\| + \left\| {\theta _i^*(x) - \theta _i^*(x')} \right\|} \right) \hfill \\
   \leqslant &\underbrace{( {{L_{f,x}} + {C_{g,x\theta }}{L_v} + \frac{{{C_{f,\theta }}{L_{g,x\theta }}}}{{{\mu _g}}}} )\left( {1 + {L_{{\theta ^*}}}} \right)}_{\triangleq L}\left\| {x - x'} \right\|,
\end{aligned}
\end{equation}
where the last step uses the Lipschitz  continuity of $v_i^*(x)$ as shown earlier. Combining the fact that $\Phi(x)=\frac{1}{m}\sum_{i=1}^{m} \Phi_i(x)$, we can derive that $\nabla \Phi(x)$ is also $L$-Lipschitz  continuous.  This completes the proof.   {\hfill $\blacksquare$}
\subsection{Proof of Proposition \ref{PR-boundness} }
\nyc{
By the strong convexity of $g_i$ in $\theta$ and the boundedness of $\nabla f_{\theta}(x,\theta_i^*(x))$, it follows that
$$\|{v_i}\left( {x,\theta^*_i(x) } \right)\| =\|{[ {\nabla _{\theta \theta }^2{g_i}\left( {x,\theta^*_i(x)  } \right)} ]^{ - 1}}{\nabla _\theta }{f_i}\left( {x,\theta^*_i(x) } \right)\| \leqslant \frac{C_{f,\theta}}{\mu_g}.$$}
This completes the proof.   {\hfill $\blacksquare$}

\section{Proof of Supporting Lemmas} \label{sec-supporting-lemma}

\subsection{Proof  of Lemma \ref{LE-hyper-heterogeneity}} \label{sec-poof-heterogeneity}
By the definition of $\Phi_i(x)$, we can compute its gradient as:
\[\nabla {\Phi _i}\left( x \right) = {\nabla _x}{f_i}(x,\theta _i^*(x)) + \nabla \theta _i^*(x){\nabla _\theta }{f_i}(x,\theta _i^*(x)).\]
Then the term $\sum\nolimits_{i = 1}^m {{{\left\| {\nabla {\Phi _i}\left( x \right) - \nabla \Phi \left( x \right)} \right\|}^2}}$ can be bounded as:
\begin{equation}\label{EQ-phi-phi}
\begin{aligned}
   \sum\limits_{i = 1}^m {{{\left\| {\nabla {\Phi _i}\left( x \right) - \nabla \Phi \left( x \right)} \right\|}^2}}   \leqslant
    & 2\sum\limits_{i = 1}^m {{{\big\| {\nabla \theta _i^*(x){\nabla _{\theta}}{f_i}(x,\theta _i^*(x)) - \frac{1}{m}\sum\limits_{j = 1}^m {\nabla \theta _j^*(x){\nabla _\theta }{f_j}(x,\theta _j^*(x))} } \big\|^2}}} \\
     &+ 2\sum\limits_{i = 1}^m {{{\big\| {{\nabla _x}{f_i}(x,\theta _i^*(x)) - \frac{1}{m}\sum\limits_{j = 1}^m {{\nabla _x}{f_j}(x,\theta _j^*(x))} } \big\|^2}}} .
\end{aligned}
\end{equation}
For the first term on the right hand of \eqref{EQ-phi-phi}, it follows that

\begin{align}
  &\sum\limits_{i = 1}^m {{{\big\| {\nabla \theta _i^*(x){\nabla _\theta }{f_i}(x,\theta _i^*(x)) - \frac{1}{m}\sum\limits_{j = 1}^m {\nabla \theta _j^*(x){\nabla _\theta }{f_j}(x,\theta _j^*(x))} } \big\|^2}}} \nonumber \hfill \\
    \mathop  \leqslant \limits^{({\rm a})} & 2\sum\limits_{i = 1}^m {{{\big\| {\nabla \theta _i^*(x)}{{\nabla _\theta }{f_i}(x,\theta _i^*(x)) - \frac{1}{m}\sum\limits_{j = 1}^m {\nabla \theta _i^*(x)}{{\nabla _\theta }{f_j}(x,\theta _i^*(x))} } \big\|^2}}} \nonumber \hfill \\
   &+ 2\sum\limits_{i = 1}^m {{{\big\| {\frac{1}{m}\sum\limits_{j = 1}^m {\left[ {\nabla \theta _j^*(x){\nabla _\theta }{f_j}(x,\theta _j^*(x)) - \nabla \theta _i^*(x){\nabla _\theta }{f_j}(x,\theta _i^*(x))} \right]} } \big\|^2}}}   \label{EQ-phi-phi-1}\hfill  \\
    \mathop  \leqslant \limits^{({\rm b})} & 2\max _i \left\{ {{{\left\| {\nabla \theta _i^*(x)} \right\|^2}}} \right\}b_{{f}}^2
   + 2\sum\limits_{i = 1}^m {{{\big\| {\frac{1}{m}\sum\limits_{j = 1}^m {\left[ {\nabla \theta _j^*(x){\nabla _\theta }{f_j}(x,\theta _j^*(x)) - \nabla \theta _i^*(x){\nabla _\theta }{f_j}(x,\theta _i^*(x))} \right]} } \big\|^2}}}, \nonumber
\end{align}
where step (a) first introduces the term $\frac{1}{m}\sum\nolimits_{i = 1}^m{\nabla \theta _i^*(x){\nabla _\theta }{f_j}(x,\theta _i^*(x))}$ and uses Young's inequality and  Cauchy-Schwartz inequality;   step (b) use the bounded data heterogeneity   in Assumption \ref{ASS-heterogeneity}. In what follows, we proceed in providing the upper bound for the last term on the right hand of \eqref{EQ-phi-phi-1} as follows:
\begin{align}
  &\sum\limits_{i = 1}^m {{{\big\| {\frac{1}{m}\sum\limits_{j = 1}^m {\left[ {\nabla \theta _j^*(x){\nabla _\theta }{f_j}(x,\theta _j^*(x)) - \nabla \theta _i^*(x){\nabla _\theta }{f_j}(x,\theta _i^*(x))} \right]} } \big\|^2}}}  \hfill \nonumber\\
    \mathop  \leqslant \limits^{({\rm a})} & 2\frac{1}{m}\sum\limits_{i = 1}^m {\sum\limits_{j = 1}^m {{{\left\| {\nabla \theta _i^*(x)} \right\|}^2}{{\left\| {{\nabla _\theta }{f_j}(x,\theta _j^*(x)) - {\nabla _\theta }{f_j}(x,\theta _i^*(x))} \right\|}^2}} }  \nonumber \hfill \\
   &+ 2\frac{1}{m}\sum\limits_{i = 1}^m {\sum\limits_{j = 1}^m {{{\left\| {{\nabla _\theta }{f_j}(x,\theta _j^*(x))} \right\|}^2}{{\left\| {\nabla \theta _j^*(x) - \nabla \theta _i^*(x)} \right\|}^2}} } \label{EQ-phi-phi-1-2}  \hfill \\
   \mathop  \leqslant \limits^{({\rm b})} &2\frac{1}{m}\sum\limits_{i = 1}^m {\sum\limits_{j = 1}^m {L_{{f},\theta }^2{{\left\| {\nabla \theta _i^*(x)} \right\|}^2}{{\left\| {\theta _j^*(x) - \theta _i^*(x)} \right\|}^2}} } + 2\frac{1}{m}\sum\limits_{i = 1}^m {\sum\limits_{j = 1}^m {C_{{f},\theta }^2{{\left\| {\nabla \theta _j^*(x) - \nabla \theta _i^*(x)} \right\|}^2}} }  \nonumber \hfill \\
   \leqslant&2{\max _i}\left\{ {L_{{f},\theta }^2{{\left\| {\nabla \theta _i^*(x)} \right\|}^2} } \right\}\frac{1}{m}\sum\limits_{i = 1}^m {\sum\limits_{j = 1}^m {{{\left\| {\theta _j^*(x) - \theta _i^*(x)} \right\|}^2}} } +  2\frac{1}{m}\sum\limits_{i = 1}^m {\sum\limits_{j = 1}^m {C_{{f},\theta }^2{{\left\| {\nabla \theta _j^*(x) - \nabla \theta _i^*(x)} \right\|}^2}} }, \nonumber
\end{align}
where step (a) first combines the term $\nabla \theta _i^*(x){\nabla _\theta }{f_j}(x,\theta _j^*(x))$ and Young’s inequality, and then follows from Jensen inequality and Cauchy-Schwartz inequality; step (b) employs Lipschitz continuity of $\nabla _{\theta} f_i$  and the boundness of $\nabla _{\theta} f_i(x,\theta_i^*(x))$ in Assumption \ref{ASS-OUTLEVEL}. For the last term of \eqref{EQ-phi-phi-1-2}, Noting that $\nabla \theta _i^*(x) =  - \nabla _{x\theta }^2{g_i}\left( {x,\theta _i^*(x)} \right){\left[ {{\nabla _{\theta \theta }}{g_i}\left( {x,\theta _i^*(x)} \right)} \right]^{ - 1}},$ we arrive at
\begin{align}
 & {\| {\nabla \theta _j^*(x) - \nabla \theta _i^*(x)} \|^2} \nonumber \\
   \leqslant& 2{\| {\nabla _{x\theta }^2{g_j}\!\left( {x,\theta _j^*(x)} \right)\!{{[ {{\nabla _{\theta \theta }}{g_j}\!( {x,\theta _j^*(x)} )} ]}^{ - 1}}\!\left[ {\nabla _{\theta \theta }^2{g_i}\left( {x,\theta _i^*(x)} \right) \!- \!\nabla _{\theta \theta }^2{g_j}\!\left( {x,\theta _j^*(x)} \right)} \right]\!{{[ {\nabla _{\theta \theta }^2{g_i}\!( {x,\theta _i^*(x)} )}\! ]}^{ - 1}}} \!\|^2} \nonumber \\
   &+ 2{\| {[ {\nabla _{x\theta }^2{g_j}\left( {x,\theta _j^*(x)} \right) - \nabla _{x\theta }^2{g_i}\left( {x,\theta _i^*(x)} \right)} ]{{[ {\nabla _{\theta \theta }^2{g_i}\left( {x,\theta _i^*(x)} \right)} ]}^{ - 1}}} \|^2} \nonumber \\
   \leqslant&  \yc {2\frac{{C_{g,x\theta }^2}}{{\mu _g^4}}}{\| {\nabla _{\theta \theta }^2{g_i}\left( {x,\theta _i^*(x)} \right) - \nabla _{\theta \theta }^2{g_j}\left( {x,\theta _j^*(x)} \right)} \|^2} + 2\frac{1}{{\mu _g^2}}{\| {\nabla _{x\theta }^2{g_i}\left( {x,\theta _i^*(x)} \right) - \nabla _{x\theta }^2{g_j}\left( {x,\theta _j^*(x)} \right)} \|^2} \nonumber \\
  \leqslant &( {4\frac{{C_{g,x\theta }^2L_{g,\theta \theta }^2}}{{\mu _g^4}} + 4\frac{{L_{g,x\theta }^2}}{{\mu _g^2}}} ){\| {\theta _i^*(x) - \theta _j^*(x)} \|^2} + 4\frac{{C_{g,x\theta }^2}}{{\mu _g^4}}{\| {\nabla _{\theta \theta }^2{g_i}( {x,\theta _j^*(x)} ) - \nabla _{\theta \theta }^2{g_j}( {x,\theta _j^*(x)} )} \|^2} \nonumber \\
   &+ 4\frac{1}{{\mu _g^2}}{\| {\nabla _{x\theta }^2{g_i}( {x,\theta _j^*(x)} ) - \nabla _{x\theta }^2{g_j}( {x,\theta _j^*(x)} )} \|^2}, \label{EQ-nablaphi-nablaphi}
\end{align}
which implies that  the last term of \eqref{EQ-phi-phi-1-2} can be bounded by:
\begin{align}
  &2C_{f,\theta }^2\frac{1}{m}\sum\limits_{i = 1}^m {\sum\limits_{j = 1}^m {{{\left\| {\nabla \theta _j^*(x) - \nabla \theta _i^*(x)} \right\|}^2}} }  \label{EQ-nablaphi-nablaphi-1}  \\
   \leqslant& ( {8\frac{{C_{f,\theta }^2C_{g,x\theta }^2L_{g,\theta \theta }^2}}{{\mu _g^4}} + 8\frac{{C_{f,\theta }^2L_{g,x\theta }^2}}{{\mu _g^2}}} )\frac{1}{m}\sum\limits_{i = 1}^m {\sum\limits_{j = 1}^m { {{{\left\| {\theta _i^*(x) - \theta _j^*(x)} \right\|}^2}} } }  + 8\frac{{C_{f,\theta }^2C_{g,x\theta }^2}}{{\mu _g^4}}b_{g}^2 + 8\frac{{C_{f,\theta }^2}}{{\mu _g^2}}b_{g}^2.  \nonumber
\end{align}
As for the second term on the right hand of \eqref{EQ-phi-phi}, we bound it by:
\begin{align}
 &\sum\limits_{i = 1}^m {{{\| {{\nabla _x}{f_i}(x,\theta _i^*(x)) \!- \!\frac{1}{m}\sum\limits_{j = 1}^m {{\nabla _x}{f_j}(x,\theta _j^*(x))} } \|^2}}} \nonumber \\
    \mathop  \leqslant \limits^{({\rm a})}\! & 2b_{{f}}^2 + 2\frac{1}{m}\sum\limits_{i = 1}^m {\sum\limits_{j = 1}^m {{{\left\| {{\nabla _x}{f_j}(x,\theta _j^*(x)) \! - \! {\nabla _x}{f_j}(x,\theta _i^*(x))} \right\|}^2}} }  \nonumber  \\
    \mathop  \leqslant \limits^{({\rm b})} & 2b_{{f}}^2 + 2{L_{{f},x}^2}\frac{1}{m}\sum\limits_{i = 1}^m {\sum\limits_{i = j}^m {{{\left\| {\theta _j^*(x) - \theta _i^*(x)} \right\|}^2}} },  \label{EQ-phi-phi-2}
\end{align}
where step (a) introduces the term $ \frac{1}{m}\sum\nolimits_{j = 1}^m {{\nabla _x}{f_j}(x,\theta _i^*(x))}$ and uses Young’s inequality and the bounded data heterogeneity  in Assumption \ref{ASS-heterogeneity};  step (b) comes from Lipschitz continuity  of $\nabla _{x} f_i$. In addition, utilizing  the strong convexity of $g_i$ in $\theta$, we further derive that:
\begin{align}
 \frac{1}{m}\sum\limits_{i = 1}^m {\sum\limits_{j = 1}^m {{{\left\| {\theta _j^*(x) - \theta _i^*(x)} \right\|}^2}} }
\leqslant & \frac{1}{m}\sum\limits_{i = 1}^m {\sum\limits_{j = 1}^m {\frac{1}{{\mu _{{g}}^2}}{{\left\| {{\nabla _\theta }{g_i}(x,\theta _j^*(x)) - {\nabla _\theta }{g_i}(x,\theta _i^*(x))} \right\|}^2}} }\label{EQ-theta-theta}  \\
 \mathop  = \limits^{({\rm a})} & \frac{1}{m}\sum\limits_{i = 1}^m {\sum\limits_{j = 1}^m {\frac{1}{{\mu _{{g}}^2}}{{\left\| {{\nabla _\theta }{g_i}(x,\theta _j^*(x)) - {\nabla _\theta }{g_j}(x,\theta _j^*(x))} \right\|}^2}} }
\leqslant &  \!{\frac{1}{{\mu _{{g}}^2}}}b_{{g}}^2, \nonumber
\end{align}
where   step (a) uses the fact that $\nabla _{\theta} g_i(x,\theta _i^*(x))=0, i \in \mathcal{V}$.

Then, by substituting the results \eqref{EQ-phi-phi-1}, \eqref{EQ-phi-phi-1-2}, \eqref{EQ-nablaphi-nablaphi}, \eqref{EQ-nablaphi-nablaphi-1},  \eqref{EQ-phi-phi-2}, \eqref{EQ-theta-theta} into \eqref{EQ-phi-phi} and rearranging the terms, we reach the following inequality:
\begin{align}
  &\sum\limits_{i = 1}^m {{{\left\| {\nabla {\Phi _i}\left( x \right) - \nabla \Phi \left( x \right)} \right\|}^2}}  \nonumber \\
   \leqslant& 4b_{{f}}^2 + 4{\max _i}\{ {{{\left\| {\nabla \theta _i^*(x)} \right\|}^2}} \}b_{{f}}^2 +16\frac{{C_{f,\theta }^2C_{g,x\theta }^2}}{{\mu _g^4}}b_{g }^2 + 16\frac{{C_{f,\theta }^2}}{{\mu _g^2}}b_{g }^2   \nonumber\\
   &+ 4\Big( {L_{f,x}^2 + {\max _i}\{ {L_{f,\theta }^2{{\left\| {\nabla \theta _i^*(x)} \right\|}^2}} \} + { {4\frac{{C_{f,\theta }^2C_{g,x\theta }^2L_{g,\theta \theta }^2}}{{\mu _g^4}} + 4\frac{{C_{f,\theta }^2L_{g,x\theta }^2}}{{\mu _g^2}}}}} \Big)\frac{1}{{\mu _g^2}}b_{{g}}^2 \\
    \leqslant& \!\underbrace{(4\!+\!4\frac{{C_{g,x\theta }^2}}{{\mu _g^2}})}_{\triangleq C_1(\mu_g,C_{g,x\theta})}\!b_{{f}}^2 \!+\! \underbrace{4\Big(\! 4\frac{{C_{f,\theta }^2C_{g,x\theta }^2}}{{\mu _g^4}}\!+\!4\frac{{C_{f,\theta }^2}}{{\mu _g^2}}+ {\frac{L_{f,x}^2}{\mu_g^2} \!+\! \frac{{L_{f,\theta }^2C_{g,x\theta }^2}}{{\mu _g^4}} \!+\! { {4\frac{{C_{f,\theta }^2C_{g,x\theta }^2L_{g,\theta \theta }^2}}{{\mu _g^6}} \!+ \!4\frac{{C_{f,\theta }^2L_{g,x\theta }^2}}{{\mu _g^4}}} }} \!\Big)}_{\triangleq C_2(\mu_g,L_{f,x},L_{f,\theta},L_{g,x\theta},L_{g,\theta\theta},C_{f,\theta},C_{g,x\theta})}\!b_{{g}}^2,  \nonumber
\end{align}
where the last step uses the following bound of $\nabla\theta _i^*(x)$:
\begin{equation}\label{EQ-nabla-theta}
\begin{aligned}
{\left\| {\nabla \theta _i^*(x)} \right\|^2} = {\| {{\nabla _{x\theta }}{g_i}( {x,\theta _i^*(x)} ){{[ {{\nabla _{\theta \theta }}{g_i}\left( {x,\theta _i^*(x)} \right)} ]}^{ - 1}}} \|^2} \leqslant \frac{{C_{g,x\theta }^2}}{{\mu _g^2}}.
\end{aligned}
\end{equation}
This completes the proof.   {\hfill $\blacksquare$}
\subsection{Proof  of Lemma \ref{LE-descent}} \label{sec-proof-LE-descent}
For ease of presentation, we recall  that the update of  $y^{k+1}$ in the case with  gradient tracking scheme in Algorithm \ref{alg:1} is given by:
\[{y^{k + 1}} = \mathcal{W}{y^k} + {z^{k + 1}} - {z^k}. \]
Let the matrix $J \triangleq \frac{{1_m^{\rm T} }}{m} \otimes {I_n}$ denote the average operator among nodes.  Since the weighted matrix $\mathcal{W}$ is doubly stochastic, multiplying the matrix $J$ in both sides of above equality yields:
\begin{equation}
\begin{aligned}
{{\bar y}^{k + 1}} = {\bar y}^k + {{\bar z}^{k + 1}} - {{\bar z}^k}.
\end{aligned}
\end{equation}
By applying induction, we can establish that ${{\bar y}^k}={{\bar z}^k}$ when the initial condition is $y^0=z^0$. This implies that each node is capable of tracking the full gradient ${{\bar z}^k}$. On the other hand, considering the case with the local gradient scheme in Algorithm \ref{alg:1}, we can directly observe that $y^k=z^k$ based on \eqref{ALG-LPDBO-h1}. Thus, in both the local gradient and tracking gradient schemes, it can be deduced that $\bar y^k=\bar z^k$. Furthermore, by multiplying the matrix $J$ on both sides of (\ref{ALG-LPDBO}a), we can determine the average state of the update (\ref{ALG-LPDBO}a) across all nodes as follows:
\begin{equation}\label{EQ-x-bar}
\begin{aligned}
{\bar x}^{k+1}={\bar x}^{k}-\tau\alpha {\bar y ^k}.
\end{aligned}
\end{equation}
Since the overall objective function $\Phi$ is smooth by Proposition \ref{PR-smooth}, it holds that:
\begin{equation}
\begin{aligned}
 \mathbb{E}[ {\Phi ({{\bar x}^{k + 1}})|{{{\mathcal{F}}}^k}} ]
   \mathop  \leqslant \limits^{(\rm a)}& \Phi ({{\bar x}^k}) - \tau \alpha \mathbb{E}[ {\langle {\nabla \Phi ({{\bar x}^k}),{{\bar y}^k}} \rangle |{{{\mathcal{F}}}^k}} ] + \frac{{\tau ^2 \alpha ^2 {L}}}{2}\mathbb{E}[ {{{\| {{{\bar y}^k}} \|}^2}|{{{\mathcal{F}}}^k}} ] \\
   \mathop \leqslant  \limits^{(\rm b)}& \Phi ({{\bar x}^k}) + \frac{\tau\alpha}{2} \mathbb{E}[ {{{\| {\nabla\Phi ({{\bar x}^k}) - {{\bar z}^k}} \|}^2}|{{{\mathcal{F}}}^k}} ] \\
   &- \frac{\tau\alpha }{2}{\| \nabla{\Phi ({{\bar x}^k})} \|^2} - \frac{\tau\alpha }{2}( {1 - \tau\alpha L} )\mathbb{E}[ {{{\| {{{\bar y}^k}} \|}^2}|{{{\mathcal{F}}}^k}} ], \\
\end{aligned}
\end{equation}
where    step (a) uses the recursion \eqref{EQ-x-bar}; step (b) holds due to the equation $2a^{\rm T}b=\|a\|^2+\|b\|^2-\|a-b\|^2$; and Young's inequality as well as the fact that  $\bar y^k=\bar z^k$. Then taking the total expectation yields the desired result. This completes the proof. {\hfill $\blacksquare$}
\subsection{Proof  of Lemma \ref{LE-HypergradientE}}
\label{sec-proof-LE-HypergradientE}
Recall that the update of ${{z}^{k+1}}$  in \eqref{ALG-LPDBO-g} is given by: ${z^{k}} = {s^{k}} + (1 - \gamma )( {{z^k} - {s^k}} ). $
Then considering the update \eqref{ALG-LPDBO-g}, we can obtain a recursive expression for the term $  \nabla \Phi ({{\bar x}^{k + 1}}) - {\text{ }}{{\bar z}^{k + 1}}$ by introducing the terms $\mathbb{E}[{{\bar s}^k}|{{{\mathcal{F}}}^k}]$ and $\nabla \Phi (\bar{x}^k)$ as follows:
\begin{equation}\label{EQ-s-z-bar}
\begin{aligned}
  \nabla \Phi ({{\bar x}^{k + 1}}) - {\text{ }}{{\bar z}^{k + 1}} = &\nabla \Phi ({{\bar x}^{k + 1}}) - ({{\bar s}^k} + (1 - \gamma )({{\bar z}^k} - {{\bar s}^k})) \\
   = &(1 - \gamma )(\nabla \Phi ({{\bar x}^k}) - {{\bar z}^k}) + \gamma (\nabla \Phi ({{\bar x}^k}) - \mathbb{E}[{{\bar s}^k}|{{{\mathcal{F}}}^k}]) \\
   &+ \gamma (\mathbb{E}[{{\bar s}^k}|{{{\mathcal{F}}}^k}] - {{\bar s}^k}) + \nabla \Phi ({{\bar x}^{k + 1}}) - \nabla \Phi ({{\bar x}^k}). \\
\end{aligned}
\end{equation}
Taking the square norm on both sides under the conditional expectation of $\mathcal{F}^k$, it follows from \eqref{EQ-s-z-bar} that:
\begin{equation}\label{EQ-reduce-variance}
\begin{aligned}
  &\mathbb{E}[ {{{\| {\nabla \Phi ({{\bar x}^{k + 1}}) - {{\bar z}^{k + 1}}} \|}^2}|{{{\mathcal{F}}}^k}} ] \\
   \mathop  = \limits^{({\rm a})} &[ {{{\| {(1 - \gamma )(\nabla \Phi ({{\bar x}^k}) - {{\bar z}^k}) + \gamma (\nabla \Phi ({{\bar x}^k}) - \mathbb{E}[{{\bar s}^k}|{{{\mathcal{F}}}^k}]) + \nabla \Phi ({{\bar x}^{k + 1}}) - \nabla \Phi ({{\bar x}^k})} \|}^2}|{{{\mathcal{F}}}^k}} ] \\
   &+ {\gamma ^2}\mathbb{E}[ {{{\| {\mathbb{E}[{{\bar s}^k}|{{{\mathcal{F}}}^k}] - {{\bar s}^k}} \|}^2}|{{{\mathcal{F}}}^k}} ] \\
    \mathop  \leqslant \limits^{({\rm b})} &(1 - \gamma ){\| {\nabla \Phi ({{\bar x}^k}) - {{\bar z}^k}} \|^2} + 2\gamma {\left\| {\nabla \Phi ({{\bar x}^k}) - \mathbb{E}[{{\bar s}^k}|{{{\mathcal{F}}}^k}]} \right\|^2} + \frac{2}{\gamma }{\| {\nabla \Phi ({{\bar x}^{k + 1}}) - \nabla \Phi ({{\bar x}^k})} \|^2} \\
    &+ {\gamma ^2}\mathbb{E}[ {{{\| {\mathbb{E}[{{\bar s}^k}|{{{\mathcal{F}}}^k}] - {{\bar s}^k}} \|}^2}|{{{\mathcal{F}}}^k}} ] \\
    \mathop  \leqslant \limits^{({\rm c})} & (1 - \gamma ){\| {\nabla \Phi ({{\bar x}^k}) - {{\bar z}^k}} \|^2} \!+\! 2\gamma {\| {\nabla \Phi ({{\bar x}^k}) \!- \!\mathbb{E}[{{\bar s}^k}|{{{\mathcal{F}}}^k}]} \|^2} \\
    &+ \frac{{2{L^2}}}{\gamma }{\tau ^2}{\alpha ^2}{\| {{{\bar y}^k}} \|^2}\! +\! {\gamma ^2}\mathbb{E}[ {{{\| {\mathbb{E}[{{\bar s}^k}|{{{\mathcal{F}}}^k}]\! - \!{{\bar s}^k}} \|}^2}|{{{\mathcal{F}}}^k}}],
\end{aligned}
\end{equation}
where   step (a) holds due to the fact that $\bar s^{k}$ is  independent of $\bar z^k$ and  $\mathbb{E}[\bar s^{k}|\mathcal{F}^k]$ is an unbiased estimate of ${\bar s^{k}}$, with  the additional assurance that   the samples in each node  are  independent of each other;   step (b) is derived by the convexity of $l_2$-norm and Young's inequality with the condition  $0<\gamma<1$; step (c) follows from  Lipschitz continuity of $\nabla \Phi$ and the recursion \eqref{ALG-LPDBO-b}. Next, we will bound  the last term  in \eqref{EQ-reduce-variance}. Note that  from the recursion \eqref{ALG-LPDBO-f} that:
\begin{equation}\label{EQ-sk1-sk1-zeta}
\begin{aligned}
\mathbb{E}[\bar s^{k}|\mathcal{F}^k] - {{\bar s}^{k}} =& J {\nabla _x}F( {{x^{k}},{\theta ^{k}}} ) - J \nabla _{x\theta }^2G( {{x^{k}},{\theta ^{k}}} ){v^{k}} \\
   &- J {\nabla _x}\hat F( {{x^{k}},{\theta ^{k+1}};\varsigma _2^{k}} ) + J \nabla _{x\theta }^2\hat G( {{x^{k}},{\theta ^{k}};\xi _3^{k}} ){v^{k}}. \\
\end{aligned}
\end{equation}
Taking the square norm on both sides of \eqref{EQ-sk1-sk1-zeta} under the total expectation, it follows that
\begin{equation}\label{EQ-variance-s}
\begin{aligned}
  \mathbb{E}[ {{{\| \mathbb{E}[\bar s^{k}|\mathcal{F}^k] - {{\bar s}^{k}} \|}^2}} ]
    \mathop  = \limits^{({\rm a})} &\mathbb{E}[ {{{\| J {{\nabla _x}F( {{x^{k}},{\theta ^{k}}} ) - J {\nabla _x}\hat F( {{x^{k}},{\theta ^{k}};\varsigma _2^{k}} )} \|}^2}} ] \\
   &+\mathbb{E}[ {{{\|J {\nabla _{x\theta }^2G( {{x^{k}},{\theta ^{k}}} ){v^{k}} - J\nabla _{x\theta }^2\hat G( {{x^{k}},{\theta ^{k}};\xi _3^{k}} ){v^{k}}} \|}^2}} ]
   \\
   \mathop  \leqslant \limits^{({\rm b})} & \frac{1}{m^2}\res{(m{\sigma _{f,x}^2}+\|v^k\|^2{\sigma _{g,x\theta}^2})},\\
     \mathop  \leqslant \limits^{({\rm c})} & \frac{1}{m}\res{({\sigma _{f,x}^2}+\nyc{\hat M^2}{\sigma _{g,x\theta}^2})}, \\
\end{aligned}
\end{equation}
where  step (a) uses the independence between the samples $\varsigma _2^{k}$ and $\xi _3^{k}$ and unbiased estimates of stochastic gradients in Assumption \ref{ASS-STOCHASTIC};   step (b) follows from the bounded variances in Assumption \ref{ASS-STOCHASTIC}; step (c) is derived by the following fact:
\begin{equation}\label{EQ-v-mm}
 \frac{1}{m}\|v^k\|^2 \leqslant 2 \frac{1}{m}\|v^*(\bar x^k)\|^2+  2\frac{1}{m}\|v^k-v^*(\bar x^k)\|^2  \leqslant 2M^2+ 2 \frac{1}{m}\|v^k-v^*(\bar x^k)\|^2 \triangleq \nyc{\hat M^2}.
\end{equation}
with the last step following the result that $\|v_i^{*}(\bar{x}^k)\| \le M$ in Proposition \ref{PR-boundness}. Then, taking total expectation and combining the upper bound \eqref{EQ-variance-s}, we reach the result \eqref{LE-stochastic-error} as follows:
\begin{equation}
\begin{aligned}
\mathbb{E}[ {{{\| \nabla \Phi ({{\bar x}^{k + 1}}) - {\text{ }}{{\bar z}^{k + 1}} \|}^2}} ]
\leqslant &
(1 - \gamma )\mathbb{E}[{\| {\nabla \Phi ({{\bar x}^k}) - {{\bar z}^k}} \|^2}] \!+\!\underbrace{2\frac{\gamma}{\alpha}}_{\triangleq r_z} \alpha \mathbb{E}[{\| {\nabla \Phi ({{\bar x}^k}) \!- \!\mathbb{E}[{{\bar s}^k}|{{{\mathcal{F}}}^k}]} \|^2}] \\
&+ \underbrace{\frac{{2{L^2}}}{\gamma }}_{\triangleq r_y}{\tau ^2}{\alpha ^2}\mathbb{E}[{\| {{{\bar y}^k}} \|^2}]\! + \frac{1}{m}\underbrace {\res{({\sigma _{f,x}^2}+\nyc{\hat M^2}{\sigma _{g,x\theta}^2})}{\frac{\gamma^2}{\alpha^2}}}_{\triangleq \sigma _{\bar z}^2}{\alpha ^2}.
\end{aligned}
\end{equation}
In what follows, we will derive the remaining result \eqref{EQ-hypergradient-error}. To this end, by the definition of $\nabla \Phi(\bar x^k)$, we first have:
\[\nabla \Phi ({{\bar x}^k}) = J{\nabla _x}F( {{{{1}}_m} \otimes {{\bar x}^k},{\theta ^*}({{\bar x}^k})} ) - J\nabla _{x\theta }^2G( {{{{1}}_m} \otimes {{\bar x}^k},{\theta ^*}({{\bar x}^k})} )( {{v^*}({{\bar x}^k})} ).\]
Additionally, it follows from  the recursion  \eqref{ALG-LPDBO-d} that:
\[ \mathbb{E}[ \bar{s}^k|\mathcal{F}^k]=J\nabla_xF(x^k,\theta^k)-J\nabla_{x\theta}^2G(x^k,\theta^k)v^k.\]
Thus, the term ${\nabla \Phi ({{\bar x}^k}) - {{\bar s}^k}}$ can be expressed as:
\begin{align}
\nabla \Phi ({{\bar x}^k}) - \mathbb{E}[ \bar{s}^k|\mathcal{F}^k] =& J{\nabla _x}F( {{{{1}}_m}\otimes{{\bar x}^k},{\theta ^*}({{\bar x}^k})} ) - J{\nabla _x}F( {{x^k},{\theta ^k}} ) +J\nabla _{x\theta }^2G({x}^k,{\theta ^k})(v^*(\bar x^k)-v^k)  \nonumber\\
&\nyc{+ J( {\nabla _{x\theta }^2G( {{x^k},{\theta ^k}} ) - \nabla _{x\theta }^2G( {{{{1}}_m}\otimes{{\bar x}^k},{\theta ^*}({{\bar x}^k})} )} ){v^*(\bar x ^k)}.}
\label{EQ-nablaphi-s}
\end{align}
Then taking the square norm on both sides of \eqref{EQ-nablaphi-s} under the total expectation and employing the  inequality $\|a+b\|^2 \leqslant 2\|a\|^2+2\|b\|^2$ twice and Jensen inequality,  we get:
\begin{align}
  &\mathbb{E}[ {{{\| {\nabla \Phi ({{\bar x}^k}) - \mathbb{E}[ \bar{s}^k|\mathcal{F}^k]} \|}^2}} ] \nonumber \hfill \\
   \leqslant &2\frac{1}{m}\mathbb{E}[ {{{\| {{\nabla _x}F( {{{{1}}_m}\otimes{{\bar x}^k},{\theta ^*}({{\bar x}^k})} ) - {\nabla _x}F( {{x^k},{\theta ^k}} )} \|}^2}} ] \nonumber \hfill \\
   + &4\frac{1}{m}\mathbb{E}[ {{{\| { - \nabla _{x\theta }^2G({x}^k,{\theta ^k}){v^*}({{\bar x}^k}) + \nabla _{x\theta }^2G({x}^k,{\theta ^k}){v^k}} \|}^2}} ] \hfill \\
   + &4\frac{1}{m}\mathbb{E}[ {{{\| { - \nabla _{x\theta }^2G( {{{{1}}_m}\otimes{{\bar x}^k},{\theta ^*}({{\bar x}^k})} ){v^*}({{\bar x}^k}) + \nabla _{x\theta }^2G( {{x^k},{\theta ^k}} ){v^*}({{\bar x}^k})} \|}^2}} ] \nonumber \hfill \\
   \leqslant & \underbrace {( {2L_{f,x }^2 \!+\! 4{M^2}L_{g,x\theta }^2} )}_{ \triangleq \yc{\Ld}}\frac{1}{m}\mathbb{E}[ {{{\| {{x^k}\! - \!{{{1}}_m}\otimes{{\bar x}^k}} \|}^2}\! +\! {{\| {{\theta ^k} \!-\! \theta^*({{\bar x}^k})} \|}^2}} ] \!+ \!4C_{g,x\theta }^2\frac{1}{m}\mathbb{E}[ {{{\| {{v^k} \!-\! {v^*}({{\bar x}^k})} \|}^2}} ]. \hfill \nonumber
\end{align}
To obtain the last step, we use Lipschitz continuity of $\nabla_{x}f_i$ and $\nabla_{x\theta}^2g_i$ in Assumptions \ref{ASS-OUTLEVEL} and \ref{ASS-INNERLEVEL}  as well as the boundedness of $\|v_i^*(\bar x^k)\|$ in Proposition \ref{PR-boundness}. This completes the proof.   {\hfill $\blacksquare$}

\subsection{Proof  of Lemma \ref{LE-Vstar}}
\label{sec-proof-LE-Vstar}
First note that the term $\mathbb{E}[\| {{v^{k + 1}} - {v^*}({{\bar x}^{k + 1}})} \|^2|\mathcal{F}^k]$ can be expanded as:
\begin{equation}\label{EQ-vk-vstar}
\begin{aligned}
\mathbb{E}[{\| {{v^{k + 1}} - {v^*}({{\bar x}^{k + 1}})} \|^2}| \mathcal{F}^k] = &  \underbrace {\mathbb{E}[{{\| {{v^{k + 1}} - {v^*}({{\bar x}^k})} \|}^2}| \mathcal{F}^k]}_{ \triangleq A_1^v} + {\mathbb{E}[\| {{v^*}({{\bar x}^k}) - {v^*}({{\bar x}^{k + 1}})} \|^2| \mathcal{F}^k]} \\
    & + \underbrace {\mathbb{E}[2\langle {{v^{k + 1}} - {v^*}({{\bar x}^k}),{v^*}({{\bar x}^k}) - {v^*}({{\bar x}^{k + 1}})} \rangle| \mathcal{F}^k] }_{ \triangleq A_2^v}.
\end{aligned}
\end{equation}
We first bound the term $A_1^v$. To this end, for reader's convenience,   we repeat the argument on the recursion of $v^{k+1}$ in \eqref{ALG-LPDBO-b} by combining \eqref{ALG-LPDBO-e} as follows:
\[{v^{k + 1}} = (I-\lambda \nabla _{\theta \theta }^2\hat G({x^k},{\theta ^k};\xi _2^{k }) ){v^k} +  \lambda \hat F( {{x^k},{\theta ^k};\varsigma _1^{k }} ).\]
Substituting the above expression into the term $A_1^v$, we get:
\begin{align}
  A_1^v
   =& \mathbb{E}[ {{{\| {( {I - \lambda \nabla _{\theta \theta }^2\hat G({x^k},{\theta ^k};\xi _2^{k })} ){v^k} + \lambda {\nabla _\theta }\hat F( {{x^k},{\theta ^k};\varsigma _1^{k}} ) - {v^*}({{\bar x}^k})} \|}^2}} | \mathcal{F}^k] \nonumber \\
   = & {{{\| {( {I - \lambda \nabla _{\theta \theta }^2G({x^k},{\theta ^k})} ){v^k} + \lambda {\nabla _\theta }F( {{x^k},{\theta ^k}} ) - {v^*}({{\bar x}^k})} \|}^2}} \label{EQ-H1v} \\
   &\!\!+ \!{\lambda ^2}\mathbb{E}[ \!{{{\| {\nabla _{\theta \theta }^2\hat G({x^k},{\theta ^k};\xi _2^{k })v^k \!- \!\nabla _{\theta \theta }^2G({x^k},{\theta ^k})v^k \!+\! {\nabla _\theta }F( {{x^k},{\theta ^k}} )\! -\! {\nabla _\theta }\hat F( {{x^k},{\theta ^k};\varsigma _1^{k}} )} \|^2}}} | \mathcal{F}^k ] \nonumber  \\
   \leqslant&  {{{\| {( {I - \lambda \nabla _{\theta \theta }^2G({x^k},{\theta ^k})} ){v^k} + \lambda {\nabla _\theta }F( {{x^k},{\theta ^k}} ) - {v^*}({{\bar x}^k})} \|}^2}}  +m\res{(\sigma _{f,\theta}^2+\hat
 M^2\sigma _{g,\theta\theta}^2)}{\lambda ^2}, \nonumber
\end{align}
where the last step follows from the bounded variances in Assumption \ref{ASS-STOCHASTIC} and the inequality \eqref{EQ-v-mm}.
Noting that ${v^*}({{\bar x}^k})$ admits the following expression:
\begin{equation}\label{EQ-Vstar}
\begin{aligned}
{v^*}({{\bar x}^k}) = {[ {\nabla _{\theta \theta }^2G({{{1}}_m} \otimes {{\bar x}^k},{\theta ^*}({{\bar x}^k}))} ]^{ - 1}}{\nabla _\theta }F({{{1}}_m} \otimes {{\bar x}^k},{\theta ^*}({{\bar x}^k})),
\end{aligned}
\end{equation}
then we rewrite that:
\begin{equation}
\begin{aligned}
 & ( {I - \lambda \nabla _{\theta \theta }^2G({x^k},{\theta ^k})} ){v^k} + \lambda {\nabla _\theta }F( {{x^k},{\theta ^k}} ) - {v^*}({{\bar x}^k}) \\
   =& \lambda ( {{\nabla _\theta }F( {{x^k},{\theta ^k}} ) - {\nabla _\theta }F({{{1}}_m} \otimes {{\bar x}^k},{\theta ^*}({{\bar x}^k}))} ) + ( {I - \lambda \nabla _{\theta \theta }^2G({x^k},{\theta ^k})} )( {{v^k} - {v^*}({{\bar x}^k})} ) \\
   &+ \lambda( {\nabla _{\theta \theta }^2G({{{1}}_m} \otimes {{\bar x}^k},{\theta ^*}({{\bar x}^k})) - \nabla _{\theta \theta }^2G({x^k},{\theta ^k})} ){v^*}({{\bar x}^k}). \\
\end{aligned}
\end{equation}
With the above expression,  the first term on the right hand of \eqref{EQ-H1v} can be bounded by:
\begin{equation}\label{EQ-H-v0}
\begin{aligned}
 &  {{{\| {( {I - \lambda \nabla _{\theta \theta }^2G({x^k},{\theta ^k})} ){v^k} + \lambda {\nabla _\theta }F( {{x^k},{\theta ^k}} ) - {v^*}({{\bar x}^k})} \|}^2}}  \\
   \mathop  \leqslant \limits^{({\rm a})}  &{( {1 + \frac{{{\mu _g}\lambda }}{3}} )^2} {{{\| {( {I - \lambda \nabla _{\theta \theta }^2G({x^k},{\theta ^k})} )( {{v^k} - {v^*}({{\bar x}^k})} )} \|}^2}}  \\
   &+ ( {1 + \frac{3}{{{\mu _g}\lambda }}} ){\lambda ^2}{{\| {{\nabla _\theta }F( {{x^k},{\theta ^k}} ) - {\nabla _\theta }F({{{1}}_m} \otimes {{\bar x}^k},{\theta ^*}({{\bar x}^k}))} \|}^2} \\
   &+ ( {1 + \frac{{{\mu _g}\lambda }}{3}} )( {1 + \frac{3}{{{\mu _g}\lambda }}} ){\lambda ^2} {{{\| {( {\nabla _{\theta \theta }^2G({{{1}}_m} \otimes {{\bar x}^k},{\theta ^*}({{\bar x}^k})) - \nabla _{\theta \theta }^2G({x^k},{\theta ^k})} ){v^*}({{\bar x}^k})} \|}^2}}  \\
    \mathop  \leqslant \limits^{({\rm b})} & {( {1 + \frac{{{\mu _g}\lambda }}{3}} )^2}{( {1 - {\mu _g}\lambda } )^2} {{{\| {{v^k} - {v^*}({{\bar x}^k})} \|}^2}} \\
   &+ ( {1 + \frac{3}{{{\mu _g}\lambda }}} )L_{f,\theta }^2{\lambda ^2} [ {{{\| {{x^k} - {{{1}}_m} \otimes {{\bar x}^k}} \|}^2} + {{\| {{\theta ^k} - {\theta ^*}({{\bar x}^k})} \|}^2}}] \\
   &+ ( {1 + \frac{3}{{{\mu _g}\lambda }}} )( {1 + \frac{{{\mu _g}\lambda }}{3}} )L_{g,\theta \theta }^2{M^2}{\lambda ^2}[ {{{\| {{x^k} - {{{1}}_m} \otimes {{\bar x}^k}} \|}^2} + {{\| {{\theta ^k} - {\theta ^*}({{\bar x}^k})} \|}^2}} ], \\
\end{aligned}
\end{equation}
where   step (a) uses Young's inequality twice and step (b) follows from the strong convexity of $g_i$, Lipschitz continuity of $\nabla_{\theta}f_i$ and $\nabla_{\theta\theta}^2g_i$, and the boundedness of $\|v_i^*(\bar x^k)\|$. By considering the condition $1-\lambda \mu_g>0$ and rearranging the terms, we can further derive from the inequality \eqref{EQ-H-v0} the following results:
\begin{align} &  {{{\| {( {I - \lambda \nabla _{\theta \theta }^2G({x^k},{\theta ^k})} ){v^k} + \lambda {\nabla _\theta }F( {{x^k},{\theta ^k}} ) - {v^*}({{\bar x}^k})} \|}^2}} \label{EQ-H-v} \\
   \leqslant& {( {1 \!- \!\frac{\mu_g \lambda}{3}} )} ( {1 \!-\! {\mu _g}\lambda } ) {{{\| {{v^k}  \!- \! {v^*}({{\bar x}^k})} \|}^2}}  \!+ \!\frac{{2\lambda }}{{{\mu _g}}}\underbrace {\left( {2L_{f,\theta }^2 \!+ \!4{M^2}L_{g,\theta \theta }^2} \right)}_{ \triangleq \yc{\Lb}}[ {{{\| {{x^k} - \!{{{1}}_m} \otimes {{\bar x}^k}} \|}^2}\! + \!{{\| {{\theta ^k} - {\theta ^*}({{\bar x}^k})} \|}^2}} ]. \nonumber
\end{align}
Besides, substituting  \eqref{EQ-H-v} into  \eqref{EQ-H1v}, we reach an upper bound for the term $ A_1^v$ in \eqref{EQ-vk-vstar} as follows:
\begin{equation}\label{EQ-V-boundness}
\begin{aligned}
 A_1^v\!\leqslant& {( {1\!-\!\frac{\mu_g \lambda}{3}} )}\!\left( {1\! -\! {\mu _g}\lambda } \right)\! {{{\| {{v^k} \!-\! {v^*}({{\bar x}^k})} \|}^2}}  \!  +\! m(1 \!+ \!{M^2}){\lambda ^2}{\sigma ^2} \! \\
 &+ \frac{{2\lambda }}{{{\mu _g}}}\yc{\Lb}[ {{{\| {{x^k}\! -\! {{{1}}_m} \otimes {{\bar x}^k}} \|}^2} \!\!+ \!{{\| {{\theta ^k}\! - \!{\theta ^*}({{\bar x}^k})}\! \|}^2}} ].
\end{aligned}
\end{equation}
Now, we analyze the term ${\mathbb{E}[\| {{v^*}({{\bar x}^k}) - {v^*}({{\bar x}^{k + 1}})} \|^2]}$ in \eqref{EQ-vk-vstar}. Employing Lipschitz continuity of $v_i^*(x)$ in Proposition \ref{PR-smooth}, the recursion \eqref{ALG-LPDBO-a} and the upper bound of the variances in \eqref{LE-stochastic-error}, we have that:
\begin{equation}\label{EQ-vk-vk}
\begin{aligned}
  \mathbb{E}[ {{{\| {{v^*}({{\bar x}^k}) - {v^*}({{\bar x}^{k + 1}})} \|}^2}} | \mathcal{F}^k ] \leqslant mL_{{v^*}}^2\mathbb{E}[ {{{\| {{{\bar x}^k} - {{\bar x}^{k + 1}}} \|}^2}} | \mathcal{F}^k ] \leqslant  m\tau^2{\alpha ^2}L_{{v^*}}^2{{{\| {{{\bar y}^k}} \|}^2}}. \\
\end{aligned}
\end{equation}
It  remains to analyze the term $A_2^v$ in \eqref{EQ-vk-vstar}, which can be rewritten as the following expression by  leveraging Cauchy-Schwartz inequality:
\begin{equation}\label{EQ-Z}
\begin{aligned}
  A_2^v \leqslant \frac{{{\mu _g}\lambda }}{3}A_1^v + \frac{3}{{{\mu _g}\lambda }}\mathbb{E}[ {{{\| {{v^*}({{\bar x}^k}) - {v^*}({{\bar x}^{k + 1}})} \|}^2}|{{{\mathcal{F}}}^k}} ]. \\
\end{aligned}
\end{equation}

Furthermore, in combination with the results \eqref{EQ-vk-vstar}, \eqref{EQ-V-boundness},  \eqref{EQ-vk-vk}, \eqref{EQ-Z},   we reach an evolution for $\mathbb{E}[\| {{v^{k + 1}} - {v^*}({{\bar x}^{k + 1}})} \|^2 | \mathcal{F}^k]$ as follows:
\begin{equation}
\begin{aligned}
  &\mathbb{E}[ {{{\| {{v^{k + 1}} - {v^*}({{\bar x}^{k + 1}})} \|}^2}|{{{\mathcal{F}}}^k}} ] \\
   \mathop  \leqslant \limits^{({\rm a})} & ( {1 + \frac{{{\mu _g}\lambda }}{3}} )( {1 - \frac{{{\mu _g}\lambda }}{3}} )( {1 - {\mu _g}\lambda } ){\| {{v^k} - {v^*}({{\bar x}^k})} \|^2}\\
  & + ( {1 + \frac{{{\mu _g}\lambda }}{3}} )\frac{{2\lambda }}{{{\mu _g}}}\yc{\Lb}[ {{{\| {{x^k} - {1_m} \otimes {{\bar x}^k}} \|}^2} + {{\| {{\theta ^k} - {\theta ^*}({{\bar x}^k})} \|}^2}} ] \\
   &+ \frac{{2L_{{v^*}}^2{\tau^2\alpha ^2}}}{{\varpi \lambda }}m{\| {{{\bar y}^k}} \|^2} + ( {1 + \frac{{{\mu _g}\lambda }}{3}} )m(1 + \nyc{\hat M^2}){\lambda ^2}{\sigma ^2} \\
    \mathop  \leqslant \limits^{({\rm b})} & ( {1 - {\mu _g}\lambda} ){\| {{v^k} - {v^*}({{\bar x}^k})} \|^2} + \underbrace {\frac{{4\yc{\Lb}{\lambda }}}{{{\mu _g}\alpha}}}_{ \triangleq {q_x}}\alpha [ {{{\| {{x^k} - {1_m} \otimes {{\bar x}^k}} \|}^2} + {{\| {{\theta ^k} - {\theta ^*}({{\bar x}^k})} \|}^2}} ] \\
   &+ \underbrace {\frac{{2L_{{v^*}}^2}}{{\varpi {\lambda } }}}_{ \triangleq {q_s}}m{\tau^2\alpha ^2}{\| {{{\bar y}^k}} \|^2} + m\underbrace {2\res{(\sigma _{f, \theta}^2+\nyc{\hat
 M^2}\sigma _{g, \theta\theta}^2)}\frac{\lambda ^2}{\alpha ^2} }_{ \triangleq {\sigma_v ^2 }}{\alpha ^2}, \\
\end{aligned}
\end{equation}
where in   step (a) we denote $\varpi  \triangleq \frac{{{\mu _g}}}{3}$, and in   step (b) we use  the fact that $\mu_g\lambda<1$.  This completes the proof.   {\hfill $\blacksquare$}

\subsection{Proof  of Lemma \ref{LE-ThetaStar}}
\label{sec-proof-LE-ThetaStar}
The main idea  to prove the evolution of the inner-level errors $\mathbb{E}[ {{{\| {{\theta ^{k + 1}} - {\theta ^*}({{\bar x}^{k + 1}})} \|^2}}} |\mathcal{F}^k ]$ is similar to the one used for the Hv errors.   We start by decomposing the inner-level errors  as  follows:
\begin{equation}\label{EQ-theta-start}
\begin{aligned}
  \mathbb{E}[ {{{\| {{\theta ^{k + 1}} - {\theta ^*}({{\bar x}^{k + 1}})} \|}^2}} |\mathcal{F}^k ] &= \underbrace {\mathbb{E}[ {{{\| {{\theta ^{k + 1}} - {\theta ^*}({{\bar x}^k})} \|}^2}}|\mathcal{F}^k ]}_{ \triangleq  A_1^{\theta}} + \mathbb{E}[ {{{\left\| {{\theta ^*}({{\bar x}^k}) - {\theta ^*}({{\bar x}^{k + 1}})} \right\|}^2}} |\mathcal{F}^k] \\
   &+ \underbrace {2\mathbb{E}[ {\langle {{\theta ^{k + 1}} - {\theta ^*}({{\bar x}^k}),{\theta ^*}({{\bar x}^k}) - {\theta ^*}({{\bar x}^{k + 1}})} \rangle }|\mathcal{F}^k ]}_{ \triangleq A_2^{\theta}}. \\
\end{aligned}
\end{equation}
As for the term $A_1^{\theta}$, we have that:
\begin{equation}
\begin{aligned}
  {{\theta ^{k + 1}} \!-\! {\theta ^*}({{\bar x}^{k}})}= {\theta ^k} - \beta {\nabla _\theta }G({{{1}}_m} \otimes {{\bar x}^k},{\theta ^k}) - {\theta ^*}({{\bar x}^k})\! +\! \beta  {( {{\nabla _\theta }G({{{1}}_m} \!\otimes\! {{\bar x}^k},{\theta ^k})\! -\! {\nabla _\theta }\hat G({x^k},{\theta ^k};\xi _1^{k})} )}. \\
\end{aligned}
\end{equation}
Taking square norm on both sides under the conditional expectation $\mathcal{F}^k$, we have:
\begin{equation}\label{EQ-Q}
\begin{aligned}
  A_1^{\theta} =& \mathbb{E}[ {{{\| {{\theta ^k} - \beta {\nabla _\theta }\hat G({x^k},{\theta ^k};\xi _1^{k}) - {\theta ^*}({{\bar x}^k})} \|}^2}} |\mathcal{F}^k ] \\
   =&   {{{\| {{\theta ^k} - \!\beta {\nabla _\theta }G({{{1}}_m} \otimes {{\bar x}^k},{\theta ^k}) - {\theta ^*}({{\bar x}^k})} \|}^2}}  + {\beta ^2}\mathbb{E}[  {{\nabla _\theta }G({{{1}}_m} \otimes {{\bar x}^k},{\theta ^k}) - {\nabla _\theta }\hat G({x^k},{\theta ^k};\xi _1^{k})}  |\mathcal{F}^k ] \\
   & + 2\beta {\langle {{\theta ^k} - \beta {\nabla _\theta }G({{{1}}_m} \otimes {{\bar x}^k},{\theta ^k}) - {\theta ^*}({{\bar x}^k}),\mathbb{E}[ {{\nabla _\theta }G({{{1}}_m} \otimes {{\bar x}^k},{\theta ^k}) \!- \!{\nabla _\theta }\hat G({x^k},{\theta ^k};\xi _1^{k})}  |\mathcal{F}^k \! ]} \rangle }  \\
    \mathop  \leqslant \limits^{({\rm a})} & \!  ( {1 \!+ \!\beta {\omega _\theta }} ){{{\| {{\theta ^k} \!- \!\beta {\nabla _\theta }G({{{1}}_m} \otimes {{\bar x}^k},{\theta ^k}) \! -\!  {\theta ^*}({{\bar x}^k})} \|}^2}}\! +\! ( {\beta \! + \!\frac{1}{{{\omega _\theta }}}} )\beta L_{g,\theta }^2{{{\| {{x^k}\! -\! {{{1}}_m} \otimes {{\bar x}^k}} \|}^2}}  + m{\beta ^2}{\sigma _{g,\theta}^2} \\
       \mathop  \leqslant \limits^{({\rm b})}&  \left( {1 + \beta {\omega _\theta }} \right) {{{\| {{\theta ^k} - \beta {\nabla _\theta }G({{{1}}_m} \otimes {{\bar x}^k},{\theta ^k}) - {\theta ^*}({{\bar x}^k})} \|}^2}}   + \frac{2}{{{\omega _\theta }}}\beta L_{g,\theta }^2{{{\| {{x^k} - {{{1}}_m} \otimes {{\bar x}^k}} \|}^2}}  + m{\beta ^2}\res{\sigma _{g,\theta}^2}, \\
\end{aligned}
\end{equation}
where step (a)  uses the variance decomposition
and Cauchy-Schwartz inequality with parameter $\omega_{\theta}=\frac{{{\mu _g}{L_{g,\theta }}}}{ 2({{\mu _g} + {L_{g,\theta }}})}$, and the  step (b) comes from Lipschitz continuity of $\nabla_{\theta}g_i$ and the condition that $\beta<\frac{1}{\omega_{\theta}}$ in \eqref{EQ-beta}. Next, we proceed in providing an upper bound for the first term on the right hand of \eqref{EQ-Q} as follows:
\begin{align}
  & {{{\| {{\theta ^k} - \beta {\nabla _\theta }G({{{1}}_m} \otimes {{\bar x}^k},{\theta ^k}) - {\theta ^*}({{\bar x}^k})} \|}^2}} \nonumber \\
    \mathop  = \limits^{({\rm a})} & {{{\| {{\theta ^k} - {\theta ^*}({{\bar x}^k})} \|}^2}}  + {\beta ^2} {{{\| {{\nabla _\theta }G({{{1}}_m} \otimes {{\bar x}^k},{\theta ^k}) - {\nabla _\theta }G({{{1}}_m} \otimes {{\bar x}^k},{\theta ^*}({{\bar x}^k}))} \|}^2}}  \nonumber \\
   &- 2\beta \langle {{\nabla _\theta }G({{{1}}_m} \otimes {{\bar x}^k},{\theta ^k}) - {\nabla _\theta }G({{{1}}_m} \otimes {{\bar x}^k},{\theta ^*}({{\bar x}^k})),{\theta ^k} - {\theta ^*}({{\bar x}^k})} \rangle \nonumber \\
   \mathop  \leqslant \limits^{({\rm b})} & \! ( \!{1 \!- \!2\beta {\frac{{{\mu _g}{L_{g,\theta }}}}{{{\mu _g}\! + \!{L_{g,\theta }}}}}} \! ) {{{\| {{\theta ^k} \!-\! {\theta ^*}({{\bar x}^k})} \|}^2}} \! +\! ( {\beta  -\! \frac{2}{{{\mu _g}\! +\! {L_{g,\theta }}}}} )\beta  {{{\| {{\nabla _\theta }G({{{1}}_m}\! \otimes\! {{\bar x}^k},{\theta ^k}) \!-\! {\nabla _\theta }G({{{1}}_m} \! \otimes \! {{\bar x}^k},{\theta ^*}({{\bar x}^k}))} \|}^2}} \nonumber  \\
    \mathop  \leqslant \limits^{({\rm c})} &( {1 - 2\beta {\frac{{{\mu _g}{L_{g,\theta }}}}{{{\mu _g} + {L_{g,\theta }}}}}} ) {{{\| {{\theta ^k} - {\theta ^*}({{\bar x}^k})} \|}^2}}, \label{EQ-theta-sc-sm}
\end{align}
where   step (a) uses the fact that ${{\nabla _\theta }G({{{1}}_m} \otimes {{\bar x}^k},{\theta ^*}({{\bar x}^k}))}=0$;   step (b) come from the strong convexity and smoothness of $g_i$; step (c) holds due to the step-size condition $\beta=c_{\beta}{\alpha}  < \frac{2}{{{\mu _g} + {L_{g,\theta }}}}$. Then, plugging \eqref{EQ-theta-sc-sm} into  \eqref{EQ-Q} yields
\begin{align}
  A_1^{\theta} \leqslant & ( {1 + \beta {\omega _\theta }} )( {1 - 2\beta \frac{{{\mu _g}{L_{g,\theta }}}}{{{\mu _g} + {L_{g,\theta }}}}} ) {{{\| {{\theta ^k} - {\theta ^*}({{\bar x}^k})} \|}^2}} + ( {\beta  + \frac{1}{\omega_{\theta} }} )\beta L_{g,\theta }^2 {{{\| {{x^k} - {{{1}}_m} \otimes {{\bar x}^k}} \|}^2}}  + m{\beta ^2}\res{\sigma _{g,\theta}^2} \nonumber \\
   \leqslant& ( {1 -\frac{3}{2}  \beta \frac{{{\mu _g}{L_{g,\theta }}}}{{{\mu _g} + {L_{g,\theta }}}}} ) {{{\| {{\theta ^k} - {\theta ^*}({{\bar x}^k})} \|}^2}} + \frac{2}{{{\omega _\theta }}}\beta L_{g,\theta }^2 {{{\| {{x^k} - {{{1}}_m} \otimes {{\bar x}^k}} \|}^2}}  + m{\beta ^2}\res{\sigma _{g,\theta}^2}. \label{EQ-Q-boundness}
\end{align}
In addition, note that the inner-product term in \eqref{EQ-theta-start} follows that
\begin{equation}\label{EQ-ip}
\begin{aligned}
\mathbb{E}[ {{{\| {{\theta ^*}({{\bar x}^k}) - {\theta ^*}({{\bar x}^{k + 1}})} \|}^2}} |\mathcal{F}^k ] \leqslant mL_{{\theta ^*}}^2{\tau^2 \alpha ^2}{{{\| {{{\bar y}^k}} \|}^2}}.
\end{aligned}
\end{equation}
Next, we will control the term $A_2^{\theta}$ in \eqref{EQ-theta-start} by  Cauchy-Schwartz inequality as follows:
\begin{equation}\label{EQ-A}
\begin{aligned}
A_2^{\theta} \leqslant {\omega _\theta }\beta A_1^{\theta} + \frac{1}{{{\omega _\theta }\beta }} \mathbb{E}[ {{{\| {{\theta ^*}({{\bar x}^k}) - {\theta ^*}({{\bar x}^{k + 1}})} \|}^2}} |\mathcal{F}^k ].
\end{aligned}
\end{equation}
In what follows, by leveraging the results \eqref{EQ-Q-boundness}, \eqref{EQ-ip}, \eqref{EQ-A}, we can control the inner-level errors $  \mathbb{E}[ {{{\| {{\theta ^{k + 1}} - {\theta ^*}({{\bar x}^{k + 1}})} \|^2}}} |\mathcal{F}^k ] $  as follows:
\begin{align}
   &\mathbb{E}[ {{{\| {{\theta ^{k + 1}} - {\theta ^*}({{\bar x}^{k + 1}})} \|}^2}} |\mathcal{F}^k ] \label{EQ-vk-}\\
   \leqslant &( {1 + {\omega _\theta }\beta } )( {1 - \frac{3}{2}\beta \frac{{{\mu _g}{L_{g,\theta }}}}{{{\mu _g} + {L_{g,\theta }}}}} ){\| {{\theta ^k} - {\theta ^*}({{\bar x}^k})} \|^2}  \nonumber \\
   &+ ( {1 + {\omega _\theta }\beta } )\frac{2}{{{\omega _\theta }}}\beta L_{g,\theta }^2{\| {{x^k} - {1_m} \otimes {{\bar x}^k}} \|^2} + ( {1 + {\omega _\theta }\beta } ){\beta ^2}{\sigma ^2} + ( {1 + \frac{1}{{{\omega _\theta }\beta }}} )L_{{\theta ^*}}^2{\tau^2 \alpha ^2} {{{\| {{{\bar y}^k}} \|}^2}}  \nonumber \\
   \leqslant& ( {1 - \frac{{{\mu _g}{L_{g,\theta }}}}{{{\mu _g} + {L_{g,\theta }}}}} {\beta }){\| {{\theta ^k} - {\theta ^*}({{\bar x}^k})} \|^2} \!+\! m\underbrace {\frac{{2L_{{\theta ^*}}^2}}{{{\omega _\theta }{\beta }}}}_{ \triangleq {q_s}}{\tau^2 \alpha ^2} {{{\| {{{\bar y}^k}} \|}^2}} \\
   &+ \underbrace {\frac{4L_{g,\theta }^2\beta}{{{\omega _\theta }}\alpha}}_{ \triangleq {q_x}}\alpha  {\| {{x^k} - {1_m} \otimes {{\bar x}^k}} \|^2}\! +\! m\underbrace {2 \sigma _{g,\theta}^2\frac{\beta ^2}{\alpha^2}}_{ \triangleq {{\sigma _{\theta}^2} }}{\alpha ^2}, \nonumber
\end{align}
where the last inequality is derived by  the fact $\omega_{\theta}\beta<1$ and  $c_{\beta}=\frac{\beta}{\alpha}$.  This completes the proof.   {\hfill $\blacksquare$}
\subsection{Proof  of Lemma \ref{LE-CE}}
\label{sec-proof-CE}
First,  recall the recursion of $x^{k+1}$ as follows:
 \[{x^{k + 1}} =(1-\tau)x^k+ \tau(\mathcal{W}{x^k} - \alpha {y^k}),\]
by which, we have
\begin{equation}\label{EQ-x-average}
\begin{aligned}
 & {x^{k + 1}} - {1_m} \otimes {{\bar x}^{k + 1}} \\
 = &( {1 - \tau } )( {{x^k} - {1_m} \otimes {{\bar x}^k}} ) + \tau ( {( {{{\mathcal{W}}} - {{\mathcal{J}}}} )( {{x^k} - {1_m} \otimes {{\bar x}^k}} ) - \alpha ( {{y^k} - {1_m} \otimes {{\bar y}^k}} )} ) \hfill
\end{aligned}
\end{equation}
where $\mathcal{J}=\frac{1_m 1_m^{\rm T}}{m} \otimes  I_n$.
Then employing Young's inequality with  parameter $\eta>0$ yields
\begin{align}
 & \mathbb{E}[ {{{\| {{x^{k + 1}} - {1_m} \otimes {{\bar x}^{k + 1}}} \|}^2}} {|{{{\mathcal{F}}}^k}} ] \nonumber \\
   \leqslant & {(1 - \tau )^2}(1 + \frac{\tau }{{1 - \tau }}){\| {{x^k} - {1_m} \otimes {{\bar x}^k}} \|^2} + {\tau ^2}(1 + \frac{{1 - \tau }}{\tau }){\| {W{x^k} - \alpha {y^k} - ( {{1_m} \otimes {{\bar x}^k} - \alpha {{\bar y}^k}} )} \|^2} \nonumber  \\
   \leqslant &(1 - \tau ){\| {{x^k} - {1_m} \otimes {{\bar x}^k}} \|^2} + \tau {\| {\left( {{{\mathcal{W}}} - {{\mathcal{J}}}} \right)( {{x^k} - {1_m} \otimes {{\bar x}^k}} ) - \alpha ( {{y^k} - {1_m} \otimes {{\bar y}^k}} )} \|^2} \nonumber  \\
   \leqslant &(1 - \tau ){\| {{x^k} - {1_m} \otimes {{\bar x}^k}} \|^2} + \tau (1 + \eta ){\| {{{\mathcal{W}}} - {{\mathcal{J}}}} \|^2}{\| {{x^k} - {1_m} \otimes {{\bar x}^k}} \|^2} + \tau (1 + \frac{1}{\eta }){\| {{y^k} - {1_m} \otimes {{\bar y}^k}} \|^2} \nonumber \\
   \leqslant &(1 - \tau \frac{{1 - \rho }}{2}){\| {{x^k} - {1_m} \otimes {{\bar x}^k}} \|^2} + \frac{{2\tau {\alpha ^2}}}{{1 - \rho }}{\| {{y^k} - {1_m} \otimes {{\bar y}^k}} \|^2} \label{EQ-y}
\end{align}
where the last step takes $\eta  = \frac{{1 - \rho }}{{2\rho }}$ and uses the fact that ${\left\| {\mathcal{W} - \mathcal{J}} \right\|^2} = \rho $. 
Taking the total expectation, we get the desired result. This completes the proof.   {\hfill $\blacksquare$}

\subsection{Proof of Lemma \ref{LE-heterogeity}} \label{sec-proof-LE-heterogeity}
For the term ${\| {{y^k} - {1_m} \otimes {{\bar y}^k}} \|^2}$, we  can further split it into:
\begin{equation}
\begin{aligned}
  {\| {{y^k} - {1_m} \otimes {{\bar y}^k}} \|^2} \leqslant &{\| {{y^k}} \|^2} \\
   =& \sum\limits_{i = 1}^m {{{\| {y_i^k - \nabla {\Phi _i}( {{{\bar x}^k}} ) + \nabla {\Phi _i}( {{{\bar x}^k}} ) - \nabla \Phi ( {{{\bar x}^k}} ) + \nabla \Phi ( {{{\bar x}^k}} )} \|}^2}}  \\
   \leqslant& 3\sum\limits_{i = 1}^m {{{\| {y_i^k - \nabla {\Phi _i}( {{{\bar x}^k}} )} \|}^2}} \! +\! 3 {\sum\limits_{i = 1}^m {{{\| {\nabla {\Phi _i}( {{{\bar x}^k}} ) \!-\! \nabla \Phi ( {{{\bar x}^k}} )} \|}^2}} }\! + \!3m{\| {\nabla \Phi ( {{{\bar x}^k}} )} \|^2} \\
   \leqslant& 3\sum\limits_{i = 1}^m {{{\| {y_i^k - \nabla {\Phi _i}( {{{\bar x}^k}} )} \|}^2}}  + 3{b^2} + 3m{\| {\nabla \Phi ( {{{\bar x}^k}} )} \|^2} \\
   =& \nyc{3{{{\| {\nabla \tilde \Phi ({{\bar x}^k}) - {z^k}} \|}^2}} + 3{b^2} + 3m{\| {\nabla \Phi ( {{{\bar x}^k}} )} \|^2},}
\end{aligned}
\end{equation}
where the second inequality follows from Lemma \ref{LE-hyper-heterogeneity} and the last step holds due to the fact that  $y _i^k =z _i ^k$ under local gradient scheme \eqref{ALG-LPDBO-h1} and $\nabla \tilde{\Phi} (\bar{x}^k)={\rm col}\{\nabla \Phi _i(\bar{x}^k)\}_{i=1}^{m}$.
Then, taking the total expectation,
we get the desired result. This completes the proof.   {\hfill $\blacksquare$}

\subsection{Proof  of Lemma \ref{LE-TE}}
\label{sec-proof-LE-TE}
According to the recursion \eqref{ALG-LPDBO-h2}, we known that the  update of  $y^{k+1}$ can be derived as:
\[{y^{k + 1}} = \mathcal{W}{y^k} + {z^{k + 1}} - {z^k},\]
which further gives that
\begin{equation}
\begin{aligned}
  {y^{k + 1}} - {1_m} \otimes {{\bar y}^{k + 1}}  = & \mathcal{W}{y^k} + {z^{k + 1}} - {z^k} - {1_m} \otimes ( {{{\bar y}^k} + {{\bar z}^{k + 1}} - {{\bar z}^k}} ) \\
   =& ( {\mathcal{W} - \mathcal{J}} )( {{y^k} - {1_m} \otimes {\bar y^k}} ) + ( {I - \mathcal{J}} )( {{z^{k + 1}} - {z^k}} ). \\
\end{aligned}
\end{equation}
Then taking the square norm on both sides of the above expression under the conditional expectation $\mathcal{F}^k$ and using Young's inequality with the parameter $\eta=\frac{1-\rho}{2\rho}$, we get that:
\begin{equation}\label{EQ-y-y}
\begin{aligned}
  &\mathbb{E}[ {{{\| {{y^{k + 1}}- {1_m} \otimes {{\bar y}^{k + 1}}} \|}^2}} |\mathcal{F}^k ] \\
  \leqslant& ( {1 + \eta } ){\|  {\mathcal{W} - \mathcal{J}} \|^2} {{{\| {{y^k} - {1_m} \otimes {\bar y^k}} \|}^2}} + ( {1 + \frac{1}{\eta }} ){\| {I - \mathcal{J}} \|^2}\mathbb{E}[ {{{\| {{z^{k + 1}} - {z^k}} \|}^2}} |\mathcal{F}^k ] \\
   \leqslant &\frac{{1 + \rho }}{2} {{{\| {{y^k}- {1_m} \otimes {\bar y^k}} \|}^2}}  + \frac{2}{{1 - \rho }}\mathbb{E}[ {{{\| {{z^{k + 1}} - {z^k}} \|}^2}}  |\mathcal{F}^k  ],
\end{aligned}
\end{equation}
where the last inequality uses the fact that $\| {\mathcal{W}} - {\mathcal{J}} \|^2 = \rho $  and  $\left\| {I - \mathcal{J}} \right\|^2 \leqslant 1$. For the last term in \eqref{EQ-y-y}, it follows from  \eqref{ALG-LPDBO-g} that:
\begin{equation}
\begin{aligned}
  {z^{k + 1}} - {z^k}
   = \gamma (\nabla \tilde \Phi ({{\bar x}^k}) - {z^k}) + \gamma ({\mathbb{E}}[{s^k}|{\mathcal{F}^k}] - \nabla \tilde \Phi ({{\bar x}^k})) + \gamma ({s^k} - {\mathbb{E}}[{s^k}|{\mathcal{F}^k}]), \\
\end{aligned}
\end{equation}
which further implies that:
\begin{equation}\label{EQ-z-s}
\begin{aligned}
{\mathbb{E}}[ {{\| z^{k+1}-z^k\|}^2} |{{{\mathcal{F}}}^k} ] \leqslant & 2{\gamma ^2}{{\| \nabla \tilde \Phi ({{\bar x}^k}) - {z^k}\|}^2}  \\
 &+ 2\gamma ^2{\| {\nabla \tilde \Phi ({{\bar x}^k}) - {\mathbb{E}}[{s^k}|{\mathcal{F}^k}]} \|^2} + m{({\sigma _{f,x}^2}+\nyc{\hat M^2}{\sigma _{g,x\theta}^2})}{\gamma ^2}.
\end{aligned}
\end{equation}
Substituting the above inequality into \eqref{EQ-y-y}, we reach
\begin{equation}
\begin{aligned}
 & {\mathbb{E}}[ {{{\| {{y^{k + 1}} - {1_m} \otimes {{\bar y}^{k + 1}}} \|}^2}|{{{\mathcal{F}}}^k}} ] \\
   \leqslant & \frac{{1 + \rho }}{2}{\| {{y^k} - {1_m} \otimes {{\bar y}^k}} \|^2} + \frac{4}{{1 - \rho }}\gamma ^2{\| {\nabla \tilde \Phi ({{\bar x}^k}) - {\mathbb{E}}[{s^k}|{\mathcal{F}^k}]} \|^2} \\
   &+ \frac{4}{{1 - \rho }}{{\gamma^2}}{{\| \nabla \tilde \Phi ({{\bar x}^k}) - {z^k}\|}^2}  + \frac{2}{{1 - \rho }}m\underbrace{\res{({\sigma _{f,x}^2}+\nyc{\hat M^2}{\sigma _{g,x\theta}^2})}\frac{\gamma^2}{\alpha^2}}_{ \triangleq \sigma_y^2}{\alpha ^2},
\end{aligned}
\end{equation}
This completes the proof.   {\hfill $\blacksquare$}



\end{document}